\documentclass{article}
\usepackage{preamble}
\usepackage{preamble_local}

\title{Mean-field variational inference with the TAP free energy: Geometric and statistical properties in linear models}

\author{Michael Celentano\thanks{Department of Statistics, University of California, Berkeley. E-mail: \tt{mcelentano@berkeley.edu}} \and Zhou Fan\thanks{Department of Statistics and Data Science, Yale University. E-mail: \tt{zhou.fan@yale.edu}} \and Licong Lin\thanks{Department of Statistics, University of California, Berkeley. E-mail: \tt{liconglin@berkeley.edu}} \and Song Mei\thanks{Department of Statistics and Department of EECS, University of California, Berkeley. E-mail: \tt{songmei@berkeley.edu}}}

\date{}

\begin{document}

\maketitle

\begin{abstract}

We study mean-field variational inference in a Bayesian linear model when the
sample size $n$ is comparable to the dimension $p$. In high dimensions, the
common approach of minimizing a Kullback-Leibler divergence from the posterior
distribution, or maximizing an evidence lower bound, may deviate from the true
posterior mean and underestimate posterior uncertainty.
We study instead minimization of the TAP free energy, showing in
a high-dimensional asymptotic framework that it has a local minimizer which
provides a consistent estimate of
the posterior marginals and may be used for correctly calibrated
posterior inference. Geometrically, we show that the landscape of the TAP free
energy is strongly convex in an extensive neighborhood of this local minimizer,
which under certain general conditions can be found by an Approximate Message
Passing (AMP) algorithm. We then exhibit an efficient algorithm that linearly
converges to the minimizer within this local neighborhood. In settings where it
is conjectured that no efficient algorithm can find this local neighborhood, we
prove analogous geometric properties for a local minimizer of the TAP free
energy reachable by AMP, and show that posterior inference based on
this minimizer remains correctly calibrated.
\end{abstract}

\tableofcontents

\section{Introduction}

Approximating expectations under high-dimensional posterior
probability distributions is a central goal in Bayesian inference.
For large-scale models, variational inference methods provide a popular
optimization-based approach to perform such approximations.
These methods typically start from a variational representation of the marginal
log-likelihood or model evidence,
\begin{equation}\label{eq:marginalll}
\log \sP(\y)=-\inf_\sQ \int
\bigg[\log \frac{\sQ(\bbeta)}{\sP(\y|\bbeta)\sP_0(\bbeta)}\bigg]\sQ(\de \bbeta)
\end{equation}
where $\sP(\y|\bbeta)$ is the data likelihood,
$\sP_0(\bbeta)$ is the prior density for parameters $\bbeta$,
and $\inf_\sQ$ expresses a minimization over all distributions $\sQ$ for
$\bbeta$. Variational methods then proceed by
minimizing an approximation to (\ref{eq:marginalll}), often restricting to a
computationally tractable sub-class of distributions $\sQ \in \cQ$, and using
the optimizer $\sQ_* \in \cQ$ as an approximate posterior law for $\bbeta$. We
refer readers to \cite{blei2017variational} for a recent review.

In this paper, we study mean-field variational inference for
high-dimensional linear regression models
\[\y=\X\bbeta+\beps \in \R^n\]
where the Bayesian prior for $\bbeta \in \R^p$ is a product distribution
specifying independent and identically distributed coordinates.
By ``mean-field'',
we refer to a choice of sub-class $\cQ$ comprised of product laws.
Mean-field methods for linear regression have found
particular application in statistical genetics, where they have been used to
infer genetic associations and estimate heritability of complex traits
\cite{logsdon2010variational,carbonetto2012scalable} and underlie 
popular linear-mixed-modeling software packages \cite{loh2015efficient}.

A body of recent work in statistical theory has studied the accuracy of
mean-field variational posteriors for linear models in classical regimes of
fixed dimension $p$ \cite{ormerod2017variational}, in posterior-contraction
regimes of high dimension $p$ and strong sparsity
$\|\bbeta\|_0 \ll n/(\log p)$ \cite{yang2020alpha,ray2022variational}, and in
low-complexity regimes encompassing $p \ll n$ or approximately low-rank designs
\cite{mukherjee2022variational}. However, these regimes are arguably far from
the setting of many applications, in which $p$ may be comparable to or 
larger than $n$, all or many variables may each explain an ``infinitesimal'' fraction of the total variance of $\y$, and the posterior law of $\bbeta$
may not contract strongly around the true regression vector. In such settings,
obtaining accurate quantifications of posterior uncertainty remains an important
goal.

Motivated by such applications, we study here the linear model in a
high-dimensional asymptotic framework where $n,p \to \infty$ proportionally,
and all or a fixed proportion of variables
contribute to the variance of $\y$. In this setting, we provide guarantees for
variational inference based upon a conjectured ``TAP approximation'' of the
evidence (\ref{eq:marginalll}),
\[\log \sP(\y) \approx -\inf_{\sQ \in \cQ} \cF_\TAP(\sQ)
\equiv -\inf_{\bbm,\bs} \cF_\TAP(\bbm,\bs)\]
for a non-convex free energy function $\cF_\TAP(\sQ) \equiv \cF_\TAP(\bbm,\bs)$
that depends on $\sQ$ via its marginal first and second moment vectors
$\bbm,\bs \in \R^p$, as defined in (\ref{eqn:TAP}) to follow.
This approximation stems from the work of \cite{thouless1977solution}, was
proposed as an optimization objective for variational inference in linear
models in \cite{krzakala2014variational}, and underlies also the class of
Approximate Message Passing (AMP) algorithms for Bayesian linear regression
developed in \cite{kabashima2003cdma,donoho2009message}. We will study the
setting of i.i.d.\ Gaussian design, which is representative of an ideal
scenario where the posterior correlation between any two variables of $\bbeta$
is weak, and mean-field methods should work
well. We expect our results to hold universally for random designs with
independent and standardized variables, and we discuss this further in
Section \ref{sec:literature} below.

Our work establishes the following main results, with probability
approaching 1 as $n,p \to \infty$ with $n/p \to \delta \in (0,\infty)$:
\begin{enumerate}
\item There exists a local minimizer $(\bbm_\star,\bs_\star)$ of
$\cF_\TAP(\bbm,\bs)$ that is consistent for the posterior first and second
moments of $\bbeta$ and yields a consistent approximation of
the model evidence, in the sense
\[\begin{gathered}
\frac{1}{p}\sum_{j=1}^p \big(m_{\star,j}-\E[\beta_j \mid \X,\y]\big)^2 \gotop
0,\quad
\frac{1}{p}\sum_{j=1}^p \big(s_{\star,j}-\E[\beta_j^2 \mid \X,\y]\big)^2 \gotop 0,\\
\frac{1}{p}\Big({-}\cF_\TAP(\bbm_\star,\bs_\star)-\log \sP(\by)\Big)
\gotop 0.
\end{gathered}\]

We show this by deriving and analyzing a lower bound for
$\cF_\TAP(\bbm,\bs)$ obtained via Gordon's comparison inequality. Interestingly,
this lower bound also verifies that $\cF_\TAP(\bbm,\bs)$ is bounded below
by the replica-symmetric potential when restricted to $(\bbm,\bs)$
satisfying a Nishimori-type condition.

\item Corresponding to this local minimizer $(\bbm_\star,\bs_\star)$ are
dual vectors $(\blambda_\star,\bgamma_\star)$ that characterize the
posterior marginal distributions of $\bbeta$, in the sense that the posterior
law of each variable $\beta_j$ is well-approximated by the exponential-family
law
\[\sP(\de\beta_j) \propto
e^{-(\gamma_{\star,j}/2)\beta_j^2+\lambda_{\star,j}\beta_j}\sP_0(\de \beta_j)\]
where $\sP_0$ is the prior. These laws $\{\sP(\de\beta_j)\}_{j=1,\ldots,p}$
may be used for asymptotically calibrated posterior inference.

\item For the conjectured region of noise variance parameter
$\sigma^2$ and dimension ratio $\delta=\lim_{n,p \to \infty} n/p$ where
Bayes-optimal inference is computationally feasible,
$\cF_\TAP$ is strongly convex in a radius-$O(\sqrt{p})$ local neighborhood of
$(\bbm_\star,\bs_\star)$. (We note that $\cF_\TAP$ is in general not globally convex.)

We propose and analyze a Natural Gradient Descent (NGD) algorithm
that exhibits linear
convergence within this neighborhood, and hence efficiently computes
$(\bbm_\star,\bs_\star,\blambda_\star,\bgamma_\star)$ from a local
initialization. Such a local initialization may be obtained via a finite
number of iterations of AMP.

\item More generally, for any $(\sigma^2,\delta)$, $\cF_\TAP$ is
convex in a radius-$O(\sqrt{p})$ neighborhood of the local minimizer
$(\bbm_\alg,\bs_\alg)$ that is reachable by AMP, and NGD exhibits linear
convergence to $(\bbm_\alg,\bs_\alg,\blambda_\alg,\bgamma_\alg)$.

We show these statements of local convexity by developing a version of Gordon's
comparison inequality conditional on the filtration generated by the iterates
of AMP, which may be of independent interest.

When $(\bbm_\alg,\bs_\alg) \neq (\bbm_\star,\bs_\star)$, the corresponding
laws $\sP(\de\beta_j) \propto
e^{-(\gamma_{\alg,j}/2)\beta_j^2+\lambda_{\alg,j}\beta_j}\sP_0(\de \beta_j)$ do
not consistently approximate the true posterior marginals.
However, we argue that posterior inference based upon these laws remains
well-calibrated.
\end{enumerate}

It is conjectured that the statements of (1.)\ and (2.)\ above may in fact hold
for $(\bbm_\star,\bs_\star)$ being the \emph{global} minimizer of $\cF_\TAP$.
Our results imply this conjecture under either sufficiently large noise
variance $\sigma^2$ or sufficiently small dimension ratio $\delta=\lim_{n,p \to
\infty} n/p$, because in these settings it may be checked that
$\cF_\TAP$ is globally convex and has
only a single local minimizer. When the prior distribution of $\bbeta$ is
uniform on a high-dimensional sphere,
a version of this conjecture for sufficiently large noise variance
$\sigma^2$ has also been shown previously in \cite{qiu2022tap}.
The techniques of \cite{qiu2022tap} are
specific to the spherical symmetry of the prior, and different from our
analyses here for product priors.

We provide more formal statements and further discussion of these results
in Section \ref{sec:results}.

\subsection{Further related literature}\label{sec:literature}

Our results build upon a series of works
\cite{montanari2006analysis,reeves2016replica,barbier2016mutual,barbier2019adaptive,barbier2019optimal,barbier2020mutual}
concerning Bayes-optimal inference for the high-dimensional linear model with
i.i.d.\ designs, which have made rigorous pioneering insights of
Tanaka \cite{tanaka2002statistical}
derived initially using
statistical mechanics techniques. Among other results, these works showed that
the asymptotic values of the model evidence and squared-error Bayes risk for
estimating the regression vector $\bbeta$ are
determined by the minimization of a scalar replica-symmetric potential.
If this potential has a unique critical point or, more generally, a global
minimizer coinciding with its closest local minimizer to 0, then an
AMP algorithm succeeds in computing an approximate posterior mean vector that
asymptotically attains near-optimal Bayes risk. We review some of these results
relevant to our work in Section \ref{sec:model}, and refer readers
to \cite{barbier2019optimal} for further details.

Stable fixed points of this AMP algorithm correspond to local minimizers of
the TAP free energy, although current AMP theory does not guarantee the
existence of, or convergence to, such fixed points for any finite $n$ and $p$.
We refer to recent results of \cite{li2022non,li2023approximate}
that make progress in this direction.
Our results complement the state-evolution theory of AMP by
rigorously establishing that $\cF_\TAP$ has a local minimizer that is asymptotically consistent
for the posterior marginals, and around which the landscape is locally convex.
We show that local convexity enables the convergence of alternative optimization procedures. This
perspective follows that of \cite{krzakala2014variational}, who
proposed the direct minimization of a version of the TAP free energy
(a.k.a.\ Bethe free energy) as an approach to variational inference,
and \cite{fan2021tap,celentano2023local} who studied analogous questions
regarding the TAP free energy in the $\Z_2$-synchronization spiked matrix model.
Recently, \cite{qiu2022tap} studied a
version of the TAP free energy in the linear model for a uniform prior on the
sphere, showing that its global minimizer has similar statistical properties
via different geometric techniques.

We analyze both the minimum value of $\cF_\TAP$ and the
smallest eigenvalue of its Hessian as min-max optimizations of
Gaussian processes defined by the design matrix $\X \in \R^{n \times p}$. We
apply Gordon's comparison inequality to bound these values via
auxiliary processes defined by i.i.d.\ Gaussian vectors $\g \in \R^n$ and
$\h \in \R^p$. This strategy is closest to that
of \cite{celentano2023local}, which applied related ideas around the
Sudakov-Fernique inequality for symmetric matrices to analyze
$\Z_2$-synchronization.
To establish local convexity of $\cF_\TAP$ in a neighborhood of
the (random) point $(\bbm_\star,\bs_\star)$, we circumvent the Kac-Rice analyses
of \cite{fan2021tap,celentano2023local} and instead extend the approach of
\cite{celentano2022sudakov} to an asymmetric setting, proving a version of
Gordon's inequality conditional on the filtration generated by AMP.
We remark that, in contrast to applications of the Convex Gaussian
Minmax Theorem (CGMT) \cite{stojnic2013framework,thrampoulidis2018precise}, the
Gaussian processes we analyze do not in general have a globally convex-concave
structure. Interestingly, our analyses show that a specialization of the Gordon
inequality lower bound to a domain of $(\bbm,\bs)$ that obeys a Nishimori-type
property recovers the replica-symmetric potential.

We expect the main results of our work to hold universally for random
designs $\X \in \R^{n \times p}$
having independent entries of mean 0, common variance, and sufficiently fast
tail decay. This universality class may be the pertinent
one for genetic association analyses of common variants at unlinked loci,
where genotypes are nearly independent due to recombination
\cite{carbonetto2012scalable}. Universality for some of our results, e.g.\ the
validity of our Gordon lower bound for $\cF_\TAP$ and the existence of a local
minimizer $(\bbm_\star,\bs_\star)$ that yields consistent approximations for the
posterior marginals, are readily obtained by combining our current arguments
with existing universality results for Bayes-optimal estimation in the linear
model \cite{barbier2019optimal}, Gordon-type min-max optimization problems
\cite{han2022universality}, and state evolutions of AMP algorithms
\cite{bayati2015universality,chen2021universality}. Verifying universality of
the local convexity of $\cF_\TAP$ near $(\bbm_\star,\bs_\star)$ seems more
challenging, and we leave this as an interesting mathematical question for
future work.

\section{Background on the Bayesian linear model}\label{sec:model}

In this section, which is mostly expository,
we review relevant background concerning Bayesian
inference and the TAP free energy for high-dimensional linear models.

We consider the linear model
\begin{equation}\label{eq:linearmodel}
\y=\X\bbeta+\beps
\end{equation}
with i.i.d.\ Gaussian design and Gaussian noise.
We reserve the notation $\bbeta_0 \in \R^p$ for the true regression vector,
and assume throughout the following conditions.

\begin{assumption}\label{ass:Bayesian_linear_model}
The design matrix $\X \in \R^{n \times p}$, true regression vector $\bbeta_0 \in
\R^p$, and residual error $\beps \in \R^n$ are independent, with entries
\[x_{ij} \overset{iid}{\sim} \normal(0,\tfrac{1}{p}),
\qquad \beta_{0,j} \overset{iid}{\sim} \sP_0, \qquad
\eps_i \overset{iid}{\sim} \normal(0,\sigma^2)\]
where $\sigma^2>0$. 
Here, $\sP_0$ is a prior distribution on $\R$ having compact support with at
least three distinct values. As $n,p \to \infty$,
we have $n / p \to \delta \in (0,\infty)$, and $(\sigma^2,\delta)$
and $\sP_0$ are fixed independently of $n,p$.
\end{assumption}

We allow $\sP_0$ to have a delta mass at 0 to model sparsity for a constant
fraction of all variables. We caution the reader that we scale $x_{ij}$ to
have variance $1/p$, so that the variance explained by
$\X\bbeta$ for each sample is approximately $\E_{\beta_0 \sim \sP_0}[\beta_0^2]$ and
independent of $\delta$. A rescaling of the prior by $\sqrt{\delta}$ would be
needed to translate our formulae and results to a setting where $x_{ij} \sim
\normal(0, 1 / n)$. 

\subsection{The asymptotic evidence and Bayes risk}

We review here several results of \cite{barbier2019optimal,barbier2020mutual},
borrowing also from the notational conventions and presentation in
\cite{celentano2022fundamental}.

In the linear model (\ref{eq:linearmodel}),
define the (normalized) evidence or marginal log-likelihood for $\y$,
\begin{equation}\label{eq:evidence}
F_p(\y):=\frac{1}{p}\log \sP(\y)=\frac{1}{p}
\log \int (2\pi\sigma^2)^{-n/2}\exp\Big({-}\frac{1}{2\sigma^2}\big\|\y-\X\bbeta\big\|_2^2\Big)
\prod_{j=1}^p \sP_0(\de\beta_j).
\end{equation}
For any function
$f:\R^p \to \R^k$, we will denote the posterior expectation of $f(\bbeta)$ by
\[\langle f(\bbeta) \rangle_{\X,\y}=\E[f(\bbeta) \mid \X,\y].\]
In particular, $\langle \bbeta \rangle_{\X,\y}$ is the posterior-mean estimate
of $\bbeta_0$. Define the per-coordinate squared-error Bayes risk (MMSE) for
estimating $\bbeta$,
\[
\MMSE_p := \frac{1}{p}\,\E\big[\big\| \bbeta_0 - \langle \bbeta
\rangle_{\X,\y} \big\|_2^2\big].
\]

The asymptotic evidence and MMSE are related to inference in the following
scalar channel: Fixing a given signal-to-noise parameter $\gamma>0$, 
consider the model
\begin{equation}\label{eq:scalarmodel}
\lambda=\gamma \beta_0+\sqrt{\gamma} z,
\qquad (\beta_0,z) \sim \sP_0 \times \normal(0,1).
\end{equation}
Under this model, the Bayes risk for estimating $\beta_0$ from the observation
$\lambda$ is
\begin{align}
\mmse(\gamma) &= \E_{(\beta_0, z) \sim \sP_0 \times \normal(0, 1)}\big[\big(\beta_0
- \E[\beta_0 \mid \gamma \beta_0 + \sqrt{\gamma} z]\big)^2\big].
\label{eqn:mmse_scalar}
\end{align}
Define the \emph{replica-symmetric potential}
\begin{equation}\label{eqn:phi_potential}
\phi(\gamma)=\frac{\sigma^2
\gamma}{2}-\frac{\delta}{2}\log \frac{\gamma}{2\pi\delta} + i(\gamma)
\end{equation}
where $i(\gamma)$
is the mutual information between $\beta_0$ and $\lambda$. Then by the
I-MMSE relation
$\frac{\de }{\de \gamma}i(\gamma) = \frac{1}{2} \mmse(\gamma)$
\cite[Corollary 1]{guo2011estimation}, the critical points of $\phi(\gamma)$
are the roots of the fixed-point equation
\begin{equation}\label{eq:phifixedpoint}
\mmse(\gamma)=\delta/\gamma-\sigma^2.
\end{equation}
We will assume for most results the following additional condition.
\begin{assumption}\label{ass:uniquemin}
The global minimizer
\begin{equation}\label{eq:gammastat}
\gamma_\stat:=\argmin_{\gamma>0} \phi(\gamma)
\end{equation}
exists and is unique, and $\phi''(\gamma_\stat)>0$ strictly.
\end{assumption}

This assumption imposes a mild genericity condition for $(\sigma^2,\delta)$:
Fixing any $\sigma^2>0$, the global minimizer
$\gamma_\stat$ exists and is unique for Lebesgue-a.e.\ $\delta>0$
(c.f.\ \cite[Proposition 1]{barbier2019optimal}). Furthermore, 
fixing $\sigma^2>0$, Sard's theorem
(c.f.\ \cite[Chapter 1.7]{guillemin2010differential}) 
implies that $\phi(\gamma)$ is Morse for Lebesgue-a.e.\ $\delta>0$,
i.e.\ $\phi''(\gamma) \neq 0$ whenever $\phi'(\gamma)=0$. Thus Assumption
\ref{ass:uniquemin} holds for all $\sigma^2>0$ and Lebesgue-a.e.\ $\delta>0$.

The following theorem is a direct consequence of the results of
\cite[Theorems 1, 2, 6]{barbier2019optimal} specialized to the linear model;
see also \cite{barbier2020mutual} for this specialization.
\begin{theorem}[\cite{barbier2019optimal}]\label{thm:freeenergy}
Suppose Assumption \ref{ass:Bayesian_linear_model} holds. Then
\[\lim_{n,p \to \infty} \E[F_p(\y)]=-\inf_{\gamma>0} \phi(\gamma)
=-\phi(\gamma_\stat),
\qquad \Var[F_p(\y)] \leq C/p\]
for a constant $C:=C(\sigma^2,\delta,\sP_0)>0$ and all sufficiently
large $n,p$. If also Assumption \ref{ass:uniquemin} holds, then
\[
\lim_{n, p \to \infty} \MMSE_p = \mmse(\gamma_{\stat}). 
\]
\end{theorem} 



\subsection{Mean-field approximation and the TAP free energy}

For $(\lambda,\gamma) \in \R^2$ and $\sP_0$ the prior distribution of
coordinates of $\bbeta$, consider the two-parameter exponential family laws
\begin{equation}\label{eq:scalarmean}
\sP_{\lambda,\gamma}(\de\beta) \propto e^{-(\gamma/2)\beta^2+\lambda
\beta}\sP_0(\de\beta), \qquad
\langle f(\beta) \rangle_{\lambda,\gamma}=\int f(\beta)\,
\sP_{\lambda,\gamma}(\de\beta). 
\end{equation}
Note that for $\gamma>0$, this is the posterior distribution of
$\beta_0$ and its associated posterior expectation in the scalar channel
model (\ref{eq:scalarmodel}) with observation $\lambda=\gamma
\beta_0+\sqrt{\gamma} z$.

This exponential family is minimal under Assumption
\ref{ass:Bayesian_linear_model} that $\sP_0$ has at least three points of
support, implying the following statements
(c.f.\ \cite[Proposition 3.2, Theorem 3.3]{wainwright2008graphical}):
The moment map $(\lambda,\gamma) \mapsto
(m,s)=(\langle \beta \rangle_{\lambda,\gamma},
\langle \beta^2 \rangle_{\lambda,\gamma})$ is bijective from
$\R^2$ onto its image
\begin{equation}\label{eq:momentspace}
\Gamma=\Big\{(m,s) \in \R^2: \text{ there exist } (\lambda,\gamma) \in
\R^2 \text{ such that } m=\langle \beta \rangle_{\lambda,\gamma}
\text{ and } s=\langle \beta^2 \rangle_{\lambda,\gamma}\Big\},
\end{equation}
with inverse function over $(m,s) \in \Gamma$ given by
\begin{equation}\label{eq:lambdagammastar}
\big(\lambda(m,s),\,\gamma(m,s)\big)=\argmax_{(\lambda,\gamma) \in \R^2}
-\frac{1}{2} \gamma s + \lambda m -  \log \E_{\beta \sim \sP_0}
\big[e^{-(\gamma/2) \beta^2 + \lambda \beta}\big].
\end{equation}
The domain $\Gamma$ in (\ref{eq:momentspace}) is convex and open in $\R^2$,
and for each $(m,s) \in \Gamma$ the maximizer $(\lambda(m,s),\gamma(m,s))$
in (\ref{eq:lambdagammastar}) is unique.
In Appendix \ref{sec:domain}, we explicitly characterize the set $\Gamma$.

Restricting the variational representation of the evidence
(\ref{eq:marginalll}) to distributions $\sQ=\prod_{j=1}^p
\sP_{\lambda(m_j,s_j),\gamma(m_j,s_j)}$ comprised of products
of such exponential family laws, direct calculation then gives
\begin{equation}\label{eq:naiveMF}
\int \bigg[\log \frac{\sQ(\bbeta)}{\sP(\y|\bbeta)\sP_0(\bbeta)}\bigg]\sQ(\de
\bbeta)= \underbrace{\frac{n}{2}\log
2\pi\sigma^2+D_0(\bbm,\bs)+\frac{1}{2\sigma^2}\big\|\y-\X\bbm\big\|_2^2
+\frac{n}{2\sigma^2}\big(S(\bs)-Q(\bbm)\big)}_{\cF_{\rm MF}(\bbm, \bs)}+o_\P(n),
\end{equation}
where
\[D_0(\bbm,\bbs)=\sum_{j=1}^p -\sh(m_j,s_j),
\qquad S(\bbs) = \frac{1}{p} \sum_{j = 1}^p s_j, \qquad
Q(\bbm) = \frac{1}{p} \sum_{j = 1}^p m_j^2,\]
and $-\sh(m,s)=D_{\mathrm{KL}}\big(\sP_{\lambda(m,s),\gamma(m,s)}\|\sP_0\big)$
is the relative entropy or Kullback-Leibler divergence from the prior $\sP_0$ to
the above exponential family law, given explicitly by
\begin{equation}\label{eq:hdef}
-\sh(m,s)=\sup_{(\lambda,\gamma) \in \R^2} {-}\frac{1}{2} \gamma s
+\lambda m -\log \E_{\beta \sim \sP_0}
\big[e^{-(\gamma/2) \beta^2 + \lambda \beta}\big].
\end{equation}
The free energy $\cF_{\rm MF}$ at the right side of (\ref{eq:naiveMF}) is sometimes referred to as the
``na\"ive mean-field free energy''. It has often served as the optimization
objective for variational inference in applications
\cite{carbonetto2012scalable}, and has been analyzed theoretically in
low-complexity and strong posterior-contraction regimes
\cite{mukherjee2022variational,ray2022variational}.

We add to (\ref{eq:naiveMF}) an Onsager correction term
\[\frac{n}{2}\left[\log\left(1+\frac{S(\bs)-Q(\bbm)}{\sigma^2}\right)
-\frac{S(\bs)-Q(\bbm)}{\sigma^2}\right]\]
that, in the asymptotic setting of Assumption \ref{ass:Bayesian_linear_model},
accounts for a difference between the value of (\ref{eq:naiveMF})
and the true model evidence that is given by
optimizing over all (non-product) distributions $\sQ$ in (\ref{eq:marginalll}).
This yields the \emph{TAP free energy},
defined over the domain $(\bbm,\bs) \in \Gamma^p$,
\begin{equation}\label{eqn:TAP}
\cF_{\TAP}(\bbm, \bbs) = \frac{n}{2} \log 2\pi\sigma^2 +
D_0(\bbm, \bbs) + \frac{1}{2\sigma^2} \| \y - \X \bbm
\|_2^2 + \frac{n}{2}\log\left(1+\frac{S(\bbs) -
Q(\bbm)}{\sigma^2}\right). 
\end{equation}

We illustrate the role of the Onsager correction for a simple example with
Gaussian prior in Section \ref{sec:gauss-prior}, and for more complex priors
in the simulations of Section \ref{sec:simulations} to follow.
For heuristic derivations of the TAP free energy,
we refer readers to \cite[Section 3.4.2]{maillard2019high} for an
approach from high-temperature expansions, and
\cite[Section III.B]{krzakala2012probabilistic} for an approach from belief
propagation and AMP.

\begin{remark}[Versions of the TAP free energy]
The TAP free energy (\ref{eqn:TAP}) coincides with the form of the Bethe free
energy in \cite{krzakala2014variational} upon replacing several instances of
$1/p$ by $x_{ij}^2$ and reparametrizing (\ref{eqn:TAP}) by the variables
$(\lambda_j,\gamma_j)=(\lambda(m_j,s_j),\gamma(m_j,s_j))$. We expect
these forms to have similar asymptotic properties as $n,p \to \infty$.

One may also study a reduced version of the TAP free energy, by observing that
at any critical point $(\bbm,\bs)$ of $\cF_\TAP$, the stationarity condition
for $\bbs$ yields $\gamma(m_j,s_j)=\delta/(\sigma^2+S(\bs)-Q(\bbm))$ which is
constant across coordinates $j=1,\ldots,p$. Then, identifying $\gamma \equiv
\gamma(m_j,s_j) \in \R$, the critical points of $\cF_\TAP$ are in
correspondence with those of
\begin{align}
\cF_\TAP^{\text{reduced}}(\bbm,S)
&=\max_{\blambda \in \R^p,\,\gamma \in \R} \sum_{j=1}^p
\bigg(-\frac{1}{2}\gamma S+\lambda_j m_j-\log \E_{\beta_j \sim \sP_0}
\big[e^{-(\gamma/2)\beta_j^2+\lambda_j \beta_j}\big]\bigg)\nonumber\\
&\hspace{1in}+\frac{1}{2\sigma^2}\big\|\y-\X\bbm\big\|_2^2+\frac{n}{2}\log 2\pi\sigma^2
+\frac{n}{2}\log\left(1+\frac{S-Q(\bbm)}{\sigma^2}\right)\label{eq:TAPreduced}
\end{align}
which replaces $\bbs \in \R^p$ by a scalar second-moment parameter
$S \equiv S(\bbs) \in \R$. In this work, we will study the optimization
landscape of the non-reduced free energy function (\ref{eqn:TAP}) over
$(\bbm,\bs)$, rather than of the reduced form (\ref{eq:TAPreduced}).
\end{remark}


\subsection{Approximate Message Passing}

We review an iterative Approximate Message Passing (AMP) algorithm for
Bayes posterior-mean inference in the linear model (\ref{eq:linearmodel}),
as described in \cite{donoho2010message,celentano2022fundamental}.
For $\x \in \R^p$, define
\begin{equation}\label{eq:denoiser}
\denoiser(\x,\gamma)=(\langle \beta \rangle_{\gamma x_j,\gamma})_{j=1}^p \in
\R^p,\qquad
\sS(\x,\gamma)=(\langle \beta^2 \rangle_{\gamma x_j,\gamma})_{j=1}^p \in \R^p
\end{equation}
where $\langle \cdot \rangle_{\lambda,\gamma}$ denotes the mean under
the exponential family model (\ref{eq:scalarmean}).
The AMP algorithm takes the form, with initializations $\bz^0=0$, $\bbm^1=0$,
and $\gamma_1=\delta(\sigma^2+\E_{\beta \sim \sP_0}[\beta^2])^{-1}$,
\begin{equation}\label{eq:AMPalg}
\begin{aligned}
\z^k&=\y-\X\bbm^k+\frac{\gamma_{k-1}\mmse(\gamma_{k-1})}{\delta}\,\z^{k-1},\\
\bbm^{k+1}&=\denoiser\left(\bbm^k+\frac{1}{\delta}\X^\top
\z^k,\,\gamma_k\right),
\quad \bs^{k+1} =
		\sS\left(
			\bbm^k + \frac1\delta \bX^\top \bz^k,\,\gamma_k
		\right),\\
\gamma_{k+1}&=\delta[\sigma^2+\mmse(\gamma_k)]^{-1}.
\end{aligned}
\end{equation}
It may be checked that fixed points $(\bbm,\bs)$ of this AMP
algorithm correspond approximately to stationary points of $\cF_\TAP(\bbm,\bs)$.


A rigorous state evolution characterization of AMP is available from \cite[Theorem~1]{bayati2011dynamics}, summarized in the theorem below. 

\begin{theorem}[\cite{bayati2011dynamics}]
Under Assumption \ref{ass:Bayesian_linear_model}, for each fixed iteration $k
\geq 1$ and
any test function $\psi:\R^2 \to \R$ satisfying
$|\psi(x)-\psi(x')| \leq C\|x-x'\|_2(1+\|x\|_2+\|x'\|_2)$ for some constant
$C>0$, the iterates $\bbm^{k+1},\gamma_k$ of (\ref{eq:AMPalg}) and the true
regression vector $\bbeta_0$ satisfy
\[\lim_{n,p \to \infty} \frac{1}{p}\sum_{j=1}^p \psi(m^{k+1}_j,\beta_{0,j})
=\E_{(\beta_0,z) \sim \sP_0 \times \normal(0,1)}\left[\psi\Big(\langle \beta
\rangle_{\gamma_k \beta_0+\sqrt{\gamma_k} z},\beta_0\Big)\right].\]
\end{theorem}

Applying this result with $\psi(m,\beta_0)=(m-\beta_0)^2$ shows that
$\lim_{n,p \to \infty} p^{-1}\|\bbm^{k+1}-\bbeta_0\|_2^2=\mmse(\gamma_k)$,
the scalar channel MMSE from (\ref{eqn:mmse_scalar}). Furthermore
(c.f.\ \cite[Proposition 2.4]{celentano2022fundamental}),
\begin{equation}\label{eq:gammaalg}
\lim_{k \to \infty}
\gamma_k=\gamma_\alg:=\inf\big\{\gamma>0:\delta/\gamma-\sigma^2 \leq
\mmse(\gamma)\big\},
\end{equation}
where, in light of (\ref{eq:phifixedpoint}), $\gamma_\alg$ describes the
closest local minimizer to 0 of the replica-symmetric potential $\phi(\gamma)$.
Thus when
$\gamma_\alg=\gamma_\stat$ from (\ref{eq:gammastat}), AMP asymptotically attains
the optimal squared-error Bayes risk, in the sense
$\lim_{k \to \infty} \lim_{n,p \to \infty}
p^{-1}\|\bbm^{k+1}-\bbeta_0\|_2^2=\mmse(\gamma_\stat)$. For later reference,
we state the equality of $\gamma_\stat$ and $\gamma_\alg$ here as a final
condition.

\begin{assumption}[The ``easy'' regime]\label{ass:bayesoptimal}
The prior $\sP_0$ and parameters $(\sigma^2,\delta)$ are such that
$\gamma_\stat$ in (\ref{eq:gammastat}) and $\gamma_\alg$ in
(\ref{eq:gammaalg}) coincide.
\end{assumption}

For bounded prior, Assumption \ref{ass:bayesoptimal} will hold for large enough $\sigma^2$ or small enough $\delta$. We refer to \
\cite{celentano2022fundamental} for plots illustrating examples in which this assumption holds and does not hold. For $\sP_0$ and $(\sigma^2,\delta)$ where
Assumption \ref{ass:bayesoptimal} does not hold, it is conjectured (c.f.\
\cite[Conjecture 1.1]{celentano2022fundamental}) that no polynomial-time
algorithm can achieve Bayes risk asymptotically smaller than
$\mmse(\gamma_\alg)$. In the following, we will establish results both
when Assumption \ref{ass:bayesoptimal} does and does not hold.

\subsection{An example with Gaussian prior}
\label{sec:gauss-prior}

For expositional purposes, we illustrate here a simple example with a Gaussian
prior $\sP_0=\normal(0,\tau^2)$, where exact posterior calculations may be
explicitly performed. This example demonstrates the underestimation of posterior
variance that is exhibited by na\"ive mean-field variational methods, and the
role of the Onsager correction in the TAP approach. (This section is not needed
for the rest of the paper, and can be skipped by the impatient reader.)

The overconfidence of na\"ive mean-field has been observed in numerous prior studies \cite{WangTitterington2005,carbonetto2012scalable,GiordanoBroderickJordan2018}.
The paper \cite{turner_sahani_2011} shows that the mean-field approximation of
correlated Gaussians underestimates the marginal variances. Our example here
is similar, but pertains to a high-dimensional Gaussian distribution in which
pairwise correlations are vanishingly small.

Consider the Gaussian prior $\sP_0=\normal(0,\tau^2)$. For this prior, a simple
computation shows that the model evidence is
\[\log \sP(\y)=-\frac{1}{2}\left[n \log 2\pi+\log \det
\left(\tau^2 \X\X^\top+\sigma^2 \I\right)+\y^\top
\left(\tau^2 \X\X^\top+\sigma^2 \I\right)^{-1}\y\right],\]
and the posterior distribution for $\bbeta$ is the multivariate Gaussian law 
\[\sP(\bbeta \mid \y,\X)=\N\left(\frac{\bSigma \X^\top \y}{\sigma^2},\;
\bSigma\right), \qquad
\bSigma:=\left(\frac{1}{\tau^2}\I+\frac{1}{\sigma^2}\X^\top
\X\right)^{-1}.\]
We remark that when $x_{ij} \overset{iid}{\sim} \normal(0,1/p)$ and
$n \asymp p$, this posterior is not well-approximated
in KL-divergence by any product distribution, even though all pairwise
correlations between variables are of vanishing size $O_\P(1/\sqrt{n})$.

In the high-dimensional asymptotic setting of Assumption
\ref{ass:Bayesian_linear_model}, we have
$p^{-1}\Tr \bSigma \to v_\star$ in probability, where
this value $v_*$ is the Stieltjes
transform of a Marcenko-Pastur law describing the limit eigenvalue
distribution of $\X^\top \X$. Explicitly (c.f.\ \cite[Example
2.8]{tulino2004random}), $v_*$ is given by
the unique positive root of the quadratic equation
\begin{equation}\label{eq:MPvstar}
\frac{1}{v_\star}=\frac{1}{\tau^2}+\frac{\delta}{\sigma^2+v_\star}. 
\end{equation}
By a simple Gaussian concentration-of-measure argument, which we omit
here for brevity, each diagonal
entry $\Sigma_{jj}$ concentrates with variance $O(1/n)$ around its mean,
and hence converges also in probability to $v_\star$. Thus, for large $n,p$,
the true posterior marginals are given by
\begin{equation}\label{eq:Gaussianmarginals}
\sP(\beta_j \mid \y,\X) \approx \N\left(\frac{\bSigma_{j \cdot} \X^\top
\y}{\sigma^2},\;v_\star\right)
\end{equation}
for $j=1,\ldots,p$, where $\bSigma_{j\cdot}$ is the $j^\text{th}$ row of
$\bSigma$.

Reparametrizing with the marginal variance $v_j=s_j-m_j^2$ in place of $s_j$,
the relative entropy (\ref{eq:hdef}) has the explicit form
\[-\sh(m_j,s_j)=D_{\mathrm{KL}}\big(\normal(m_j,v_j)\,\|\,
\normal(0,\tau^2)\big)=\frac{1}{2}\left(\frac{m_j^2+v_j}{\tau^2}-\log
\frac{v_j}{\tau^2}-1\right). \]
Differentiating in $(m_j,v_j)_{j=1}^p$, the minimizer of the
na\"ive mean-field free energy (\ref{eq:naiveMF}) is then given by
\[m_j=\frac{\bSigma_{j\cdot}\X^\top \y}{\sigma^2}, \qquad
v_j \approx \left(\frac{1}{\tau^2}+\frac{\delta}{\sigma^2}\right)^{-1}\]
for each $j=1,\ldots,p$.
Comparing with (\ref{eq:MPvstar}) and (\ref{eq:Gaussianmarginals}) illustrates that this approach
recovers the correct marginal posterior means, but gives inconsistent (and overconfident) estimates
of the marginal posterior variances because $v_j<v_\star$. We remark that this
consistency of the posterior mean estimate is special to the Gaussian
prior, and does not hold more generally as shown in simulation in
Section \ref{sec:simulations}.

In contrast, differentiating the TAP free energy (\ref{eqn:TAP})
shows that it has minimizer $m_j=\bSigma_{j\cdot}\X^\top\y/\sigma^2$ 
and $v_j=v_\star$ for all $j=1,\ldots,p$, thus consistently recovering both the
marginal means and variances. Furthermore, letting $(\bbm_\star,\bs_\star)$
denote this minimizer of $\cF_\TAP$ where $s_{\star,j}=m_{\star,j}^2+v_\star$,
we have the identity
\[\frac{1}{2}\,\y^\top\big(\tau^2\X\X^\top+\sigma^2\I\big)^{-1}\y
=\frac{1}{2\sigma^2}\|\y-\X\bbm_\star\|_2^2+\frac{1}{2\tau^2}
\|\bbm_\star\|_2^2.\]
We also have the convergence in probability (c.f.\ \cite[Examples 2.14,
2.10]{tulino2004random})
\begin{align*}
\frac{1}{p}\log \det\big(\tau^2\X\X^\top+\sigma^2\I)
&=\frac{n}{p}\log \sigma^2+\frac{1}{p}\log
\det\left(\frac{\tau^2}{\sigma^2}\X^\top\X+\I\right)\\
&\to \delta\log
\sigma^2+\log\left(1+\frac{v_\star+\tau^2(\delta-1)}{\sigma^2}\right)
+\delta \log\left(1+\frac{v_\star}{\sigma^2}\right)
-\left(1-\frac{v_*}{\tau^2}\right)\\
&=\delta \log(\sigma^2+v_\star)-\log\frac{v_\star}{\tau^2}
+\left(\frac{v_*}{\tau^2}-1\right).
\end{align*}
Thus $\frac{1}{p}\log \sP(\y)+\frac{1}{p}\cF_\TAP(\bbm_\star,\bs_\star) \to 0$,
so that to leading order in $n,p$ we have
$\log \sP(\y) \approx -\cF_\TAP(\bbm_\star,\bs_\star)$.

\section{Main results}\label{sec:results}

In this section, we present our main results.
Section \ref{sec:Bayes-optimal-local-minimizer} establishes that in all regimes
of $(\sigma^2,\delta)$,
the TAP free energy has a local minimizer which consistently approximates the true posterior marginals (Theorems \ref{thm:bayes_consistence} and \ref{thm:marginalposterior}). 

The existence of a local minimizer does not guarantee that it can found by efficient algorithms.
In Section \ref{sec:local-convexity},
we take a step towards developing such algorithms by showing that the TAP local minimizer is contained in a neighborhood of strong convexity with radius $O(\sqrt{p})$.
In the regime of $(\sigma^2,\delta)$ where $\gamma_\stat=\gamma_\alg$
and Bayes-optimal inference is computationally ``easy'',
such a neighborhood can be reached by AMP. Thus,
following a strategy similar to that of \cite{celentano2023local},
we describe in Section \ref{sec:NGD} a Natural Gradient Descent (NGD) algorithm,
initialized with a constant number of iterations of AMP, that exhibits linear
convergence to this Bayes-optimal local minimizer.

Finally, in Section \ref{sec:hard-regime},
we consider the ``hard'' regime of $(\sigma^2,\delta)$ where
$\gamma_\alg \neq \gamma_\stat$, and it is conjectured that no polynomial-time
algorithm can reach a $\eps\sqrt{p}$-neighborhood of the Bayes posterior-mean
estimate for $\eps>0$ a sufficiently small constant.
Even in this hard regime, we show that AMP after a constant number of iterations arrives
in a region of local strong convexity of the TAP free energy, having radius
$O(\sqrt{p})$ and containing a local minimizer. Thus,
the algorithm described above converges instead to this local minimizer.
Although this minimizer does not correspond to the true Bayes posterior marginals,
we nevertheless show that it provides asymptotically calibrated statements about
posterior uncertainty, and can thus serve as the basis for valid posterior inference.

\subsection{The Bayes-optimal local minimizer} 
\label{sec:Bayes-optimal-local-minimizer}

Our first main result concerns the existence of a local minimizer of the TAP
free energy near the true posterior mean and marginal second moments of
$\bbeta$. Denote these marginal first and second moments by
\begin{equation}\label{eq:marginalmoments}
	\bbmB = \Big(\langle \beta_j \rangle_{\X,\y}\Big)_{j=1}^p,
	\qquad
	\bsB = \Big(\langle \beta_j^2 \rangle_{\X,\y}\Big)_{j=1}^p.
\end{equation}
With high probability,
the TAP free energy has a local minimizer $\bbm_\star,\bs_\star$ approximating
these marginal moments.
\begin{theorem}\label{thm:bayes_consistence}
Let Assumptions \ref{ass:Bayesian_linear_model} and \ref{ass:uniquemin} hold.
Then with probability approaching $1$ as $n, p \to \infty$,
there exists a local minimizer $(\bbm_\star,\bbs_\star) \in \Gamma^p$
of $\cF_\TAP(\bbm,\bs)$ (where $\Gamma$ is the moment space
(\ref{eq:momentspace})) such that 
\begin{equation}\label{eq:Bayesoptimalcrit}
 p^{-1} \Big[\big\| \bbm_\star - \bbmB \big\|_2^2 +
\big\| \bs_\star - \bsB \big\|_1\Big]
\gotop 0.
\end{equation}
Furthermore,
\begin{equation}\label{eqn:TAP_energy_convergence}
\frac{1}{p}\cF_{\TAP} (\bbm_\star, \bs_\star) \gotop \inf_{\gamma>0}
\phi(\gamma) = {-}\lim_{n, p \to \infty} \frac{1}{p} \log \sP(\y). 
\end{equation}
\end{theorem}

Although $\bbm_\star$ and $\bs_\star$ only describe the first and second moments of the posterior,
they can be used as the basis for a much richer description of the posterior
marginal laws.
In particular,
the next theorem shows that,
for any fixed variable $\beta_j$,
its posterior law is well-approximated by $\sP_{\lambda,\gamma}$ as defined in \eqref{eq:scalarmean},
where $\lambda = \lambda(m_{\star,j},s_{\star,j})$ and $\gamma =
\gamma(m_{\star,j},s_{\star,j})$ are given by the duality relations \eqref{eq:lambdagammastar}.
\begin{theorem}\label{thm:marginalposterior}
Suppose Assumptions \ref{ass:Bayesian_linear_model} and \ref{ass:uniquemin}
hold. Let $(\bbm_\star,\bbs_\star) \in \Gamma^p$ be any local minimizer of
$\cF_\TAP(\bbm,\bbs)$ satisfying (\ref{eq:Bayesoptimalcrit}) as
$n,p \to \infty$. Then for any Lipschitz function $f:\supp(\sP_0) \to \R$
and any $j \in \{1,\ldots,p\}$,
\[\E\Big[\Big(\langle f(\beta_j) \rangle_{\X,\y}-\langle
f(\beta)
\rangle_{\lambda(m_{\star,j},s_{\star,j}),\gamma(m_{\star,j},s_{\star,j})}
\Big)^2 \Big] \to 0.\]
\end{theorem}
\noindent The proofs of Theorems \ref{thm:bayes_consistence}  and
\ref{thm:marginalposterior} are
contained in Appendix \ref{app:proof_bayes_consistence}.

We remark that by symmetry, the expectation in Theorem
\ref{thm:marginalposterior} is the same for all coordinates $j \in
\{1,\ldots,p\}$.
By Markov's inequaility,
Theorem \ref{thm:marginalposterior} then implies
for any (sequence of) non-random indices $j:=j(p) \in \{1,\ldots,p\}$,
\[\langle f(\beta_j) \rangle_{\X,\y}-\langle
f(\beta)
\rangle_{\lambda(m_{\star,j},s_{\star,j}),\gamma(m_{\star,j},s_{\star,j})}
\gotop 0.\]
Here we may choose any Lipschitz function $f$, so this provides a description of
$\beta_j$ at the level of its full posterior marginal law. Specializing to
$f(\beta)=\beta$ and $f(\beta)=\beta^2$ recovers the statements about
its marginal first and second moments in Theorem \ref{thm:bayes_consistence}.

\subsection{Convexity of the TAP free energy}
\label{sec:local-convexity}

%

For large enough $\sigma^2$ or small enough $\delta=\lim_{n,p \to \infty} n/p$,
the following verifies that
the TAP free energy is globally strongly convex. In these settings,
$(\bbm_\star,\bs_\star)$ described by Theorem \ref{thm:bayes_consistence}
must be the global minimizer of $\cF_\TAP$, so in particular
\[-\frac{1}{p}\log \sP(\y)-\inf_{(\bbm,\bs) \in \Gamma^p}
\frac{1}{p}\,\cF_\TAP(\bbm,\bs) \gotop 0.\]
This global convexity is summarized by the following proposition.
\begin{proposition}\label{prop:global_convexity}
Let Assumption \ref{ass:Bayesian_linear_model} hold, where the support of
$\sP_0$ is contained in $[-M,M]$. Then there exist constants $\kappa,c_0>0$
depending only on $M$ such that if $(n/p)/\sigma^2<c_0$, then
$\nabla^2 \cF_{\TAP}(\bbm, \bs) \succeq \kappa \id_{2p}$ for all
$(\bbm,\bs) \in \Gamma$.
\end{proposition}
\noindent
The proof of Proposition \ref{prop:global_convexity} is straightforward and contained in Appendix \ref{app:proof_global_convexity}.

A more difficult optimization scenario is one in which the TAP free energy is
not globally convex and may have multiple local minimizers.
Note that this can occur even in the ``easy'' regime of
Assumption \ref{ass:bayesoptimal} where $\gamma_\stat$ and $\gamma_\alg$
coincide. In this regime, 
the following result establishes that with high probability,
the TAP free energy is strongly convex in
a $O(\sqrt{p})$-radius local neighborhood of the Bayes-optimal local minimizer
described by Theorem \ref{thm:bayes_consistence}.
This implies, in particular, that this Bayes-optimal local minimizer is unique.

We denote by
\[\ball((\bbmB,\bsB),r)=\big\{(\bbm,\bs) \in
\Gamma^p:\|(\bbm,\bs)-(\bbmB,\bsB)\|_2<r\big\}\]
the subset of the radius-$r$ ball around $(\bbmB,\bsB)$ that belongs to the
parameter space $\Gamma^p$, and by
$\lambda_{\min}(\cdot)$ the smallest eigenvalue of a symmetric matrix.

\begin{theorem}\label{thm:local_convexity}
Let Assumptions \ref{ass:Bayesian_linear_model}, \ref{ass:uniquemin}, and \ref{ass:bayesoptimal} hold. 
Then there exist constants $\eps, \kappa > 0$ such that with probability
approaching $1$ as $n, p \to \infty$, $\cF_\TAP$ has a unique local minimizer $(\bbm_\star,\bs_\star)$ in $\ball((\bbmB, \bsB), \eps\sqrt{p} )$, and 
\[
\inf_{(\bbm, \bs) \in \ball((\bbmB, \bsB), \eps\sqrt{p} ) } \lambda_{\min}\Big(
\nabla^2 \cF_{\TAP}(\bbm, \bs) \Big) \ge \kappa.
\]
\end{theorem}
We prove Theorem \ref{thm:local_convexity} in Appendix \ref{app:proof_local_convexity}.
We remark that 
it may be of interest to study the local landscape of the TAP free
energy near the Bayes posterior marginals $(\bbmB,\bsB)$ even in the hard regime
$\gamma_\alg \neq \gamma_\stat$, and that this local convexity property
may continue to hold true.
As it is conjectured that no polynomial-time algorithm based on the data
$(\X,\y)$ alone can reach
this local neighborhood of $(\bbmB,\bsB)$, in Section \ref{sec:hard-regime}
we will study instead the local convexity of the TAP free energy around a
different local minimizer that is reachable by an AMP algorithm.

\subsection{Convergence of natural gradient descent}
\label{sec:NGD}

In this section, we discuss an algorithm which, for fixed (and large) $n,p$, converges linearly to the local minimizer $(\bbm_\star,\bs_\star)$ described in Theorem \ref{thm:local_convexity}.

This algorithm has two stages, the first being an AMP algorithm that
successfully navigates the global landscape of the TAP free energy, and the
second being a variant of gradient descent that optimizes the TAP free energy
within a locally convex neighborhood of a minimizer. We will refer to this second stage
as ``natural gradient descent'' (NGD), and it is given by the iteration
\begin{equation}\label{eqn:NGD}
\begin{aligned}
	\Big(\blambda_{\NGD}^k, -\frac12 \bgamma_{\NGD}^k \Big) =&~ \Big(\blambda_{\NGD}^{k-1}, - \frac12 \bgamma_{\NGD}^{k-1} \Big) - \eta \nabla \cF_{\TAP}(\bbm_{\NGD}^k, \bs_{\NGD}^k), \\
	(\bbm_{\NGD}^{k+1}, \bs_{\NGD}^{k+1}) 
		=&~ 
		\Big((\<\beta\>_{\lambda_{\NGD,j}^k,\gamma_{\NGD,j}^k})_{j=1}^p, (\<\beta^2\>_{\lambda_{\NGD,j}^k,\gamma_{\NGD,j}^k})_{j=1}^p\Big), 
\end{aligned}
\end{equation}
where $\eta>0$ is a step size parameter. A related
NGD method was used in \cite{celentano2023local} to optimize the TAP free energy
for a different model of $\Z_2$-synchronization.

This iteration can be viewed a preconditioned form of
gradient descent \cite{Amari1998} on $(\blambda,\bgamma)$, noting that
\begin{equation}
\begin{gathered}
	\nabla_{\bbm,\bs} \cF_{\TAP}(\bbm, \bs)
		=
		\bI(\blambda,\bgamma)^{-1}\nabla_{\blambda,\bgamma} \cF_{\TAP}\Big((\<\beta\>_{\lambda_j,\gamma_j})_{j=1}^p, (\<\beta^2\>_{\lambda_j,\gamma_j})_{j=1}^p\Big)
\end{gathered}
\end{equation}
where the Jacobian for the change-of-variables $(\blambda,\bgamma) \mapsto
(\bbm,\bs)$ takes the form
\[\bI(\blambda,\bgamma)^{-1}=\diag\left(\begin{pmatrix}
	\Var_{\lambda_j,\gamma_j}[\beta] &
\Cov_{\lambda_j,\gamma_j}[\beta,\beta^2]
	\\
	\Cov_{\lambda_j,\gamma_j}[\beta,\beta^2] &
\Var_{\lambda_j,\gamma_j}[\beta^2]
\end{pmatrix}\right)_{j=1}^p.\]
Here, $\Var_{\lambda,\gamma}[\,\cdot\,]$ and $\Cov_{\lambda,\gamma}[\,\cdot\,]$
are the variances and covariances under $\sP_{\lambda,\gamma}$.
Alternatively, as we will show in Appendix \ref{app:proof_NGD_convergence},
this iteration can be viewed as a mirror descent or Bregman gradient method
\cite{Blair1985,BECK2003167}
for the Bregman divergence associated to the relative entropy function
$D_0(\bbm,\bs)$.
We use NGD in place of ordinary gradient descent to adapt to the divergence of
the gradient $\nabla_{\bbm,\bs} \cF_{\TAP}$ and Hessian $\nabla_{\bbm,\bs}^2
\cF_\TAP$ at the boundaries of the parameter domain $\Gamma^p$.

Using the local convexity established in Theorem \ref{thm:local_convexity},
we can show that NGD with appropriate initialization converges linearly to $(\bbm_\star,\bs_\star)$.
\begin{theorem}\label{thm:NGD_convergence}
	Let Assumptions \ref{ass:Bayesian_linear_model}, \ref{ass:uniquemin},
and \ref{ass:bayesoptimal} hold. There exist constants $C,c,\mu,\eta_0,\eps >0 $ depending only $(\delta,\sigma^2,\sP_0)$ such that, with probability going to 1 as $n,p \rightarrow \infty$,
	if NGD is initialized at $(\bbm_{\NGD}^0,\bs_{\NGD}^0)$ satisfying
	\begin{equation}
		\frac1p \Big[ \| \bbm_{\NGD}^0 - \bbm_\star \|_2^2 +  \| \bs_{\NGD}^0 - \bs_\star \|_2^2\Big] \le \eps^2,
		\qquad
		\frac1p\cF_\TAP(\bbm_{\NGD}^0,\bs_{\NGD}^0) - \frac1p\cF_\TAP(\bbm_\star,\bs_\star)
		\leq 
		c\eps^2,
	\end{equation}
	then for all $k \geq 1$
	\begin{equation}
		\frac1p \Big[ \| \bbm_{\NGD}^k - \bbm_\star \|_2^2 +  \| \bs_{\NGD}^k - \bs_\star \|_2^2\Big]
		\le
		C(1-\mu\eta)^k.
	\end{equation}
\end{theorem}
An initialization satisfying the conditions of Theorem \ref{thm:NGD_convergence}
may be obtained by running a constant number of iterations of AMP:
Fixing some large $T_0 \geq 1$ not depending on $n,p$,
let $(\bbm_{\AMP}^{T_0},\bs_{\AMP}^{T_0})$ be the
$T_0^\text{th}$ iterates of the AMP iterations (\ref{eq:AMPalg}).
For $k > T_0$, let $(\bbm_{\NGD}^k,\bs_{\NGD}^k)$ be the NGD iterates
\eqref{eqn:NGD} initialized at $(\bbm_{\NGD}^0,\bs_{\NGD}^0) = (\bbm_{\AMP}^{T_0},\bs_{\AMP}^{T_0})$.
\begin{corollary}\label{cor:AMP+NGD-convergence}
	Under the conditions of Theorem \ref{thm:NGD_convergence},
	there exists constants $\overline{T}, C, \kappa > 0$ depending on $(\sigma^2, \delta, \sP_0)$ such that, for any fixed $T_0 \geq \overline{T}$,
	with probability going to $1$ as $n, p \to \infty$,
	the following occurs:
	Taking $(\bbm_{\NGD}^0,\bs_{\NGD}^0) = (\bbm_{\AMP}^{T_0},\bs_{\AMP}^{T_0})$,
	for all $k \geq 0$
	\begin{equation}\label{eqn:NGD_convergence-statement}
		\frac1p \Big[ \| \bbm_{\NGD}^k - \bbm_\star \|_2^2 +  \| \bs_{\NGD}^k - \bs_\star \|_2^2\Big] \le C e^{- \kappa k},
	\end{equation}
	for $(\bbm_{\NGD}^k,\bs_{\NGD}^k)$ the iterates of NGD initialized at the $T_0^\text{th}$ iterate of AMP.
\end{corollary}
The proof of Corollary \ref{cor:AMP+NGD-convergence} is contained in Appendix \ref{app:AMP+NGD-convergence}. We remark it is believed that the iterates of
AMP alone (without switching to NGD) may also satisfy
\eqref{eqn:NGD_convergence-statement}, but establishing such convergence remains
an open problem.

\subsection{Calibrated inference in the hard regime}
\label{sec:hard-regime}

Theorems \ref{thm:local_convexity} and \ref{thm:NGD_convergence} are stated for
the ``easy'' regime described by Assumption \ref{ass:bayesoptimal},
in which efficient algorithms can achieve the asymptotically Bayes-optimal
squared-error risk. In this section we show that,
even outside of this easy regime,
AMP+NGD converges to a local minimizer contained in a region of local strong convexity.

\begin{theorem}\label{thm:local_convexity_AMP_fixed_point}
	Let Assumption \ref{ass:Bayesian_linear_model} hold, and assume that $\phi''(\gamma_{\alg}) > 0$. 
	Then there exist constants $\overline{T},C,\kappa,\eps > 0$ depending on
$(\sigma^2, \delta, \sP_0)$ such that  for any fixed $T_0 \ge \overline{T}$, with probability approaching 1 as $n,p \to
\infty$, the following holds: 
	\begin{enumerate}[label={(\alph*)},ref={(\alph*)}]

		\item \label{item:local-convexity}
		If $(\bbm_{\AMP}^{T_0},\bs_{\AMP}^{T_0})$ denotes the $T_0^\text{th}$  iterate of AMP, as in Corollary \ref{cor:AMP+NGD-convergence},
		then $\cF_\TAP$ has a unique local minimizer $(\bbm_\star,\bs_\star)$ in $\ball((\bbm_{\AMP}^{T_0}, \bs_{\AMP}^{T_0}), \eps\sqrt{p} )$ with $p^{-1}\big[\| \bbm_{\AMP}^{T_0} - \bbm_\star \|_2^2 + \| \bs_{\AMP}^{T_0} - \bs_\star \|_2^2\big] \leq Ce^{-\kappa T_0}$ and
		\begin{equation}
			\inf_{(\bbm, \bs) \in \ball((\bbm_\star, \bs_\star), \eps\sqrt{p} ) } \lambda_{\min}\Big( \nabla^2 \cF_{\TAP}(\bbm, \bs) \Big) \ge \kappa.
		\end{equation}

		\item \label{item:NGD-convergence}
		 Consider the AMP+NGD algorithm, with the initialization of NGD
as $(\bbm_{\NGD}^0,\bs_{\NGD}^0) = (\bbm_{\AMP}^{T_0},\bs_{\AMP}^{T_0})$.
		Then for all $k \geq 0$, $p^{-1}\big[\| \bbm_{\NGD}^{k} - \bbm_\star \|_2^2 + \| \bs_{\NGD}^{k} - \bs_\star \|_2^2\big] \leq Ce^{-\kappa k}$.

	\end{enumerate}

\end{theorem}
\noindent We prove part \ref{item:local-convexity} of Theorem \ref{thm:local_convexity_AMP_fixed_point} in Appendix \ref{app:proof_local_convexity}.
We prove part \ref{item:NGD-convergence} in Appendix \ref{app:AMP+NGD-convergence}.

In the hard regime $\gamma_\alg \neq \gamma_\stat$, the local minimizer
$(\bbm_\star,\bs_\star)$ of Theorem \ref{thm:local_convexity_AMP_fixed_point}
does not correspond to the marginal first and second moments of the true Bayes
posterior law, as was the case in Theorem \ref{thm:bayes_consistence}.
Nevertheless, Theorem \ref{thm:local_convexity_AMP_fixed_point} shows that
AMP+NGD exhibits linear convergence to some local minimizer of the TAP free
energy.
We show in the next theorem that this local minimizer achieves the squared-error
Bayes risk that is conjecturally optimal among polynomial-time algorithms,
and furthermore that it can serve as the basis for correctly calibrated posterior inference.
\begin{theorem}
\label{thm:error-and-calibration}
	Under the conditions of Theorem \ref{thm:local_convexity_AMP_fixed_point},
	the local minimizer $(\bbm_\star,\bs_\star)$ satisfies
	\begin{equation}\label{eq:localminerror}
		\frac1p \| \bbm_\star - \bbeta_0 \|_2^2  \gotop \mmse(\gamma_\alg).
	\end{equation}
	Furthermore, for any non-empty open set $A \subset \reals$, Lipschitz and bounded function $f:\reals \rightarrow \reals$, and (sequence of) non-random indices $j:=j(p) \in \{1,\ldots,p\}$,
	\begin{equation}
		\E[f(\beta_{0,j}) \mid \lambda(m_{\star,j},s_{\star,j}) \in A]
		\rightarrow
		\E_{(\beta_0,z) \sim \sP_0 \times \normal(0,1)}\big[f(\beta_0)
\mid \gamma_{\alg}  \beta_0 + \sqrt{\gamma_{\alg}} \, z \in A\big],
	\end{equation}
and $\gamma(m_{\star,j},s_{\star,j}) \gotop \gamma_{\alg}$.
\end{theorem}
\noindent We prove Theorem \ref{thm:error-and-calibration} in Section \ref{app:calibration}.

This theorem shows that, for the purposes of posterior inference, one may view
the parameters $\lambda_{\star,j}:=\lambda(m_{\star,j},s_{\star,j})$ computed
from the TAP local minimizer $(\bbm_\star,\bs_\star)$ 
as observations from a scalar sequence model
$\lambda_{\star,j}=\gamma_\alg \beta_{0,j}+\sqrt{\gamma_\alg}\,z_j$. 
Here $\gamma_\alg$ is a deterministic parameter that may be computed from the
replica-symmetric potential (\ref{eqn:phi_potential}), or alternatively,
approximated by $\gamma_{\star,j}:=\gamma(m_{\star,j},s_{\star,j})$.
Theorem \ref{thm:error-and-calibration} states that posterior
inference for $\beta_{0,j}$ based upon $\lambda_{\star,j}$ and this sequence
model will be correctly calibrated in the asymptotic limit as $n,p \to \infty$,
even though the posterior laws $\sP_{\lambda_{\star,j},\gamma_{\star,j}}$ in
this sequence model are inaccurate approximations for the true
posterior marginals.


It is worth contrasting the behavior described in Theorem
\ref{thm:error-and-calibration} with the na\"ive mean field approximation. We
note that under the asymptotic setting of Assumption
\ref{ass:Bayesian_linear_model},
$\blambda_\star$ computed analogously from $(\bbm_\star,\bs_\star)$
that optimizes the naive mean-field free energy (\ref{eq:naiveMF}) will, in general, not be
well-approximated by a Gaussian sequence model around $\bbeta_0$,
and posterior inferences based on the approximate laws
$\sP_{\lambda_{\star,j},\gamma_{\star,j}}$ will not have a calibration
guarantee.
Indeed, as we saw in Section \ref{sec:gauss-prior} for even a simple Gaussian
prior, na\"ive mean field can lead to overconfident estimates of posterior
variance, and thus inflated error rates in statistical applications. We
provide further examinations of posterior calibration in simulation for other
priors in Section \ref{sec:exp_marginal}.

\section{Numerical simulations}\label{sec:simulations}

%
%

\subsection{Mean squared errors of MF and TAP estimators} \label{sec:exp_mse}

We perform numerical simulations to compare the mean squared errors (MSE)
$p^{-1}\|
\bbm_\star - \bbeta_0 \|_2^2$ achieved by the minimizers
$(\bbm_\star,\bs_\star)$ of the mean field free
energy $\cF_{\rm MF}$ in \eqref{eq:naiveMF} and the TAP free energy $\cF_{\TAP}$
in \eqref{eqn:TAP}. 
We consider two particular choices of the prior $\sP_0$:
\begin{enumerate}
\item A three-point distribution $\sP_0=(1/3) \cdot\delta_{1}+ (1/3) \cdot\delta_{-1}+(1/3)\cdot\delta_{0}$. 
\item A  Bernoulli-Gaussian distribution $\sP_0=(1/2)\cdot\normal(0, 1)+(1/2)\cdot\delta_0$. 
\end{enumerate}

\begin{figure}[htp]
\includegraphics[width=0.45\linewidth]{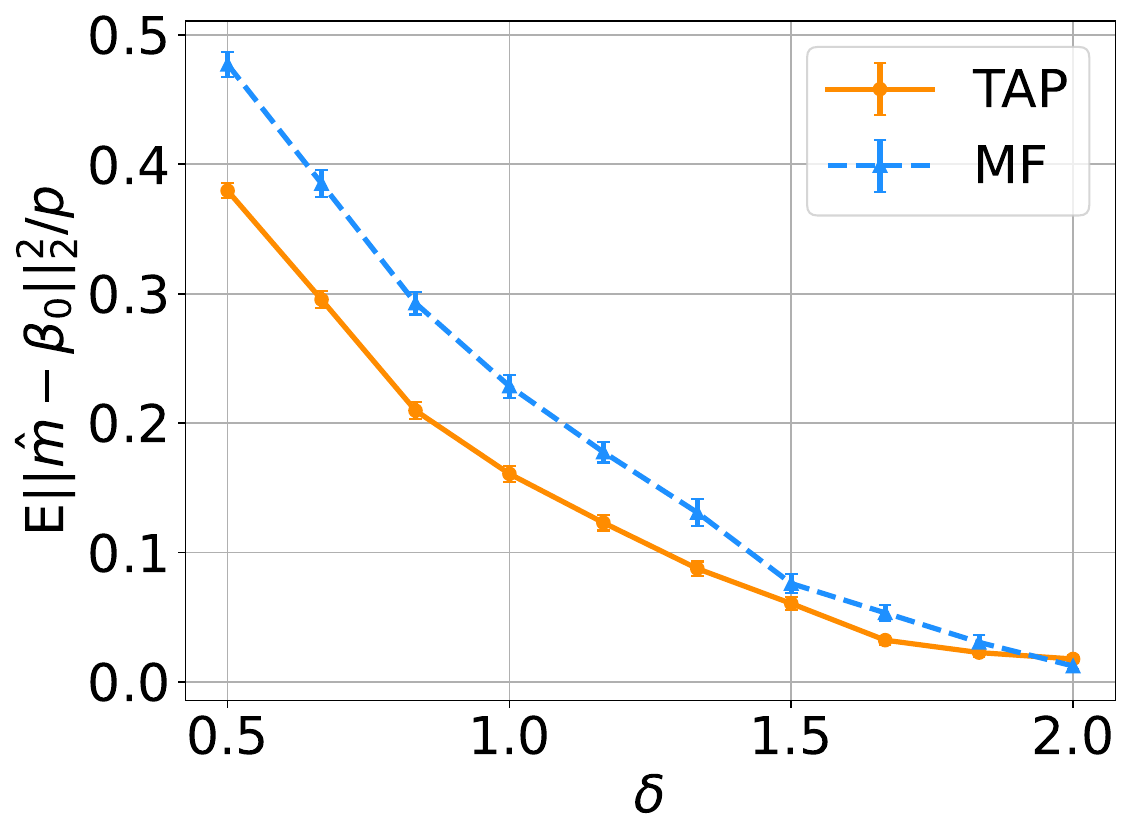}
\includegraphics[width=0.45\linewidth]{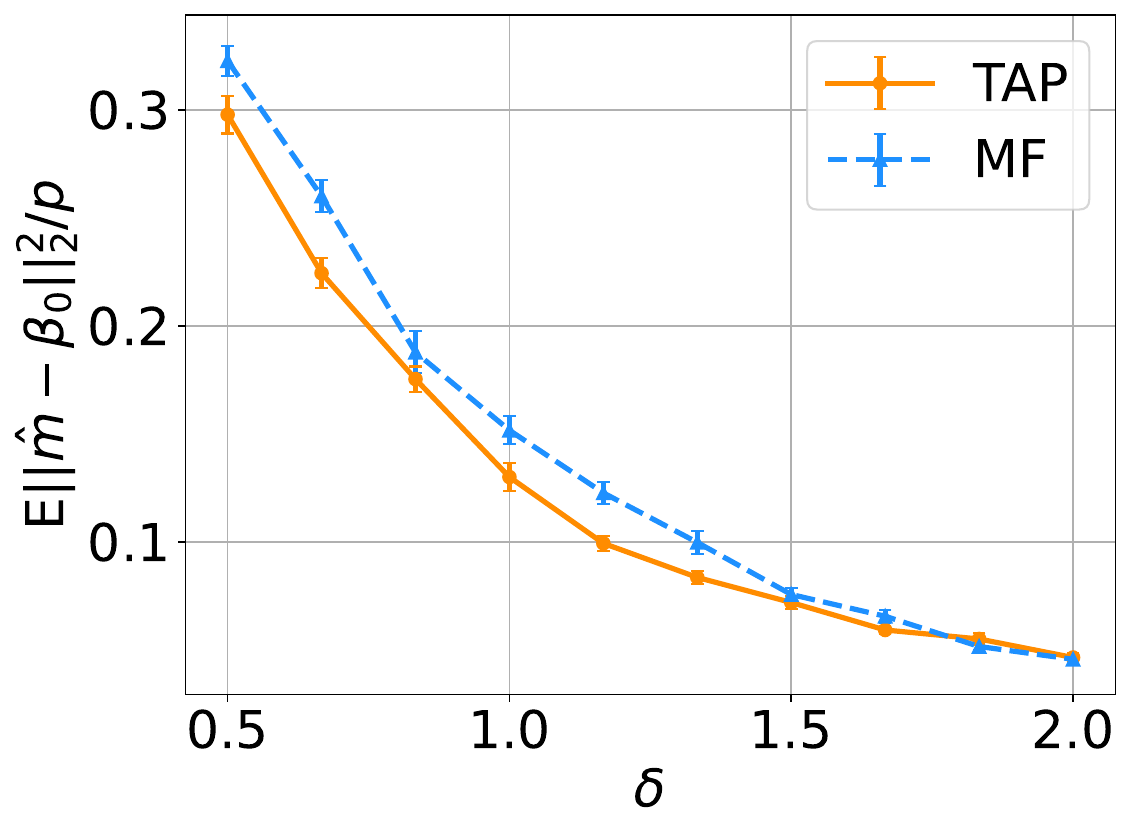}
\caption{MSE of TAP and MF posterior-mean estimators for
$\bbeta_0$ under different values of
$\delta = n / p$. Parameters: $\sigma=0.3$, $n=300$, $p=\lfloor
n/\delta\rfloor$. Left: Three-point prior $\sP_0=(1/3)\cdot
\delta_{-1}+(1/3)\cdot \delta_{0}+(1/3)\cdot \delta_{1}$. Right:
Bernoulli-Gaussian prior $\sP_0=(1/2)\cdot \normal(0,1)+(1/2)\cdot\delta_{0}$. Error
bars denote the standard deviation of the MSE estimates across 20 independent
simulations.}\label{fig:mmse_1}
\end{figure}

We (approximately) minimize both free energies using the natural gradient
descent (NGD) algorithm of \eqref{eqn:NGD}. We observe that NGD typically
converged within $20000$ iterations (in the sense of achieving a small
gradient),
despite the lack of a theoretical convergence guarantee in certain settings.
In Figure~\ref{fig:mmse_1}, we report the MSE of both the MF and TAP estimators
for the two choices of prior, fixed $\sigma^2$, and varying dimension ratios
$\delta=n/p$. As anticipated by our theory, the MSE of the TAP estimator is less
than that of the MF estimator across all settings of $\delta$ and both choices
of prior. We remark that we have chosen to illustrate a setting of $(\sigma^2, \delta)$ in which there is a more significant difference between the MSE of TAP and MF. 

\subsection{Calibration of posterior marginals}\label{sec:exp_marginal}


Let $\lambda_{\star,j}=\lambda(m_{\star,j},s_{\star,j})$ and
$\gamma_{\star,j}=\gamma(m_{\star,j},s_{\star,j})$ be defined by the duality
relations \eqref{eq:lambdagammastar}, where $(\bbm_{\star}, \bs_{\star})$ is the
(approximate) minimizer of either the TAP free energy \eqref{eqn:TAP}
or the naive mean field free energy \eqref{eq:naiveMF}.
Theorem \ref{thm:marginalposterior} implies for the TAP estimates that posterior
inferences based on the approximate posterior laws
$\beta_{0,j} \sim \sP_{\lambda_{\star,j},\gamma_{\star,j}}$
are asymptotically well-calibrated.

\begin{figure}[htp]
\includegraphics[width=0.33\linewidth]{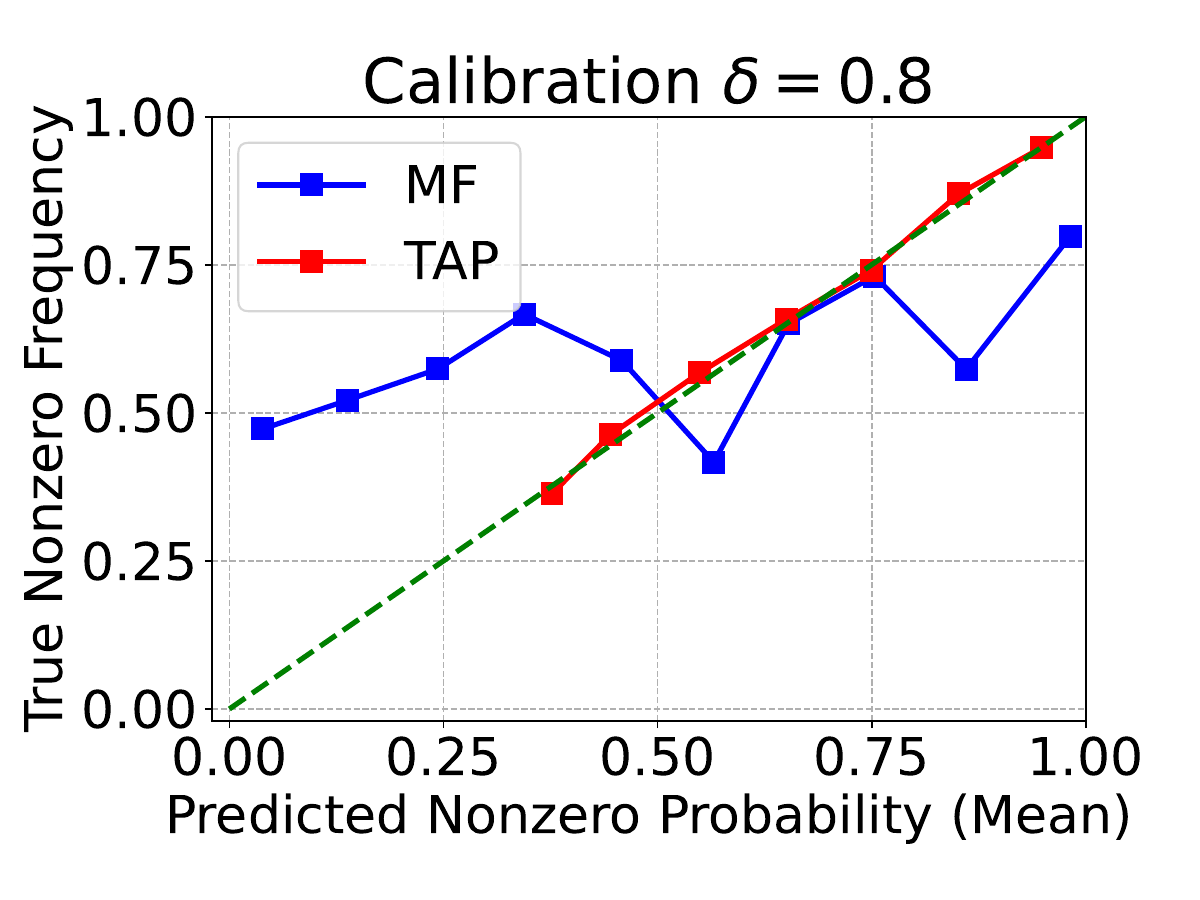}
\includegraphics[width=0.33\linewidth]{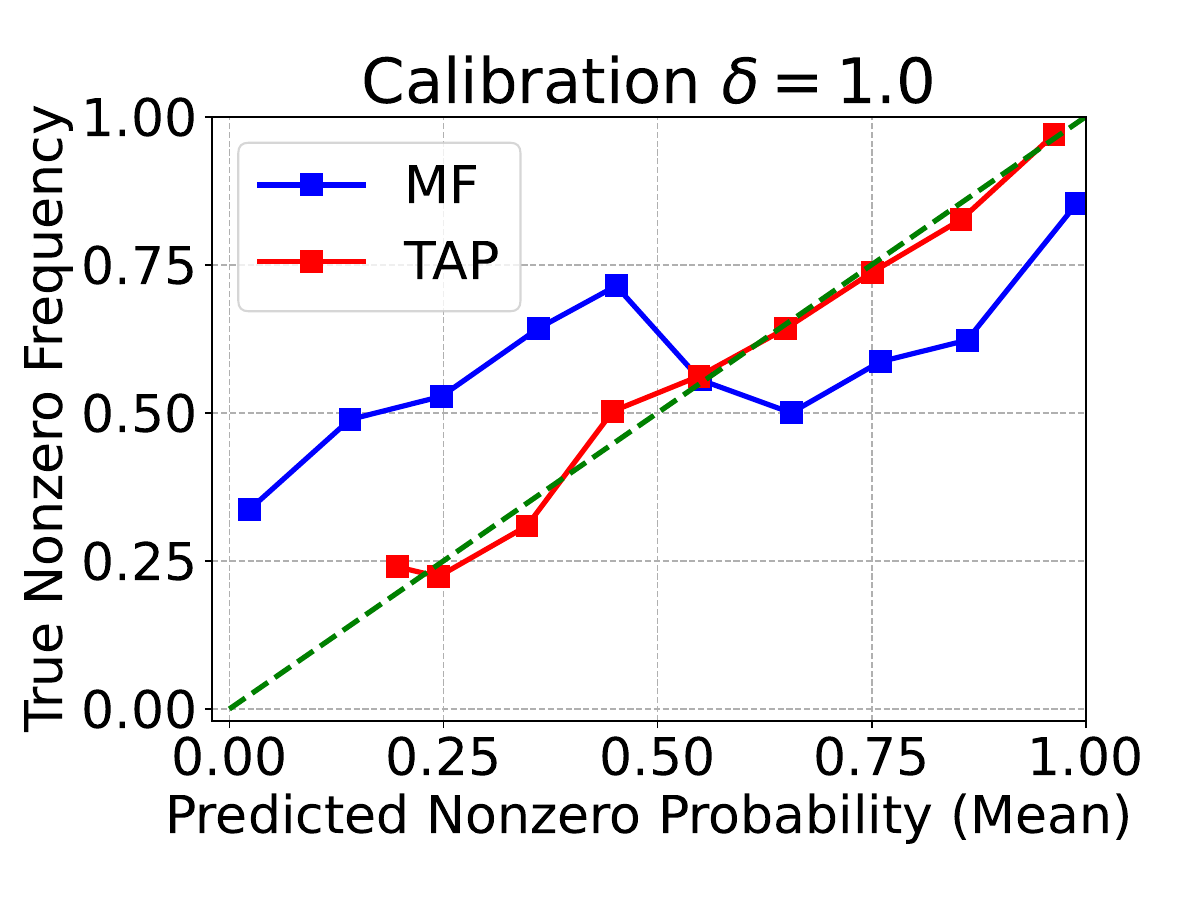}
\includegraphics[width=0.33\linewidth]{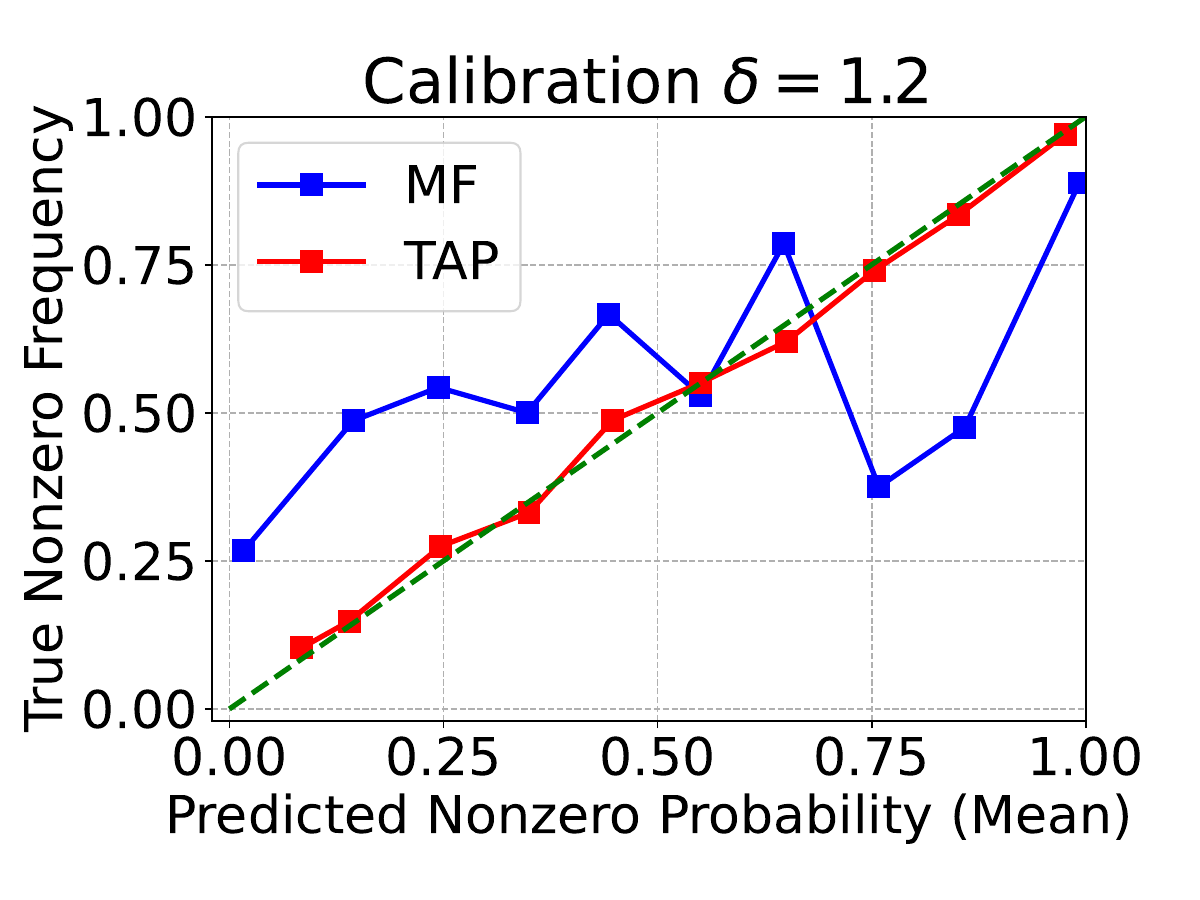}
\caption{Calibration of the true non-zero frequencies versus the estimated
Posterior Inclusion Probabilities (PIPs).
Coordinates $\beta_{0,j}$ across 10 independent simulations are binned by their
estimated PIPs into 10 bins. The true non-zero frequency for each bin with
center $p$ is given by $\frac{\#\{j : \beta_{0,j} \neq 0,
\<\ones_{\beta_{0,j}\neq0}\>_{\lambda_{\star,j},\gamma_{\star,j}} \in [p-0.05,
p+0.05] \}}{\#\{j :
\<\ones_{\beta_{0,j}\neq0}\>_{\lambda_{\star,j},\gamma_{\star,j}} \in [p-0.05,
p+0.05] \}}$. 
Parameters: $\sigma=0.3$, $n=500$, and three-point prior $\sP_0=(1/3) \cdot \delta_{-1}+(1/3) \cdot \delta_{0}+(1/3) \cdot \delta_{1}$.
}\label{fig:cali_1}
\end{figure}

Here, for both priors and both methods, we estimate from these variational
approximations $\sP_{\lambda_{\star,j},\gamma_{\star,j}}$ the Posterior
Inclusion Probabilities (PIPs)
\[\<\ones_{\beta_{0,j}\neq 0}\>_{\lambda_{\star,j},\gamma_{\star,j}},\]
and we assess the calibration of these estimates in simulation.
Figures~\ref{fig:cali_1} and \ref{fig:uni_1} show calibration plots, where the x-axis bins
coordinates $j \in \{1,\ldots,p\}$ by their estimated PIP
$\<\ones_{\beta_{0,j}\neq 0}\>_{\lambda_{\star,j},\gamma_{\star,j}}$
into ten bins $(0,0.1),(0.1,0.2),\ldots,(0.9,1)$, and the y-axis plots the true
fraction of coefficients $\beta_{0,j}$ that are non-zero within each bin.
We see that the estimated PIPs from the TAP approach are well-calibrated, in the
sense that the true fraction of non-zero coefficients is close to the estimated
PIP value within each bin. In contrast, the PIP estimates from the
naive mean field approximation exhibit varying degrees of miscalibration.

\subsection{Universality of TAP free energy}\label{sec:exp_univ}

\begin{figure}[htp]
\includegraphics[width=0.435\linewidth]{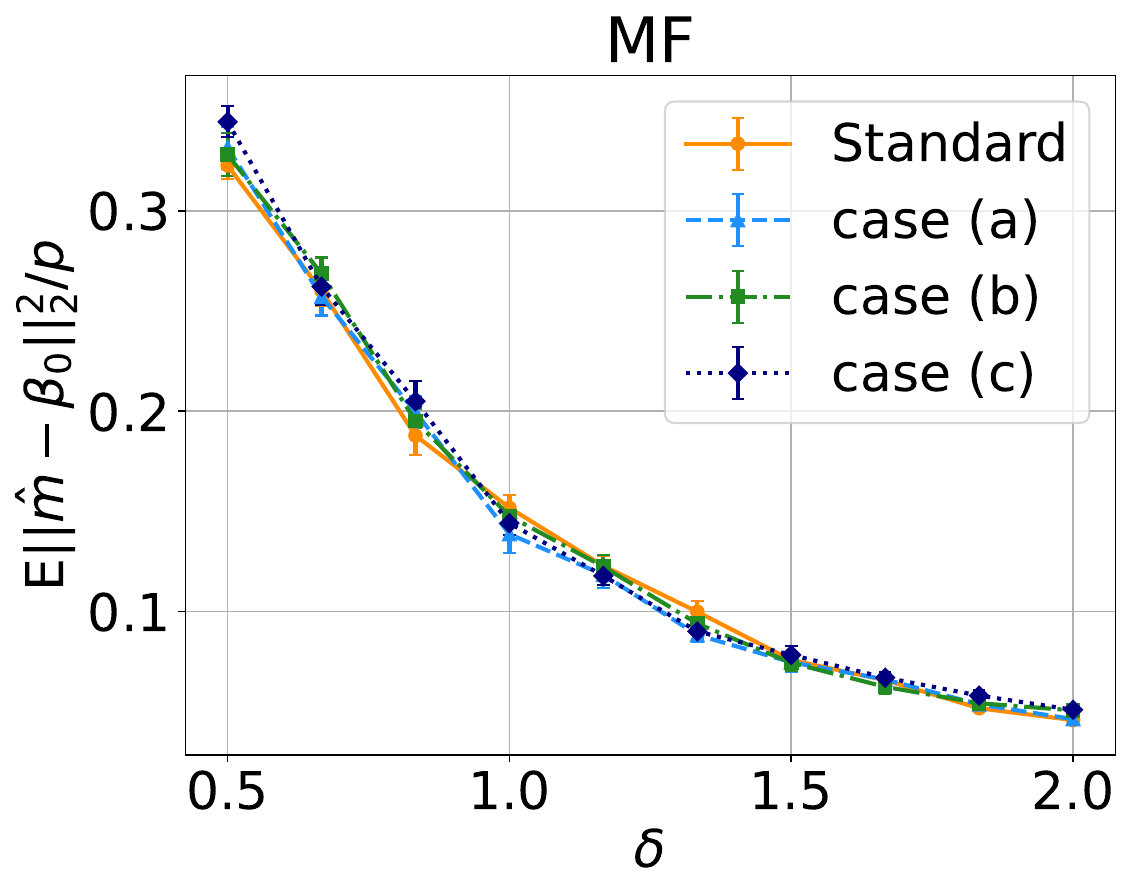}
\includegraphics[width=0.45\linewidth]{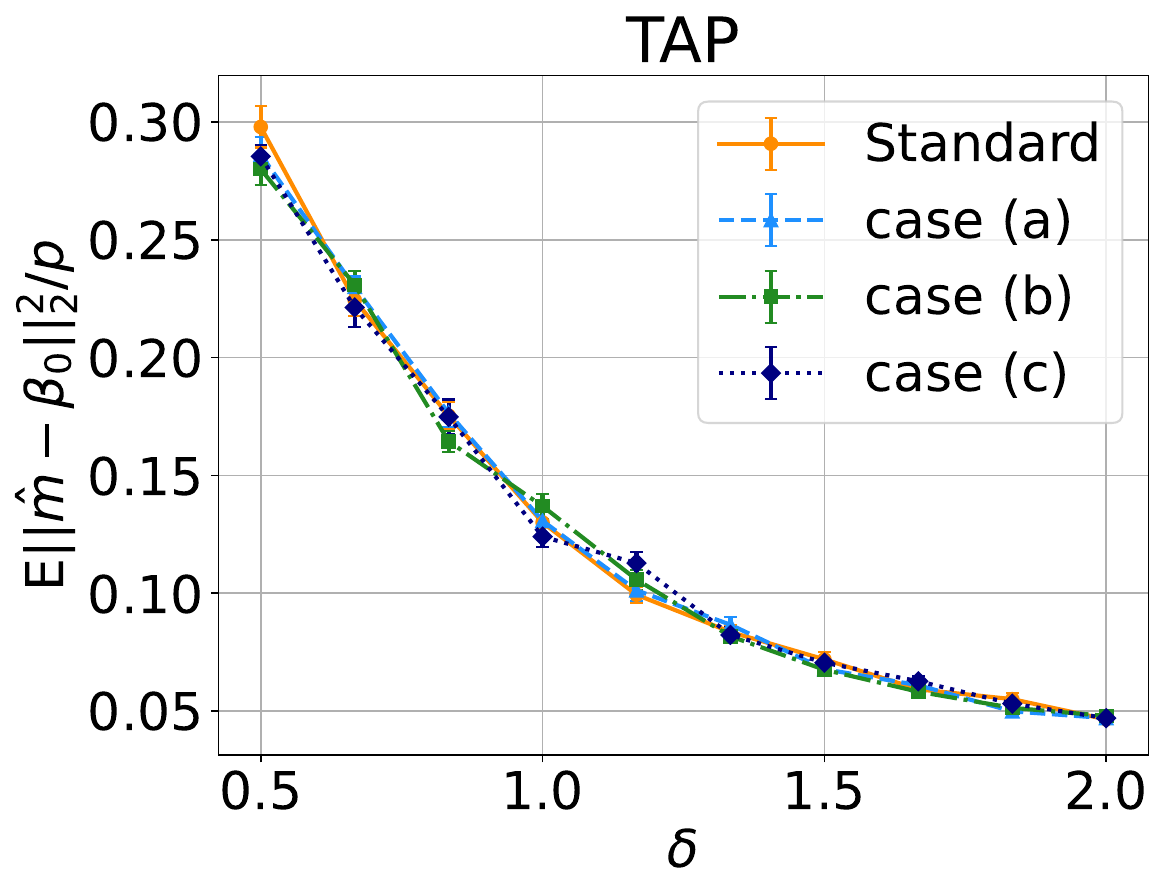}

\caption{Comparison of MSE of the MF and TAP posterior-mean estimators for
$\bbeta_0$ under three types of misspecifications. Parameters: $n=300$,
$\sigma=0.3$, $\sP_0=(1/2) \cdot \delta_{0}+ (1/2) \cdot \normal(0,1)$. Standard:
Gaussian $\X$ and $\beps$. Case (a): Radamacher $\X$. Case (b): Radamacher
$\beps$. Case (c): Bernoulli $\X$ with heterogeneous sparsity across columns.}
\label{fig:uni_05}
\end{figure}

\begin{figure}[htp]
\includegraphics[width=0.34\linewidth]{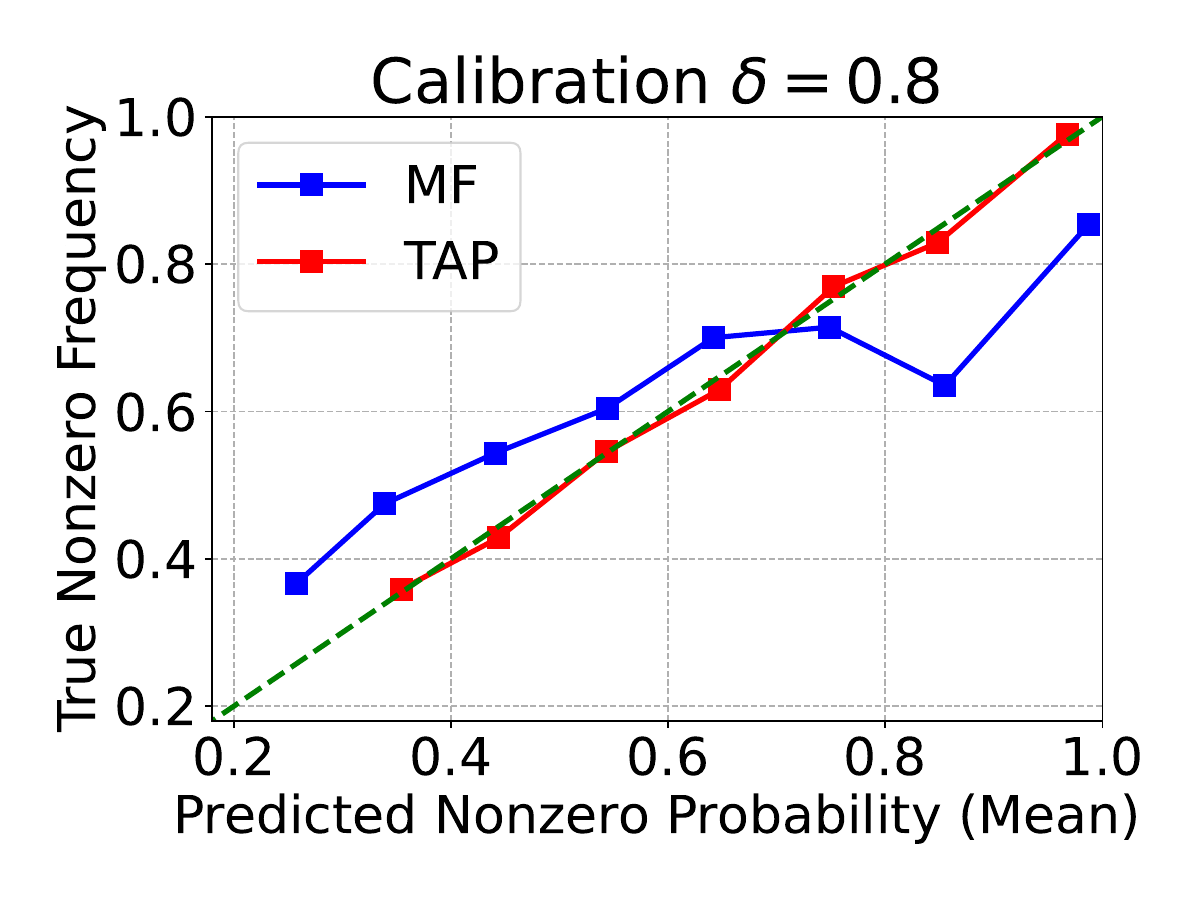}
\includegraphics[width=0.34\linewidth]{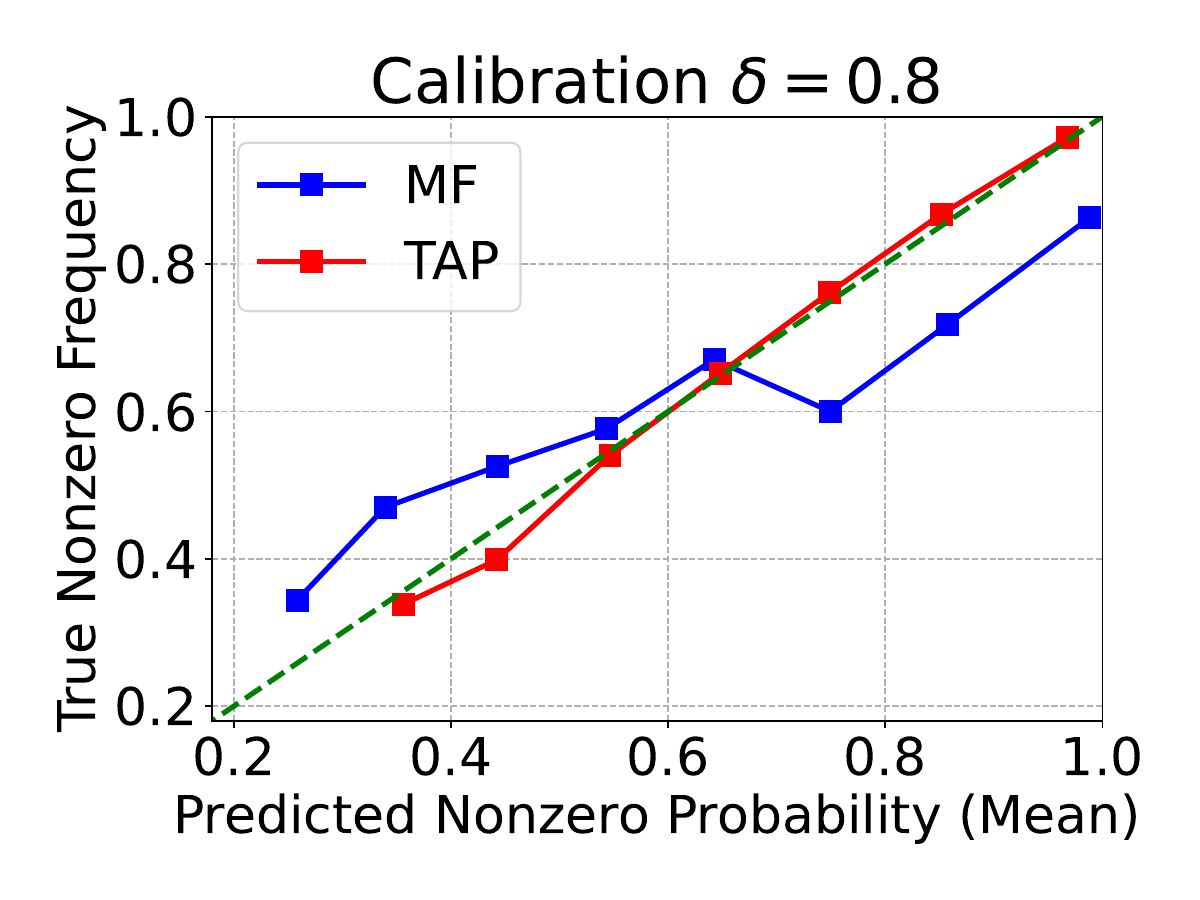}
\includegraphics[width=0.34\linewidth]{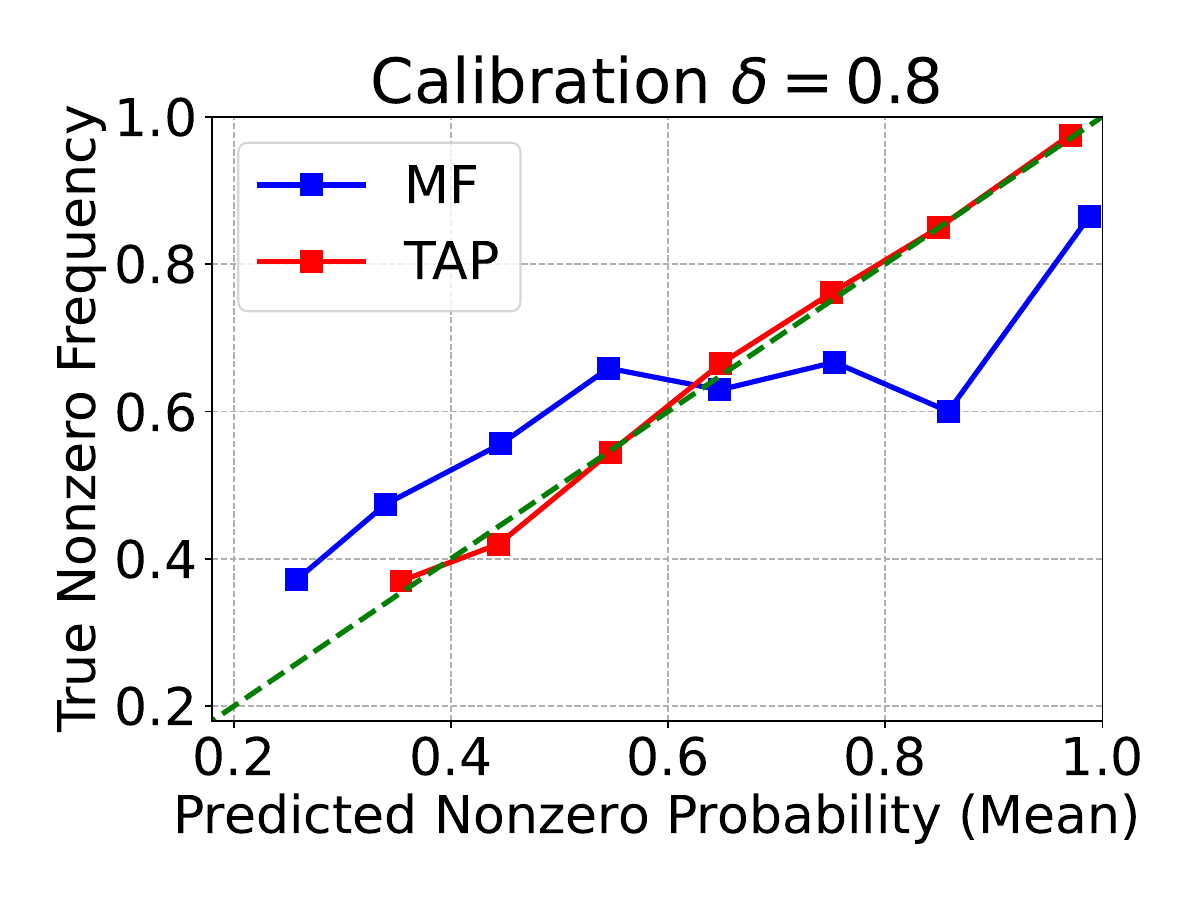}
\includegraphics[width=0.34\linewidth]{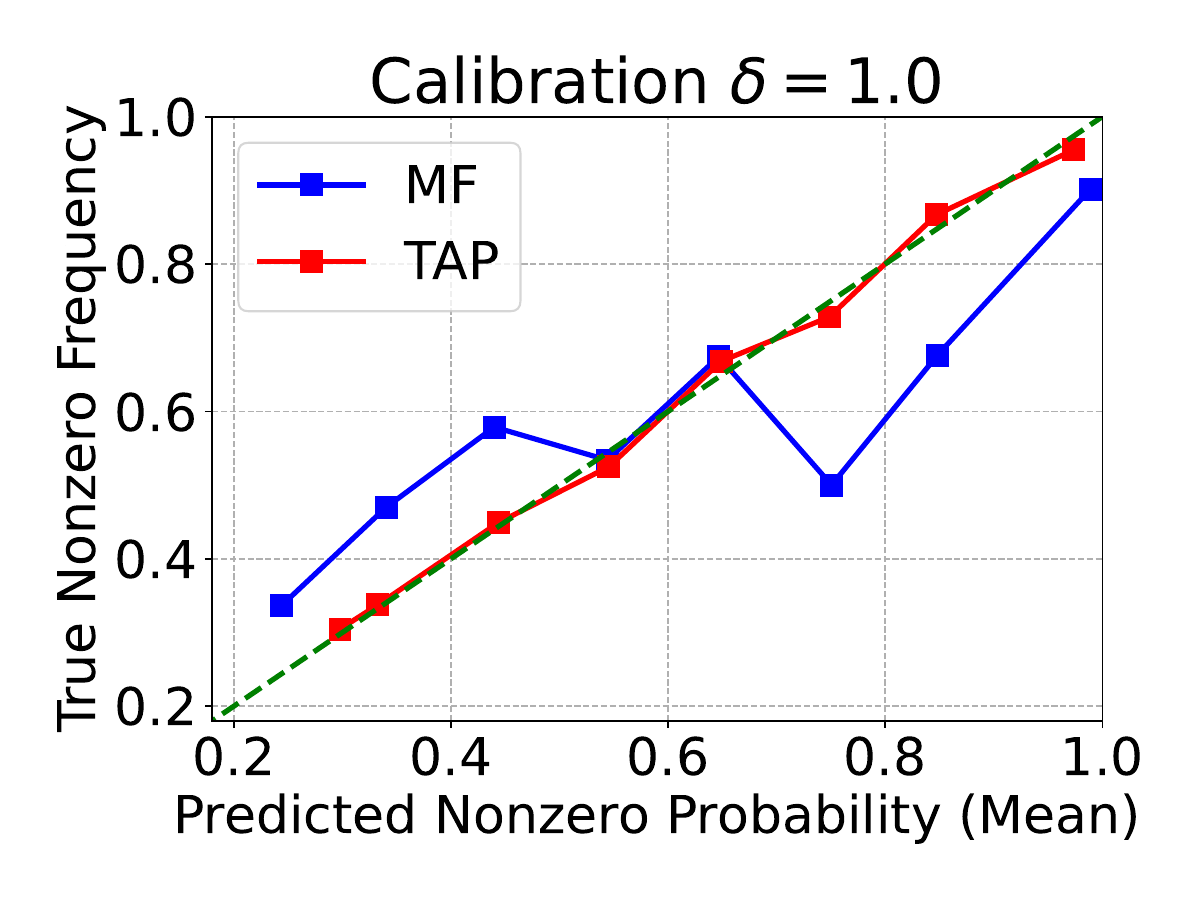}
\includegraphics[width=0.34\linewidth]{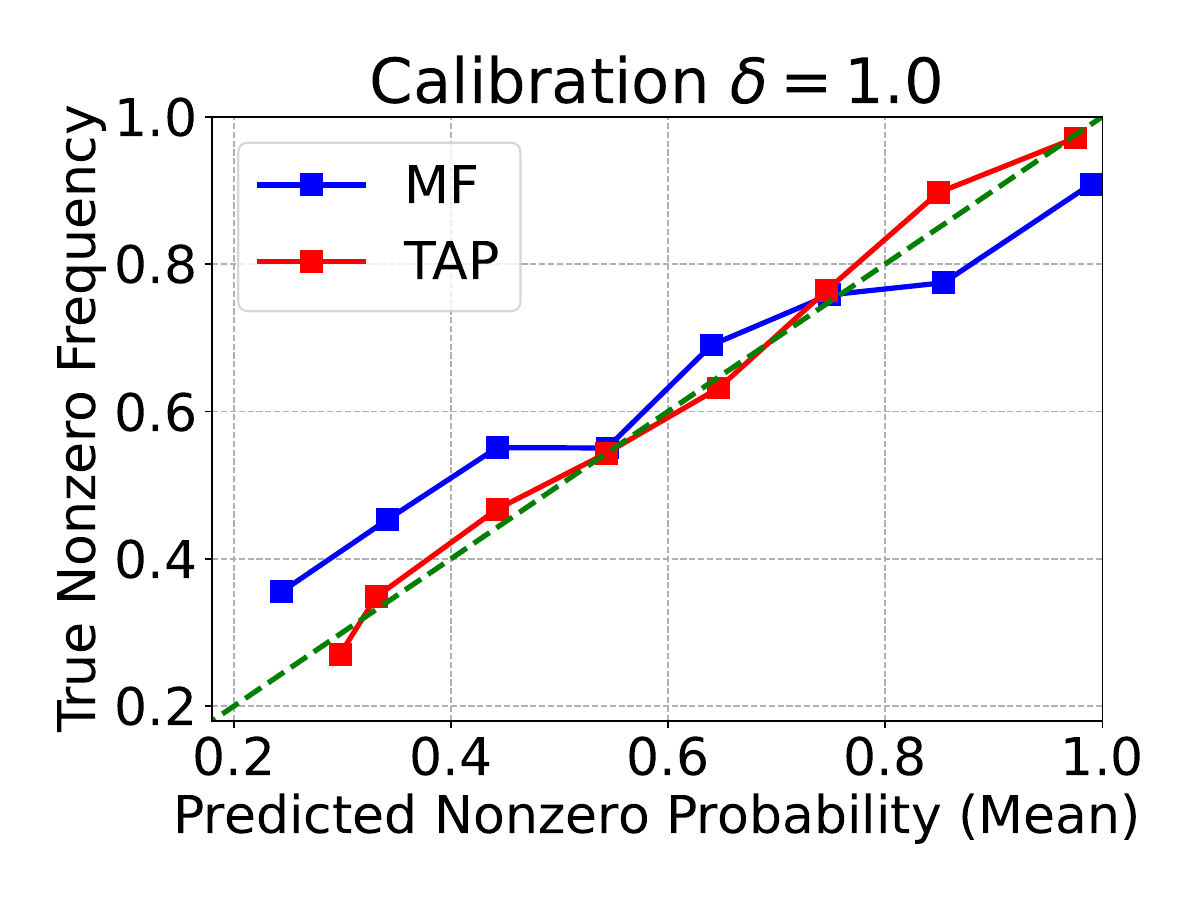}
\includegraphics[width=0.34\linewidth]{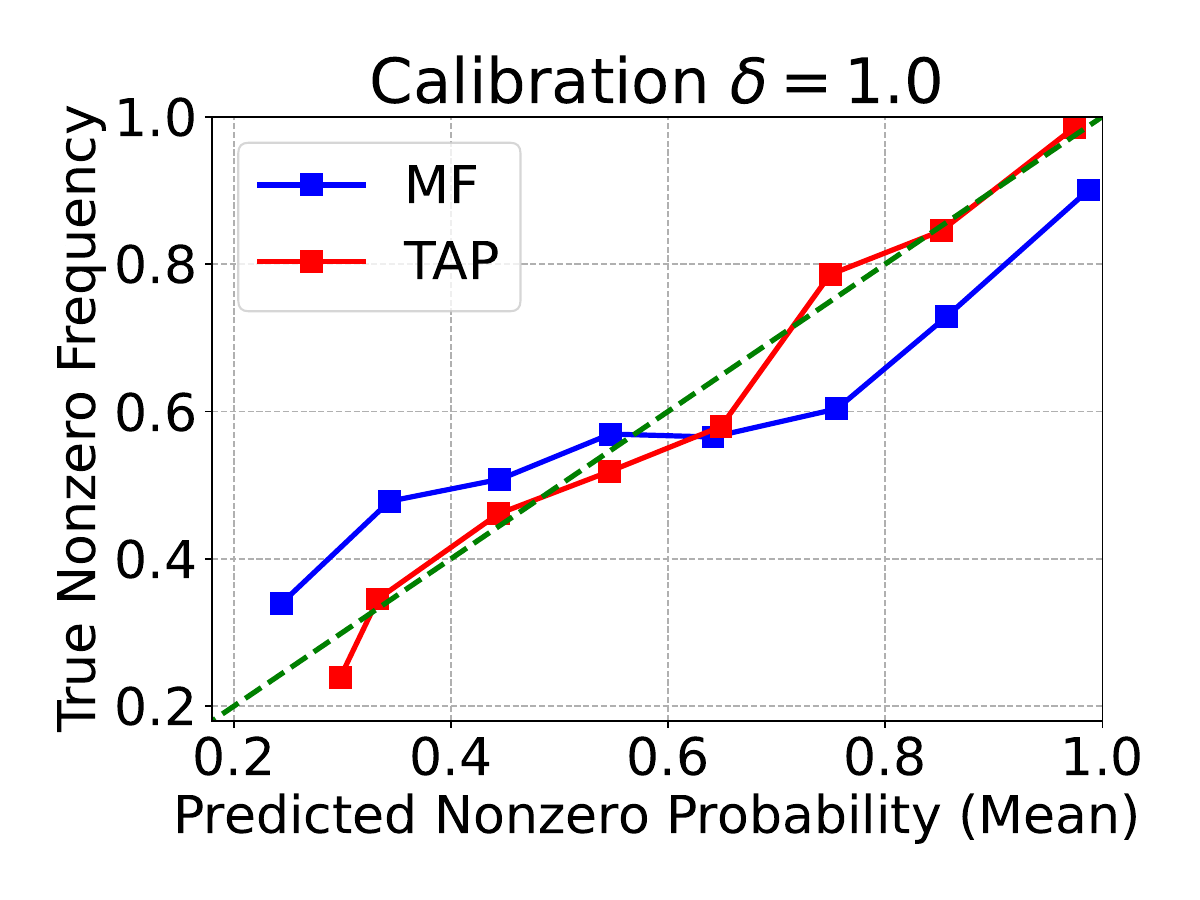}
\caption{Calibration of the estimated PIPs under three types of
misspecifications. Parameters: $n=500$, $\sigma=0.3$, Bernoulli-Gaussian prior
$\sP_0=(1/2) \cdot \delta_{0}+ (1/2) \cdot \normal(0,1)$. Left: Radamacher $\X$.
Middle: Radamacher
$\beps$. Right: Bernoulli $\X$ with heterogeneous sparsity across columns.
The plots are generated similarly to Figure~\ref{fig:cali_1}, except under a
Bernoulli-Gaussian rather than three-point prior.}
\label{fig:uni_1}
\end{figure}

Although our theoretical results rely on the Gaussian assumptions of
$(\X,\beps)$, we expect these results to be robust under sufficiently
light-tailed distributions for the random design and additive noise. Here, we verify this numerically in three scenarios: 
\begin{itemize} 
	\item[(a)] The design $\X$ has i.i.d. entries generated from  $\mathrm{Unif}\{\pm 1/\sqrt{p}\}$. The noise remains Gaussian.
	\item[(b)] The noise $\beps$ has i.i.d. entries generated from $\mathrm{Unif}\{\sigma,-\sigma\}$.
	The design remains Gaussian.
	\item[(c)] For $j\in \{1,\ldots,p\}$, the $j^\text{th}$ column of the
design $\X$ has i.i.d.\ Bernoulli entries with parameter
$0.1+0.8\cdot(j-1)/(p-1)$. Each column is then standardized to mean 0 and
variance $1/p$, and the noise remains Gaussian.
\end{itemize}
The simulations in Sections \ref{sec:exp_mse} and \ref{sec:exp_marginal} are
repeated for the three misspecification scenarios described above, with results
reported in Figures \ref{fig:uni_05} and \ref{fig:uni_1}. From
Figure~\ref{fig:uni_05}, we see that the MSE values of the MF and TAP estimators
are both universal across the three distributional misspecifications. Moreover,
we observe in Figure \ref{fig:uni_1} that the PIPs estimated from the TAP
approximation remain correctly calibrated under all three misspecification scenarios. 

\begin{figure}[htp]
\includegraphics[width=0.44\linewidth]{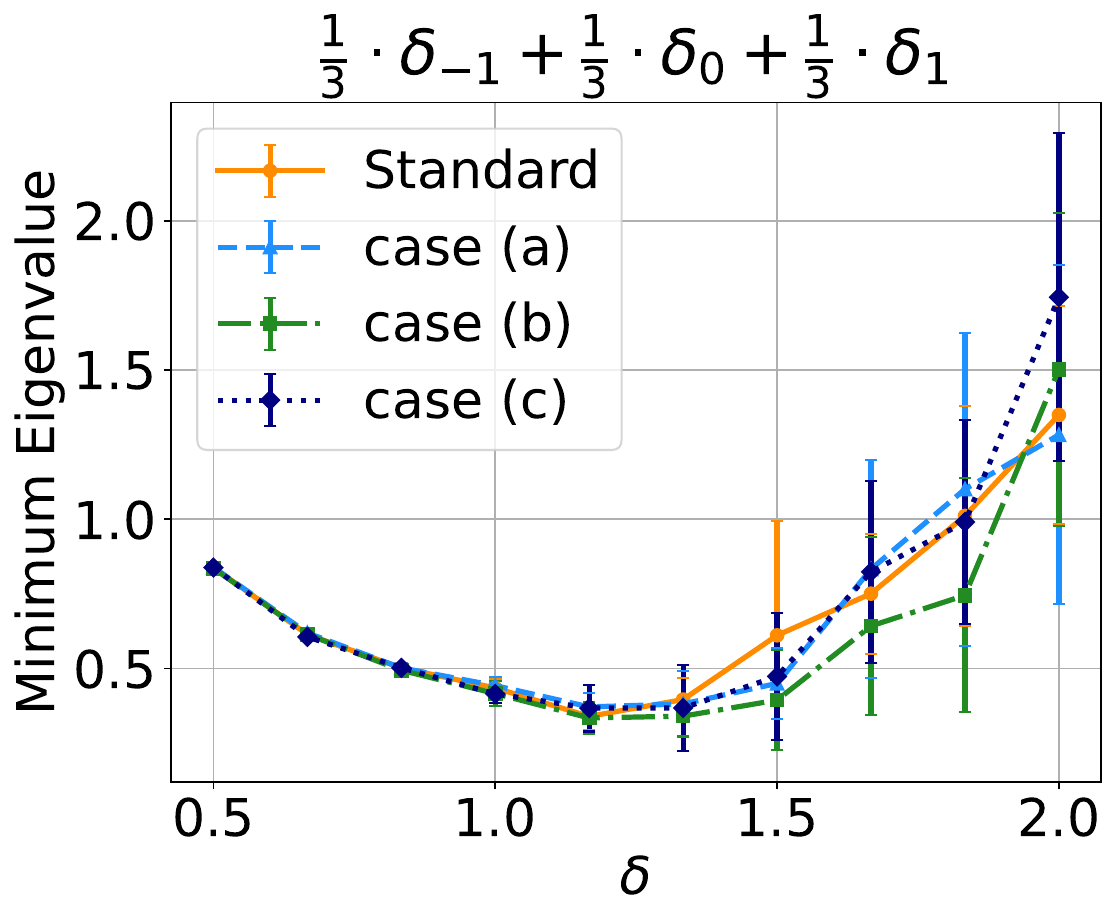}
\includegraphics[width=0.44\linewidth]{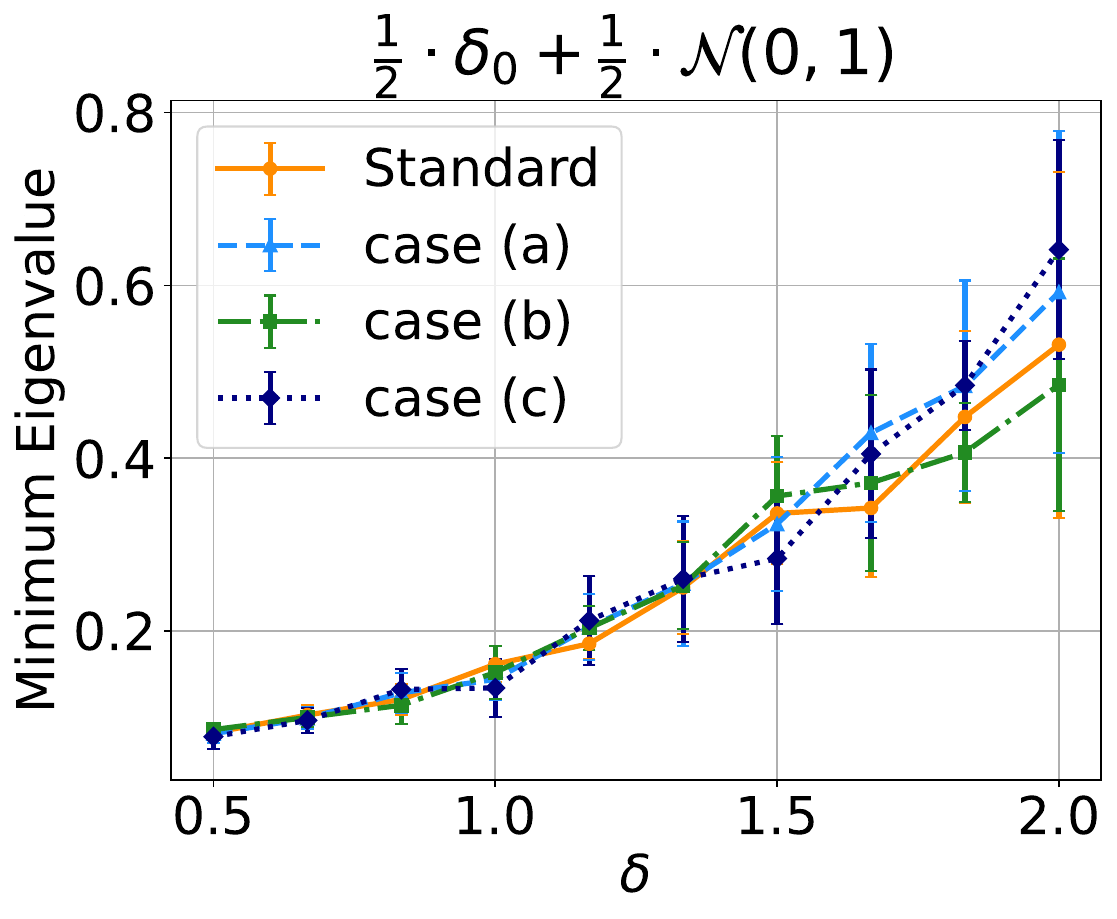}
\caption{
Minimum eigenvalue of the Hessian $\nabla^2\cF_{\TAP}(\bbm_\star,\bbs_\star)$ versus $\delta=n/p$ under different distributional misspecifications. Parameters: $n=300$, $\sigma=0.3$. Standard: Gaussian $\X$ and $\beps$. Case (a): Rademacher $\X$. Case (b): Rademacher $\beps$. Case (c): Bernoulli $\X$ with heterogeneous sparsity across columns. The approximate minimizer $(\bbm_\star,\bbs_\star)$ is computed by NGD. Left: Three-point prior $\sP_0=(1/3)\cdot \delta_{-1}+(1/3)\cdot \delta_{0}+(1/3)\cdot \delta_{1}$. Right: Bernoulli-Gaussian prior $\sP_0=(1/2)\cdot \normal(0,1)+(1/2)\cdot\delta_{0}$. Error bars show standard deviation over 10 simulations.
%
}
\label{fig:uni_2}
\end{figure}

Finally, Figure \ref{fig:uni_2} illustrates the minimum eigenvalue of the
Hessian $\nabla^2\cF_{\TAP}(\bbm_\star,\bbs_\star)$ at the  (approximate)
minimizer $(\bbm_\star,\bbs_\star)$ computed by NGD, under different scenarios
of misspecification. Under all levels of $\delta$ and across the three
distributional misspecifications, the minimum eigenvalue remains strictly
positive. The results also suggest a certain universality of this minimum eigenvalue value across these different types of misspecification.


\section{Proof ideas}

In this section, we describe some aspects of the proofs of our main results.
Our goal is not to provide a complete proof outline, but rather highlight some of the key ideas, and in particular, those which are particularly novel.
Complete proof details can be found in the appendices.

\subsection{The TAP lower bound via Gordon's comparison inequality}

The proofs of Theorems \ref{thm:bayes_consistence} and \ref{thm:marginalposterior} require (1) characterizing the posterior expectation of Lipschitz functions of $\beta_j$, and (2) establishing the existence of and properties of the Bayes-optimal local minimizer $\bbm_\star$, $\bs_\star$. 
The former requires extending results of \cite{barbier2019optimal}, which applies only to the posterior variances.
To carry out this extension, we use an adaptation of the interpolation argument of \cite[Theorem 1.7.11]{talagrand2010mean} to approximate the posterior marginals via the observation of a scalar Gaussian channel. We carry out this interpolation argument in Appendix \ref{sec:cavity}.

Here, we highlight aspects of our argument regarding the existence of a Bayes-optimal local minimizer.
The key step is showing that TAP free energy contains a local minimizer $(\bbm_\star,\bs_\star)$ satisfying $\tfrac1p \sum_{j=1}^p (m_{\star,j} - \beta_{0,j})^2 \approx \tfrac1p\sum_{j=1}^p s_{\star,j} - m_{\star,j}^2 \approx \mmse(\gamma_\stat) =: q_\stat$ and $\cF_\TAP(\bbm_\star,\bs_\star)/p \approx \inf_{\gamma > 0} \phi(\gamma) = \phi(\gamma_\stat)$.
For any $K \subseteq [0,\infty) \times [0,\infty)$, define
\begin{equation}\label{eq:GammaK-body}
	\Gamma^p[K]=\left\{(\bbm,\bs) \in \Gamma^p:\left(\frac{1}{p}
	\sum_{j=1}^p (m_j-\beta_{0,j})^2,\;\frac{1}{p} \sum_{j=1}^p s_j-m_j^2\right)
	\in K\right\}, 
\end{equation}
and let $K(\rho) \subseteq \R^2$ be the closed Euclidean-ball of radius $\rho$
around $(q_\stat,q_\stat)$.
For any $\rho_0 > \rho_1 > 0$ sufficiently small,
we will find $\iota > 0$ and a point $(\hbm,\hbs)$ such that,
with high-probability,
\begin{equation}
\label{eq:local-min-key-tools}
\begin{gathered}
	\frac1p \cF_\TAP(\hbm,\hbs) < \phi(\gamma_\stat) + \iota/2,
	\qquad
	(\hbm,\hbs) \in \Gamma^p[K(\rho_1)],
	\\
	\inf_{(\bbm,\bs) \in \Gamma^p[K(\rho_0)]} \frac1p \cF_\TAP(\bbm,\bs) > \phi(\gamma_\stat) - \iota,
	\quad
	\inf_{(\bbm,\bs) \in \Gamma^p[K(\rho_0)\setminus K(\rho_1)]} \frac1p \cF_\TAP(\bbm,\bs) > \phi(\gamma_\stat) + \iota.
\end{gathered}
\end{equation}
We will show that the first line is satisfied by taking $(\hbm,\hbs)$ to be an appropriate truncation of the Bayes estimate $(\bbm_{\sB},\bs_{\sB})$.
Then, to find the desired local minimizer, we consider a descent path from $(\hbm,\hbs)$ to any local minimizer $(\bbm_\star,\bs_\star)$. We can use the second line in the preceding display to argue that this descent path must remain in the set $\Gamma^p[K(\rho_1)]$, whence we can conclude that $\tfrac1p \sum_{j=1}^p (m_{\star,j} - \beta_{0,j})^2 \approx \tfrac1p\sum_{j=1}^p s_{\star,j} - m_{\star,j}^2 \approx q_\stat$ and 
$\cF_\TAP(\bbm_\star,\bs_\star)/p \approx \phi(\gamma_\stat)$, as desired.
Complete details are carried out in Appendix \ref{app:proof_bayes_consistence}.

A key step in this argument is proving the lower bound in the preceding display, which is the focus of Appendix \ref{app:TAP-lower-bound}.
The lower bound is based on the following lemma.
\begin{lemma}\label{lm:tap_lower_bound_body}
	Suppose Assumption \ref{ass:Bayesian_linear_model} holds.
	For any $\iota>0$ and compact sets $K\subset[0,\infty)\times [0,\infty)$,
	$K'\subset [0,\infty) \times \R^2$, there exists a constant
	$c>0$ (depending on $\iota,K,K'$) such that with probability $1-e^{-cn}$ for all
	large $n$,
	\begin{equation}
	\inf_{(\bbm,\bs) \in \Gamma^p[K]} \frac{1}{p}\,\cF_\TAP(\bbm,\bs)
	\geq \inf_{(q,r)\in K}\sup_{(\alpha,\tau,\gamma)\in K'} f(q, r; \alpha,\tau, \gamma) -\iota,
	\end{equation}
	where $f$ is a deterministic variational objective defined explicitly in Appendix \ref{app:gordon-lower-bound}.
\end{lemma}
\noindent We prove Lemma \ref{lm:tap_lower_bound_body} using Gordon's Gaussian comparison inequality  \cite{Gordon1985,gordon1988,thrampoulidis2015regularized} in Appendix \ref{app:gordon-lower-bound}.
Lemma \ref{lm:tap_lower_bound_body} reduces analysis of the left-hand sides in the second line of \eqref{eq:local-min-key-tools} to the analysis of a low-dimensional variational objective.
We show that for sufficiently large $K'$, the function $(q,r) \mapsto \sup_{(\alpha,\tau,\gamma)\in K'} f(q, r; \alpha,\tau, \gamma)$ has local minimizer of $(q,r) = (q_\stat,q_\stat)$ with value $\phi(\gamma_\stat)$ and is locally strongly convex in a neighborhood of this point.
From this, we can conlude the second line of \eqref{eq:local-min-key-tools}.

It is natural to conjecture that $(q_\stat,q_\stat)$ is in fact a global minimizer of $(q,r) \mapsto \sup_{(\alpha,\tau,\gamma)\in K'} f(q, r; \alpha,\tau, \gamma)$ for sufficiently large $K'$.
If one could show this to be the case,
we could conclude Theorem \ref{thm:bayes_consistence} for the \textit{global minimizer} of the TAP free energy.
We remark that the only gap in proving Theorem \ref{thm:bayes_consistence} for the \textit{global minimizer} is establishing this property of the variational objective $f$.

Finally, we provide a result, which may be of independent interest, which reveals the relationship between the variational lower bound of Lemma \ref{lm:tap_lower_bound_body} and the replica-symmetric potential $\phi$.
In particular, the lower bound of Lemma \ref{lm:tap_lower_bound_body} is given by the replica-symmetric potential if we restrict to points $\bbm,\bs$ satisfying a Nishimori-type condition:
\begin{equation}
	\frac{1}{p}\sum_{j=1}^p (m_j-\beta_{0,j})^2 \approx 
	\frac{1}{p} \sum_{j=1}^p s_j-m_j^2.
\end{equation}
\begin{theorem}\label{thm:replicalowerbound}
Let Assumption \ref{ass:Bayesian_linear_model} hold. Define
$\gamma(q)=\delta/(q+\sigma^2)$ for $q>0$, and fix any constant
$\iota>0$. Then there exists $\eps:=\eps(\iota,\delta,\sigma^2)>0$
such that the following
holds: For any compact set $K \subset [0,\infty) \times [0,\infty)$
satisfying $\sup_{(q,r) \in K} |q-r| \leq \eps$, with probability approaching 1
as $n,p \to \infty$,
\[\inf_{(\bbm,\bs) \in \Gamma^p[K]} \frac{1}{p}\,\cF_\TAP(\bbm,\bs)
\geq \inf_{(q,r) \in K} \phi(\gamma(q))-\iota.\]
\end{theorem}
\noindent The proof of Theorem \ref{thm:replicalowerbound} is given in Appendix \ref{sec:proof-thm-replicalowerbound}.
It follows from Lemma \ref{lm:tap_lower_bound_body} by showing that for $K'$ sufficiently large, $\phi(\gamma(q)) = \sup_{(\alpha,\tau,\gamma)\in K'} f(q,q;\alpha,\tau,\gamma)$.
By the definition of $\gamma_\stat$ (see \eqref{eq:gammastat}), 
Theorem \ref{thm:replicalowerbound} shows that the local minimizer $(\bbm_\star,\bs_\star)$ can be chosen to be globally optimal among those points satisfying the above Nishimori-type condition.
Said another way, the TAP global minimizer is far from $(\bbm_{\sB},\bs_{\sB})$ only if it does not satisfy the above Nishimori-type condition.
We remark that, although the replica-symmetric potential has appeared multiple times previously, we are unaware of previous work connecting it to the Gordon lower bound.

\subsection{The TAP local convexity via Gordon post-AMP}
\label{sec:gordon-post-AMP}

A central piece of the proof of Theorem \ref{thm:local_convexity}, Corollary \ref{cor:AMP+NGD-convergence}, and Theorem \ref{thm:local_convexity_AMP_fixed_point} is establishing that the TAP free energy is, with high probability, locally convex in a neighborhood of the AMP iterates $(\bbm^k,\bs^k)$ for sufficiently large $k$. For this, we adapt to the current setting a proof technique introduced by \cite{celentano2022sudakov} to study the TAP free energy in a low-rank matrix model involving a $\mathrm{GOE}$ matrix.

The key idea is to write the minimum value of the Hessian as a min-max problem whose objective is a Gaussian process.
In particular, in Appendix \ref{app:proof_local_convexity}, we show that the minimum eigenvalue of the Hessian of $\cF_{\TAP}$ evaluated at a point $(\bbm,\bs)$ can be written as
\begin{equation}
	\frac1p \min_{\|\bl\|_2 = \sqrt{p}}
		\langle \bl \nabla^2 \cF_{\TAP}(\bbm,\bs)\bl\rangle
		=
		2\min_{\|\bl_1\|_2^2 + \|\bl_2\|_2^2 = p} \; 
		\max_{\bu \in \reals^n} 
		\Big\{
			\frac1p \bu^\top \bX \bl_1
			+
			\Theta_{\TAP}(\bu,\bl_1,\bl_2;\bbm,\bs)
		\Big\},
\end{equation}
where $\bl_1,\bl_2 \in \reals^p$ and $\Theta_{\TAP}$ is an objective, specified in the appendix, whose randomness only involves the noise $\beps$.
Conditioning on $\beps$,
the objective in the preceding display is a Gaussian process whose randomness comes from the matrix $\bX$.
To establish a high probability lower-bound on the Hessian at a fixed point $(\bbm,\bs)$ chosen \textit{a priori} (that is independently of $\bX$),
we can use again Gordon’s comparison inequality.
It implies that for any $C,t > 0$,
\begin{equation}
\label{eq:gordon-unconditional}
\begin{aligned}
	\P\Big(
		&\min_{\|\bl_1\|_2^2 + \|\bl_2\|_2^2 = p} \; 
		\max_{\|\bu\|_2 \leq C} 
		\Big\{
			\frac1p \bu^\top \bX \bl_1
			+
			\Theta_{\TAP}(\bu,\bl_1,\bl_2;\bbm,\bs)
		\Big\} \leq t
	\Big)
	\\ 
	&\qquad\leq 
	2\P\Big(
		\min_{\|\bl_1\|_2^2 + \|\bl_2\|_2^2 = p} \; 
		\max_{\|\bu\|_2 \leq C} 
		\Big\{
			-\frac1{p^{3/2}} \| \bu \|_2 \< \bg , \bl_1 \>
			\frac1{p^{3/2}} \| \bl_1 \|_2 \< \bh , \bu \>
			+
			\Theta_{\TAP}(\bu,\bl_1,\bl_2;\bbm,\bs)
		\Big\} \leq t
	\Big),
\end{aligned}
\end{equation}
where $\bg \sim \normal(0,\id_p)$ and $\bh \sim \normal(0,\id_n)$ independent of each other and everything else.
Because the min-max problem on the right-hand side involves two high-dimensonal Gaussian vectors in place of a high-dimensional Gaussian matrix,
it is substantially easier to analyze than the min-max problem on the left-hand side.

For our results, we require not a lower bound on the Hessian at a point $(\bbm,\bs)$ chosen \textit{a priori},
but rather at all points in an $\epsilon \sqrt{p}$ ball around the AMP iterates $(\bbm^k,\bs^k)$.
That is, we require a high-probability lower-bound on
\begin{equation}
	2\min_{
		\substack{
			\|\bl_1\|_2^2 + \|\bl_2\|_2^2 = p \\
			\| \bbm - \bbm^k \|_2 / \sqrt{p} \leq \epsilon \\
			\| \bs - \bs^k \|_2 / \sqrt{p} \leq \epsilon
		}
	} \; 
	\max_{\bu \in \reals^n} 
	\Big\{
		\frac1p \bu^\top \bX \bl_1
		+
		\Theta_{\TAP}(\bu,\bl_1,\bl_2;\bbm,\bs)
	\Big\}.
\end{equation}
Because the domain of minimization depends on $\bbm^k,\bs^k$, 
which in turn depend on the random matrix $\bX$ via the iteration \eqref{eq:AMPalg},
we cannot apply Gordon’s comparison inequality directly to this problem.

The paper \cite{celentano2022sudakov} faced a similar challenge in a context involving a symmetric GOE random matrix and a minimization rather than a min-max problem. In that context, the appropriate Gaussian comparison inequality is the Sudakov-Fernique inequality. By conditioning on a sequence of AMP iterates and leveraging properties of the AMP state evolution, \cite{celentano2022sudakov} introduces an asymptotic comparison inequality---called the Sudakov-Fernique post-AMP inequality---whose analysis leads to an asymptotic lower bound on the TAP Hessian locally. In this paper,
we adapt this technique to a setting in which the matrix is standard Gaussian rather than GOE and the quantity of interest is represented by a min-max problem.
We call the resulting asymptotic comparison inequality the Gordon post-AMP inequality, and state it here.

It is convenient to state our result in terms of the quantities, for all $k \geq 0$
\begin{equation}
\begin{aligned}
      \bg^k
                  &:=
                  \bbm^k + \frac1\delta \bX^\top \bz^k - \bbeta_0,
            \qquad
            &\bh^k
                  &:=
                  \beps - \bz^k,
      \\
      \bnu^k
                  &:=
                  \bbm^k - \bbeta_0,
            \qquad
            &\br^k
            	&:=
            	-\bz^k.
\end{aligned}
\end{equation}
Let $\bR_k \in \R^{n \times k}$, $\bG_k \in \R^{p \times k}$, $\bV_k \in \R^{p
\times k}$, and $\bH_k \in \R^{n \times k}$ be the matrices whose columns are
$\{\br^{k'}\}_{1\leq k' \leq k}$, $\{ \bg^{k'} \}_{1 \leq k' \leq k}$,
$\{ \bnu^{k'} \}_{1 \leq k' \leq k}$, and
$\{ \bh^{k'} \}_{1 \leq k' \leq k}$ respectively.
Let $\proj_{\bR_k} \in \R^{n \times n}$ and $\proj_{\bV_k} \in \R^{p \times p}$ 
be the projections onto the linear spans of $\br^1,\ldots,\br^k$ and
$\bnu^1,\ldots,\bnu^k$, and let $\proj_{\bR_k}^\perp,\proj_{\bV_k}^\perp$ be the
projections onto the orthogonal complements.
Define
\begin{equation}
\begin{gathered}
    \bg^*(\bu)
          :=
          \frac{1}{\sqrt{\delta}}\,\frac{1}{n}\bG_k\bK_{g,k}^{-1}\bR_k^\top \bu
          +
          \frac1{\sqrt{n}}\|\proj_{\bR_k}^\perp \bu\|_2 \bg,
    \qquad
    \bh^*(\bl_1)
          :=
          \frac1p \bH_k \bK_{h,k}^{-1} \bV_k^\top \bl_1
          +
          \frac1{\sqrt{p}} \| \proj_{\bV_k}^\perp \bl_1 \|_2 \bh.
\end{gathered}
\end{equation}
The Gordon post-AMP inequality is the following asymptotic comparison inequality, which holds for any fixed constants $C_0,\epsilon > 0$:
\begin{proposition}[Gordon post-AMP]
\label{prop:gordon-post-amp}
      We have
      \begin{equation}
      \begin{aligned}
            &\pliminf_{n \rightarrow \infty}
            \min_{
                  \substack{
                        \|\bl\|_2/\sqrt{p} = 1 \\
                        \|\bbm - \bbm^k\|_2/\sqrt{n} \leq \epsilon \\
                        \|\bs - \bs^k \|_2/\sqrt{n} \leq \epsilon
                  }
            }\;
            \max_{\|\bu\|_2 \leq C_0\sqrt{n}}
            \Big\{
                  \frac1p \bu^\top \bX \bl_1
                  +
                  \Theta_{\TAP}(\bu,\bl_1,\bl_2;\bbm,\bs)
            \Big\}
      \\
            &\qquad\geq
            \pliminf_{n \rightarrow \infty}
            \min_{
                  \substack{
                        \|\bl\|_2/\sqrt{p} = 1 \\
                        \|\bbm - \bbm^k\|_2/\sqrt{p} \leq \epsilon \\
                        \|\bs - \bs^k \|_2/\sqrt{p} \leq \epsilon
                  }
            }\;
            \max_{\|\bu\|_2 \leq C_0\sqrt{n}}
            \Big\{
            -\sqrt{\delta}\,\frac{\<\bg^*(\bu),\bl_1\>}{p}
		    +\delta\,\frac{\<\bh^*(\bl_1),\bu\>}{n}
		    +\Theta_{\TAP}(\bu,\bl_1,\bl_2;\bbm,\bs)
		    \Big\}.
      \end{aligned}
      \end{equation}
\end{proposition}
\noindent We prove Proposition \ref{prop:gordon-post-amp} in Appendix \ref{app:proof-gordon-post-amp}.
Proposition \ref{prop:gordon-post-amp} reduces a min-max problem involving the Gaussian matrix $\bX$ to a min-max problem involving two high-dimensional Gaussian vectors $\bg,\bh$ and $O(k)$ AMP iterates $\bg^{k'}$, $\bh^{k'}$, $\bnu^{k'}$, and $\br^{k'}$ for $k' \leq k$.
Establishing local convexity of the TAP free energy in a neighborhood around the AMP iterates requires carrying out an asymptotic analysis of this reduces min-max problem.
We will consider this problem for $k$ fixed (in $n$) but still arbitrarily large.
Thus, there remain substantial challenges in carrying out this analysis. The details are provided in Appendix \ref{app:proof_local_convexity}.

\begin{remark}
	As is clear from its proof,
	a form of Proposition \ref{prop:gordon-post-amp} can be established for the general class of AMP algorithms considered \cite{berthierMontanariNguyen} and a general class of objectives $\Theta$ satisfying appropriate regularity conditions.
	Because we only consider the problem of local convexity of the TAP free energy, we limit ourselves to the statement in Proposition \ref{prop:gordon-post-amp}.
	We leave a more general statement and its application to a wider range of problems to future work.
\end{remark}

\begin{remark}
	It is worth comparing Proposition \ref{prop:gordon-post-amp} to \eqref{eq:gordon-unconditional}.
	Indeed, if $\delta = n/p$, the objectives in Proposition \ref{prop:gordon-post-amp} is obtained by replacing $\| \bu \|_2 \bg $ by $ \sqrt{n} \, \bg^*(\bu)$ and $\| \bl_1 \|_2 \bh$ by $\sqrt{p} \bh^*(\bl_1)$ in \eqref{eq:gordon-unconditional}.
	We encourage the reader to adopt the following intuition for why these replacements might be reasonable.
	Speaking imprecisely, $\| \bu \|_2 \bg$ can be thought of as a Gaussian vector whose norm is determined implicitly by $\bu$.
	On the other hand, $\sqrt{n}\bg^*(\bu)$ can be thought of as a Gaussian vector whose norm and direction is determined implicitly by $\bu$.
	Indeed, according to the AMP state evolution, stated in Appendix \ref{app:proof_local_convexity}, the iterates $\{ \bg^{k'} \}_{k' \leq k}$ behave, in a certain sense, like correlated high-dimensional Gaussian vectors.
	Because $k$ is fixed as $n \rightarrow \infty$ and $\bg^*(\bu)$ is in the linear span of $\{ \bg^{k'} \}_{k' \leq k}$ and $\bg$, it too behaves, in a certain sense, like a high-dimensional Gaussian vector, but both its direction and norm depend implicitly on $\bu$. 
	In fact, using the the AMP state evolution,
	stated in Appendix \ref{app:proof_global_convexity},
	one can show that the norm of $\sqrt{n} \bg^*(\bu)$ approximately agrees with the norm of $\| \bu \|_2 \bg$.
	Similar remarks apply to the replacement of $\| \bl_1 \|_2 \bh$ by $\sqrt{p} \bh^*(\bl_1)$.
	Thus, the structure of the Gordon post-AMP inequality agrees with that of the standard Gordon inequality, except that the effective Gaussian noise is replaced by effective Gaussian noise that is correlated with the AMP algorithm.   
\end{remark} 

\begin{remark}
	The fact that one can apply Gordon's inequality conditional on a sequence of AMP iterates is not novel.
	Thus, our main contribution is not to observe that an asymptotic comparison inequality like that in Proposition \ref{prop:gordon-post-amp} is possible, but rather to identify structure in the resulting min-max problem that facilitates its analysis.
	Just as Lemma \ref{lm:tap_lower_bound_body} provides a lower bound on the TAP free energy in terms of a low-dimensional variational problem,
	we can lower bound right-hand side in Proposition \ref{prop:gordon-post-amp} by a low-dimensional variational problem.
	As carried out in Appendix \ref{app:proof_local_convexity},
	the main steps are to (1) reduce this problem further to a min-max problem defined on Wasserstein space involving $O(k)$ scalar random variables, and (2) reduce the problem involving $O(k)$ random variables to a problem involving $O(1)$ scalar random variables. Thus, we arrive at an explicit variational problem which does not depend on $k$, and can explicitly analyze this reduced problem.
	Although, at a high level, \cite{celentano2022sudakov} followed a similar sequence of steps, 
	carrying out each step requires substantial novelty in the present setting due to the min-max structure of the problem.
\end{remark}

\section{Conclusion and discussion}

In this paper, we studied variational inference for high-dimensional Bayesian
linear models, and showed the existence of local minimizers of the TAP free
energy that consistently approximate the true marginal posterior laws. We proved
the local convexity of the TAP landscape and finite-sample convergence of
natural gradient descent, to this Bayes-optimal local minimizer in the
computationally ``easy'' regime and to a surrogate local minimizer in the
``hard'' regime. 
In both regimes, the local minimizer can be used for correctly calibrated
posterior inference.
Numerical simulations confirm that the TAP free energy can be efficiently
optimized, and that properties of its minimizers exhibit some robustness to
model misspecification. Together, these results provide theoretical justification for using the TAP free energy to perform variational inference in Bayesian linear models.

Our proof of local convexity employed a novel technique, utilizing Gordon's
inequality conditioned on the iterates of Approximate Message Passing to lower
bound the Hessian around the Bayes-optimal local minimizer. This generalizes the
technique of the Sudakov-Fernique inequality after AMP in
\cite{celentano2022sudakov} for optimization objectives defined instead by
symmetric Gaussian matrices. This proof technique can be readily generalized to
other statistical models, for example generalized linear models. 

Finally, recent work has demonstrated close connections between variational
inference for computing posterior expectations and sampling from the posterior
distribution through stochastic localization and diffusion-based methods
\cite{el2022sampling, el2023sampling, montanari2023sampling,
montanari2023posterior, ghio2023sampling, mei2023deep}. For example, such
variational approaches to sampling were studied for spin-glass models in
\cite{el2022sampling, el2023sampling} and low-rank matrix denoising models in
\cite{montanari2023posterior}. Based on this connection, we believe our results
will be of interest for developing analogous sampling algorithms for supervised
learning models such as linear regression.



%

\section*{Acknowledgement}

Michael Celentano is supported
by the Miller Institute for Basic Research in Science, University of California, Berkeley.
Zhou Fan is supported by NSF DMS-2142476.
Song Mei is supported by NSF DMS-2210827 and NSF CCF-2315725.

\clearpage

\appendix

\section{Moment space of the exponential family}\label{sec:domain}

Recall from (\ref{eq:momentspace}) the space $\Gamma \subset \R^2$ of possible
first and second moments for the exponential family (\ref{eq:scalarmean}).
We provide here an explicit characterization of this set $\Gamma$ and some
properties of the relative entropy function ${-}\sh(\cdot)$ in (\ref{eq:hdef}).

Define the support of $\sP_0$ by
\[\supp(\sP_0)=\{x \in \R:\;\sP_0((x-\eps,x+\eps))>0 \text{ for all } \eps>0\}.\]
Equivalently, this is the smallest closed set $K$ for which $\sP_0(K)=1$. Denote
the lower and upper endpoints of this support by
\begin{equation}\label{eq:supportendpoints}
a(\sP_0)=\inf\{x \in \R:x \in \supp(\sP_0)\},
\qquad b(\sP_0)=\sup\{x \in \R:x \in \supp(\sP_0)\}.
\end{equation}
The following elementary proposition first characterizes the domain of the
first moment $m=\langle \beta \rangle_{\lambda,\gamma}$,
fixing any $\gamma \in \R$.

\begin{proposition}\label{prop:oneparamexpfam}
Suppose $\sP_0$ has compact support containing at least two distinct values.
Fix any $\gamma \in \R$. Then
\[\Big\{m \in \R:\text{ there exists } \lambda \in \R \text{ such that }
m=\langle \beta \rangle_{\lambda,\gamma}\Big\}=(a(\sP_0),b(\sP_0)).\]
For any $m \in (a(\sP_0),b(\sP_0))$,
this value $\lambda:=\lambda_\gamma(m)$ for which
$m=\langle \beta \rangle_{\lambda,\gamma}$ is unique, and $\lambda_\gamma(m)$
is continuously differentiable in $(m,\gamma)$ and strictly increasing in $m$.
\end{proposition}
\begin{proof}
The law $\sP_{\lambda,\gamma}$ has support contained in $[a(\sP_0),b(\sP_0)]$
with at least two distinct values. Thus if $m=\langle \beta
\rangle_{\lambda,\gamma}$, then $m \in (a(\sP_0),b(\sP_0))$.

Conversely, fix $m \in (a(\sP_0),b(\sP_0))$ and
define $f(\lambda;m,\gamma)=\langle \beta \rangle_{\lambda,\gamma}-m$.
Note that $\partial_\lambda f(\lambda;m,\gamma)=\Var_{\beta \sim \sP_{\lambda,\gamma}}[\beta]>0$
strictly, because $\sP_{\lambda,\gamma}$ has at least two points of support.
Then there exists a unique value $\lambda=\lambda_\gamma(m)$ for which
$f(\lambda;m,\gamma)=0$, i.e.\ $\< \beta \>_{\lambda,\gamma} = m$,
because $\lim_{\lambda \rightarrow -\infty} \< \beta
\>_{\lambda,\gamma} = a(\sP_0)$, $\lim_{\lambda\rightarrow \infty}
\<\beta\>_{\lambda,\gamma} = b(\sP_0)$,
and $\lambda \mapsto f(\lambda;m,\gamma)$ is
continuous and strictly increasing. Since $\partial_\lambda
f(\lambda_\gamma(m);m,\gamma)>0$ and $f(\lambda;m,\gamma)$ is continuously differentiable with respect to $(\lambda, m, \gamma)$, the implicit function theorem implies that
$\lambda_\gamma(m)$
is continuously differentiable in $(m,\gamma)$, with derivative in $m$ given by
$\lambda_\gamma'(m)=1/\partial_\lambda
f(\lambda_\gamma(m);m,\gamma)>0$. Thus $\lambda_\gamma(m)$ is also strictly
increasing in $m$.
\end{proof}

Now fixing $\sP_0$, for each $m \in [a(\sP_0),b(\sP_0)]$, define
\[a(m)=\sup\{x \in \R:x \leq m,\,x \in \supp(\sP_0)\},\qquad 
b(m)=\inf\{x \in \R: x \geq m,\,x \in \supp(\sP_0)\}.\]
Thus $a(m) \leq m \leq b(m)$ with equality if and only if $m \in \supp(\sP_0)$.
%
%
Recall the entropy function ${-}\sh(m,s)$ from (\ref{eq:hdef}), interpreted as
an extended real-valued function ${-}\sh:\R^2 \to [0,\infty]$,
and define $\domain({-}\sh) = \{(m,s) \in \R^2:-\sh(m,s) < \infty\}$.
We then have the following characterization of the moment space $\Gamma$
and $\domain({-}\sh)$.

\begin{proposition}\label{prop:Gammacharacterization}
\label{prop:domain-h}
Suppose that $\sP_0$ has compact support with at least three distinct values.
Then the set $\Gamma$ in (\ref{eq:momentspace}) is given explicitly by
\begin{equation}\label{eq:momentspaceexplicit}
\Big\{(m,s) \in \R^2: a(\sP_0)<m<b(\sP_0),\;
(a(m)+b(m))m-a(m)b(m) < s < (a(\sP_0)+b(\sP_0))m-a(\sP_0)b(\sP_0) \Big\}. 
\end{equation}
Furthermore, $\Gamma \subseteq \domain({-}\sh) \subseteq \overline{\Gamma}$.
\end{proposition}

%

\begin{proof}
Note that since the maximization defining $-\sh(m,s)$ in (\ref{eq:hdef}) is
concave with stationary conditions $m=\langle \beta
\rangle_{\lambda,\gamma}$ and $s=\langle \beta^2 \rangle_{\lambda,\gamma}$,
these stationary conditions hold if and only if $(\lambda,\gamma)$ attains the
supremum in (\ref{eq:hdef}). Thus the definition of $\Gamma$ in
(\ref{eq:momentspace}) is equivalently
\begin{equation}\label{eq:Gammaequiv}
\Gamma=\Big\{(m,s) \in \R^2:\text{ the supremum defining } {-}\sh(m,s)
\text{ in (\ref{eq:hdef}) is attained}\Big\}
\end{equation}
Let us denote the set (\ref{eq:momentspaceexplicit}) by $\Gamma'$, so we wish to
show $\Gamma=\Gamma'$ and $\Gamma' \subseteq \domain({-}\sh) \subseteq
\overline{\Gamma'}$.

It will be convenient to first center $\supp(\sP_0)$:
For any constant $c \in \R$, let $\tilde\sP_0$ denote the law of $\beta+c$,
and let $\tilde\Gamma$ and $\tilde\sh(m,s)$ be the corresponding moment space
and entropy function for this prior $\tilde\sP_0$. We have
\begin{align}
&-\frac{1}{2}\gamma s+\lambda m-\log \E_{\beta \sim \tilde\sP_0}
\big[e^{-(\gamma/2)\beta^2+\lambda\beta}\big]\nonumber\\
&=-\frac{1}{2}\gamma s+\lambda m-\log \E_{\beta \sim \sP_0}
\big[e^{-(\gamma/2)(\beta+c)^2+\lambda(\beta+c)}\big]\nonumber\\
&=-\frac{1}{2}\gamma(s-2mc+c^2)+(\lambda-\gamma c)(m-c)
-\log \E_{\beta \sim \sP_0}\big[e^{-(\gamma/2)\beta^2+(\lambda-\gamma
c)\beta}\big].\label{eq:supportrecenter}
\end{align}
Then by (\ref{eq:Gammaequiv}), $(m,s) \in \tilde\Gamma$ if and only if
$(m-c,s-2mc+c^2) \in \Gamma$, and similarly $(m,s) \in \domain({-}\tilde \sh)$
if and only if $(m-c,s-2mc+c^2) \in \domain({-}\sh)$. Furthermore,
let $\tilde \Gamma'$ be the analogue of (\ref{eq:momentspaceexplicit}) for
$\tilde \sP_0$, defined by $a(\tilde\sP_0)=a(\sP_0)+c$, $b(\tilde\sP_0)=b(\sP_0)+c$, $\tilde a(m)=a(m-c)+c$, and $\tilde b(m)=b(m-c)+c$.
Then it is direct to check also that
$(m,s) \in \tilde\Gamma'$ if and only if $(m-c,s-2mc+c^2) \in
\Gamma'$. So it suffices to prove the proposition for $\tilde \sP_0$ instead of
$\sP_0$.

Choosing $c \in \R$ appropriately, we
may thus assume without loss of generality that $-a(\sP_0)=b(\sP_0):=M>0$.
In particular, $\pm M \in \supp(\sP_0)$.
Then (\ref{eq:momentspaceexplicit}) takes the form
\begin{equation}\label{eq:Gammacentered}
\Gamma'=\Big\{(m,s) \in \R^2: |m|<M,\; m(a(m)+b(m))-a(m)b(m)<s<M^2 \Big\}. 
\end{equation}
It may be checked that the function $m \mapsto m(a(m)+b(m))-a(m)b(m)$ is
continuous (equal to $m^2$ on $\supp(\sP_0)$, and linearly interpolating
between these values over each interval of $\R \setminus \supp(\sP_0)$),
so $\Gamma'$ is open. Under the given condition that $\sP_0$ has at least three
points of support, this function is not identically equal to $M^2$, so $\Gamma'$
is also non-empty, and its closure $\overline{\Gamma'}$ is defined by
replacing all three $<$ above by $\leq$.

First, we show that $-\sh(m,s)=\infty$ if $(m,s) \not \in \overline{\Gamma'}$,
and hence $\domain({-}\sh) \subseteq \overline{\Gamma'}$.
We consider several cases:
\begin{itemize}
\item Suppose $s < m(a(m)+b(m)) - a(m)b(m)$. For any $\gamma \in \R$,
	take $\lambda = \gamma (a(m)+b(m))/2$.
	Then
	\begin{equation}
	\begin{aligned}
		-\sh(m,s)
			&\geq
			\frac{\gamma}{2}
			\Big(
				m(a(m)+b(m)) + \frac{(a(m)+b(m))^2}{4}-s
			\Big)
			- \log \E_{\beta \sim \sP_0}\big[e^{-(\gamma/2)\big(\beta - \frac{a(m)+b(m)}{2}\big)^2}\big]
		\\
			&\geq 
			\frac{\gamma}{2}
			\Big(
				m(a(m)+b(m)) - \frac{(a(m)+b(m))^2}{4} + \frac{(a(m)-b(m))^2}{4} - s
			\Big) \\
			&=\frac{\gamma}{2}
			\Big(
				m(a(m)+b(m)) - a(m)b(m) - s
			\Big), 
	\end{aligned}	
	\end{equation}
the second inequality holding
	because $\E_{\beta \sim \sP_0}[e^{-(\gamma/2)(\beta -
\frac{a(m)+b(m)}{2})^2}] \leq e^{-(\gamma/2)(\frac{a(m)-b(m)}{2})^2}$ using
that $\supp(\sP_0)$ does not intersect $(a(m),b(m))$.
	Taking $\gamma \rightarrow \infty$ gives $-\sh(m,s) = \infty$.

\item Suppose $s > M^2$. Take $\lambda = 0$. Then
	\begin{equation}
	\begin{aligned}
		-\sh(m,s)
			\geq
			-\frac12 \gamma s
			- \log \E_{\beta \sim \sP_0}\big[e^{-(\gamma/2)\beta^2}\big]
			\geq 
			-\frac12 \gamma s
			+\frac12 \gamma M^2,
	\end{aligned}	
	\end{equation}
	because $\supp(\sP_0) \subseteq [-M,M]$
so $\E_{\beta \sim \sP_0}[e^{-(\gamma/2)\beta^2}] \geq
e^{-(\gamma/2)M^2}$. Taking $\gamma \rightarrow -\infty$ gives $-\sh(m,s) = \infty$.

\item Suppose $|m| > M$. Take $\gamma = 0$. Then similarly
	\begin{equation}
	\begin{aligned}
		-\sh(m,s)
			\geq
			\lambda m
			- \log \E_{\beta \sim \sP_0}\big[e^{\lambda \beta}\big]
			\geq 
			\lambda m
			- |\lambda| M,
	\end{aligned}	
	\end{equation}
	and taking $\lambda \rightarrow \infty$ if $m > M$ and $\lambda \rightarrow -\infty$ if $m < -M$ gives $-\sh(m,s) = \infty$.
\end{itemize}
	\noindent Thus we have shown that $\domain({-}\sh) \subseteq
\overline{\Gamma'}$.

	Next, we show that $\Gamma' \subseteq \Gamma$.
	To do so, fix any $(m,s) \in \Gamma'$. Then $m \in (-M,M)$,
so Proposition \ref{prop:oneparamexpfam} verifies that
for each $\gamma \in \reals$, there exists a unique value
$\lambda_\gamma(m)$ such that $\< \beta \>_{\lambda_\gamma(m),\gamma} = m$,
and furthermore $\lambda_\gamma(m)$ is continuous in $\gamma$.
	We use the following key lemma, which is proved after the completion of the proof of Proposition \ref{prop:domain-h}.
	\begin{lemma}
	\label{lem:large-gamma-moments}
		We have $\lim_{\gamma \rightarrow \infty} \< \beta^2
\>_{\lambda_\gamma(m),\gamma} = m(a(m)+b(m)) - a(m)b(m)$
		and $\lim_{\gamma \rightarrow -\infty} \< \beta^2
\>_{\lambda_\gamma(m),\gamma} = M^2$.
	\end{lemma}
\noindent Then, since $m(a(m)+b(m))-a(m)b(m)<s<M^2$ and
$\gamma \mapsto \langle \beta^2
\rangle_{\lambda_\gamma(m),\gamma}$ is continuous by the continuity of
$\lambda_\gamma(m)$, there exists $\gamma$ for which also $\langle \beta^2
\rangle_{\lambda_\gamma(m),\gamma}=s$, so $(m,s) \in \Gamma$. Thus we have shown
that $\Gamma' \subseteq \Gamma$.

	
The definition (\ref{eq:Gammaequiv}) of $\Gamma$ implies immediately $\Gamma
\subseteq \domain({-}\hs)$, so $\Gamma' \subseteq \Gamma \subseteq
\domain({-}\hs) \subseteq \overline{\Gamma'}$. Since $\Gamma,\Gamma'$ are both open
in $\R^2$, the inclusion $\Gamma' \subseteq \Gamma \subseteq \overline{\Gamma'}$
implies $\Gamma=\Gamma'$, completing the proof.
\end{proof}

\begin{proof}[Proof of Lemma \ref{lem:large-gamma-moments}]
We will fix $m$ and write as shorthand $\lambda(\gamma)=\lambda_\gamma(m)$.
Recall $\sP_{\lambda,\gamma}$ from (\ref{eq:scalarmean}), which has the
equivalent form $\sP_{\lambda,\gamma}(\de \beta) \propto
e^{-(\gamma/2)(\beta-\lambda/\gamma)^2}\sP_0(\de \beta)$.

First consider $m \in \supp(\sP_0)$, so that $a(m) = b(m) = m$.
For any $\eps>0$, we show at the conclusion of the proof the following claim:
	\begin{equation}
	\label{eq:supp-large-gamma}
	\begin{gathered}
		\lim_{\gamma \rightarrow \infty}
			\sup_{\lambda : \lambda / \gamma \geq m}
			\sP_{\lambda,\gamma}(\beta \leq m - \eps) = 0,
		\qquad
		\lim_{\gamma \rightarrow \infty}
			\sup_{\lambda : \lambda / \gamma \leq m}
			\sP_{\lambda,\gamma}(\beta \geq m + \eps) = 0.
	\end{gathered}
	\end{equation}
Fixing $\eps,\delta>0$,
we have $0 = \<\beta-m\>_{\lambda(\gamma),\gamma}  \geq
\epsilon \sP_{\lambda(\gamma), \gamma}(\beta-m \geq \epsilon) - \delta - 2M
\sP_{\lambda(\gamma), \gamma}(\beta-m \leq {-}\delta)$, which implies that 
	$
	\sP_{\lambda(\gamma),\gamma}(\beta \geq m+\epsilon)
		\leq \frac{\delta + 2M\sP_{\lambda(\gamma), \gamma}(\beta \leq
m-\delta)}{\eps}$. Likewise,
	$
	\sP_{\lambda(\gamma),\gamma}(\beta \leq m-\epsilon)
		\leq 
		\frac{\delta + 2M\sP_{\lambda(\gamma), \gamma}(\beta \geq m+\delta)}{\epsilon}
	$.
	Considering the cases
$\lambda(\gamma)/\gamma\geq m$ and $\lambda(\gamma)/\gamma\leq m$ separately,
	we conclude that 
	\begin{align*}
		\sP_{\lambda(\gamma),\gamma}(|\beta-m|\geq \epsilon)
			&\leq \sup_{\lambda : \lambda / \gamma \geq m}
			\Big\{
				\sP_{\lambda,\gamma}(\beta \leq m - \epsilon)
				+
				\frac{\delta + 2M\sP_{\lambda, \gamma}(\beta
\leq m - \delta)}{\epsilon}			\Big\}
		\\ &\hspace{1in}+
			\sup_{\lambda : \lambda / \gamma \leq m}
			\Big\{
				\sP_{\lambda,\gamma}(\beta \geq m + \epsilon)
				+
				\frac{\delta + 2M\sP_{\lambda, \gamma}(\beta
                                        \geq m + \delta)}{\epsilon}\Big\}.
	\end{align*}
	Choosing $\delta=\epsilon^2$ and
using (\ref{eq:supp-large-gamma}) gives $\lim_{\gamma \rightarrow
\infty}\sP_{\lambda(\gamma),\gamma}(|\beta-m|\geq \epsilon) \leq 2\epsilon$.
	Thus, $\lim_{\gamma \rightarrow \infty} |\< \beta^2
\>_{\lambda(\gamma),\gamma}-m^2| \leq 2M\epsilon + (2M)^2(2\eps)$.
	Since this holds for any $\epsilon>0$, we conclude that $\lim_{\gamma \rightarrow \infty} \< \beta^2 \>_{\lambda(\gamma),\gamma} = m^2$.

	Next, consider $m \not \in \supp(\sP_0)$.
For any $\eps>0$, we show at the conclusion of the proof that
	\begin{equation}
	\label{eq:off-supp-large-gamma}
		\lim_{\gamma \rightarrow \infty}
		\inf_{\lambda/\gamma \in [a(m),b(m)]}
		\sP_{\lambda,\gamma}\big(\beta \in [a(m)-\epsilon,a(m)]\cup[b(m),b(m)+\epsilon]\big)
		=1.
	\end{equation}
For $\eps>0$ small enough such that $a(m)+\eps<m-\eps$,
(\ref{eq:supp-large-gamma}) shows $\lim_{\gamma \to \infty}
\sup_{\lambda:\lambda/\gamma \leq a(m)}
\sP_{\lambda,\gamma}(\beta \geq m-\eps)=0$. This implies
$\lim_{\gamma \to \infty}
\sup_{\lambda:\lambda/\gamma \leq a(m)} \langle \beta \rangle_{\lambda,\gamma}
\leq m-\eps$, so we must have $\lambda(\gamma)/\gamma \geq a(m)$ for
sufficiently large $\gamma$. Similarly, $\lambda(\gamma)/\gamma \leq b(m)$
for sufficiently large $\gamma$.
	Then (\ref{eq:off-supp-large-gamma}) implies
$\lim_{\gamma \rightarrow \infty} \big\< (\beta - (a(m)+b(m))/2)^2
\big\>_{\lambda(\gamma),\gamma} = (b(m) - a(m))^2/4$. Together with $\langle
\beta \rangle_{\lambda(\gamma),\gamma}=m$, this  implies
	$\lim_{\gamma \rightarrow \infty}\<\beta^2\>_{\lambda(\gamma),\gamma} = m(a(m)+b(m))- a(m)b(m)$.

	Finally, consider any $m \in (-M,M)$. We show at the conclusion of the
proof that
	\begin{equation}
	\label{eq:small-gamma}
		\lim_{\gamma \rightarrow -\infty}
		\sup_{\lambda \in \reals} \sP_{\lambda,\gamma}(|\beta|\leq M-\epsilon)
		=0.
	\end{equation}
	Then, since $\< \beta^2 \>_{\lambda,\gamma} \geq
(M-\eps)^2(1-\sP_{\lambda,\gamma}(|\beta|\leq M-\epsilon))$, this shows
	$\lim_{\gamma \rightarrow -\infty} \< \beta^2 \>_{\lambda(\gamma),\gamma} = M^2$.

	It remains to show \eqref{eq:supp-large-gamma}, \eqref{eq:off-supp-large-gamma}, and \eqref{eq:small-gamma}.
If $\lambda / \gamma \geq m$,
	then $\inf_{\beta \leq m - \eps} (\beta - \lambda/\gamma)^2
\geq (\lambda/\gamma-m+\eps)^2$ and $\sup_{|\beta - m| \leq \eps/2}
(\beta - \lambda/\gamma)^2 \leq (\lambda/\gamma-m+\eps/2)^2$,
where $(\lambda/\gamma-m+\eps)^2-(\lambda/\gamma-m+\eps/2)^2
 \geq 3\eps^2/4$. Hence
	\begin{equation}
		\sP_{\lambda,\gamma}(\beta \leq m-\eps)
		\leq 
		\frac{\sP_{\lambda,\gamma}(\beta \leq
m-\eps)}{\sP_{\lambda,\gamma}(|\beta-m| \leq \eps/2)} 
			\leq
			e^{-(\gamma/2)(3\eps^2/4)}
			\frac{\sP_0(\beta \leq m-\eps)}{\sP_0(|\beta-m| \leq
\eps/2)}.
	\end{equation}
Here, $\sP_0(|\beta-m| \leq \eps/2)>0$ because $m \in \supp(\sP_0)$.
Then, taking $\gamma \to \infty$ for fixed $\eps>0$,
the first statement of \eqref{eq:supp-large-gamma} holds, and the second
statement is analogous.

	Next we prove \eqref{eq:off-supp-large-gamma}.
	If $\lambda/\gamma \in [a(m),b(m)]$,
	then $(\inf_{\beta \geq b(m)+\epsilon} (\beta - \lambda/\gamma)^2) - (\sup_{b(m) \leq \beta \leq b(m)+\epsilon/2} (\beta - \lambda/\gamma)^2) \geq \epsilon^2/2$.
	Thus,
	\begin{equation}
		\sP_{\lambda,\gamma}\big(\beta \geq b(m) + \epsilon\big)
		\leq 
		\frac{\sP_{\lambda,\gamma}\big(\beta \geq b(m) + \epsilon\big)}
		{\sP_{\lambda,\gamma}\big(b(m) \leq \beta \leq b(m) + \epsilon/2\big)}
		\leq 
		e^{-(\gamma/2)(\epsilon^2/2)}
		\frac{\sP_0\big(\beta \geq b(m) + \epsilon\big)}
		{\sP_0\big(b(m) \leq \beta \leq b(m) + \epsilon/2\big)}.
	\end{equation}
We remark that $\sP_0(m \leq \beta<b(m))=0$ and
$\sP_0(b(m) \leq \beta \leq b(m)+\eps/2)>0$ by definition of $b(m)$, so this
implies $\sP_{\lambda,\gamma}(m \leq \beta<b(m))=0$ and
	$\lim_{\gamma \rightarrow \infty}
		\sup_{\lambda/\gamma \in [a(m),b(m)]}
		\sP_{\lambda,\gamma}(\beta \geq b(m)+\epsilon) = 0$.
	Likewise, $\sP_{\lambda,\gamma}(a(m)<\beta \leq m)=0$ and
	$\lim_{\gamma \rightarrow \infty}
		\sup_{\lambda/\gamma \in [a(m),b(m)]}
		\sP_{\lambda,\gamma}(\beta \leq a(m)-\epsilon) = 0$,
so \eqref{eq:off-supp-large-gamma} follows.

	Finally, we prove \eqref{eq:small-gamma}. Suppose $\gamma<0$.
	If $\lambda/\gamma \leq 0$,
	then $\inf_{\beta \in [M-\epsilon/2,M]} (\beta - \lambda / \gamma)^2 - \sup_{|\beta| \leq M-\epsilon}(\beta - \lambda / \gamma)^2 \geq (M-\epsilon)\epsilon$.
	Thus,
	\begin{equation}
		\sP_{\lambda,\gamma}\big(|\beta|\leq M - \epsilon\big)
		\leq 
		\frac{\sP_{\lambda,\gamma}\big(|\beta|\leq M -
\epsilon\big)}{\sP_{\lambda,\gamma}\big(\beta \in [M - \epsilon/2,M]\big)}
		\leq 
		e^{(\gamma/2)(M-\epsilon)\epsilon}
		\frac{\sP_0\big(|\beta|\leq M - \epsilon\big)}{\sP_0\big(\beta \in [M - \epsilon/2,M]\big)}.
	\end{equation}
Here, $\sP_0(\beta \in [M-\eps/2,M])>0$ since $M \in \supp(\sP_0)$. Then
$\lim_{\gamma \to -\infty}
\sup_{\lambda:\lambda/\gamma \leq 0} \sP_{\lambda,\gamma}(|\beta| \leq
M-\eps)=0$. The case $\lambda/\gamma \geq 0$ is analogous,
and this shows \eqref{eq:small-gamma}.
\end{proof}

We conclude this section with the following consequence of the above analyses,
which we will use later in the main proofs.

\begin{lemma}\label{lemma:hboundary}
The function ${-}\hs:\R^2 \to [0,\infty]$ is lower semi-continuous and convex.
For each $(\bar m,\bar s)$ on the boundary of $\Gamma$,
if ${-}\hs(\bar m,\bar s)<\infty$, then
there exists a smooth path $\{(m^t,s^t)\}_{t \in (0,1)}$ in $\Gamma$ such that 
\begin{equation}\label{eq:hboundary}
\lim_{t \to 0} (m^t,s^t)=(\bar m,\bar s), \qquad
\lim_{t \to 0} \frac{{-}\hs(\bar m,\bar s)+\hs(m^t,s^t)}{\|(\bar m,\bar s)
-(m^t,s^t)\|_2}=\infty.
\end{equation}
\end{lemma}
\begin{proof}
We may again first center $\supp(\sP_0)$: Indeed, let $\tilde \sP_0$ denote
the law of $\beta+c$ for any $c \in \R$, with corresponding domain $\tilde
\Gamma$ and relative entropy function $-\tilde \hs(m,s)$.
If $\{(\tilde m^t,\tilde s^t)\}_{t \in (0,1)}$ is a
path in $\tilde \Gamma$, then (\ref{eq:supportrecenter}) implies that
$(m^t,s^t)=(\tilde m^t-c,\tilde s^t-2\tilde m^tc+c^2)$ is a corresponding
path in $\Gamma$ for which ${-}\tilde \hs(\tilde m^t,\tilde
s^t)={-}\hs(m^t,s^t)$. Hence the result for $\tilde \sP_0$ implies
that for $\sP_0$.

Thus, assume that $\sP_0$ is centered so that
$-a(\sP_0)=b(\sP_0)=M>0$, and $\Gamma$ takes the form
(\ref{eq:Gammacentered}). Let us first show that for each point $(\bar m,\bar
s)$ on the boundary $\partial \Gamma$, there exists a path $\{(m^t,s^t)\}_{t \in
(0,1)}$ in (the interior of) $\Gamma$ for which
\begin{equation}\label{eq:hderivboundary}
\lim_{t \to 0} (m^t,s^t)=(\bar m,\bar s), \qquad
\lim_{t \to 0} \frac{\de}{\de t}[{-}\hs(m^t,s^t)]=-\infty.
\end{equation}
Consider first the case $\bar m \in (-M,M)$ and $\bar s=M^2$.
We choose $m^t=\bar m$ and $s^t=M^2-t\eps$ for some small $\eps>0$.
Then, differentiating (\ref{eq:hdef}) by the envelope theorem,
$\frac{\de}{\de t}[{-}\hs(m^t,s^t)]=\frac{\eps}{2}\gamma(m^t,s^t)$.
For each $\gamma \in \R$ and $m \in (-M,M)$, let $\lambda_\gamma(m)$ be
the value in Proposition \ref{prop:oneparamexpfam} such that
$m=\langle \beta \rangle_{\lambda_\gamma(m),\gamma}$. Then,
denoting $\gamma^t=\gamma(m^t,s^t)$, we must have
$(\lambda(m^t,s^t),\gamma(m^t,s^t))=(\lambda_{\gamma^t}(\bar m),\gamma^t)$,
and Lemma \ref{lem:large-gamma-moments} implies that $\gamma^t \to -\infty$
as $s^t \to M^2$. Thus (\ref{eq:hderivboundary}) holds. For $\bar m \in (-M,M)$
and $\bar s=m(a(m)+b(m))-a(m)b(m)$, we may similarly choose
$m^t=\bar m$ and $s^t=m(a(m)+b(m))-a(m)b(m)+t\eps$ for some small $\eps>0$,
and the proof is analogous.

For $\bar m=-M$ and $\bar s=M^2$, fix $\gamma=0$ and choose
$m^t=-M+t\eps$ and $s^t=\langle \beta^2 \rangle_{\lambda_0(m^t),0}$
for some small $\eps>0$. Then
$\frac{\de}{\de t}[{-}\hs(m^t,s^t)]=\eps \lambda(m^t,s^t)=\eps \lambda_0(m^t)$.
Proposition \ref{prop:oneparamexpfam} shows $\lambda_0(m) \to
-\infty$ as $m \to -M$, so (\ref{eq:hderivboundary}) holds.
For $\bar m=M$ and $\bar s=M^2$, we may similarly choose $m^t=M-\eps t$ and
$s^t=\langle \beta^2 \rangle_{\lambda_0(m^t),0}$, and the proof is
analogous. This shows the existence of $\{(m^t,s^t)\}_{t \in (0,1)}$
satisfying (\ref{eq:hderivboundary}) for all $(\bar m, \bar s) \in \partial
\Gamma$.

We now show the second statement of (\ref{eq:hboundary})
using some convex analysis: First note that this statement holds trivially if
${-}\sh(\bar m,\bar s)=\infty$, since ${-}\sh(m^t,s^t)$ is finite for all
$(m^t,s^t) \in \Gamma$. Thus, let us assume henceforth ${-}\sh(\bar m,\bar
s)<\infty$. Next, note that ${-}\hs:\R^2
\to [0,\infty]$ defined by (\ref{eq:hdef}) is a supremum of linear functions,
and hence is lower semi-continuous and convex. Thus ${-}\hs$ is closed, and
by the Gale-Klee-Rockafellar Theorem
(c.f.\ \cite[Theorem 10.2]{rockafellar1997convex}), its restriction
to any closed line segment in $\domain({-}\hs)$ is continuous. 
If $(\bar m,\bar s) \in \partial \Gamma$ is such that $\bar m \in (-M,M)$ and
${-}\hs(\bar m,\bar s)<\infty$, then $(\bar m,\bar s) \in \domain({-}\hs)$ and
the path $(m^t,s^t)$ constructed above is linear, so this continuity implies
that ${-}\hs(\bar m,\bar s)=\lim_{t \to 0} {-}\hs(m^t,s^t)$. Then the desired
statement (\ref{eq:hboundary}) follows from integrating
(\ref{eq:hderivboundary}) to get
\[\frac{{-}\hs(\bar m,\bar s)+\hs(m^t,s^t)}{\|(\bar m,\bar s)
-(m^t,s^t)\|_2}=-\frac{\int_0^t \frac{\de}{\de\tau}[{-}\hs(m^\tau,s^\tau)]\,\de \tau}
{t\eps}\to \infty \quad \text{ as } \quad t \to 0.\]

Take $(\bar m,\bar s)=(-M,M^2)$. Consider the function
${-}\bar \hs:\R \to [0,\infty]$ defined by
\[{-}\bar \hs(m)=\sup_{\lambda \in \R} \lambda m-\log \E_{\beta \sim \sP_0}
[e^{\lambda \beta}].\]
This corresponds to restricting $\gamma=0$ in the supremum of (\ref{eq:hdef}),
and hence we have $-\bar \hs(m) \leq -\hs(m,s)$ for all $(m,s) \in \R^2$.
Here ${-}\bar \hs$ is also lower semi-continuous and convex.
If $-\hs(\bar m,\bar s)<\infty$, then also $-\bar \hs(m)<\infty$,
so $-\bar \hs$ is continuous on $[-M,M)$. Then,
letting $(m^t,s^t)=(-M+\eps t,\langle \beta^2 \rangle_{\lambda_0(m^t),0})$
be the path constructed above,
\[{-}\hs(\bar m,\bar s) \geq {-}\bar \hs(\bar m)=\lim_{t \to 0} {-}\bar \hs(m^t)
=\lim_{t \to 0} {-}\hs(m^t,s^t),\]
where the last equality holds because $\gamma=0$ in the supremum of
(\ref{eq:hdef}) defining ${-}\hs(m^t,s^t)$, by construction of $(m^t,s^t)$.
We have also $\|(\bar m,\bar s)-(m^t,s^t)\|_2 \leq t\eps \cdot \sqrt{1+(2M)^2}$,
from the fact that $(m^t,s^t) \in \Gamma$ and the shape of $\Gamma$ in
(\ref{eq:supportrecenter}). Then the desired statement (\ref{eq:hboundary})
follows again from integrating (\ref{eq:hderivboundary}) to get
\[\frac{{-}\hs(\bar m,\bar s)+\hs(m^t,s^t)}{\|(\bar m,\bar s)
-(m^t,s^t)\|_2} \geq -\frac{\int_0^t
\frac{\de}{\de\tau}[{-}\hs(m^\tau,s^\tau)]\,\de \tau}{t\eps \cdot
\sqrt{1+(2M)^2}} \to \infty \quad \text{ as } \quad t \to 0.\]
The proof for $(\bar m,\bar s)=(M,M^2)$ is analogous.
\end{proof}

\section{Global convexity in low SNR}\label{app:proof_global_convexity}

We prove Proposition \ref{prop:global_convexity} on the global convexity of
$\cF_\TAP$ for sufficiently large $\sigma^2$ or small $\delta$.

\begin{proof}[Proof of Proposition \ref{prop:global_convexity}]
The Hessian of $D_0(\bbm,\bs)$ in the expression (\ref{eqn:TAP}) for $\cF_\TAP$
has smallest eigenvalue $\lambda_{\min}(\nabla^2 D_0(\bbm,\bs))=\min_{j=1}^p
\lambda_{\min}(\nabla^2[-\sh(m_j,s_j)])$.
Differentiating $-\sh(m,s)$ using the envelope theorem, we have
\[\nabla[{-}\sh(m,s)]=(\lambda(m,s),-\tfrac{1}{2}\gamma(m,s)),
\qquad \nabla^2[{-}\sh(m,s)]=\frac{\partial
(\lambda(m,s),-\frac{1}{2}\gamma(m,s))}{\partial (m,s)}\]
The Jacobian $\frac{\partial (m,s)}{\partial (\lambda,-\frac{1}{2}\gamma)}$
of the moment map $(\lambda,-\frac{1}{2}\gamma) \mapsto (m,s)
=(\langle \beta \rangle_{\lambda,\gamma},\langle \beta^2
\rangle_{\lambda,\gamma})$ is the covariance matrix of
$(\beta,\beta^2)$ under the law $\sP_{\lambda,\gamma}$.
Denoting by $\Var_{\lambda,\gamma},\Cov_{\lambda,\gamma}$ the variance and
covariance under this law, this implies
\begin{equation}\label{eq:hesshinverse}
\nabla^2[{-}\sh(m,s)]^{-1}
=\begin{pmatrix} \Var_{\lambda(m,s),\gamma(m,s)}[\beta] &
\Cov_{\lambda(m,s),\gamma(m,s)}[\beta,\beta^2] \\
\Cov_{\lambda(m,s),\gamma(m,s)}[\beta,\beta^2] &
\Var_{\lambda(m,s),\gamma(m,s)}[\beta^2] \end{pmatrix}.
\end{equation}
This matrix has operator norm at most a constant depending only on the
size $M$ of the support of $\sP_0$, so
\begin{equation}\label{eq:hesshlowerbound}
\lambda_{\min}(D_0(\bbm,\bs)) \geq c(M)
\end{equation}
for a constant $c(M)>0$. For the remaining terms
\[
E(\bbm, \bs) = \frac{n}{2}\log 2\pi\sigma^2
+\frac{1}{2\sigma^2}\big\|\y-\X\bbm\big\|_2^2 +\frac{n}{2}\log \Big(1 + \frac{ S(\bs)-Q(\bbm)}{\sigma^2}\Big)
\]
of $\cF_\TAP$ in (\ref{eqn:TAP}),
denoting $V(\bbm,\bs)=\sigma^2+S(\bs)-Q(\bbm)$, we have
\[\nabla^2 E(\bbm, \bs)
=\begin{pmatrix} \frac{1}{\sigma^2} \X^\top \X - \frac{n}{p}
\frac{1}{V(\bbm,\bs)} \id_p - \frac{2n}{p} \frac{\bbm \bbm^\top/p}{V(\bbm,\bs)}
& \frac{n}{p} \frac{\bbm\ones^\top/p}{V(\bbm,\bs)} \\
\frac{n}{p} \frac{\ones\bbm^\top/p}{V(\bbm,\bs)} &
-\frac{n}{2p} \frac{\ones\ones^\top/p}{V(\bbm,\bs)} \end{pmatrix}\]
Then, applying $V(\bbm,\bs) \geq \sigma^2$, $\|\bbm\bbm^\top/p\| \leq M^2$,
and $\|\ones\ones^\top/p\| \leq 1$, this shows
$\lambda_{\min}(\nabla^2 E(\bbm,\bs)) \geq -\frac{nC(M)}{p\sigma^2}$
for a constant $C(M)>0$. Thus, as long as
$(n/p)/\sigma^2<c(M)/2C(M)$, we have
$\lambda_{\min}(\cF_\TAP(\bbm,\bs)) \geq c(M)/2$, as desired.
\end{proof}

\section{Cavity field approximations of the posterior marginals}\label{sec:cavity}

We prove in this appendix the following result, which formalizes an
approximation for the marginal posterior law of any variable $\beta_j$ by
a univariate cavity field variable $\hat{\lambda}_j$. Recall the mean $\langle
\cdot \rangle_{\lambda,\gamma}$ under the exponential family law
(\ref{eq:scalarmean}), and recall also
$\gamma_\stat=\argmin_{\gamma>0} \phi(\gamma)$ which is uniquely defined under
Assumption \ref{ass:uniquemin}.

For a fixed coordinate $j \in [p]$, write $\x^j \in \R^n$ and $\X^- \in \R^{n \times (p-1)}$ for the $j^\text{th}$
and all-but-$j^\text{th}$ columns of $\X$. Write
$\beta_j,\beta_{0,j}$ and $\bbeta^-,\bbeta_0^- \in \R^{p-1}$
for the $j^\text{th}$ and all-but-$j^\text{th}$ coordinates of
$\bbeta,\bbeta_0$, and denote $\y^-=\X^-\bbeta_0^-+\beps$.
Define the leave-one-out posterior measure and associated posterior mean
\begin{equation}\label{eq:LOOposterior}
\sP^-(\de\bbeta^-)=\frac{e^{-\frac{1}{2\sigma^2}\|\X^-\bbeta^--\y^-\|^2}
\prod_{k \neq j} \sP_0(\de\beta_k)}
{\int e^{-\frac{1}{2\sigma^2}\|\X^-\bbeta^--\y^-\|^2}\prod_{k \neq j}
\sP_0(\de\beta_k)},
\qquad \langle f(\bbeta^-) \rangle_- = \int f(\bbeta^-) \sP^-(\de\bbeta^-). 
\end{equation}
Then we define a random variable $\hat{\lambda}_j$ as:
\begin{equation}\label{eq:hatlambdadef}
\hat{\lambda}_j=\gamma_\stat \beta_{0,j}
-\frac{1}{\sigma^2}{\x^j}^\top(\X^-\langle \bbeta^- \rangle_--\y^-).
\end{equation}
\begin{lemma}\label{lemma:cavity}
Suppose Assumptions \ref{ass:Bayesian_linear_model} and \ref{ass:uniquemin}
hold. For any $j \in
\{1,\ldots,p\}$ and
for any bounded function $f:\supp(\sP_0) \to
\R$,
\[\lim_{n,p \to \infty} \E\Big[\Big(\langle f(\beta_j) \rangle_{\X,\y}-\langle
f(\beta)\rangle_{\hat\lambda_j,\gamma_\stat}\Big)^2\Big]=0.\]
\end{lemma}

In the remainder of this section, we prove Lemma \ref{lemma:cavity} using
an interpolation argument, adapted from the argument of \cite[Theorem
1.7.11]{talagrand2010mean} for the Sherrington-Kirkpatrick model.

\subsection{Overlap concentration}

Let us introduce the shorthands
\[\v_0=\X\bbeta_0, \quad \v=\X\bbeta, \quad
\v_k=\X\bbeta_k,~~~ k \ge 1, \]
where $\bbeta_0$ is the true parameter and
$\bbeta,\bbeta_1,\bbeta_2,\ldots$ denote independent samples (replicas) from the
posterior distribution of $\bbeta$ given $(\X,\y)$. We write as shorthand
$\langle \cdot \rangle=\langle \cdot \rangle_{\X,\y}$ for the joint posterior
expectation over $\bbeta,\bbeta_1,\bbeta_2,\ldots$ fixing $\X,\bbeta_0,\y$.
We denote
\[\bar\v=\v-\v_0, \quad \dot\v=\v-\langle \v \rangle\]
and similarly $\bar\v_k=\v_k-\v_0$, $\dot\v_k=\v_k-\langle \v \rangle$
for $k \ge 1$. We first show the following
concentration-of-overlaps result (in the space of $\v$ rather than $\bbeta$),
which will be needed for the later interpolation argument.

\begin{lemma}\label{lemma:overlaps}
Under Assumptions \ref{ass:Bayesian_linear_model} and \ref{ass:uniquemin},
set $\rho=-\gamma_\stat\sigma^4+\delta\sigma^2$. Then
\begin{align}
&\lim_{n,p \to \infty}
\E\Big\langle\Big(p^{-1}\dot \v_1^\top \dot\v_2 \Big)^2\Big\rangle
=0,\label{eq:overlapconc1}\\
&\lim_{n,p \to \infty}
\E\Big\langle\Big(p^{-1}\|\dot\v\|^2-\rho\Big)^2\Big\rangle
=0.\label{eq:overlapconc2}
\end{align}
\end{lemma}

In this section, let us
denote $\omega=\sigma^{-2}$ and write
$\y(\omega)=\X\bbeta_0+\omega^{-1/2}\z$, where
$z_i \overset{iid}{\sim} \normal(0,1)$. We write $\E=\E_{\X,\bbeta_0,\z}$ and
$\sP_0(\de\bbeta)=\prod_{j=1}^p \sP_0(\de\beta_j)$. We denote the evidence
(\ref{eq:evidence}) as a function of $\omega$ by
\begin{align*}
F_p(\omega)&=\frac{1}{p}\log \int (\omega/2\pi)^{n/2}
\exp\left(-\frac{\omega}{2} \|\X\bbeta-\y(\omega)\|^2\right)\sP_0(\de\bbeta)\\
&=\frac{n}{2p}\log (\omega/2\pi)+\frac{1}{p}\log \int 
\exp\left(-\frac{\omega}{2} \|\X(\bbeta-\bbeta_0)\|^2
+\omega^{1/2}\z^\top\X(\bbeta-\bbeta_0)
-\frac{1}{2}\|\z\|^2 \right)\sP_0(\de \bbeta). 
\end{align*}
Recalling $\bar\v=\v-\v_0=\X(\bbeta-\bbeta_0)$, let us define
\[L(\bar\v):=-\frac{1}{2}\|\bar{\v}\|^2+\frac{1}{2\omega^{1/2}} \z^\top \bar\v\]
so that $F_p'(\omega)=(n/2p)\omega^{-1}+p^{-1}\langle L(\bar\v) \rangle$.
We emphasize that the definitions of $L(\bar\v)$ and the posterior average
$\langle \cdot \rangle$ depend on $\omega$, and we will write $L_\omega(\bar\v)$
and $\langle \cdot \rangle_\omega$ if we wish to make this dependence explicit.

\begin{lemma}\label{lemma:Lconc}
Suppose Assumptions \ref{ass:Bayesian_linear_model} and \ref{ass:uniquemin}
hold. Then
\begin{align}
\lim_{n,p \to \infty} p^{-1}\,\E\big|\langle L(\bar\v) \rangle-
\E\langle L(\bar\v) \rangle\big|&=0,\label{eq:Lconc1}\\
\lim_{n,p \to \infty} p^{-1}\,\E\big\langle \big|L(\bar\v)-\langle
L(\bar\v) \rangle\big|\big\rangle&=0.\label{eq:Lconc2}
\end{align}
\end{lemma}
\begin{proof}
The argument parallels that of \cite[Section 5]{lelarge2019fundamental} in the
spiked matrix model, which in turn adapts the proof of the Ghirlanda-Guerra
identities in \cite{panchenko2013sherrington}.

To show (\ref{eq:Lconc1}), observe that
for any function $f(\bar\v,\omega)$ that is differentiable in $\omega$, we have
\begin{equation}\label{eq:derivativerule}
\frac{d}{d\omega} \langle f(\bar\v,\omega) \rangle_\omega
=\langle \partial_\omega f(\bar\v,\omega) \rangle
+\langle f(\bar\v,\omega)L(\bar\v) \rangle
-\langle f(\bar\v,\omega) \rangle \langle L(\bar\v) \rangle. 
\end{equation}
Noting that $F_p'(\omega)=(n/2p)\omega^{-1}+p^{-1}\langle L_\omega(\bar\v) \rangle$
and applying this with $f(\bar\v,\omega)=p^{-1}L_\omega(\bar\v)$,
\begin{equation}\label{eq:Fpp}
F_p''(\omega)=-\frac{n}{2p}\omega^{-2}+\frac{d}{d\omega}\,p^{-1}\langle L_\omega(\bar\v) \rangle
=-\frac{n}{2p}\omega^{-2}+p^{-1}\Big(
{-}\frac{1}{4\omega^{3/2}}\z^\top \langle \bar\v \rangle+
\big\langle L(\bar\v)^2 \big\rangle-\langle L(\bar\v) \rangle^2
\Big).
\end{equation}
Then, applying $|\z^\top \bar\v| \leq \|\z\| \cdot \|\bar\v\| \leq 
\|\X\|_{\op} \|\z\| \cdot
(\|\bbeta\|+\|\bbeta_0\|) \leq 2K\sqrt{p} \cdot \|\X\|_{\op} \|\z\|$, the
function
\begin{equation}\label{eq:convexification}
G_p(\omega):=F_p(\omega)-\frac{n}{2p}\log \omega-\omega^{1/2} \cdot \frac{2K
\cdot \|\X\|_{\op}\|\z\|}{\sqrt{p}}
\end{equation}
satisfies $G_p''(\omega) \geq 0$. It may be checked that
$\E G_p(\omega)$ is differentiable in
$\omega$ with derivative $\E G_p'(\omega)$, by the dominated convergence
theorem. Then both $G_p(\omega)$ and $\E G_p(\omega)$ are differentiable and
convex, so by \cite[Lemma 3.2]{panchenko2013sherrington},
\[|G_p'(\omega)-\E G_p'(\omega)| \leq
\E G_p'(\omega+\eps)-\E G_p'(\omega-\eps)+\Delta_p/\eps\]
for $\Delta_p=\sum_{x \in \{\omega-\eps,\omega,\omega+\eps\}} |G_p(x)-\E G_p(x)|$
and any $\eps>0$. Then
\begin{align}
p^{-1}\,\E\big|\langle L(\bar\v) \rangle-
\E\langle L(\bar\v) \rangle\big|
&=\E\left|G_p'(\omega)-\E G_p'(\omega)+\frac{K(\|\X\|_{\op}
\|\z\|-\E[\|\X\|_{\op}\|\z\|])}{\omega^{1/2}\sqrt{p}}
\right|\nonumber\\
&\leq \E G_p'(\omega+\eps)-\E G_p'(\omega-\eps)
+\frac{\E\Delta_p}{\eps}+\frac{K\,\E\big|\|\X\|_{\op}\|\z\|-\E[\|\X\|_{\op}\|\z\|]\big|}{\omega^{1/2}\sqrt{p}}. \label{eq:Lvbound1}
\end{align}
We have $\Var[F_p(\omega)] \leq C/p$ from Theorem \ref{thm:freeenergy}, while
\[\Var[\|\X\|_\op \|\z\|]
=\E \|\X\|_\op^2 \E\|\z\|^2-(\E\|\X\|_\op \E\|\z\|)^2
=\Var[\|\X\|_\op] \cdot \E\|\z\|^2
+(\E\|\X\|_\op)^2 \cdot \Var[\|\z\|] \leq C.\]
Then, applying Cauchy-Schwarz and
taking the limit $n,p \to \infty$ for fixed $\eps,K,\omega$, the last two
terms on the right side of (\ref{eq:Lvbound1})
vanish. For the first term of (\ref{eq:Lvbound1}), let us denote
$F_\RS(\omega)=-\phi(\gamma_\stat(\omega);\omega)$, making explicit the
dependence of both $\phi(\cdot)$ and its minimizer $\gamma_\stat$ on $\omega$.
Then applying Theorem \ref{thm:freeenergy} and the convergence
$\E\|\X\|_\op\|\z\|/\sqrt{p} \to 1+\sqrt{\delta}$, we have
\begin{equation}\label{eq:Gomega}
\E G_p(\omega) \to G(\omega):=F_\RS(\omega)-(\delta/2)\log \omega
-2K(1+\sqrt{\delta})\omega^{1/2}
\end{equation}
for each fixed $\omega$. This implies by convexity of $\E G_p(\omega)$ that
\begin{equation}\label{eq:derivativeconv}
\E G_p'(\omega) \to G'(\omega)=F_\RS'(\omega)-(\delta/2)\omega^{-1}-K(1+\sqrt{\delta})\omega^{-1/2}
\end{equation}
as long as
$F_\RS$ is differentiable at $\omega$. Under Assumption \ref{ass:uniquemin},
since $\phi''(\gamma_\stat;\omega)>0$ strictly, the implicit function theorem
implies that both $\gamma_\stat(\omega)=\argmin_{\gamma} \phi(\gamma;\omega)$
and $F_\RS(\omega)=-\phi(\gamma_\stat(\omega);\omega)$ are
indeed continuously differentiable in an open neighborhood of $\omega$.
Then, applying (\ref{eq:derivativeconv}) to the first term of
(\ref{eq:Lvbound1}), for all sufficiently small $\eps$,
\[\limsup_{n,p \to \infty} p^{-1}\,\E\big|\langle L(\bar\v) \rangle-
\E\langle L(\bar\v) \rangle\big|
\leq G'(\omega+\eps)-G'(\omega-\eps).\]
Taking the limit $\eps \to 0$ and applying continuous differentiability of
$G(\omega)$ shows (\ref{eq:Lconc1}).

To show (\ref{eq:Lconc2}), define
$d_p(\bar\v,\omega)=|L_\omega(\bar\v)-\langle L_\omega(\bar\v) \rangle_\omega|$
and
\[D_p(\omega)=\frac{1}{p}\langle d_p(\bar\v,\omega) \rangle_\omega, \qquad
V_p(\omega)=\frac{1}{p} \langle d_p(\bar\v,\omega)^2 \rangle_\omega
=\frac{1}{p}\Big(\langle L(\bar\v)^2 \rangle-\langle L(\bar\v) \rangle^2 \Big)\]
so that (\ref{eq:Lconc2}) is the statement $\E D_p(\omega) \to 0$.
First, by (\ref{eq:Fpp}), (\ref{eq:convexification}),
and the bound $|\z^\top \bar\v| \leq 2K\sqrt{p} \cdot \|\X\|_\op\|\z\|$,
\[V_p(\omega) \leq F_p''(\omega)+\frac{n}{2p}\omega^{-2}
+\frac{K\|\X\|_\op\|\z\|}{2\omega^{3/2}
\sqrt{p}} =G_p''(\omega).\]
Then fixing $\eps>0$
such that $G(\omega)$ in (\ref{eq:Gomega}) is continuously
differentiable in a neighborhood of $(\omega,\omega+\eps)$,
\begin{equation}\label{eq:smoothedVbound}
\limsup_{n,p \to \infty} \E \int_\omega^{\omega+\eps} V_p(x)\de x
\leq \limsup_{n,p \to \infty}\Big(\E G_p'(\omega+\eps)-\E G_p'(\omega)\Big)
=G'(\omega+\eps)-G'(\omega). 
\end{equation}

Next, applying again (\ref{eq:Fpp}) and $|\z^\top \bar\v| \leq 2K\sqrt{p}
\cdot \|\X\|_\op\|\z\|$, observe that
\begin{align*}
&\limsup_{x \to \omega} \sup_{\bar\v \in [-2K,2K]^p}
\frac{|p^{-1}d_p(\bar\v,x)
-p^{-1}d_p(\bar\v,\omega)|}{|x-\omega|}\\
&\leq \limsup_{x \to \omega} \sup_{\bar\v \in [-2K,2K]^p}
\frac{|p^{-1}[L_x(\bar\v)-\langle L_x(\bar\v)
\rangle_x]-p^{-1}[L_\omega(\bar\v)-\langle L_\omega(\bar\v) \rangle_\omega]|}
{|x-\omega|}\\
&\leq \left|\frac{d}{d\omega} \frac{1}{p}\,\frac{1}{2\omega^{1/2}}\right|
\cdot \sup_{\bar\v \in [-2K,2K]^p} |\z^\top \bar\v|
+\left|\frac{d}{d\omega} \frac{1}{p} \langle L_\omega(\bar\v)
\rangle_\omega \right| \leq V_p(\omega)+\frac{K\|\X\|_\op\|\z\|}{\omega^{3/2}\sqrt{p}}. 
\end{align*}
For any function $f(\bar\v,\omega)$, possibly non-differentiable in
$\omega$, we have analogously to (\ref{eq:derivativerule})
\begin{align*}
&\limsup_{x \to \omega} \frac{|\langle f(\bar \v,x) \rangle_x
-\langle f(\bar \v,\omega) \rangle_\omega|}{|x-\omega|}\\
&\leq \limsup_{x \to \omega} \frac{|\langle f(\bar \v,x) \rangle_x
-\langle f(\bar \v,\omega) \rangle_x|}{|x-\omega|}
+\limsup_{x \to \omega} \frac{|\langle f(\bar \v,\omega) \rangle_x
-\langle f(\bar \v,\omega) \rangle_\omega|}{|x-\omega|}\\
&\leq \limsup_{x \to \omega} \sup_{\bar\v \in [-2K,2K]^p}
\frac{|f(\bar \v,x)-f(\bar \v,\omega)|}{|x-\omega|}
+\Big|\langle f(\bar\v,\omega)L(\bar\v) \rangle_\omega
-\langle f(\bar\v,\omega)\rangle_\omega \langle L(\bar\v) \rangle_\omega\Big|. 
\end{align*}
Applying this with
$f(\bar\v,\omega)=p^{-1}d_p(\bar\v,\omega)$, we obtain
\begin{align*}
\limsup_{x \to \omega} \frac{|D_p(x)-D_p(\omega)|}
{|x-\omega|} &\leq V_p(\omega)+\frac{K\|\X\|_\op\|\z\|}{\omega^{3/2}\sqrt{p}}
+\frac{1}{p}\Big|\langle d_p(\bar{\v},\omega) L(\bar\v) \rangle
-\langle d_p(\bar{\v},\omega)\rangle \langle L(\bar\v) \rangle\Big|\\
&\leq 2V_p(\omega)+\frac{K\|\X\|_\op\|\z\|}{\omega^{3/2}\sqrt{p}}.
\end{align*}
For any $\omega' \in (\omega,\omega+\eps)$, integrating this bound from
$x=\omega$ to $x=\omega'$ implies
\[D_p(\omega) \leq D_p(\omega')
+\int_\omega^{\omega'}
\left(2V_p(x)+\frac{K\|\X\|_\op\|\z\|}{x^{3/2}\sqrt{p}}\right)\de x
\leq D_p(\omega')
+\int_\omega^{\omega+\eps}
\left(2V_p(x)+\frac{K\|\X\|_\op\|\z\|}{x^{3/2}\sqrt{p}}\right)\de x. \]
Then integrating a second time from $\omega'=\omega$ to $\omega'=\omega+\eps$
implies
\[\eps\,D_p(\omega) \leq \int_\omega^{\omega+\eps} D_p(\omega')
\de \omega'+\eps \int_\omega^{\omega+\eps}
\left(2V_p(x)+\frac{K\|\X\|_\op\|\z\|}{x^{3/2}\sqrt{p}}\right)\de x. \]
Taking expectations on both sides, and applying the Cauchy-Schwarz inequalities
$\E \int_\omega^{\omega+\eps} D_p(\omega')\de\omega'
\leq (\eps\,\E \int_\omega^{\omega+\eps} D_p(\omega')^2\de\omega')^{1/2}$ and
$D_p(\omega')^2 \leq p^{-1} V_p(\omega')$, we get
\[\eps\cdot
\E D_p(\omega) \leq \left(\frac{\eps}{p} \cdot \E\int_\omega^{\omega+\eps}
V_p(\omega')\de\omega'\right)^{1/2}+2\eps \cdot \E\int_\omega^{\omega+\eps}
V_p(x)\de x+C\eps^2\]
for a constant $C:=C(\omega,K,\delta)>0$. Then, dividing by $\eps$ and
applying (\ref{eq:smoothedVbound}),
\[\limsup_{n,p \to \infty} \E D_p(\omega)
\leq 2\Big(G'(\omega+\eps)-G'(\omega)\Big)+C\eps. \]
Taking $\eps \to 0$ and recalling that $G(\omega)$ is continuously
differentiable in an open neighborhood of $\omega$,
we get $\E D_p(\omega) \to 0$ which shows (\ref{eq:Lconc2}).
\end{proof}

\begin{proof}[Proof of Lemma \ref{lemma:overlaps}]
There exist constants $C,c>0$ depending on $(\omega,\delta,K)$
such that $\P[\|\X\|_\op>C+t]<e^{-c(1+t)p}$ for any $t \geq 0$
(c.f.\ \cite[Theorem
4.4.5]{vershynin2018high}) and $\E[\sup_{\bar \v \in [-2K,2K]^p}
L(\bar\v)^2] \leq Cp^2$ for all large $n,p$.
Then applying Lemma \ref{lemma:Lconc} on the event $\|\X\|_\op \leq C$
and Cauchy-Schwarz on the event $\|\X\|_\op>C$, we also have
\begin{align*}
\lim_{n,p \to \infty} p^{-1} \E\left[\big|\langle L(\bar\v) \rangle
-\E \langle L(\bar\v)\rangle\big| \cdot \|\X\|_\op^2\right]&=0, \\
\lim_{n,p \to \infty} p^{-1} \E\left[\big\langle\big|\langle L(\bar\v) \rangle
-\E \langle L(\bar\v)\rangle\big| \big\rangle\cdot \|\X\|_\op^2\right]&=0. 
\end{align*}
Recall that $\bar\v=\v-\v_0$, and let $\v_1,\v_2$ denote independent replicas
of $\v=\X\bbeta$ under the posterior measure. Then, applying the above and
the bounds $|\v^\top \v_0|,|\v_1^\top \v_2| \leq K^2p \cdot \|\X\|_\op^2$,
we have
\begin{align}
p^{-2}\Big(\E\big\langle L(\bar\v) \cdot \v^\top \v_0 \big\rangle
-\E\langle L(\bar\v) \rangle \langle \v^\top \v_0 \rangle\Big)
&=p^{-2}\,\E\Big\langle \big[L(\bar\v)-\langle L(\bar\v)\rangle\big] \cdot
\v^\top \v_0 \Big\rangle \to 0, \label{eq:Lvconc1}\\
p^{-2}\Big(\E\langle L(\bar\v) \cdot \v^\top \v_0
\rangle-\E\langle L(\bar\v) \rangle \E \langle \v^\top \v_0 \rangle
\Big)
&=p^{-2}\,\E\Big\langle \big[L(\bar\v)-\E\langle L(\bar\v) \rangle\big] \cdot
\v^\top \v_0 \Big\rangle \to 0, \label{eq:Lvconc2}\\
p^{-2}\Big(\E\langle L(\bar\v) \rangle \langle \v_1^\top \v_2
\rangle-\E\langle L(\bar\v) \rangle \E \langle \v_1^\top \v_2 \rangle
\Big)
&=p^{-2}\,\E\Big\langle \big[\langle L(\bar\v) \rangle-\E \langle L(\bar\v)
\rangle\big] \cdot \v_1^\top \v_2 \Big \rangle \to 0. \label{eq:Lvconc3}
\end{align}
The lemma will follow from deriving implications of these three statements.

By the tower property of conditional expectation (i.e.\ the Nishimori identity),
\[\E\langle f(\v_1,\ldots,\v_k,\v_0) \rangle=
\E\langle f(\v_1,\ldots,\v_k,\v_{k+1}) \rangle\]
so we may treat $\v_0$ corresponding to the true signal as an additional
replica inside $\E\langle \cdot \rangle$. Let us introduce the notations
\[S_{jk}=\frac{1}{p}\E \langle \v_j^\top \v_k \rangle,
\qquad S_{ijkl}=\frac{1}{p^2} \E \langle \v_i^\top \v_j \cdot \v_k^\top \v_l
\rangle. \]
Observe that for any function $f:\R^p \times \R^p \to \R^p$, by Stein's lemma,
\begin{equation}
\E \z^\top \langle f(\v_0,\v) \rangle
=\E\sum_{j=1}^p \partial_{z_j} \langle f_j(\v_0,\v) \rangle
=\omega^{1/2}\E\Big[\langle \bar\v^\top f(\v_0,\v) \rangle-
\langle \bar\v \rangle^\top \langle f(\v_0,\v) \rangle\Big]\label{eq:stein1}
\end{equation}
and similarly for any $g:(\R^p)^{k+1} \to \R$,
\begin{align}
&\E \z^\top \langle f(\v_0,\v) \rangle \langle g(\v_0,\v_1,\ldots,\v_k) \rangle
=\omega^{1/2}\E\bigg[\langle \bar\v^\top f(\v_0,\v) \rangle \langle
g(\v_0,\v_1,\ldots,\v_k) \rangle \notag\\
&\hspace{1in}+\sum_{i=1}^k \langle f(\v_0,\v) \rangle^\top \langle \bar \v_i
g(\v_0,\v_1,\ldots,\v_k) \rangle
-(k+1)\langle \bar \v \rangle^\top \langle f(\v_0,\v) \rangle
\langle g(\v_0,\v_1,\ldots,\v_k) \rangle \bigg]. \label{eq:stein2}
\end{align}
Then, applying (\ref{eq:stein1}),
\begin{align*}
\E \big\langle L(\bar\v) \cdot \v^\top \v_0 \big\rangle
&=-\frac{1}{2}\E\big\langle \|\bar\v\|^2 \cdot \v^\top \v_0\big\rangle
+\frac{1}{2\omega^{1/2}}\E \z^\top\big\langle \bar{\v} \cdot \v^\top \v_0
\big\rangle\\
&=-\frac{1}{2}\E \langle \bar\v \rangle^\top \langle \bar\v \cdot
\v^\top \v_0 \rangle
=-\frac{1}{2}\E \langle (\v_1-\v_0)^\top(\v_2-\v_0) \cdot \v_2^\top \v_0. 
\rangle
\end{align*}
Expanding this product and using symmetry between replicas gives
\begin{equation}\label{eq:Lv1}
p^{-2} \E\big\langle L(\bar\v) \cdot \v^\top \v_0 \big\rangle
=\frac{1}{2}({-}S_{1220}+S_{1020}+S_{0220}-S_{0020})
=\frac{1}{2}({-}S_{0001}+S_{0101}). 
\end{equation}
Similarly, applying (\ref{eq:stein2}) with $k=1$,
\begin{align*}
\E\langle L(\bar\v)\rangle \langle \v^\top \v_0 \rangle
&=-\frac{1}{2}\E\langle \|\bar\v\|^2 \rangle \langle \v^\top \v_0 \rangle
+\frac{1}{2\omega^{1/2}} \E \z^\top\langle \bar\v \rangle\langle \v^\top \v_0
\rangle\\
&=\frac{1}{2}\E\Big[\langle \bar\v \rangle^\top \langle \bar\v \cdot \v^\top \v_0
\rangle-2\|\langle \bar\v \rangle\|^2\langle \v^\top \v_0 \rangle \Big]\\
&=\frac{1}{2}\E\langle (\v_1-\v_0)^\top (\v_2-\v_0) \cdot \v_2^\top \v_0
-2 (\v_1-\v_0)^\top(\v_2-\v_0) \cdot \v_3^\top \v_0 \rangle
\end{align*}
and expanding and simplifying using symmetry between replicas yields
\begin{equation}\label{eq:Lv2}
p^{-2}\E\langle L(\bar\v)\rangle \langle \v^\top \v_0 \rangle
=\frac{1}{2}({-}S_{0001}-S_{0101}+4S_{0102}-2S_{0123}).
\end{equation}
We deduce from (\ref{eq:Lvconc1}), (\ref{eq:Lv1}), and (\ref{eq:Lv2}) that
\begin{equation}\label{eq:Sid1}
S_{0101}-2S_{0102}+S_{0123} \to 0.
\end{equation}

Next, applying (\ref{eq:stein1}), observe that
\[\E\langle L(\bar\v) \rangle=-\frac{1}{2}\E\langle \|\bar\v\|^2 \rangle
+\frac{1}{2\omega^{1/2}} \E \z^\top \langle \bar\v \rangle
=-\frac{1}{2}\E\|\langle \bar\v \rangle\|^2
=-\frac{1}{2}\E\langle (\v_1-\v_0)^\top(\v_2-\v_0) \rangle. \]
Then $p^{-1}\E\langle L(\bar\v)\rangle=({-}S_{00}+S_{01})/2$ so
\begin{equation}\label{eq:Lv3}
p^{-2}\,\E\langle L(\bar\v) \rangle \E\langle \v^\top \v_0 \rangle
=p^{-2}\,\E\langle L(\bar\v) \rangle \E\langle \v_1^\top \v_2 \rangle
=\frac{1}{2}({-}S_{00}S_{01}+S_{01}^2). 
\end{equation}
Note that for any $j,k$, we have
\begin{align}
|S_{00jk}-S_{00}S_{jk}|&=p^{-2}\Big|\E\Big\langle(\|\v_0\|^2-\E
\langle \|\v_0\|^2 \rangle)\v_j^\top \v_k\Big\rangle\Big|\nonumber\\
&\leq K^2p^{-1}\E\Big[\Big|\|\v_0\|^2-\E \|\v_0\|^2\Big| \cdot \|\X\|_\op^2\Big]
\leq K^2\E[\|\X\|_\op^4]^{1/2}
\cdot \Var[p^{-1}\|\v_0\|^2]^{1/2} \to 0. \label{eq:Sid00}
\end{align}
We then deduce from (\ref{eq:Lvconc2}),
(\ref{eq:Lv1}), (\ref{eq:Lv3}), and the consequence $S_{0001}-S_{00}S_{01} \to
0$ of (\ref{eq:Sid00}) that
\begin{equation}\label{eq:Sid2}
S_{0101}-S_{01}^2 \to 0.
\end{equation}

Finally, applying (\ref{eq:stein2}) with $k=2$, observe that
\begin{align*}
\E\langle L(\bar\v) \rangle \langle \v_1^\top \v_2 \rangle
&=-\frac{1}{2}\E\langle \|\bar\v\|^2 \rangle \langle \v_1^\top \v_2 \rangle
+\frac{1}{2\omega^{1/2}} \E\z^\top \langle \bar\v \rangle \langle \v_1^\top \v_2
\rangle\\
&=\frac{1}{2}\E\Big[\langle \bar\v \rangle^\top \langle \bar\v_1 \cdot \v_1^\top
\v_2 \rangle+\langle \bar\v \rangle^\top \langle \bar\v_2 \cdot \v_1^\top
\v_2 \rangle-3\|\langle \bar\v \rangle\|^2\langle \v_1^\top \v_2 \rangle\Big]\\
&=\frac{1}{2}\E\big\langle 2(\v_1-\v_0)^\top(\v_2-\v_0)\v_2^\top \v_3
-3(\v_1-\v_0)^\top(\v_2-\v_0) \v_3^\top \v_4\big\rangle
\end{align*}
and simplifying gives
\begin{equation}\label{eq:Lv4}
p^{-2}\E\langle L(\bar\v) \rangle \langle \v_1^\top \v_2 \rangle
=\frac{1}{2}({-}S_{0012}+S_{0123}).
\end{equation}
Applying $S_{0012}-S_{00}S_{12} \to 0$ from (\ref{eq:Sid00}), where
$S_{00}S_{12}=S_{00}S_{01}$ by symmetry, we then deduce from (\ref{eq:Lvconc3}),
(\ref{eq:Lv3}), and (\ref{eq:Lv4}) that
\begin{equation}\label{eq:Sid3}
S_{0123}-S_{01}^2 \to 0.
\end{equation}
To summarize,
these conclusions (\ref{eq:Sid1}), (\ref{eq:Sid2}), (\ref{eq:Sid3}), and
(\ref{eq:Sid00}) show that, up to asymptotically vanishing errors,
$S_{0101},S_{0102},S_{0123}$ all coincide with $S_{01}^2$,
$S_{0000}$ coincides with $S_{00}^2$, and $S_{0001},S_{0012}$ both coincide with
$S_{00}S_{01}$.

Then, recalling $\dot\v=\v-\langle \v \rangle$, we have
\begin{align*}
p^{-2}\,\E \langle (\dot \v_1^\top \dot \v_2)^2 \rangle
&=p^{-2}\,\E \big\langle (\v_1-\v_3)^\top (\v_2-\v_4)
\cdot (\v_1-\v_5)^\top (\v_2-\v_6) \big\rangle\\
&=S_{0101}-2S_{0102}+S_{0123} \to 0
\end{align*}
which is (\ref{eq:overlapconc1}). Furthermore,
\begin{align}
&p^{-2}\E \langle \|\dot \v\|^4 \rangle
-\big(p^{-1}\E\langle \|\dot\v\|^2 \rangle\big)^2\nonumber\\
&=p^{-2}\E\big\langle (\v_1-\v_2)^\top (\v_1-\v_3) \cdot 
(\v_1-\v_4)^\top (\v_1-\v_5) \big\rangle
-\Big(p^{-1}\E\big\langle (\v_1-\v_2)^\top (\v_1-\v_3) \rangle\Big)^2\nonumber\\
&=\big(S_{0000}-4S_{0001}+2S_{0012}+4S_{0102}-3S_{0123}\big)
-\big(S_{00}-S_{01}\big)^2 \to 0. \label{eq:overlapconc2tmp}
\end{align}
Recall that $p^{-1}\E \langle L(\bar\v) \rangle=({-}S_{00}+S_{01})/2$, so we
have
\[p^{-1}\E \langle \|\dot \v\|^2 \rangle
=S_{00}-S_{01}=-2p^{-1}\E \langle L(\bar\v) \rangle
=-2\E F_p'(\omega)+(n/p)\omega^{-1}.\]
By (\ref{eq:derivativeconv}),
with $G_p(\omega)$ as defined in (\ref{eq:convexification}),
we have $\E F_p'(\omega) \to F_{\RS}'(\omega)$. Recalling
$F_\RS(\omega)=-\phi(\gamma_\stat(\omega);\omega)$
and differentiating (\ref{eqn:phi_potential}) by the envelope theorem,
$-F_{\RS}'(\omega)=\partial_\omega \phi(\gamma_\stat(\omega);\omega)
=-\gamma_\stat(\omega)/(2\omega^2)$.
Thus $p^{-1}\E\langle \|\dot\v\|^2 \rangle \to
{-}\gamma_\stat/\omega^2+\delta/\omega$, which is the quantity $\rho$
defined in the statement of the lemma. 
Applying this to (\ref{eq:overlapconc2tmp}) gives (\ref{eq:overlapconc2}).
\end{proof}

We pause to record here the following consequence of the above result, which we
will use in later proofs.

\begin{corollary}\label{cor:ymmseconc}
Suppose Assumptions \ref{ass:Bayesian_linear_model} and \ref{ass:uniquemin}
hold. Then as $n,p \to \infty$,
\[p^{-1} \E \|\y-\X\langle \bbeta \rangle\|_2^2 \to \gamma_\stat \sigma^4,
\qquad p^{-1}\|\y-\X\langle \bbeta \rangle\|_2^2 \gotop \gamma_\stat \sigma^4.\]
\end{corollary}
\begin{proof}
Note that $\y-\X\langle \bbeta \rangle=\omega^{-1/2}\z-\langle \bar\v \rangle$.
By the Nishimori identity and Lemma~\ref{lemma:overlaps},
\begin{align*}
\E[(p^{-1}\|\langle \bar\v \rangle\|_2^2-\rho)^2]
&=\E\left[\left(p^{-1}\|\X(\bbeta_0-\langle \bbeta
\rangle)\|_2^2-\rho\right)^2\right]\\
&=\E\left\<\left(p^{-1}\|\X(\bbeta-\langle \bbeta \rangle)\|_2^2-\rho\right)^2\right\> 
=\E\left\<\left(p^{-1}\|\dot \v\|_2^2-\rho\right)^2\right\> \to
0.
\end{align*}
This implies
\begin{equation}\label{eq:ymmseconc1}
p^{-1}\E\|\langle \bar \v \rangle\|_2^2 \to \rho, \qquad
p^{-1}\|\langle \bar \v \rangle\|_2^2 \gotop \rho.
\end{equation}

Similarly, by the Nishimori identity and Lemma~\ref{lemma:overlaps},
\[\E\left\<(p^{-1}\|\bar\v\|_2^2-2\rho)^2 \right\>
=\E\left\<(p^{-1}\|\X(\bbeta-\bbeta_0)\|_2^2-2\rho)^2\right\>
=\E\left\<(p^{-1}\|\dot \v_1-\dot \v_2\|_2^2-2\rho)^2\right\> \to 0\]
so this implies
\begin{equation}\label{eq:ymmseconc2}
p^{-1}\E\langle \|\bar\v\|_2^2 \rangle \to 2\rho, \qquad
p^{-1}\langle \|\bar\v\|_2^2 \rangle \gotop 2\rho.
\end{equation}
We recall from the
preceding proof that $-2p^{-1}\E \langle L(\bar\v) \rangle=p^{-1}\E\langle
\|\dot \v\|^2
\rangle \to \rho$. Then the first statement of Lemma \ref{lemma:Lconc} implies
$\E|p^{-1}\omega^{-1/2} \z^\top \langle \bar\v\rangle
-p^{-1}\langle \|\bar \v\|^2 \rangle+\rho| \to 0$,
which together with (\ref{eq:ymmseconc2}) shows
\begin{equation}\label{eq:ymmseconc3}
p^{-1}\E[\omega^{-1/2} \z^\top \langle \bar\v\rangle] \to \rho, \qquad
p^{-1}\,\omega^{-1/2} \z^\top \langle \bar\v\rangle \gotop \rho.
\end{equation}

Combining (\ref{eq:ymmseconc1}), (\ref{eq:ymmseconc3}), and
$p^{-1}\E\|\z\|_2^2 \to \omega^{-1}\delta$ and
$p^{-1}\omega^{-1}\|\z\|_2^2 \gotop \omega^{-1}\delta$, we obtain
$p^{-1}\|\y-\X\langle \bbeta \rangle\|_2^2
=p^{-1}\|\omega^{-1/2}\z-\langle \bar \v \rangle\|_2^2
\gotop \omega^{-1}\delta-\rho=\gamma_\stat/\omega^2=\gamma_\stat\sigma^4$,
and similarly for the convergence in expectation.
\end{proof}

\subsection{Cavity method interpolation}

We recall the notation $\dot\v=\X(\bbeta-\langle \bbeta \rangle)$ from the
preceding section, where we write as shorthand
$\langle \cdot \rangle=\langle \cdot \rangle_{\X,\y}$ for the posterior mean.

\begin{lemma}\label{lemma:cavityinterpolation}
Under Assumptions \ref{ass:Bayesian_linear_model} and \ref{ass:uniquemin}, set
$\rho=-\gamma_\stat \sigma^4+\delta \sigma^2$.
Let $V:\R \to \R$ be any smooth function such that $\E_{\xi \sim \normal(0,1)}
V(\rho^{1/2}\xi)=0$ and $|V(x)|,|V'(x)|,|V''(x)| \leq C e^{C|x|}$ for a
constant $C>0$ and all $x \in \R$. Let $\x \sim \normal(0,p^{-1}I) \in \R^n$ be
independent of $\X,\bbeta_0,\beps$. Then
\[\lim_{n,p \to \infty} \E\left[\big\langle V(\x^\top \dot\v)
\big\rangle^2\right]=0.\]
\end{lemma}
\begin{proof}
Let $\xi_1,\xi_2 \overset{iid}{\sim} \normal(0,1)$ be independent of
$\x,\X,\bbeta_0,\beps$, let $\v_1,\v_2$ denote two independent replicas of
$\v=\X\bbeta$ under the posterior measure, and denote
$\dot\v_k=\v_k-\langle \v \rangle$ for $k=1,2$. Introduce
\[S_k(t)=\sqrt{t}\,\x^\top \dot\v_k+\sqrt{1-t}\,\rho^{1/2}\xi_k\]
and the interpolation
\[\varphi(t)=\E\big\langle V(S_1(t))V(S_2(t)) \big \rangle\]
where $\E$ is over $\x,\X,\bbeta_0,\beps,\xi_1,\xi_2$.
The quantity we wish to bound is
$\varphi(1)=\E\langle V(\x^\top \dot \v_1)V(\x^\top \dot \v_2)
\rangle=\E[\langle V(\x^\top \dot\v) \rangle^2]$, and we have $\varphi(0)=\E
V(\rho^{1/2}\xi_1)V(\rho^{1/2}\xi_2)=0$ by assumption on $V$.

Differentiating in $t$ and applying Stein's lemma for the expectations
over $\x$ and $\xi_1$,
\begin{align*}
\varphi'(t)&=2\E\left\langle
\left(\frac{1}{2\sqrt{t}} \x^\top \dot \v_1
-\frac{1}{2\sqrt{1-t}}\rho^{1/2}\xi_1\right) V'(S_1(t))V(S_2(t))
\right\rangle\\
&=\E \left\langle
\left(\frac{1}{p}\|\dot \v_1\|^2-\rho\right)V''(S_1(t))V(S_2(t))
+\left(\frac{1}{p}\dot\v_1^\top \dot \v_2\right)V'(S_1(t))V'(S_2(t))
\right\rangle.
\end{align*}
We apply Cauchy-Schwarz over $\E\langle \cdot \rangle$, together with Lemma
\ref{lemma:overlaps}, the bound
\begin{align*}
\E\Big\langle V''(S_1(t))^4 \Big\rangle
&\leq \E\Big\langle (Ce^{C|S_1(t)|})^4\Big\rangle
\leq C^4\E_{\X,\y}\Big\langle \E_{\x,\xi_1} e^{4C|\x^\top \dot
\v_1|+4C|\rho^{1/2}\xi_1|}\Big\rangle \leq C'
\end{align*}
for a constant $C':=C'(C,K,\rho)$, and similarly for $V(S_2(t))$, $V'(S_1(t))$,
and $V'(S_2(t))$. This gives
\[\lim_{n,p \to \infty} \sup_{t \in (0,1)}
|\varphi'(t)|=0.\]
Integrating this bound from $t=0$ to $t=1$ gives
$\lim_{n,p \to \infty} \varphi(1)=0$ as desired.
\end{proof}

\begin{proof}[Proof of Lemma \ref{lemma:cavity}]
Recall the leave-one-out posterior measure (\ref{eq:LOOposterior}), with the
notations $\X^- \in \R^{n \times (p-1)}$ and $\x^j \in \R^n$ for the columns of
$\X$ and the notations $\bbeta^-,\bbeta_0^- \in \R^{p-1}$ and
$\beta_j,\beta_{0,j}$ for the coordinates of $\bbeta,\bbeta_0$.
Setting $\y^-=\X^-\bbeta_0^-+\beps$, we have
\[\X\bbeta-\y=(\beta_j-\beta_{0,j})\x^j+\X^-\bbeta^--\y^-.\]
Then the posterior density of $\bbeta$ is
\begin{align*}
\sP(\de\bbeta \mid \X,\y) &\propto
\exp\left(-\frac{1}{2\sigma^2}\|\X\bbeta-\y\|^2\right)\sP_0(\de\bbeta)\\
&\propto
\exp\left(-\frac{(\beta_j-\beta_{0,j})^2}{2\sigma^2}\|\x^j\|^2
-\frac{\beta_j-\beta_{0,j}}{\sigma^2}{\x^j}^\top(\X^-\bbeta^--\y^-)\right)
\sP^-(\de\bbeta^-)\sP_0(\de\beta_j).
\end{align*}
We define
\begin{equation}\label{eq:hbetadef}
\dot\v^-=\X^-(\bbeta^--\langle \bbeta^- \rangle_-),
\qquad h(\beta_j)=\exp\left(-\frac{\beta_j-\beta_{0,j}}{\sigma^2}{\x^j}^\top (\X^-\langle \bbeta^-\rangle_--\y^-)\right)
\end{equation}
where $h(\beta_j)$ depends only on the coordinate $\beta_j$
and not the remaining coordinates
$\bbeta^-$. Then the posterior density may be written as
\[\sP(\de\bbeta \mid \X,\y) \propto h(\beta_j)
e^{-\frac{(\beta_j-\beta_{0,j})^2}{2\sigma^2}\|\x^j\|^2}
e^{-\frac{\beta_j-\beta_{0,j}}{\sigma^2}{\x^j}^\top \dot\v^-}
\sP^-(\de\bbeta^-)\sP_0(\de\beta_j).\]
This shows the identity
\begin{equation}\label{eq:postexpr}
\langle f(\beta_j) \rangle_{\X,\y}
=\frac{\int
f(\beta_j)h(\beta_j)e^{-\frac{(\beta_j-\beta_{0,j})^2}{2\sigma^2}\|\x^j\|^2}
\big\langle e^{-\frac{\beta_j-\beta_{0,j}}{\sigma^2}{\x^j}^\top \dot\v^-} \big\rangle_-
\sP_0(\de\beta_j)}
{\int h(\beta_j)e^{-\frac{(\beta_j-\beta_{0,j})^2}{2\sigma^2}\|\x^j\|^2}\big\langle
e^{-\frac{\beta_j-\beta_{0,j}}{\sigma^2}{\x^j}^\top \dot\v^-} \big\rangle_-
\sP_0(\de\beta_j)}.
\end{equation}

For any fixed $\beta_j,\beta_{0,j} \in \R$, we may
apply Lemma \ref{lemma:cavityinterpolation}
for the leave-one-out posterior measure
with $\langle \cdot \rangle_-$ and $\dot\v^-$ in
place of $\langle \cdot \rangle$ and $\dot\v$, and with the function
\[V(x)=e^{-\frac{\beta_j-\beta_{0,j}}{\sigma^2}x}
-e^{\frac{(\beta_j-\beta_{0,j})^2\rho}{2\sigma^4}}.\]
Then the lemma implies
\begin{equation}\label{eq:postapprox1}
\lim_{n,p \to \infty}
\E_{\x^j,\X^-,\bbeta_0^-,\beps}
\Big[\Big(\big\langle e^{-\frac{\beta_j-\beta_{0,j}}{\sigma^2}{\x^j}^\top \dot\v^-}
\big\rangle_- -e^{\frac{(\beta_j-\beta_{0,j})^2\rho}{2\sigma^4}}\Big)^2\Big]=0.
\end{equation}
As $\|\x^j\|^2 \to \delta$, the bounded convergence theorem implies also
\begin{equation}\label{eq:postapprox2}
\lim_{n,p \to \infty}
\E_{\x^j}
\Big[\Big(e^{-\frac{(\beta_j-\beta_{0,j})^2}{2\sigma^2}\|\x^j\|^2}
-e^{-\frac{(\beta_j-\beta_{0,j})^2\delta}{2\sigma^2}}\Big)^2\Big]=0.
\end{equation}
Noting that $\delta/\sigma^2-\rho/\sigma^4=\gamma_\stat$, and comparing
the definition of $h(\beta_j)$ in (\ref{eq:hbetadef}) with $\hat{\lambda}_j$
from (\ref{eq:hatlambdadef}), we have exactly
\begin{equation}\label{eq:cavityexpr}
\frac{\int f(\beta_j)h(\beta_j)e^{-\frac{(\beta_j-\beta_{0,j})^2\delta}{2\sigma^2}}
e^{\frac{(\beta_j-\beta_{0,j})^2\rho}{2\sigma^4}}\sP_0(\de\beta_j)}
{\int h(\beta_j)e^{-\frac{(\beta_j-\beta_{0,j})^2\delta}{2\sigma^2}}
e^{\frac{(\beta_j-\beta_{0,j})^2\rho}{2\sigma^4}}\sP_0(\de\beta_j)}
=\langle f(\beta) \rangle_{\hat{\lambda}_j,\gamma_\stat}
\end{equation}
where the right side is defined by (\ref{eq:scalarmean}).
By assumption, $f(\beta)$ is bounded over $\beta \in \supp(\sP_0)$.
For a sufficiently large constant $C_0>0$, the event
$\cE:=\{\|\x\|^2 \leq C_0,\;\|\y\|^2 \leq C_0,\;\|\X\|_{\op} \leq
C_0\}$ has probability approaching 1. On $\cE$, the quantities
\[h(\beta_j),\;e^{-\frac{(\beta_j-\beta_{0,j})^2}{2\sigma^2}\|\x^j\|^2},\;
e^{-\frac{\beta_j-\beta_{0,j}}{\sigma^2}{\x^j}^\top \dot\v^-},\;
e^{-\frac{(\beta_j-\beta_{0,j})^2\gamma_\stat}{2}}\]
are all bounded above and below by a constant, uniformly over
$\beta_j,\beta_{0,j} \in \supp(\sP_0)$.
Then applying (\ref{eq:postapprox1}), (\ref{eq:postapprox2}), and the bounded
convergence theorem to compare
(\ref{eq:postexpr}) with (\ref{eq:cavityexpr}), we get
\[\lim_{n,p \to \infty} \E\Big[\1\{\cE\}\big(\langle f(\beta_j) \rangle_{\X,\y}
-\langle f(\beta) \rangle_{\hat{\lambda}_j,\gamma_\stat}\big)^2\Big]=0.\]
As $f(\beta)$ is bounded and $\P[\cE] \to 1$, we have also
$\lim_{n,p \to \infty} \E[\1\{\cE^c\}(\langle f(\beta_j) \rangle_{\X,\y}
-\langle f(\beta) \rangle_{\hat{\lambda}_j,\gamma_\stat})^2]=0$,
yielding Lemma \ref{lemma:cavity}.
\end{proof}

\section{TAP lower bound via Gordon's comparison inequality}
\label{app:TAP-lower-bound}
In this section, we prove the asymptotic lower bound for the TAP free energy
$\cF_\TAP(\bbm,\bs)$ stated in Lemma \ref{lem:muhat-l-to-mu-l}.
We prove this lower bound using Gordon's comparison inequality
We further prove Theorem \ref{thm:replicalowerbound},
which specializes Lemma \ref{lm:tap_lower_bound_body} to sets on which the Nishimori-type condition
\begin{equation}\label{eq:Nishimorims}
\frac{1}{p}\sum_{j=1}^p (m_j-\beta_{0,j})^2 \approx
\frac{1}{p} \sum_{j=1}^p s_j-m_j^2,
\end{equation}
is satisfied.
For the reader's convenience,
we recall that 
for any $K \subseteq [0,\infty) \times [0,\infty)$, the set $\Gamma^p[K]$ is defined as
\begin{equation}\label{eq:GammaK}
\Gamma^p[K]=\left\{(\bbm,\bs) \in \Gamma^p:\left(\frac{1}{p}
\sum_{j=1}^p (m_j-\beta_{0,j})^2,\;\frac{1}{p} \sum_{j=1}^p s_j-m_j^2\right)
\in K\right\},
\end{equation}
and
$K(\rho) \subseteq \R^2$ denotes the closed Euclidean-ball of radius $\rho$
around $(q_\star,q_\star)$.
In this section, we will also prove the following Corollary of Lemma \ref{lm:tap_lower_bound_body},
which makes precise the lower bound in \eqref{eq:local-min-key-tools}.
\begin{corollary}\label{cor:lower_TAP_f_star}
Suppose Assumption \ref{ass:Bayesian_linear_model} holds.
Let $\gamma_\star>0$ be any local minimizer of $\phi(\cdot)$ that satisfies
$\phi''(\gamma_\star)>0$ strictly. Define
$q_\star=\delta/\gamma_\star-\sigma^2$, 
\begin{itemize}
\item[(a)] There exists $\rho_0:=\rho_0(\delta,\sigma^2,\sP_0,\gamma_\star)>0$
such that for any $\iota>0$ and some $c>0$ (depending on $\rho_0,\iota$), with
probability at least $1-e^{-cn}$ for all large $n,p$,
\begin{equation}
\inf_{(\bbm,\bs) \in \Gamma^p[K(\rho_0)]} \frac{1}{p}\,\cF_\TAP(\bbm,\bs)
\geq \phi(\gamma_\star)-\iota.
\end{equation}
\item[(b)] For any $\rho_1$ such that $0<\rho_1<\rho_0$, there exist 
$c,\iota_0>0$ (depending on $\rho_0,\rho_1$) such that with probability at
least $1-e^{-cn}$ for all large $n,p$,
\begin{equation}
\inf_{(\bbm,\bs) \in \Gamma^p[K(\rho_0) \setminus K(\rho_1)]}
\frac{1}{p}\,\cF_\TAP(\bbm,\bs) \geq \phi(\gamma_\star)+\iota_0.
\end{equation}
\end{itemize}
\end{corollary}

\noindent
The proof of Corollary \ref{cor:lower_TAP_f_star} is given at the end of Section
\ref{subsec:gordonconvexity}.

\subsection{Gordon comparison lower bound: proof of Lemma \ref{lm:tap_lower_bound_body}}
\label{app:gordon-lower-bound}

We will apply a version of Gordon's comparison inequality similar to
\cite[Theorem 3(i)]{thrampoulidis2015regularized}, which we state in the
following lemma.

\begin{lemma}\label{lemma:gordon}
Let $\G \in \R^{n \times p}$, $\bg \in \R^p$, and $\bh \in \R^n$ be independent
with i.i.d.\ $\normal(0,1)$ entries. Let $k,l \geq 0$, let
$S \subset \R^p \times \R^k$ and $T \subset \R^n \times \R^l$ be
compact, let $\psi:S \times T \to \R$ be continuous, and define
\begin{align*}
\Phi(\G)&=\inf_{(\bw,\bv) \in S} \sup_{(\bu,\bs) \in T} \bu^\top
\G\bw+\psi(\bu,\bs,\bw,\bv),\\
\phi(\bg,\bh)&=\inf_{(\bw,\bv) \in S} \sup_{(\bu,\bs) \in T}
\|\bu\|_2\bg^\top \bw +\|\bw\|_2 \bh^\top \bu+\psi(\bu,\bs,\bw,\bv).
\end{align*}
Then for any $c \in \R$, we have
$\P[\Phi(\G) \leq c] \leq 2\,\P[\phi(\bg,\bh) \leq c]$.
\end{lemma}
\begin{proof}
The proof is the same as that of
\cite[Theorem 3(i)]{thrampoulidis2015regularized}, which shows this result for
$k=l=0$: Let $g \sim \normal(0,1)$ be independent of $\G$, and define
\[\Phi(\G,g)=\inf_{(\bw,\bv) \in S} \sup_{(\bu,\bs) \in T} \bu^\top
\G\bw+g\|\bw\|_2\|\bu\|_2+\psi(\bu,\bs,\bw,\bv).\]
When $S$ and $T$ are finite sets, the classical Gordon comparison theorem
(c.f.\ \cite[Theorem A.0.1]{thrampoulidis2015regularized}) gives
$\P[\Phi(\G,g) \leq c] \leq \P[\phi(\bg,\bh) \leq c]$. For $S$ and $T$ compact,
$\psi$ is uniformly continuous on $S \times T$, so the same inequality holds via
a covering net argument---we refer to \cite[Proof of Theorem
1]{thrampoulidis2015regularized} for details. Finally, we have
$\Phi(\G) \geq \Phi(\G,g)$ when $g<0$, so that
\[\P[\Phi(\G) \leq c]=\P[\Phi(\G) \leq c \mid g<0]
\leq \P[\Phi(\G,g) \leq c \mid g<0] \leq \frac{\P[\Phi(\G,g) \leq c]}{\P[g<0]}
\leq 2\,\P[\phi(\bg,\bh) \leq c].\]
\end{proof}

We now define the variational objective $f: \R^5 \to \R$ appearing in the statement of Lemma \ref{lm:tap_lower_bound_body}. Let $a(\sP_0),b(\sP_0)$
be the lower and upper endpoints of $\supp(\sP_0)$, as defined in
(\ref{eq:supportendpoints}). Define $e:\R^5 \to \R$ by
\begin{equation}\label{eq:def_f_2}
e(\beta_0,z;\alpha,\tau,\gamma)=\alpha z\beta_0+\frac{\tau}{2}\beta_0^2+
\inf_{m \in (a(\sP_0),b(\sP_0))} \sup_{\lambda \in \R}
\Big[\lambda m-\log \E_{\beta \sim\sP_0} \big[e^{-(\gamma/2)
\beta^2+\lambda \beta}\big]+\frac{\tau-\gamma}{2}m^2-\alpha zm -\tau \beta_0 m
\Big].
\end{equation}
Then our variational objective is
\begin{equation}\label{eq:def_f_1}
f(q,r;\alpha,\tau,\gamma) = \alpha
\sqrt{\delta(q+\sigma^2)}-\frac{\alpha^2\sigma^2}{2}+
\frac{\delta}{2}\log 2\pi(\sigma^2+r)-\frac{\gamma r}{2}-\frac{\tau q}{2}
+\E_{(\beta_0,z) \sim \sP_0 \times \normal(0, 1)}[e(\beta_0, z; \alpha,\tau,\gamma)]. 
\end{equation}
This is the $f$ appearing in Lemma \ref{lm:tap_lower_bound_body}.

\begin{proof}[Proof of Lemma \ref{lm:tap_lower_bound_body}]
We denote $\dotpd{\bbm}$ as the entrywise square of $\bbm$,
and $\ones$ as the all-1's vector in $\R^p$.
Let us change variables to $\bw=\bbm-\bbeta_0$ and
$\bbv=\bbs-\dotpd{\bbm}$, and define the
inverse map $\bbm(\bw)=\bw+\bbeta_0$ and
$\bs(\bw,\bv)=\bv+\dotpd{\bbm(\bw)}=\bv+\dotpd{(\bw+\bbeta_0)}$.
Recall $\Gamma^p[K]$ from (\ref{eq:GammaK}), and
define the ($\bbeta_0$-dependent) domains
\[\Omega^p=\Big\{(\bw,\bbv):(\bbm(\bw),\bs(\bw,\bv)) \in \Gamma^p\Big\},
\qquad \Omega^p[K]=\Big\{(\bw,\bbv) \in \Gamma^p:
\big(p^{-1}\|\bw\|_2^2,p^{-1}\ones^\top \bv \big) \in K \Big\}\]
so that $(\bw,\bv) \in \Omega^p[K]$ if and only if $(\bbm,\bs) \in \Gamma^p[K]$.
Recall the definition of $\cF_\TAP$ in (\ref{eqn:TAP}), and note that in this
definition, $S(\s)-Q(\bbm)=p^{-1}\ones^\top \v$.
Then, applying the variational representation
$\|\y-\X\bbm\|_2^2/(2\sigma^2)=\sup_{\bu\in\R^{n}}[\bu^{\top}(\X\bbm-\y)-\sigma^2\|\bu\|_2^2/2]=\sup_{\bu\in\R^{n}}[\bu^{\top}\X\bw-\bu^\top\beps-\sigma^2\|\bu\|_2^2/2]$, we obtain
\begin{align}
&\inf_{(\bbm, \bbs) \in \Gamma^p[K]}  \cF_{\TAP}(\bbm,\bbs)\nonumber\\
&=\inf_{(\bw,\bbv) \in \Omega^p[K]}
\sup_{\bu\in\R^{n}} \bu^{\top}\X \bw
\underbrace{\phantom{}-\bu^{\top}\beps-\frac{\sigma^2}{2}\|\bu\|_2^2 +
D_0(\bbm(\bw),\bs(\bw,\bv)) + \frac{n}{2} \log
2\pi(\sigma^2+p^{-1}\ones^\top \v)}_{:=\psi(\bu,\bw,\bv)}. \label{eq:TAPvariational}
\end{align}

\noindent
{\bf Step 1. Comparison using Gordon's inequality.}
We proceed to lower bound (\ref{eq:TAPvariational}) using Lemma
\ref{lemma:gordon}. For any $M \in (0, \infty]$, define
\[\Gamma_M=\Big\{(m,s):m=\langle \beta \rangle_{\lambda,\gamma},\;
s=\langle \beta^2 \rangle_{\lambda,\gamma} \text{ for some }
(\lambda,\gamma) \in [-M,M]^2\Big\}. \]
Note that $\Gamma_M$ is a continuous image of a compact set, and hence is
compact. Then by compactness of $K$, the domain
\[\Omega_M^p[K]:=\Big\{(\bw,\bv):(\bbm(\bw),\bs(\bw,\bv)) \in \Gamma_M^p,\;
\big(p^{-1}\|\bw\|_2^2,p^{-1}\ones^\top \bv \big) \in K \Big\}\]
is also compact. For any $M,M' \in (0,\infty]$, define
\begin{align*}
\Phi_{M,M'}(\X) &=
\inf_{(\bw,\bv) \in \Omega_M^p[K]} \sup_{\|\bu\|_2\leq M'\sqrt{n}} 
\bu^{\top}\X \bw+\psi(\bu,\bw,\bv), \\
\phi_{M,M'}(\bg,\bh) &=
\inf_{(\bw,\bv) \in \Omega_M^p[K]} \sup_{\|\bu\|_2\leq M'\sqrt{n}} 
\|\bu\|_2\,\frac{\bg^\top \bw}{\sqrt{p}}
+\frac{\|\bw\|_2}{\sqrt{p}}\,\bh^\top \bu+\psi(\bu,\bw,\bv),
\end{align*}
where $\psi(\bu,\bw,\bv)$ is as defined in (\ref{eq:TAPvariational}),
and $\bg \sim \normal(\bzero, \id_p)$ and $\bh \sim \normal(\bzero, \id_{n})$ are
independent of each other and of $\beps,\bbeta_0$.
Then Lemma \ref{lemma:gordon} applied with $\bW=\sqrt{p}\,\X$
gives, for any $c \in \R$ (conditionally on $\bbeta_0,\beps$, and hence also
unconditionally)
\begin{equation}\label{eq:gordondirect}
\P[p^{-1}\Phi_{M,M'}(\X) \leq c] \leq 2\,\P[p^{-1}\phi_{M,M'}(\bg,\bh) \leq c]. 
\end{equation}

Note that taking $M'=\infty$, the inner suprema over $\bu$ in $\Phi_{M,\infty}$
and $\phi_{M,\infty}$ are attained at, respectively, 
\[
\bu_\Phi=(\X\bw-\beps)/\sigma^2, ~~~~~ \bu_\phi= \bigg(1+\frac{\bg^\top
\bw}{\sqrt{p}}\bigg/ \bigg\|\frac{\|\bw\|_2}{\sqrt{p}}\bh-\beps \bigg\|_2 \bigg)
\times \bigg(\frac{\|\bw\|_2}{\sqrt{p}}\bh-\beps \bigg)\bigg/\sigma^2. 
\]
Applying that $\|\bw\|_2$ is bounded over $(\bw,\bv) \in
\Omega^p[K]$ by compactness of $K$, this implies that $\|\bu_\Phi\|_2$ and
$\|\bu_\phi\|_2$ are also bounded over $(\bw,\bv) \in
\Omega^p[K]$ for any fixed realizations of $\X,\beps,\bbeta_0,\bg,\bh$, so
$\Phi_{M,\infty}(\X)=\lim_{M' \to \infty} \Phi_{M,M'}(\X)$
and $\phi_{M,\infty}(\bg,\bh)=\lim_{M' \to \infty} \Phi_{M,M'}(\bg,\bh)$.
Then, since each
point $(\bw,\bv) \in \Omega^p[K]$ belongs to $\Omega_M^p[K]$ for sufficiently
large $M$, we have
\[\Phi_{\infty,\infty}(\X)=\lim_{M \to \infty}
\lim_{M' \to \infty} \Phi_{M,M'}(\X), \qquad
\phi_{\infty,\infty}(\bg,\bh)=\lim_{M \to \infty}
\lim_{M' \to \infty} \phi_{M,M'}(\bg,\bh).\]
Taking these limits in (\ref{eq:gordondirect}) and applying
$\inf_{(\bbm,\bs) \in \Gamma^p[K]} \cF_\TAP(\bbm,\bs)=\Phi_{\infty,\infty}(\X)$
by (\ref{eq:TAPvariational}), we have for any $c \in \R$ and $\iota>0$,
\begin{equation}\label{eq:gordoninfty}
\P\left[p^{-1}\inf_{(\bbm,\bs) \in \Gamma^p[K]} \cF_\TAP(\bbm,\bs) \leq c-\iota
\right]
\leq 2\,\P[p^{-1}\phi_{\infty,\infty}(\bg,\bh) \leq c+\iota]. 
\end{equation}\\



\noindent
{\bf Step 2. Asymptotics of the comparison process.} We now derive an asymptotic
lower bound for $\phi_{\infty,\infty}$. We have
\begin{align*}
\phi_{\infty,\infty}(\bg,\bh)
   &=
    \inf_{(q,r)\in K} \mathop{\inf_{(\bw,\bv) \in \Omega^p}}_{p^{-1}\|\bw
\|_2^2=q,\,p^{-1}\ones^\top \bbv=r} \sup_{\bu\in\R^{n}}
   \|\bu\|_2\,\frac{\bg^\top\bw}{\sqrt{p}}
+\frac{\|\bw\|_2}{\sqrt{p}}\,\bh^\top\bu+\psi(\bu,\bw,\bv)\\
   &\geq
 \inf_{(q,r)\in K}\inf_{(\bw,\bv) \in \Omega^p}
\sup_{(\tau,\gamma) \in \R^2} \sup_{\bu\in\R^{n}}
   \|\bu\|_2 \frac{\bg^\top \bw}{\sqrt{p}}
+(\sqrt{{q}}\bh-\beps)^\top \bu-\frac{\sigma^2}{2}\|\bu\|_2^2+
D_0(\bbm(\bw),\bbs(\bv,\bw))\\
&\hspace{1in} + \frac{n}{2} \log 2\pi(\sigma^2+ r)
   +\frac{\tau}{2} \Big (\|\bw\|_2^2-pq \Big)+\frac{\gamma}{2} \Big(\ones^\top \bbv-pr\Big), 
   \end{align*}
where the inequality recalls the definition of $\psi(\cdot)$ from
(\ref{eq:TAPvariational}) and
introduces $(\tau,\gamma)$ as Lagrangian multipliers
for the constraints $\|\bw\|_2^2 = pq$ and $\ones^\top \bbv=pr$. We
observe that the supremum over $\bu$ is attained at
$\bu=\alpha\sqrt{p}(\sqrt{{q}}\bh-\beps)/\|\sqrt{{q}}\bh-\beps\|_2$ for some
$\alpha\geq0$. Applying this and the definition $D_0(\bbm,\bs)=-\sum_j
\sh(m_j,s_j)$ with $\sh(\cdot)$ as in (\ref{eq:hdef}), we have
\begin{align*}
\phi_{\infty,\infty}(\bg,\bh)
& \geq 
 \inf_{(q,r)\in K} \inf_{(\bw,\bv) \in \Omega^p}
 \sup_{(\alpha,\tau,\gamma) \in [0,\infty) \times \R^2}
\sup_{(\blambda,\bgamma) \in\R^p \times \R^p}
   \alpha\,\bg^\top \bw+\alpha\sqrt{p}\|\sqrt{q}\bh-\beps\|_2-\frac{\sigma^2}{2}\alpha^2p\\
&\hspace{1in} +\blambda^\top \bbm(\bw)
   -\frac{1}{2}\bgamma^\top \bbs(\bv,\bw)
- \sum_{j=1}^p \log \E_{\beta \sim \sP_0}[e^{-(\gamma_j/2) \beta^2 + \lambda_j
\beta}]\\
&\hspace{1in}
+ \frac{n}{2} \log 2\pi(\sigma^2+r)+\frac{\tau}{2}\Big(\|\bw\|_2^2-pq\Big)
+\frac{\gamma}{2} \Big(\ones^\top \bv-pr \Big). 
\end{align*}
%

Note that if $(\bw,\bv) \in \Omega^p$, i.e.\ $(\bbm(\bw),\bs(\bw,\bv)) \in
\Gamma^p$, then in particular $\bbm(\bw)=\bw+\bbeta_0 \in (a(\sP_0),b(\sP_0))^p$
by the characterization of $\Gamma$ in Proposition
\ref{prop:Gammacharacterization}. Thus, we may
lower bound this quantity by relaxing $\inf_{(\bw,\bv) \in \Omega^p}$
to $\inf_{(\bw,\bv):\bw+\bbeta_0 \in (a(\sP_0),b(\sP_0))^p}$, and also
restricting $\sup_{(\alpha,\tau,\gamma) \in [0,\infty) \times \R^2}$
to $\sup_{(\alpha,\tau,\gamma) \in K'}$ for any compact domain $K' \subset
[0,\infty) \times \R^2$. We obtain a further lower bound by
exchanging these $\inf_{(\bw,\bv)}$ and $\sup_{(\alpha,\tau,\gamma)}$,
and specializing the inner supremum over
$\bgamma \in \R^p$ to $\gamma_j=\gamma$ for all $j=1,\ldots,p$. Applying this
specialization and recalling $\bs(\bv,\bw)=\bv+\dotpd{(\bw+\bbeta_0)}$,
the dependence on $\bv$ cancels, and we obtain
\begin{align*}
&\phi_{\infty,\infty}(\bg,\bh)\\
&\geq
 \inf_{(q,r)\in K}
\sup_{(\alpha,\tau,\gamma) \in K'}
\inf_{\bw:\bw+\bbeta_0 \in (a(\sP_0),b(\sP_0))^p}
\sup_{\blambda\in\R^p}\;
\alpha\,\bg^\top
\bw+\alpha\sqrt{p}\|\sqrt{q}\bh-\beps\|_2-\frac{\sigma^2}{2}\alpha^2p+
\blambda^\top(\bw+\bbeta_0) \\
   &\qquad - \frac{1}{2}\sum_{j=1}^p\gamma(w_j+\beta_{0,j})^2- \sum_{j=1}^p \log
\E_{\beta \sim \sP_0}[e^{-(\gamma/2) \beta^2 + \lambda_j \beta}]  + \frac{n}{2}
\log 2\pi(\sigma^2 + r)+\frac{\tau}{2} \Big(\|\bw\|_2^2-pq \Big)-\frac{\gamma pr}{2}\\ 
   &=
 \inf_{(q,r)\in K}
\sup_{(\alpha,\tau,\gamma) \in [0,\infty) \times \R^2}
p\Bigg[{\alpha}\frac{\|\sqrt{q}\bh-\beps\|_2}{\sqrt{p}}-\frac{\alpha^2\sigma^2}{2}+\frac{n}{2p}\log
2\pi(\sigma^2+r)-\frac{\tau q}{2}-\frac{\gamma r}{2}\Bigg]\\
&\qquad+\sum_{j=1}^p \underbrace{
\inf_{w_j \in (a(\sP_0)-\beta_{0,j},b(\sP_0)-\beta_{0,j})}
\sup_{\lambda_j}\Big[\alpha g_jw_j+\lambda_j(w_j+\beta_{0,j})
-\frac{\gamma}{2}(w_j+\beta_{0,j})^2-\log \E_{\beta \sim \sP_0}
[e^{-(\gamma/2)\beta^2+\lambda_j \beta}]
+\frac{\tau}{2}w_j^2\Big]}_{:=e(\beta_{0,j},-g_j;\alpha,\tau,\gamma)}. 
\end{align*}
Making the change of variables $w_j=m_j-\beta_{0,j}$ and $z_j=-g_j$, we see that
the above quantity $e(\cdot)$ coincides with the definition (\ref{eq:def_f_2}).

Finally, defining the errors
\begin{align*}
E_1 &= \sup_{(q,r)\in K} \left| \sqrt{\delta(q+\sigma^2)}
-\frac{\|\sqrt{q}\bh-\beps\|_2}{\sqrt{p}} \right|, \\
E_2 &= \sup_{(\alpha,\tau,\gamma)\in K'} \left|\frac{1}{p}\sum_{j=1}^p
e(\beta_{0,j},-g_j;\alpha,\tau,\gamma)-\E_{(\beta_0,z)\sim \sP_0 \times
\normal(0,1)} e(\beta_0,z;\alpha,\tau,\gamma) \right|,
\end{align*}
and recalling the definition of $f(\cdot)$ in (\ref{eq:def_f_1}), we get
 \begin{equation}\label{eqn:lower_bound_phi_infty}
    \phi_{\infty,\infty}(\bg,\bh)
 \geq p \inf_{(q,r)\in K} \sup_{(\alpha, \tau, \gamma)\in K'}
f(q,r;\alpha,\tau,\gamma)-p\left(\sup_{(\alpha,\tau,\gamma)\in K'} \alpha\right)
E_1 - p\,E_2.
\end{equation}
Since $\bh,\beps$ are independent Gaussian vectors, the coordinates of
$\sqrt{q}\bh-\beps$ are i.i.d.\ and equal in law to $\normal(0,q+\sigma^2)$.
Then, by the convergence $n/p \to \delta$ and a standard chi-squared tail bound,
we have for any $\iota>0$, a constant $c>0$, and all large $n,p$,
$\P[E_1 \leq \iota] \geq 1-e^{-cn}$.
The error term $E_2$ is controlled by the following lemma, which we
prove below.
\begin{lemma}\label{lem:bound_E2}
For any $\iota>0$, there exists a constant $c>0$ depending only on
$(\iota,\sP_0,K')$ such that for all large $p$,
$\P[E_2 \le \iota] \ge 1 - e^{-cp}$.
\end{lemma}
\noindent Applying these bounds for $E_1$ and $E_2$ to (\ref{eq:gordoninfty})
and (\ref{eqn:lower_bound_phi_infty}), for any $\iota>0$, we obtain
\[
\begin{aligned}
 &~\P\Big[p^{-1}\inf_{(\bbm,\bs) \in \Gamma^p[K]}
\cF_\TAP(\bbm,\bs) \leq \inf_{(q,r)\in K}\sup_{(\alpha,\tau,\gamma)\in K'}
f(q,r;\alpha,\tau,\gamma)-\iota \Big]\\
&\leq
 2\,\P\Big[p^{-1}\phi_{\infty,\infty}(\bg,\bh) \leq \inf_{(q,r)\in K}
\sup_{(\alpha,\tau,\gamma)\in K'} f(q,r;\alpha,\tau,\gamma)-\iota/2\Big]
\leq e^{-cn}
 \end{aligned}
 \]
for a constant $c>0$ and all large $n,p$. This completes the proof of the lemma.
\end{proof}

\begin{proof}[Proof of Lemma~\ref{lem:bound_E2}]
We first show concentration for fixed $(\alpha,\tau,\gamma) \in K'$.
Proposition \ref{prop:oneparamexpfam} verifies that in the definition
(\ref{eq:def_f_2}) of $e(\cdot)$, for any $m \in
(a(\sP_0),b(\sP_0))$, the supremum over $\lambda$ is attained
at some value $\lambda(m) \in \R$. Hence
\begin{align*}
e(\beta_0,z;\alpha,\tau,\gamma)&=\alpha z\beta_0+\frac{\tau}{2}\beta_0^2+
e_1(\beta_0,g;\alpha,\tau,\gamma), \\
e_1(\beta_0,g;\alpha,\tau,\gamma)&:=\inf_{m \in (a(\sP_0),b(\sP_0))}
\Big[\lambda(m)m-\log \E_{\beta \sim\sP_0} \big[e^{-(\gamma/2)
\beta^2+\lambda(m)\beta}\big]+\frac{\tau-\gamma}{2}m^2-\alpha zm -\tau m\beta_0
\Big]. 
\end{align*}
This function $e_1(\cdot)$ is an infimum of linear
$L$-Lipschitz functions of $\beta_0,z$ for a
constant $L:=L(\sP_0,K')>0$, and hence $e_1(\cdot)$ is concave
and $L$-Lipschitz in $(\beta_0,z)$. Then the average of
$e_1(\beta_{0,j},z_j;\alpha,\tau,\gamma)$ over coordinates $j=1,\ldots,p$ (with
$z_j=-g_j$) is
also concave and $L/\sqrt{p}$-Lipschitz in $(\bbeta_0,\bz)$. Here $\bbeta_0$ has
bounded entries, and $\bz$ is standard Gaussian, so
for constants $C,c>0$ (depending on $\sP_0,K'$) and any $t \geq 0$,
\[\P\left[\left|\frac{1}{p}\sum_{j=1}^p e_1(\beta_{0,j},z_j;\alpha,\tau,\gamma)
-\E_{(\beta_0,z) \sim \sP_0 \times \normal(0,1)} e_1(\beta_0,z;\alpha,\tau,\gamma)
\right|>t\right] \leq Ce^{-cpt^2}\]
by the Talagrand and Borell-TIS concentration inequalities (c.f.\ \cite[Theorem
5.2.2 and 5.2.16]{vershynin2018high}). The first term $\alpha
z_j\beta_{0,j}+(\tau/2)\beta_{0,j}^2$ of $e(\beta_{0,j},z_j;\alpha,\tau,\gamma)$
is a sub-exponential random variable, hence its average over $j=1,\ldots,p$
also concentrates around its mean by Bernstein's inequality (c.f.\ \cite[Theorem
2.8.1]{vershynin2018high}). Putting this together, for any $\iota>0$, there
exist constants $C,c>0$ depending on $\iota,\sP_0,K'$ such that for any
fixed $(\alpha,\tau,\gamma) \in K'$,
\[\P\left[\left|\frac{1}{p}\sum_{j=1}^p e(\beta_{0,j},z_j;\alpha,\tau,\gamma)
-\E_{(\beta_0,z) \sim \sP_0 \times \normal(0,1)} e(\beta_0,z;\alpha,\tau,\gamma)
\right|>\iota/2\right] \leq Ce^{-cp}.\]

We now apply a covering net argument over $(\alpha,\tau,\gamma)$.
Let $\mathcal N$ be a covering net of $K'$ of cardinality $|\mathcal N|\leq
Cp^3$, such that each $(\alpha,\tau,\gamma) \in K'$ has
$(\alpha',\tau',\gamma') \in \normal$ with
$|\alpha-\alpha'|+|\tau-\tau'|+|\gamma-\gamma'|<1/p$.
Note that for fixed $\bbeta_0,\bz$, the function
$p^{-1}\sum_j e(\beta_{0,j},z_j;\alpha,\tau,\gamma)$ is
$L\|\bz\|_1/p$-Lipschitz in $\alpha$ and $L$-Lipschitz in $\tau$ and $\gamma$,
for a constant $L:=L(\sP_0)>0$. Then for all sufficiently large $p$,
 \begin{align*}
\P[E_2>\iota]
&\leq 
  \P\left[\sup_{(\tau,\alpha,\gamma) \in \mathcal N}
\left|\frac{1}{p}\sum_{j=1}^p e(\beta_{0,j},z_j;\alpha,\tau,\gamma)
-\E_{(\beta_0,z) \sim \sP_0 \times \normal(0,1)} e(\beta_0,z;\alpha,\tau,\gamma)
\right|>\iota/2\right]+\P[\|\bz\|_1/p>2]\\
 &\leq Cp^3e^{-cp}+e^{-cp} \leq e^{-c'p}
 \end{align*}
for constants $C,c,c'>0$ depending only on $\iota,\sP_0,K'$, as desired.
\end{proof}

\subsection{Lower bound by the replica-symmetric potential}
\label{sec:proof-thm-replicalowerbound}

We now use Lemma \ref{lm:tap_lower_bound_body} to 
prove Theorem \ref{thm:replicalowerbound}. The proof follows from the next
proposition, which verifies that the variational objective
$f(q,r;\alpha,\tau,\gamma)$ coincides with the
replica-symmetric potential $\phi(\gamma)$ upon a specialization of its
parameters.

\begin{proposition}\label{prop:specialization}
Let $f(\cdot)$ be as defined in (\ref{eq:def_f_1}), and let $\phi(\gamma)$ be
the potential (\ref{eqn:phi_potential}). Fix any $\gamma>0$, and set
\[q=r=\delta/\gamma-\sigma^2, \qquad \tau=\gamma, \qquad \alpha=\sqrt{\gamma}. \]
Then $f(q,r;\alpha,\tau,\gamma)=\phi(\gamma)$.
\end{proposition}
\begin{proof}
Specializing (\ref{eq:def_f_2}) to $\tau=\gamma$ and $\alpha=\sqrt{\gamma}$,
for any $m \in (a(\sP_0),b(\sP_0))$,
Proposition \ref{prop:oneparamexpfam} shows that
the supremum over $\lambda$ in (\ref{eq:def_f_2}) is attained at some
$\lambda(m)$ which is differentiable and strictly increasing in $m$.
Furthermore, the derivative of
\[h_1(m):=\lambda(m)m-\log \E_{\beta \sim
\sP_0}[e^{-(\gamma/2)\beta^2+\lambda(m)\beta}]-\sqrt{\gamma}zm-\gamma\beta_0m\]
is $h_1'(m)=\lambda(m)-\sqrt{\gamma}z-\gamma \beta_0$ by the envelope theorem.
Since $\lambda(m)$ is strictly increasing with $\lambda(m) \to
-\infty$ as $m \to a(\sP_0)$ and $\lambda(m) \to \infty$ as $m \to b(\sP_0)$,
this shows that the supremum of $h_1(m)$ is also attained a unique value
$m \in (a(\sP_0),b(\sP_0))$. Thus, there are unique values
$(\lambda_\star(\beta_0,g),m_\star(\beta_0,g))$
defined by the stationary conditions
\[\lambda_\star(\beta_0,z)=\gamma \beta_0+\sqrt{\gamma} z,
\qquad m_\star(\beta_0,z)=\langle \beta \rangle_{\lambda_\star(\beta_0,z),\gamma}\]
of (\ref{eq:def_f_2}),
and applying these definitions in (\ref{eq:def_f_2}) gives
\[\E_{(\beta_0,z) \sim \sP_0 \times \normal(0,1)}[e(\beta_0,z;\sqrt{\gamma},\gamma,\gamma)]
=\E_{(\beta_0,z) \sim \sP_0 \times \normal(0,1)}
\left[\frac{\tau}{2}\beta_0^2-\log \E_{\beta \sim \sP_0}
[e^{-(\gamma/2)\beta^2+\lambda_\star(\beta_0,z) \beta}]\right].\]
It is direct to check that this is exactly the mutual information
$i(\gamma)=\E[\log \frac{\sP(\lambda|\beta_0)}{\sP(\lambda)}]$ in
(\ref{eqn:phi_potential}). Then, specializing (\ref{eq:def_f_1})
further to $q=r=\delta/\gamma-\sigma^2$, we obtain
\[f\left(\frac{\delta}{\gamma}-\sigma^2,\frac{\delta}{\gamma}-\sigma^2;
\sqrt{\gamma},\gamma,\gamma\right)
=\frac{\sigma^2\gamma}{2}+\frac{\delta}{2}\log \frac{2\pi \delta}{\gamma}
+i(\gamma)=\phi(\gamma).\]
\end{proof}

\begin{proof}[Proof of Theorem \ref{thm:replicalowerbound}]
Fix any $\eps>0$, and suppose $K \subset [0,\infty) \times [0,\infty)$ is a
compact set such that $\sup_{(q,r) \in K} |q-r| \leq \eps$. Define
$\gamma(q)=\delta/(q+\sigma^2)$ as in the theorem statement, and take $K'
\subset [0,\infty) \times \R^2$ to be any compact set containing all points
$\{(\sqrt{\gamma(q)},\gamma(q),\gamma(q)):(q,r) \in K\}$. Applying
Lemma \ref{lm:tap_lower_bound_body} with $\iota=\eps$ and
with these choices of $K$ and $K'$, and further
lower bounding the supremum over $(\alpha,\tau,\gamma) \in K'$ by the
specialization $(\alpha,\tau,\gamma)=(\sqrt{\gamma(q)},\gamma(q),\gamma(q))$,
we have
\[\inf_{(\bbm,\bs) \in \Gamma^p[K]}
\frac{1}{p}\,\cF_\TAP(\bbm,\bs)
\geq \inf_{(q,r) \in K} f(q,r;\sqrt{\gamma(q)},\gamma(q),\gamma(q))-\eps\]
with probability approaching 1. Observe that $\gamma(q) \leq \delta/\sigma^2$
for any $q \geq 0$. Then, applying Lipschitz continuity in $r$ of
$f(q,r;\alpha,\tau,\gamma)$ as defined by (\ref{eq:def_f_1}), we have
\[f(q,r;\sqrt{\gamma(q)},\gamma(q),\gamma(q))
\geq f(q,q;\sqrt{\gamma(q)},\gamma(q),\gamma(q))-\frac{\delta}{\sigma^2}\eps
=\psi(\gamma(q))-\frac{\delta}{\sigma^2}\eps,\]
where the last equality applies Proposition \ref{prop:specialization}. Thus,
with probability approaching 1,
\[\inf_{(\bbm,\bs) \in \Gamma^p[K]}
\frac{1}{p}\,\cF_\TAP(\bbm,\bs)
\geq \inf_{(q,r) \in K} \psi(\gamma(q))
-\left(1+\frac{\delta}{\sigma^2}\right)\eps,\]
implying Theorem \ref{thm:replicalowerbound}.
\end{proof}

\subsection{Local convexity of the lower bound}\label{subsec:gordonconvexity}

We now show that each local minimizer of the replica-symmetric potential
$\phi(\cdot)$ corresponds to a local minimizer of the variational objective
$f(\cdot)$, and furthermore the lower bound $\sup_{(\alpha,\tau,\gamma)}
f(q,r;\alpha,\tau,\gamma)$ for $\cF_\TAP$ implied by Lemma
\ref{lm:tap_lower_bound_body} is strongly convex in $(q,r)$ at each such minimizer.

\begin{lemma}\label{lm:tap_lower_bound_function}
Recall the potential $\phi(\cdot)$ from (\ref{eqn:phi_potential}), and the
variational objective $f(\cdot)$ from (\ref{eq:def_f_1}).
\begin{enumerate}
\item[(a)] If $\gamma_\star>0$ is any critical point of $\phi$,
i.e.\ $\phi'(\gamma_\star)=0$, then
\begin{equation}\label{eqn:stationary_in_lemma_tap_lower_bound_function}
q_\star=r_\star=\delta/\gamma_\star-\sigma^2, \qquad \tau_\star=\gamma_\star,
\qquad \alpha=\sqrt{\gamma_\star}
\end{equation}
is a critical point of $f$, i.e.\ $\nabla
f(q_\star,r_\star;\alpha_\star,\tau_\star,\gamma_\star)=0$.
\item[(b)] Suppose $\gamma_\star>0$ is a local minimizer of $\phi$ with
$\phi''(\gamma_\star)>0$ strictly, and define
$q_\star, r_\star, \tau_\star, \alpha_\star$ by
\eqref{eqn:stationary_in_lemma_tap_lower_bound_function}. Then for
any compact subset $K'\subseteq [0,\infty) \times \R^2$ containing
$(\alpha_\star,\tau_\star,\gamma_\star)$ in its interior,
\[\sup_{(\alpha,\tau,\gamma)\in K'}       
f(q_\star,r_\star;\alpha,\tau,\gamma)=f(q_\star,r_\star;\alpha_\star,\tau_\star,\gamma_\star)=\phi(\gamma_\star)\]
Moreover, define
\begin{equation}\label{eqn:bar_f_in_lemma_calculation}
\bar f(q,r)=\sup_{(\alpha,\tau,\gamma) \in K'} f(q,r;\alpha,\tau,\gamma)
\end{equation}
Then for some constants $\rho,c>0$ depending on
$(\delta,\sigma^2,\sP_0,\gamma_\star,K')$, we have
\begin{equation}
\bar f(q,r) \geq \phi(\gamma_\star)+c \cdot [ (r-r_\star)^2+(q-q_\star)^2 ]
\text{ for all } (q,r) \text{ such that } |q-q_\star|,|r-r_\star|\leq\rho.
\end{equation}
\end{enumerate}
\end{lemma}

\begin{proof}
The stationary conditions for $\inf_m$ and $\sup_\lambda$ in (\ref{eq:def_f_2})
are
\begin{align}
0&=\lambda_\star+(\tau-\gamma)m_\star-\alpha z-\tau \beta_0 \label{eq:stationary1}\\
0&=m_\star-\langle \beta \rangle_{\lambda_\star,\gamma}\label{eq:stationary2}
\end{align}
Let $\gamma_\star>0$ be any critical point of $\phi(\cdot)$.
We have shown in the proof of Proposition \ref{prop:specialization} that at
$(\alpha,\tau,\gamma)=(\alpha_\star,\tau_\star,\gamma_\star)=(\sqrt{\gamma_\star},\gamma_\star,\gamma_\star)$,
these equations (\ref{eq:stationary1}--\ref{eq:stationary2})
have unique solutions $m_\star \in (a(\sP_0),b(\sP_0))$ and
$\lambda_\star \in \R$, which realize the supremum and infimum in
(\ref{eq:def_f_2}). Then the implicit function theorem implies that these
solutions extend smoothly to 
$\lambda_\star=\lambda_\star(\beta_0,z;\alpha,\tau,\gamma)$ and
$w_\star=w_\star(\beta_0,z;\alpha,\tau,\gamma)$ solving
(\ref{eq:stationary1}--\ref{eq:stationary2}) in an open neighborhood
of $(\alpha_\star,\tau_\star,\gamma_\star)$, and continuity of $e(\cdot)$
implies that they also realize the supremum in
(\ref{eq:def_f_2}) for $(\alpha,\tau,\gamma)$ within a sufficiently small
such neighborhood.\\

\noindent {\bf Proof of part (a).} 
Computing the gradient of $f$ in this neighborhood, we obtain
\begin{subequations}
\begin{align}
\partial_{q}f=&~\frac{1}{2}\left(\alpha\sqrt{\frac{\delta}{q +
\sigma^2}}-\tau\right)\label{eq:grad_q_f}\\
\partial_{r}f=&~\frac{1}{2}\left(\frac{\delta}{\sigma^2+r}-\gamma\right)\label{eq:grad_v_f}\\
\partial_{\tau}f =&~\frac{1}{2} \Big(\E_{\beta_0,z}
[(m_\star-\beta_0)^2]-q\Big)=\frac{1}{2}\Big(\E_{\beta_0,z}[(\langle \beta
\rangle_{\lambda_\star,\gamma}-\beta_0)^2]-q\Big)\label{eq:grad_tau_f}\\
\partial_{\gamma}f=&~\frac{1}{2}\Big(\E_{\beta_0,z}[\langle \beta^2
\rangle_{\lambda_\star,\gamma}]-m_\star^2-r\Big)=\frac{1}{2}\Big(\E_{\beta_0,z}[\Var_{\beta
\sim \sP_{\lambda_\star,\gamma}}[\beta]]-r\Big)\label{eq:grad_gamma_f}\\
\partial_{\alpha}f=&~
\sqrt{\delta(q + \sigma^2)}-\sigma^2\alpha-\E_{\beta_0,z}[zm_\star]. 
\label{eq:grad_alpha_f}
\end{align}
\end{subequations}
We have exchanged differentiation in $(\alpha,\tau,\gamma)$ with
expectation in $(\beta_0,g)$ of $e(\beta_0,g;\alpha,\tau,\gamma)$ using
the dominated convergence theorem, and applied
the stationary condition (\ref{eq:stationary2})
in the second equalities of (\ref{eq:grad_tau_f}--\ref{eq:grad_gamma_f}).
The last derivative (\ref{eq:grad_alpha_f})
may be simplified using Stein's lemma and (\ref{eq:stationary2}),
\[\E[zm_\star]=\E[\partial_z m_\star]
=\E[\partial_z \langle \beta \rangle_{\lambda_\star,\gamma}]
=\E[\Var_{\beta \sim \sP_{\lambda_\star,\gamma}}[\beta] \cdot \partial_z
\lambda_\star]. \]
Differentiating (\ref{eq:stationary1}--\ref{eq:stationary2}) both in $z$, we
have $0=\partial_z \lambda_\star+(\tau-\gamma)\partial_z m_\star-\alpha$ and
$0=\partial_z m_\star-\Var_{\beta \sim \sP_{\lambda_\star,\gamma}}[\beta]
\cdot \partial_z \lambda_\star$, so $\partial_z
\lambda_\star=\alpha/[1+\Var_{\beta \sim
\sP_{\lambda_\star,\gamma}}[\beta](\tau-\gamma)]$ and hence
\begin{equation}\label{eq:grad_alpha_f_final}
\partial_\alpha f=\sqrt{\delta(q+\sigma^2)}-\sigma^2 \alpha
-\alpha\,\E_{\beta_0,g}\left[\frac{\Var_{\beta \sim
\sP_{\lambda_\star,\gamma}}[\beta]}{1+\Var_{\beta \sim
\sP_{\lambda_\star,\gamma}}[\beta](\tau-\gamma)}\right]. 
\end{equation}

Specializing to $(q_\star,r_\star;\alpha_\star,\tau_\star,\gamma_\star)$,
the stationary condition (\ref{eq:stationary1}) gives
$\lambda_\star=\gamma_\star \beta_0+\sqrt{\gamma_\star}z$, so the law
$\sP_{\lambda_\star,\gamma_\star}$ defined by (\ref{eq:scalarmean})
is precisely the posterior distribution for
$\beta_0$ given $\lambda_\star$. Then we have
\[\E_{\beta_0,z}[(\langle \beta \rangle_{\lambda_\star,\gamma_\star}-\beta_0)^2]
=\E_{\beta_0,z}[\Var_{\beta \sim \sP_{\lambda_\star,\gamma_\star}}[\beta]]
=\mmse(\gamma_\star). \]
We recall from (\ref{eq:phifixedpoint})
that $\mmse(\gamma_\star)=\delta/\gamma_\star-\sigma^2$,
because $\phi'(\gamma_\star)=0$. Then it is easily checked that
(\ref{eq:grad_q_f}--\ref{eq:grad_gamma_f}) and (\ref{eq:grad_alpha_f_final})
all vanish at $(q_\star,r_\star;\alpha_\star,\tau_\star,\gamma_\star)$.\\


\noindent
{\bf Proof of part (b).} The following lemmas record the
Hessians of both $\phi(\cdot)$ and $f(\cdot)$; we defer their proofs to after
the completion of the proof of Lemma \ref{lm:tap_lower_bound_function}.

\begin{lemma}\label{lem:second_derivative_phi}
For any $\gamma>0$,
\[\phi''(\gamma)=\frac{1}{2}\left(
\frac{\delta}{\gamma^2}-\E_{(\beta_0,z) \sim \sP_0 \times
\normal(0,1)}[\Var[\beta_0 \mid \gamma \beta_0+\sqrt{\gamma} z]^2]\right). 
\]
\end{lemma}

\begin{lemma}\label{lemma:hessian_computation}
Let $\lambda_\star=\gamma_\star\beta_0+\sqrt{\gamma_\star}z$, where $\beta_0
\sim \sP_0$ and $z \sim \normal(0,1)$ are independent. Set
\[b_\star=\E_{\beta_0,z}[\Var_{\beta \sim
\sP_{\lambda_\star,\gamma_\star}}[\beta]^2], \qquad
k_\star=\E_{\beta_0,z}[(\beta_0-\langle
\beta\rangle_{\lambda_\star,\gamma_\star})^4]-3b_\star. \]
Then the Hessian of $f(\cdot)$ at
$(q_\star,r_\star;\alpha_\star,\tau_\star,\gamma_\star)$ is
\[\nabla^2 f(q_\star,r_\star;\alpha_\star,\tau_\star,\gamma_\star)
=\begin{pmatrix}
   -\frac{\gamma_\star^2}{4\delta} & 0 & \frac{\sqrt{\gamma_\star}}{2} & -\frac{1}{2} & 0\\
   0 & -\frac{\gamma_\star^2}{2\delta} & 0 & 0 &-\frac{1}{2}\\
   \frac{\sqrt{\gamma_\star}}{2} & 0 & -\frac{\delta}{\gamma_\star}-k_\star\gamma_\star & b_\star\sqrt{\gamma_\star} &\frac{k_\star\sqrt{\gamma_\star}}{2}\\
   -\frac{1}{2} & 0 & b_\star\sqrt{\gamma_\star} & -b_\star & 0\\
   0 & -\frac{1}{2} & \frac{k_\star\sqrt{\gamma_\star}}{2} & 0 & -\frac{b_\star}{2}-\frac{k_\star}{4}
\end{pmatrix}.\]
\end{lemma}

Consider $g(q,r;\tau,\gamma)=f(q,r;\sqrt{\gamma},\tau,\gamma)$,
which specializes the variational objective $f(\cdot)$ to
$\alpha=\sqrt{\gamma}$. By Proposition
\ref{prop:specialization}, part (a) of Lemma
\ref{lm:tap_lower_bound_function} already proven, 
and Lemma \ref{lemma:hessian_computation}, we have
\[g(q_\star,r_\star;\tau_\star,\gamma_\star)=\phi(\gamma_\star),
\qquad \nabla g(q_\star,r_\star;\tau_\star,\gamma_\star)=0,
\qquad
\nabla^2 g(q_\star,r_\star;\tau_\star,\gamma_\star)=D^\top \nabla^2
f(q_\star,r_\star;\alpha_\star,\tau_\star,\gamma_\star)D\]
where
\[ D=\begin{pmatrix}
     1 & 0 & 0 & 0\\
      0 & 1 & 0 & 0\\
       0& 0 & 0 & \frac{1}{2\sqrt{\gamma_\star}}\\
       0& 0 & 1 & 0\\
       0& 0 & 0 & 1
 \end{pmatrix}, \quad\text{ so } \quad
\nabla^2 g(q_\star,r_\star;\tau_\star,\gamma_\star)
=\underbrace{\begin{pmatrix}
     -\frac{\gamma_\star^2}{4\delta} & 0& -\frac{1}{2} &\frac{1}{4} \\
   0&-\frac{\gamma_\star^2}{2\delta} & 0 &-\frac{1}{2} \\
   -\frac{1}{2} & 0 & -b_\star  &\frac{b_\star}{2}\\
   \frac{1}{4} & -\frac{1}{2} & \frac{b_\star}{2} &
-\frac{\delta}{4\gamma_\star^2}-\frac{b_\star}{2}
\end{pmatrix}}_{:=\begin{pmatrix} H_{qr,qr} & H_{qr,\tau\gamma} \\
H_{\tau\gamma,qr} & H_{\tau\gamma,\tau\gamma} \end{pmatrix}}.\]
Observe that by Lemma \ref{lem:second_derivative_phi} and the given condition
$\phi''(\gamma_\star)>0$, we have
\[\det H_{\tau\gamma,\tau\gamma}:=\det \nabla_{\tau,\gamma}^2 g(q_\star,r_\star;\tau_\star,\gamma_\star)
=\frac{b_\star^2}{4}+\frac{\delta b_\star}{4\gamma_\star^2}
=\frac{b_\star}{2}\phi''(\gamma_\star)+\frac{b_\star^2}{4}>0.\]
Also $\Tr H_{\tau\gamma,\tau\gamma}<0$, so this implies
$H_{\tau\gamma,\tau\gamma} \prec 0$. We may compute also the Schur complement
\begin{equation}\label{eq:schur}
S:=H_{qr,qr}-H_{qr,\tau\gamma}H_{\tau\gamma,\tau\gamma}^{-1}H_{\tau\gamma,qr}
=\begin{pmatrix}
-\frac{\gamma_\star^2}{4\delta}+\frac{1}{4b_\star} & 0 \\ 0 &
   -\frac{\gamma_\star^2}{2\delta}+\left(b_\star+\frac{\delta}{\gamma_\star^2}\right)^{-1}
\end{pmatrix}. 
\end{equation}
The given condition $\phi''(\gamma_\star)>0$ is, by Lemma
\ref{lem:second_derivative_phi}, equivalent to
$\delta-b_\star\gamma_\star^2>0$, from which it is easily verified that
$S \succ 0$.

For the given compact set $K' \subset [0,\infty) \times \R^2$, define
\[\bar g(q,r)=\sup_{(\tau,\gamma):(\sqrt{\gamma},\tau,\gamma) \in K'}
g(q,r;\tau,\gamma).\]
Then we have $\bar f(q,r) \geq \bar g(q,r)$, because the supremum defining $\bar
g(q,r)$ is taken over a smaller domain. Since
$\det H_{\tau\gamma,\tau\gamma} \neq 0$,
the implicit function theorem shows that in a neighborhood
$\{(q,r):|q-q_\star|,|r-r_\star| \leq \rho\}$
of $(q_\star,r_\star)$,
there exist analytic functions $\tau=\tau(q,r)$ and $\gamma=\gamma(q,r)$
where
$(\tau_\star,\gamma_\star)=(\tau(q_\star,r_\star),\gamma(q_\star,r_\star))$
and $0=\nabla_{\tau,\gamma} g(q,r;\tau(q,r),\gamma(q,r))$.
For sufficiently small $\rho$, we must have
$(\sqrt{\gamma(q,r)},\tau(q,r),\gamma(q,r)) \in K'$ because
$(\alpha_\star,\tau_\star,\gamma_\star) \in K'$ by assumption. Then
$\bar{g}(q,r)=g(q,r;\tau(q,r),\gamma(q,r))$ by the concavity
$H_{\tau\gamma,\tau\gamma} \prec 0$. Differentiating in $(q,r)$, where the
derivatives of $(\tau(q,r),\gamma(q,r))$ are obtained by implicitly
differentiating $0=\nabla_{\tau,\gamma} g(q,r;\tau(q,r),\gamma(q,r))$, we get
$\nabla^2 \bar{g}(q_\star,r_\star)=S$, the Schur complement matrix in
(\ref{eq:schur}). Since $S \succ 0$, this shows that
\begin{align*}
\bar f(q,r) \geq \bar g(q,r) \geq \bar g(q_\star,r_\star)
+c(q-q_\star)^2+c(r-r_\star)^2
=\phi(\gamma_\star)+c(q-q_\star)^2+c(r-r_\star)^2
\end{align*}
for $|q-q_\star|,|r-r_\star|\leq\rho$, where $\rho,c>0$ are some small constants
depending on $\gamma_\star$ and the function $\bar g$, and hence on
$(\delta,\sigma^2,\sP_0,\gamma_\star,K')$. This completes the proof of part (b).
\end{proof}

\begin{proof}[Proof of Lemma \ref{lem:second_derivative_phi}]
This result is standard, see e.g.\ \cite[Theorem 2]{payaro2009hessian},
but for convenience we include a proof here.

By the I-MMSE relationship $\frac{\de }{\de \gamma}i(\gamma) = \frac{1}{2} \mmse
(\gamma)$, we have
\begin{equation}\label{eq:phisecondder}
\phi''(\gamma)=\frac{\delta}{2\gamma^2}+\frac{1}{2}\,\mmse'(\gamma). 
\end{equation}
Let $\lambda=\gamma \beta_0+\sqrt{\gamma} z$, and write 
$\kappa_j[f_1(\beta),\ldots,f_j(\beta) \mid \lambda]$ 
and $\kappa_j[f(\beta) \mid \lambda]=\kappa_j[f(\beta),\ldots,f(\beta) \mid
\lambda]$ for the $j^\text{th}$ pure and
mixed cumulants under the posterior law $\sP_{\lambda,\gamma}(\de\beta) \propto
e^{-(\gamma/2)\beta^2+\lambda \beta}\sP_0(\de\beta)$.
For example,
\begin{align*}
\kappa_1[f(\beta) \mid \lambda]&=\langle f(\beta) \rangle_{\lambda,\gamma}
=\E[f(\beta_0) \mid \gamma\beta_0+\sqrt{\gamma} z]\\
\kappa_2[f(\beta) \mid \lambda]
&=\langle f(\beta)^2 \rangle_{\lambda,\gamma}
-\langle f(\beta) \rangle_{\lambda,\gamma}^2
=\Var[f(\beta_0) \mid \gamma\beta_0+\sqrt{\gamma} z]\\
\kappa_2[f(\beta),g(\beta) \mid \lambda]
&=\langle f(\beta)g(\beta) \rangle_{\lambda,\gamma}
-\langle f(\beta) \rangle_{\lambda,\gamma}\langle
g(\beta) \rangle_{\lambda,\gamma}=\Cov[f(\beta_0),g(\beta_0) \mid \gamma
\beta_0+\sqrt{\gamma}z]. 
\end{align*}
Then by the chain rule, we have
\[\frac{d}{d\gamma} \kappa_j[f_1(\beta),\ldots,f_j(\beta) \mid \lambda]
=\left(\beta_0+\frac{z}{2\sqrt{\gamma}}\right)
\kappa_{j+1}[f_1(\beta),\ldots,f_j(\beta),\beta]
-\frac{1}{2}\kappa_{j+1}[f_1(\beta),\ldots,f_j(\beta),\beta^2]\]
and by Stein's lemma, we have
\[\E[z\kappa_j[f_1(\beta),\ldots,f_j(\beta) \mid \lambda]]
=\E\left[\frac{d}{dz}\kappa_j[f_1(\beta),\ldots,f_j(\beta) \mid \lambda]\right]
=\sqrt{\gamma}\,\E[\kappa_{j+1}[f_1(\beta),\ldots,f_j(\beta),\beta \mid
\lambda]].\]

Applying these identities,
 \begin{align}
 \mmse'(\gamma)
 =&~ \frac{d}{d\gamma}\,\E\big[\kappa_2[\beta \mid \lambda]\big]\nonumber\\
 =&~ \E\left[\left(\beta_0+\frac{z}{2\sqrt{\gamma}}\right)\kappa_3[\beta \mid
\lambda]-\frac{1}{2}\kappa_3[\beta,\beta,\beta^2\mid \lambda]\right]\nonumber\\
 =&~ \E\left[\kappa_1[\beta \mid \lambda]
\kappa_3[\beta \mid \lambda]
+\frac{1}{2}\kappa_4[\beta \mid \lambda]
-\frac{1}{2}\kappa_3[\beta,\beta,\beta^2 \mid
\lambda]\right]\label{eq:mmseprimetmp}
 \end{align}
where we have used the tower property of conditional expectation
for the first term of the second line. Then, expanding in moments
\begin{equation}\label{eq:cumulants}
\begin{aligned}
\kappa_1[\beta \mid \lambda]&=\langle \beta \rangle_{\lambda,\gamma}\\
\kappa_3[\beta \mid \lambda]
&=\langle \beta^3 \rangle_{\lambda,\gamma}
-3\langle \beta^2 \rangle_{\lambda,\gamma}\langle \beta \rangle_{\lambda,\gamma}
+2\langle \beta \rangle_{\lambda,\gamma}^3\\
\kappa_3[\beta,\beta,\beta^2 \mid \lambda]
&=\langle \beta^4 \rangle_{\lambda,\gamma}
-\langle \beta^2 \rangle_{\lambda,\gamma}^2
-2\langle \beta^3 \rangle_{\lambda,\gamma}\langle \beta \rangle_{\lambda,\gamma}
+2\langle \beta^2 \rangle_{\lambda,\gamma}
\langle \beta \rangle_{\lambda,\gamma}^2\\
\kappa_4[\beta \mid \lambda]
&=\langle \beta^4 \rangle_{\lambda,\gamma}
-4\langle \beta^3 \rangle_{\lambda,\gamma}
\langle \beta \rangle_{\lambda,\gamma}
-3\langle \beta^2 \rangle_{\lambda,\gamma}^2
+12\langle \beta^2 \rangle_{\lambda,\gamma}
\langle \beta \rangle_{\lambda,\gamma}^2
-6\langle \beta \rangle_{\lambda,\gamma}^4
\end{aligned}
\end{equation}
and cancelling terms, we arrive at
\begin{equation}\label{eq:mmseprime}
\mmse'(\gamma)=\E[{-}\langle \beta^2 \rangle_{\lambda,\gamma}^2
+2\langle \beta^2 \rangle\langle \beta \rangle_{\lambda,\gamma}^2
-\langle \beta \rangle_{\lambda,\gamma}^4]={-}\E[\Var[\beta_0 \mid
\gamma\beta_0+\sqrt{\gamma}z]^2]. 
\end{equation}
Applying this to (\ref{eq:phisecondder}) completes the proof.
\end{proof}


\begin{proof}[Proof of Lemma~\ref{lemma:hessian_computation}]
The expressions for second-order derivatives involving $(q,r)$ follow directly
from differentiating (\ref{eq:grad_q_f}--\ref{eq:grad_v_f}) and specializing to
$(q_\star,r_\star;\alpha_\star,\tau_\star,\gamma_\star)$. In the remainder of
the proof, we describe the calculation of the lower-right $3 \times 3$ submatrix
$\nabla^2_{\alpha,\tau,\gamma} f$.

Recall that in an open neighborhood of $(\alpha_\star,\tau_\star,\gamma_\star)$,
the infimum and supremum in (\ref{eq:def_f_2}) are realized at solutions of the
equations (\ref{eq:stationary1}--\ref{eq:stationary2}).
As in the proof of Lemma \ref{lem:second_derivative_phi},
write $\kappa_j[\cdot \mid \lambda_\star]$ for the $j^\text{th}$ mixed and pure
cumulants under $\sP_{\lambda_\star,\gamma_\star}$. Then,
differentiating (\ref{eq:stationary1}--\ref{eq:stationary2})
implicitly in $(\alpha,\tau,\gamma)$
and specializing to $\tau_\star=\gamma_\star$ and
$\alpha_\star=\sqrt{\gamma_\star}$, we have
\[\partial_\alpha\lambda_\star=z, \quad
\partial_\tau\lambda_\star=-m_\star+\beta_0, \quad
\partial_\gamma\lambda_\star=m_\star, \quad
\partial_\alpha m_\star=\kappa_2[\beta \mid \lambda_\star] \cdot z,
\quad \partial_\tau m_\star=\kappa_2[\beta \mid \lambda_\star]
\cdot (-m_\star+\beta_0),\]
\[\partial_\gamma m_\star=\kappa_2[\beta \mid \lambda_\star] \cdot m_\star
-\frac{1}{2}\kappa_2[\beta,\beta^2 \mid \lambda_\star]
=-\frac{1}{2}\kappa_3[\beta \mid \lambda_\star].\]
Here, this second equality for $\partial_\gamma m_\star$ follows from applying
$m_\star=\langle \beta \rangle_{\lambda_\star,\gamma_\star}$ from
(\ref{eq:stationary2}) and expanding also
$\kappa_2[\cdot \mid \lambda_\star]$ and $\kappa_3[\cdot \mid \lambda_\star]$
in moments. Taking derivatives of (\ref{eq:grad_tau_f}--\ref{eq:grad_alpha_f})
and applying these formulas, we get
\begin{align}\label{eqn:3by3matrixhessianf}
\nabla^2_{\alpha,\tau,\gamma}f(\alpha_\star,\tau_\star,\gamma_\star)
=\E\begin{pmatrix}
-z^2\kappa_2[\beta \mid \lambda^*]-\sigma^2 &
z(m_\star-\beta_0)\kappa_2[\beta \mid \lambda_\star] & 
\frac{1}{2}z \kappa_3[\beta \mid \lambda_\star]\\
z(m_\star-\beta_0)\kappa_2[\beta \mid \lambda_\star] &
-(m_\star-\beta_0)^2\kappa_2[\beta \mid \lambda_\star] & 
 -\frac{1}{2}(m_\star-\beta_0) \kappa_3[\beta \mid \lambda_\star]\\
\frac{1}{2}z \kappa_3[\beta \mid \lambda_\star] &
 -\frac{1}{2}(m_\star-\beta_0) \kappa_3[\beta \mid \lambda_\star] &  
\partial_\gamma^2 f
\end{pmatrix}.
\end{align}
Differentiating (\ref{eq:stationary1}--\ref{eq:stationary2}) in $z$, we have
$\partial_z \lambda_\star=\sqrt{\gamma_\star}$ and $\partial_z
m_\star=\sqrt{\gamma_\star}\kappa_2[\beta \mid \lambda_\star]$. Then
applying $m_\star=\langle \beta \rangle_{\lambda_\star,\gamma_\star}$
from (\ref{eq:stationary2}) together with
the tower property of conditional expectation and Stein's lemma, we get
\begin{equation}\label{eq:hessiancalculations}
\begin{aligned}
\E[(m_\star-\beta_0)^2\kappa_2[\beta \mid \lambda_\star]]&=\E\big[\kappa_2[\beta
\mid \lambda_\star\big]^2]\\
\E[(m_\star-\beta_0)\kappa_3[\beta \mid \lambda_\star]]&=\E\big[
\E[m_\star-\beta_0 \mid \lambda_\star]\kappa_3[\beta \mid \lambda_\star]\big] =0\\
\E[z \kappa_3[\beta \mid \lambda_\star]]&=\E[\partial_z \kappa_3[\beta \mid
\lambda_\star]]=\sqrt{\gamma_\star}\,\E\big[\kappa_4[\beta \mid
\lambda_\star]\big]\\
\E[z(m_\star-\beta_0)\kappa_2[\beta \mid \lambda_\star]]&=\E[(\partial_z
m_\star) \kappa_2[\beta \mid \lambda_\star] + (m_\star-\beta_0)\partial_z \kappa_2[\beta
\mid \lambda_\star]]\\
&=\sqrt{\gamma_\star}\,\E\big[\kappa_2[\beta \mid
\lambda_\star]^2+(m_\star-\beta_0)\kappa_3[\beta \mid \lambda_\star]\big]
=\sqrt{\gamma_\star}\,\E\big[\kappa_2[\beta \mid \lambda_\star]^2\big]\\
\E[z^2 \kappa_2[\beta \mid \lambda_\star]]&=\E[\kappa_2[\beta \mid
\lambda_\star]+\partial_z^2 \kappa_2[\beta \mid \lambda_\star]]
=\E\big[\kappa_2[\beta \mid \lambda_\star] +\gamma_\star \kappa_4[\beta \mid
\lambda_\star] \big]. 
\end{aligned}
\end{equation}
Identifying $\E[\kappa_2[\beta \mid \lambda_\star]]+\sigma^2
=q_\star+\sigma^2=\delta/\gamma_\star$, and
identifying the definitions of $b_\star,k_\star$ in Lemma
\ref{lemma:hessian_computation} as
$b_\star=\E[\kappa_2[\beta \mid \lambda_\star]^2]$ and
$k_\star=\E[\kappa_4[\beta \mid \lambda_\star]]$, we obtain the desired form
for all but the lower right entry of (\ref{eqn:3by3matrixhessianf}).
For this lower right entry, by (\ref{eq:grad_gamma_f}) and the above identity
$\partial_\gamma \lambda_\star=m_\star$,
\begin{align*}
\partial_\gamma^2 f
=\frac{1}{2}\partial_\gamma \E[\kappa_2[\beta \mid \lambda_\star,\gamma_\star]]
&=\E\left[\frac{1}{2}\,m_\star \kappa_3[\beta \mid \lambda_\star]
-\frac{1}{4}\kappa_3[\beta,\beta,\beta^2 \mid \lambda_\star]\right]
=\E\left[{-}\frac{1}{4}\kappa_4[\beta \mid \lambda_\star]
-\frac{1}{2}\kappa_2[\beta \mid \lambda_\star]^2\right]
\end{align*}
where this last equality applies (\ref{eq:mmseprimetmp}) and the final form of
$\mmse'(\gamma)$ computed in (\ref{eq:mmseprime}). This is $-b_\star/2-k_\star/4$, completing the proof.
\end{proof}

Finally, let us apply Lemma \ref{lm:tap_lower_bound_function} to show
Corollary \ref{cor:lower_TAP_f_star} stated at the start of this section.

\begin{proof}[Proof of Corollary~\ref{cor:lower_TAP_f_star}]
Let $K' \subset [0,\infty) \times \R^2$ be a compact set containing
$(\alpha_\star,\tau_\star,\gamma_\star)=(\sqrt{\gamma_\star},\gamma_\star,\gamma_\star)$
in its interior.
By Lemma~\ref{lm:tap_lower_bound_function}(b), there exists $\rho_0>0$
such that 
\[
\sup_{(\alpha,\tau,\gamma) \in K'} f(q,r;\alpha,\tau,\gamma) \geq \phi(\gamma_\star)
~~~~ \text{for any } (q,r) \in K(\rho_0).
\]
Part (a) follows from Lemma \ref{lm:tap_lower_bound_body}, applied with this set $K'$
and with $K=K(\rho_0)$.

Also by Lemma~\ref{lm:tap_lower_bound_function}(b),
there exist $\rho_0>\rho_1>0$ and $\iota_0>0$ such that 
\[
\sup_{(\alpha,\tau,\gamma) \in K'} f(q,r; \alpha, \tau, \gamma) \geq
\phi(\gamma_\star) + 2 \iota_0, ~~~~ \text{for any } (q,r) \in
\overline{K(\rho_0) \setminus K(\rho_1)}.
\]
Applying Lemma \ref{lm:tap_lower_bound_body} with this $K'$ and with
$K=\overline{K(\rho_0) \setminus K(\rho_1)}$, with probability at
least $1-e^{-cn}$, $\inf_{(\bbm,\bs) \in \Gamma^p[K]} \cF_\TAP(\bbm,\bs)
\geq \phi(\gamma_\star) + 2\iota_0 - \iota_0$
which proves part (b). 
\end{proof}

\section{The Bayes-optimal TAP local minimizer}\label{app:proof_bayes_consistence}

We prove Theorems \ref{thm:bayes_consistence} and \ref{thm:marginalposterior},
which rely on the next three lemmas. Recall the marginal first and second
moment vectors $\bbm_\sB$ and $\bs_\sB$
of the posterior law from (\ref{eq:marginalmoments}). We will denote
\[\gamma_\star=\gamma_\stat, \qquad
q_\star=\mmse(\gamma_\star)=\gamma_\star/\delta-\sigma^2\]
and write $\bbm^2$ for the entrywise square of a vector $\bbm \in \R^p$.

\begin{lemma}[Concentration for Bayes posterior marginals]\label{lm:concen_bayes} 
Under Assumptions \ref{ass:Bayesian_linear_model} and \ref{ass:uniquemin},
\begin{align}
   p^{-1} \|\bbeta_0-\bbmB\|_2^2 \gotop&~ q_\star, \label{eqn:lm:concen_bayes1}\\
     Q(\bbmB)=p^{-1}\|\bbmB \|_2^2 \gotop&~ \E[\beta_0^2]-q_\star, \label{eqn:lm:concen_bayes2} \\
  S(\bsB)=p^{-1} \|\bsB\|_1 \gotop&~ \E[\beta_0^2].\label{eqn:lm:concen_bayes3}
  \end{align}
\end{lemma}

\begin{lemma}[Concentration for TAP stationary point]\label{lm:stationary_concen} 
Let Assumptions \ref{ass:Bayesian_linear_model} and \ref{ass:uniquemin} hold.
Then there exists $(\bbm_\star,\bs_\star) \in \Gamma^p$ that is a (measurable)
function of $(\X,\y)$, such that $(\bbm_\star,\bs_\star)$ is a local minimizer
of $\cF_\TAP(\bbm,\bs)$ with probability approaching 1, and
\begin{align}
p^{-1} \|\bbeta_0-\bbm_\star\|_2^2 \gotop&~ q_\star, \label{eqn:stationary_concen_lemma_eq1}\\
S(\bbs_\star)-Q(\bbm_\star)=p^{-1} \|\bbs_\star-\bbm_\star^2\|_1 \gotop&~ q_\star, \label{eqn:stationary_concen_lemma_eq2}\\
p^{-1}\cF_{\TAP} (\bbm_\star, \bs_\star) \gotop&~ \phi(\gamma_\star),
\end{align}
where $\phi$ is the replica-symmetric potential in (\ref{eqn:phi_potential}).
\end{lemma}

\begin{lemma}\label{lem:marginalposterior-in-proof}
Let Assumptions \ref{ass:Bayesian_linear_model} and \ref{ass:uniquemin}
hold. Suppose that with probability approaching 1, there exists a local
minimizer $(\bbm_\star,\bbs_\star) \in \Gamma^p$ of $\cF_\TAP(\bbm,\bbs)$
satisfying 
\begin{align}
p^{-1}  \|\bbm_\star-  \bbmB \|_2^2 \gotop&~ 0, \label{eqn:marginalposterior-in-proof-cd1}\\
p^{-1} \|\bbs_\star-\bbm_*^2\|_1 \gotop&~ q_\star. \label{eqn:marginalposterior-in-proof-cd2}
\end{align}
Then for any Lipschitz function $f:\supp(\sP_0) \to \R$
and any index $j \in \{1,\ldots,p\}$,
\[\E\Big[\Big(\langle f(\beta_j) \rangle_{\X,\y}-\langle
f(\beta)
\rangle_{\lambda(m_{\star,j},s_{\star,j}),\gamma(m_{\star,j},s_{\star,j})}
\Big)^2 \Big] \to 0,
\]
where $\< \cdot \>_{\X, \y}$ is the posterior average of $\bbeta$ given $(\X, \y)$, and $\< \cdot \>_{\lambda, \gamma}$ is the average under the law (\ref{eq:scalarmean}). 
\end{lemma}

The proofs of Lemmas \ref{lm:concen_bayes}, \ref{lm:stationary_concen}, and
\ref{lem:marginalposterior-in-proof} are contained in Section
\ref{sec:pf_lm:concen_bayes},  \ref{sec:pf_lm:stationary_concen}, and
\ref{app:proof_marginalposterior-in-proof}, respectively. Using these, we first
show Theorems \ref{thm:bayes_consistence} and \ref{thm:marginalposterior}.

\begin{proof}[Proof of Theorem~\ref{thm:bayes_consistence}]
Let $(\bbm_\star,\bs_\star) \in \Gamma^p$ be as specified in
Lemma \ref{lm:stationary_concen}, which is a TAP local
minimizer satisfying (\ref{eqn:TAP_energy_convergence})
with probability approaching 1. Moreover,
since $\bbm_\sB=\E[\bbeta_0 \mid \X,\y]$ and $\bbm_\star$ is a 
function of $(\X,\y)$, we have $\E\|\bbm_\star -\bbeta_0\|_2^2 = \E\|\bbm_\star -
\bbmB\|_2^2 + \E\|\bbmB-\bbeta_0\|_2^2$. Therefore, by Lemmas
\ref{lm:concen_bayes} and \ref{lm:stationary_concen} and the bounded convergence
theorem, we have $p^{-1} \|\bbs_\star-\bbm_*^2\|_1 \gotop q_\star$ and
\begin{align*}
   p^{-1}  \E\|\bbm_\star-\bbmB\|_2^2= p^{-1}
\big(\E\|\bbm_\star-\bbeta_0\|_2^2-\E\|\bbmB-\bbeta_0\|_2^2\big) \to q_\star - q_\star = 0.
\end{align*} 
By Markov's inequality, this verifies the needed conditions (\ref{eqn:marginalposterior-in-proof-cd1})
and (\ref{eqn:marginalposterior-in-proof-cd2}) of Lemma
\ref{lem:marginalposterior-in-proof}. Then
applying Lemma \ref{lem:marginalposterior-in-proof}
with $f (\beta) = \beta^2$, we obtain 
\[
\E \big[ \big| (\bsB)_j - s_{\star,j}\big| \big]
\leq \E \big[ \big( (\bsB)_j - s_{\star,j}\big)^2 \big]^{1/2}
= \E \big[ \big(\langle \beta_j^2 \rangle_{\X,\y}-\langle \beta^2 \rangle_{\lambda(m_{\star,j},s_{\star,j}),\gamma(m_{\star,j},s_{\star,j})}
\big)^2 \big]^{1/2} \to 0,
\] 
where these expectations are the same for every index
$j \in \{1,\ldots,p\}$ by symmetry of the linear model across coordinates.
Thus also $p^{-1} \E[\|\bbs_\star-\bsB\|_1] \to0$, which implies
(\ref{eq:Bayesoptimalcrit}) by Markov's inequality.
\end{proof}

\begin{proof}[Proof of Theorem \ref{thm:marginalposterior}]
We check the conditions (\ref{eqn:marginalposterior-in-proof-cd1}) and
(\ref{eqn:marginalposterior-in-proof-cd2}) of Lemma
\ref{lem:marginalposterior-in-proof}. Note that
(\ref{eqn:marginalposterior-in-proof-cd1}) is implied by the assumption that
$(\bbm_\star, \bs_\star)$ satisfies (\ref{eq:Bayesoptimalcrit}). Furthermore,
applying (\ref{eqn:lm:concen_bayes2}) and (\ref{eqn:lm:concen_bayes3}) in Lemma
\ref{lm:concen_bayes}, the assumption that $(\bbm_\star,
\bs_\star)$ satisfies (\ref{eq:Bayesoptimalcrit}), and the condition 
$\|\bs_\star-\bbm_\star^2\|_1=S(\bs_\star)-Q(\bbm_\star)$ because
$\bs \geq \bbm^2 \geq 0$ entrywise for any $(\bbm,\bs) \in \Gamma$, we get
\[\Big\vert p^{-1} \| \bs_\star - \bbm_\star^2 \|_1 - q_\star \Big\vert
\le \Big\vert S(\bsB) - \E[\beta_0^2] \Big\vert + \Big\vert 
Q(\bbmB) - (\E[\beta_0^2] - q_\star) \Big\vert + |S(\bs_\star)-S(\bsB)|
+ |Q(\bbm_\star)-Q(\bbmB)| \to 0.\]
This proves that $(\bbm_\star, \bs_\star)$ satisfies
(\ref{eqn:marginalposterior-in-proof-cd2}). Then the desired result follows
from Lemma \ref{lem:marginalposterior-in-proof}.
\end{proof}

\subsection{Proof of Lemma~\ref{lm:concen_bayes}}\label{sec:pf_lm:concen_bayes}

This lemma is proved by combining and adapting several results of
\cite{barbier2019optimal}. We write as shorthand $\langle \cdot \rangle=\langle
\cdot \rangle_{\X,\y}$ for the joint posterior expectation over independent
replicas $\bbeta,\bbeta_1,\bbeta_2$ fixing $\X,\bbeta_0,\y$. Thus $\bbmB=\langle
\bbeta \rangle$ and $\bsB=\langle \bbeta^2 \rangle$.\\

\noindent
{\bf Step 1. Convergence of the overlap.} We first show that
\begin{align}\label{eq:pf_lm:concen_bayes_claim1}
\big\< \big( p^{-1} \bbeta_0^\top \bbeta  - ( \E[\beta_0^2]-q_\star) \big)^2 \big\> \gotop 0. 
\end{align}
Indeed, let $Q:=Q(\bbeta_0,\bbeta)=p^{-1} \bbeta_0^\top\bbeta$.
By Theorem 2 of~\cite{barbier2019optimal} (where it is straightforward
to check that conditions (h1)-(h5) therein are satisfied), we have
\begin{align}
\big\< \big( |Q|  - ( \E[\beta_0^2]-q_\star) \big)^2
\big\> \gotop 0.  \label{eq:pf_lm:concen_bayes_claim1_pf1}
\end{align}
This equation is slightly different from (\ref{eq:pf_lm:concen_bayes_claim1}) by
the absolute value operator inside $\langle \cdot \rangle$. We may remove this
absolute value as follows: Observe that for a new sample
$(y_\mathrm{new},\x_\mathrm{new})$, the Bayes-optimal generalization error is
$\E[(y_\mathrm{new}-\E[y_\mathrm{new} \mid \X,\y,\x_\mathrm{new}])^2]
=\E[(\x_\mathrm{new}^\top(\bbeta_0-\langle \bbeta \rangle))^2]+\sigma^2
=p^{-1}\E[\|\bbeta_0-\langle \bbeta \rangle\|_2^2]+\sigma^2$. Then
Theorem 3 of~\cite{barbier2019optimal} specialized to the linear model
implies that $p^{-1} \E\|\bbeta_0 -\<\bbeta\> \|_2^2 \to  q_\star$.
Using Nishimori's identity $\E \langle \bbeta_0^\top \bbeta \rangle
=\E \langle \bbeta_1^\top \bbeta_2 \rangle=\E\|\langle \bbeta \rangle\|_2^2$,
we have
\begin{align}
\E\<Q\>=p^{-1}\E\|\<\bbeta\> \|_2^2=p^{-1}\E\| \bbeta_0
\|_2^2-p^{-1}\E\|\bbeta_0-\<\bbeta\> \|_2^2 \to \E[\beta_0^2]-q_\star\label{eq:pf_lm:concen_bayes_claim1_pf2}.
\end{align}
Note that $Q$ is uniformly bounded, by boundedness of the supports of
$\bbeta_0,\bbeta$. Then,
taking expectation in \eqref{eq:pf_lm:concen_bayes_claim1_pf1} using the bounded
convergence theorem and combining with
\eqref{eq:pf_lm:concen_bayes_claim1_pf2}, we must have
$\E\<||Q|-Q|\>=\E\<|Q|\>-\E\<Q\>
 \to 0$, where we have used that $||Q| - Q| = |Q| - Q$ because $Q \leq |Q|$. Applying this and boundedness of $Q$
back to \eqref{eq:pf_lm:concen_bayes_claim1_pf1}
gives $\E \langle (Q-(\E[\beta_0^2]-q_\star))^2
\rangle \to 0$, which shows (\ref{eq:pf_lm:concen_bayes_claim1}) by Markov's
inequality.\\

\noindent
{\bf Step 2. Proof of (\ref{eqn:lm:concen_bayes1}) and
(\ref{eqn:lm:concen_bayes2}).} Note that \eqref{eq:pf_lm:concen_bayes_claim1},
\eqref{eq:pf_lm:concen_bayes_claim1_pf2}, and boundedness of $Q$
imply $\E\<|Q-\E\<Q\>|\> \leq \E[\< (Q-\E\< Q\>)^2 \>^{1/2}] \to 0$.
As a consequence, applying $\<Q\>=p^{-1} \bbeta_0^\top \<\bbeta\>$,
\begin{align*} 
p^{-1}  \E\Big| \|\<\bbeta \>\|_2^2 -\E\|\<\bbeta \>\|_2^2\Big|
&=p^{-1}\E \Big| \< \bbeta_1^\top \bbeta_2\> - \E\|\<\bbeta\>\|_2^2 \Big|
\leq 
p^{-1}\E \Big\<\Big|\bbeta_1^\top \bbeta_2-\E\|\<\bbeta\>\|_2^2 \Big|\Big\>\\
&=p^{-1} \E \Big\<\Big| \bbeta_0^\top \bbeta-\E\bbeta_0^\top \<\bbeta\> \Big|
\Big\>=\E\<|Q -\E \<Q\>|\>\to0,
\end{align*}
where the inequality is Jensen's inequality for $\langle \cdot \rangle$,
and the first equality on the second line uses Nishimori's identity. This implies the
concentration of $p^{-1}\|\bbmB\|_2^2=p^{-1}\|\<\bbeta\>\|^2_2$, and combining with
(\ref{eq:pf_lm:concen_bayes_claim1_pf2}) proves (\ref{eqn:lm:concen_bayes2}).
Furthermore, combining (\ref{eqn:lm:concen_bayes2}), 
$\<Q\>-\E\<Q\> \gotop 0$, (\ref{eq:pf_lm:concen_bayes_claim1_pf2}),
and $p^{-1} \| \bbeta_0 \|_2^2 \gotop \E[\beta_0^2]$ proves
(\ref{eqn:lm:concen_bayes1}).\\

\noindent
{\bf Step 3. Proof of (\ref{eqn:lm:concen_bayes3}).} Note that by the entrywise
positivity $\bbeta^2 \geq 0$, we have
$\|\bsB\|_1=\langle \|\bbeta^2 \|_1 \rangle$.
Then, applying again Jensen's inequality and Nishimori's identity, we obtain
\[\E \Big| p^{-1}   \|\bsB\|_1  -\E[ \beta_0^2 ] \Big|
 \le \E  \Big\<  \Big| p^{-1}  \| \bbeta^2  \|_1  -\E[ \beta_0^2 ] \Big| \Big\>     = \E \Big| p^{-1} \|\bbeta_0^2\|_1 -\E[ \beta_0^2  ] \Big|  \to0.
\]
This proves (\ref{eqn:lm:concen_bayes3}) and thus finishes the proof of Lemma~\ref{lm:concen_bayes}.

\subsection{Proof of Lemma~\ref{lm:stationary_concen}}\label{sec:pf_lm:stationary_concen}

The proof will combine the lower bounds for $\cF_\TAP$ obtained in
Corollary \ref{cor:lower_TAP_f_star} with an upper bound for $\cF_\TAP$
evaluated at a point near the Bayes posterior marginal vectors $(\bbmB,\bsB)$.

To control the entropy term $D_0(\bbm,\bs)$ of $\cF_\TAP$, we use the following
truncation of $(\bbmB,\bsB)$: Define the domain
\[\Gamma_M=\Big\{(m,s):m=\langle \beta \rangle_{\lambda,\gamma},\;
s=\langle \beta^2 \rangle_{\lambda,\gamma} \text{ for some }
(\lambda,\gamma) \in [-M,M]^2\Big\}.\]
Let $\Proj_M(m,s)$ be the projection in Euclidean distance onto
$\Conv(\Gamma_M)$, the convex hull of $\Gamma_M$. Note that $\Gamma_M$ is the
continuous image of a compact set, so both $\Gamma_M$ and
$\Conv(\Gamma_M)$ are compact, and $\Proj_M$ is uniquely defined.
We write as shorthand also $\Proj_M(\bbm,\bs)$
for the application of $\Proj_M$ to each coordinate pair $(m_j,s_j)$, and denote
\[(\bbmB^M,\bsB^M)=\Proj_M(\bbmB,\bsB).\]
Observe that $\{\Conv(\Gamma_M)\}_{M \geq 0}$ is an increasing family such
that $\cup_M \Conv(\Gamma_M) = \Gamma$. The closure $\widebar\Gamma$ is compact,
and therefore
\begin{align}
\sup_{(m,s)\in \widebar\trunset} \|(m,s) - \ProjM(m,s) \|_2 \to 0 ~~\text{as
}M\to\infty.\label{eq:Projapprox}
\end{align}
This provides a uniform approximation of $(m,s) \in \Gamma$ by $\Proj_M(m,s)$,
which we will use throughout the proof.

\begin{lemma}\label{lm:stationary_concen_tech}
Suppose Assumptions \ref{ass:Bayesian_linear_model} and \ref{ass:uniquemin}
hold. Then for any $\eps>0$, there exists $M>0$ such that for all large $n,p$,
\[p^{-1}\E[\cF_\TAP(\bbmB^M,\bsB^M)] \leq \phi(\gamma_\star)+\eps.\]
\end{lemma}
\begin{proof}
We write $\eps(M)$ for a positive constant satisfying
$\eps(M) \to 0$ as $M \to \infty$ and changing from instance to instance.

Recall the form of $\cF_\TAP$ from (\ref{eqn:TAP}).
By Lemma~\ref{lm:concen_bayes}, $Q(\bbmB) \gotop \E[\beta_0^2]-q_\star$
and $S(\bbmB) \gotop \E[\beta_0^2]$. Then by the bounded convergence theorem
and the approximation (\ref{eq:Projapprox}), also
\begin{equation}\label{eq:Onsagerlimit}
\limsup_{n,p \to \infty} \frac{n}{2p}\,\E \log(\sigma^2+S(\bsB^M)-Q(\bbmB^M))
\leq \frac{\delta}{2} \log (\sigma^2+q_\star)+\eps(M)
=\frac{\delta}{2}\log \frac{\delta}{\gamma_\star} + \eps(M).
\end{equation}
By Corollary \ref{cor:ymmseconc}, the approximation (\ref{eq:Projapprox}),
and the bound $\E\|\X\|_{\op}^2 \leq C$ for a constant $C>0$, also
\begin{equation}\label{eq:ymmselimit}
\limsup_{n,p \to \infty}
\frac{1}{2\sigma^2} \cdot \frac{1}{p}\,\E\|\y-\X\bbmB^M\|_2^2 
\leq \frac{\gamma_\star\sigma^2}{2}+\eps(M).
\end{equation}

It remains to evaluate the limit of $p^{-1}\E[D_0(\bbmB^M,\bsB^M)]$. For this,
we apply the cavity approximation of Lemma \ref{lemma:cavity} and a truncation
of the entropy function. Let $\hat{\lambda}_j$ be as defined in
(\ref{eq:hatlambdadef}), copied here for readers' convenience:
\begin{equation*}
\hat{\lambda}_j=\gamma_\star \beta_{0,j}
-\frac{1}{\sigma^2}{\x^j}^\top(\X^{-j} \langle \bbeta^{-j} \rangle_{-j}
-\y^{-j}).
\end{equation*}
Here, $\x^j,\X^{-j}$ are the $j^\text{th}$ and all-but-$j^\text{th}$ columns of
$\X$, similarly for $\beta_j,\bbeta^{-j}$ and $\beta_{0,j},\bbeta_0^{-j}$,
and $\y^{-j}=\X^{-j}\bbeta_0^{-j}+\beps$ and
$\<f(\bbeta^{-j})\>_{-j}=\E[\bbeta_0^{-j} \mid \y^{-j},\X^{-j}]$.
Define in addition
\[\bar{\lambda}_j=\gamma_\star \beta_{0,j}-\frac{\sqrt{p \cdot
\gamma_\star}}{\|(\X^{-j}\<\bbeta^{-j}\>_{-j}
-\y^{-j})\|_2}{\x^j}^\top(\X^{-j}\<\bbeta^{-j}\>_{-j} - \y^{-j}).\]
We denote
\[\hat{m}_j=\<\beta\>_{\hat\lambda_j,\gamma_\star},
\quad \hat{s}_j=\<\beta^2\>_{\hat\lambda_j,\gamma_\star},
\quad \bar{m}_j=\<\beta\>_{\bar\lambda_j,\gamma_\star},
\quad \bar{s}_j=\<\beta^2\>_{\bar\lambda_j,\gamma_\star},\]
and $(\hat m_j^M,\hat s_j^M)=\Proj_M(\hat m_j,\hat s_j)$ and
$(\bar m_j^M,\bar s_j^M)=\Proj_M(\bar m_j,\bar s_j)$.

By these definitions and by Lemma \ref{lemma:cavity} applied
with $f(\beta)=\beta$ and $f(\beta)=\beta^2$, we have
\[\E[((\bbmB)_j-\hat m_j)^2] \to 0, \qquad \E[((\bsB)_j-\hat s_j)^2] \to 0,\]
where these and all subsequence expectations are the same for every $j \in
\{1,\ldots,p\}$ by coordinate symmetry of the model. Note that
$\|\nabla_{\lambda}(\<\beta\>_{\lambda,\gamma_\star},\<\beta^2\>_{\lambda,\gamma_\star})\|_2\leq
C$ for a constant $C:=C(\sP_0)>0$, since $\sP_0$ is compactly supported.
Then
\begin{align*} 
\E\Big\|(\hat m_j, \hat s_j)-(\bar m_j,\bar s_j)\Big\|_2^2
&\leq C \cdot \E(\widehat\lambda_j-\widebar\lambda_j)^2\\
&= C \cdot \E\left[ \left\vert
\left(\frac{1}{\sigma^2}-\frac{\sqrt{p \cdot \gamma_\star}}{\|(\X^{-j}\<\bbeta^{-j}\>_{-j}
-\y^{-j})\|_2}\right){\xb^j}^\top(\X^{-j}\<\bbeta^{-j}\>_{-j} -\y^{-j})
\right\vert^2 \right] \\
&= C \cdot \E\left[\left|\frac{\|\X^{-j}\<\bbeta^{-j}\>_{-j}
-\y^{-j}\|_2}{\sigma^2 \sqrt{p}}-{\sqrt{\gamma_\star}}\right|^2\right]
\end{align*}
where the last line first evaluates the expectation over $\xb^j \sim
\mathcal{N}(0,p^{-1}\bI)$, noting that it is independent of all other variables.
This upper bound vanishes as $n,p \to \infty$, by Corollary \ref{cor:ymmseconc}
applied to the leave-one-out model $\y^{-j} = \X^{-j} \bbeta^{-j} + \beps$. Then
\begin{equation}\label{eq:msapprox}
\E[\|((\bbmB)_j,(\bsB)_j)-(\bar m_j,\bar s_j)\|_2^2] \to 0.
\end{equation}

Observe that since $\Conv(\Gamma_M)$ is compact and
${-}\hs(m,s)$ is continuous on $\Conv(\Gamma_M)$, it must be uniformly
continuous on $\Conv(\Gamma_M)$, so there exists a decreasing bounded
function $f:[0,\infty) \to \R$ with $f(0)=0$ and continuous at 0 such that
\[\Big|{-}\hs((\bbmB)_j^M,(\bsB)_j^M)+\hs(\bar m_j^M,\bar s_j^M)\Big|
\leq f(\|((\bbmB)_j^M,(\bsB)_j^M)-(\bar m_j^M,\bar s_j^M)\|_2).\]
Since $\Conv(\Gamma_M)$ is also convex and thus $\Proj_M$ is 1-Lipschitz,
we have
\[\|((\bbmB)_j^M,(\bsB)_j^M)-(\bar m_j^M,\bar s_j^M)\|_2
\leq \|((\bbmB)_j,(\bsB)_j)-(\bar m_j,\bar s_j)\|_2.\]
Combining these statements, and applying (\ref{eq:msapprox}) and the bounded
convergence theorem,
\[\lim_{n,p \to \infty}
\E\Big|{-}\hs((\bbmB)_j^M,(\bsB)_j^M)+\hs(\bar m_j^M,\bar s_j^M)\Big|=0.\]
Now observe that $\bar{\lambda}_j$ is equal in law to
$\lambda_\star=\gamma_\star \beta_0+\sqrt{\gamma_\star} z$ for
$(\beta_0,z) \sim \sP_0 \times \mathcal{N}(0,1)$, because
$\xb^j \sim \mathcal{N}(0,p^{-1}\bI)$ is independent of the remaining variables
defining $\bar \lambda_j$. Set $(m_\star,s_\star)=(\langle \beta
\rangle_{\lambda_\star,\gamma_\star},\langle \beta^2
\rangle_{\lambda_\star,\gamma_\star})$ as functions of $(\beta_0,z)$,
and let $(m_\star^M,s_\star^M)=\Proj_M(m_\star,s_\star)$. Then
$\E[{-}\hs(\bar m_j^M,\bar s_j^M)]=\E[{-}\hs(m_\star^M,s_\star^M)]$,
where this quantity does not depend on $n$ and $p$. Thus, the above implies
\begin{equation}\label{eq:happrox}
\lim_{n,p \to \infty}
\Big|\E[{-}\hs((\bbmB)_j^M,(\bsB)_j^M)]-\E[{-}\hs(m_\star^M,m_\star^M)]\Big|=0.
\end{equation}

Finally, observe that
\begin{align*}
\Big|\E[{-}\hs(m_\star^M,s_\star^M)]+\E[\hs(m_\star,s_\star)]\Big|
&\leq \E\Big[\big(|\hs(m_\star^M,s_\star^M)|+|\hs(m_\star,s_\star)|\big)
\cdot \1\{(m_\star,s_\star) \notin \Conv(\Gamma_M)\}\Big]\\
&\leq \E\Big[\big(|\hs(m_\star^M,s_\star^M)|+|\hs(m_\star,s_\star)|\big)
\cdot \1\{\lambda_\star \notin [-M,M]\}\Big].
\end{align*}
Consider the map
$f(\lambda,\gamma)={-}\sh(\langle \beta \rangle_{\lambda,\gamma},\langle
\beta^2 \rangle_{\lambda,\gamma})$ whose gradient is, by the chain rule,
\[\nabla f(\lambda,\gamma)
=\frac{\partial(m,s)}{\partial(\lambda,\gamma)}^\top
\nabla[{-}\hs(m,s)]\]
where $(m,s)=(\langle \beta \rangle_{\lambda,\gamma},\langle
\beta^2 \rangle_{\lambda,\gamma})$. For $(\lambda,\gamma) \in [-M,M]^2$, we have
$\|\nabla[{-}\hs(m,s)]\|_2=\|(\lambda,-\frac{1}{2}\gamma)\|_2 \leq 2M$.
Since $\|\frac{\partial(m,s)}{\partial(\lambda,\gamma)}\|$ is also bounded
because $\sP_0$ has bounded support, and since $f(0,0)=0$,
for a constant $C>0$ we obtain $|f(\lambda,\gamma)| \leq CM$ for all
$(\lambda,\gamma) \in [-M,M]^2$. Thus,
\[|\hs(m_\star^M,s_\star^M)| \leq CM, \qquad
|\hs(m_\star,s_\star)| \leq C|\lambda_\star|.\]
Applying this above yields
\[\lim_{M \to \infty}
\Big|\E[{-}\hs(m_\star^M,s_\star^M)]-\E[{-}\hs(m_\star,s_\star)]\Big|
\leq \lim_{M \to \infty}
\E\Big[C(M+\lambda_\star) \cdot \1\{\lambda_\star \notin [-M,M]\}\Big]=0,\]
where the last statement holds under the
law $\lambda_\star=\gamma_\star\beta_0+\sqrt{\gamma_\star} z$. Then we get
\[\lim_{M \to \infty} \E[{-}\hs(m_\star^M,s_\star^M)]
=\E[{-}\hs(m_\star,s_\star)]=\E_{\beta_0,z}[D_{\mathrm{KL}}(\sP_{\lambda_\star,\gamma_\star}\|\sP_0)]=i(\gamma_\star),\]
where we have recalled that $\sP_{\lambda_\star,\gamma_\star}$ is the
posterior distribution of $\beta_0$ given
$\lambda_\star=\gamma_\star\beta_0+\sqrt{\gamma_\star} z$, and thus
this limit is the mutual information $i(\gamma_\star)$ between $\beta_0$
and $\lambda_\star$. Then (\ref{eq:happrox}) implies
\begin{equation}\label{eq:hlimit}
\limsup_{n,p \to \infty} p^{-1}\E[D_0(\bbmB^M,\bsB^M)]
=\limsup_{n,p \to \infty} \E[{-}\hs((\bbmB)_j^M,(\bsB)_j^M)]
\leq \E[{-}\hs(m_\star^M,s_\star^M)] \leq i(\gamma_\star)+\eps(M).
\end{equation}

Combining (\ref{eq:Onsagerlimit}), (\ref{eq:ymmselimit}), and (\ref{eq:hlimit}),
\[\limsup_{n,p \to \infty} p^{-1}\E[\cF_\TAP(\bbmB^M,\bsB^M)]
\leq {-}\frac{\delta}{2}\log \frac{\gamma_\star}{2\pi\delta}
+\frac{\gamma_\star
\sigma^2}{2}+i(\gamma_\star)+\eps(M)=\phi(\gamma_\star)+\eps(M),\]
and taking $M>0$ sufficiently large proves Lemma~\ref{lm:stationary_concen_tech}.
\end{proof}

\begin{corollary}\label{cor:tap_upperbound_proj_bayes}
Let Assumptions \ref{ass:Bayesian_linear_model} and \ref{ass:uniquemin} hold.
Then for any $\iota,\eps>0$, there exists $M=M(\iota,\eps)>0$ such that for all
large $n,p$,
\[\P\Big[p^{-1}\cF_\TAP(\bbmB^M,\bsB^M) \geq \phi(\gamma_\star)+\iota\Big]<\eps.\]
\end{corollary}
\begin{proof}
By Corollary \ref{cor:lower_TAP_f_star}, there exists $\rho_0>0$ such that
for any $\eps_0>0$, with probability approaching 1,
\[
\inf_{(\bbm, \bs) \in \Gamma^p[K(\rho_0)]} p^{-1} \cF_{\TAP}(\bbm, \bs) \ge
\phi(\gamma_\star)-\eps_0.
\]
Here, $K(\rho_0)$ is the ball of radius $\rho_0$ around
$(q_\star,q_\star)$, and $\Gamma^p[K]$ is the domain defined in
(\ref{eq:GammaK}). By Lemma~\ref{lm:concen_bayes} and (\ref{eq:Projapprox}), for
all sufficiently large $M$, we have $\P ((\bbmB^M,\bsB^M) \in \Gamma^p[K(\rho_0)] )
\to 1$ as $n,p \to \infty$. This implies that for any $\eps_0>0$ and all
sufficiently large $M>0$, the event
\[
\cE = \Big\{ p^{-1} \cF_{\TAP}(\bbmB^M,\bsB^M) \geq \phi(\gamma_\star) - \eps_0
\Big\}
\]
has probability approaching 1 as $n,p \to \infty$. Then by Lemma
\ref{lm:stationary_concen_tech}, for sufficiently large $M:=M(\eps_0)$,
\begin{align*}
\limsup_{n,p \to \infty} \E[ p^{-1} \cF_{\TAP}(\bbmB^M,\bsB^M) \mid \cE] \le
\limsup_{n,p \to \infty} \frac{\E[p^{-1} \cF_{\TAP}(\bbmB^M,\bsB^M)]}{\P[\cE]}
\leq \phi(\gamma_\star)+\eps_0.
\end{align*}
As a consequence, it follows from Markov's inequality and the definition of
$\cE$ that for any $\iota>0$ and all large $n,p$,
\begin{align}
    \P \Big[ p^{-1} \cF_{\TAP}(\bbmB^M,\bsB^M) \geq \phi(\gamma_\star) +
\iota \;\Big|\; \cE \Big] \leq \frac{\E[p^{-1}
\cF_{\TAP}(\bbmB^M,\bsB^M)-\phi(\gamma_\star)+\eps_0 \mid \cE]}{\iota+\eps_0}
\leq \frac{2\eps_0}{\iota+\eps_0}.
\end{align}
For any $\eps>0$, choosing $\eps_0:=\eps_0(\iota,\eps)$ sufficiently small
ensures that this probability is at most $\eps/2$. Together with the statement
$\P[\cE^c] \leq \eps/2$ for all large $n,p$, this proves the lemma.
\end{proof}

\begin{proof}[Proof of Lemma~\ref{lm:stationary_concen}]
Let $K(\rho)$ be the ball of radius $\rho$ around $(q_\star,q_\star)$ as
defined in Corollary \ref{cor:lower_TAP_f_star}, and
let $\cE:=\cE(\rho_0,\rho_1,\iota,M)$ denote the event on which the
following conditions hold:
\begin{itemize}
\item[(a)] We have the TAP free energy lower bound
\begin{align}
\inf_{(\bbm,\bbs) \in \Gamma^p[K(\rho_0) \setminus K(\rho_1)]} 
p^{-1} \cF_{\TAP}(\bbm, \bbs) &\ge \phi(\gamma_\star)+\iota,
\label{eq:TAPlowerstatproof1}\\
\inf_{(\bbm,\bbs) \in \Gamma^p[K(\rho_1)]}
p^{-1}\cF_{\TAP}(\bbm, \bbs) &\ge \phi(\gamma_\star)-
\iota.\label{eq:TAPlowerstatproof2}
\end{align}
\item[(b)] We have the TAP free energy upper bound
\[
p^{-1} \cF_{\TAP}(\bbmB^M,\bsB^M) \leq \phi(\gamma_\star)+\iota/2.
\]
\item[(c)] We have $(\bbmB^M,\bsB^M) \in \Gamma^p[K(\rho_1)]$, i.e.\
$(p^{-1}\|\bbmB^M-\bbeta_0\|_2^2,p^{-1}\|\bsB^M-(\bbmB^M)^2\|_1)
\in K(\rho_1)$.
\end{itemize}
For any sufficiently small $\rho_0>\rho_1>0$, some
$\iota_0:=\iota_0(\rho_0,\rho_1)>0$, and any $\iota \in (0,\iota_0)$,
Corollary \ref{cor:lower_TAP_f_star} ensures that condition (a) holds with
probability approaching 1 as $n,p \to \infty$. Then for any $\eps>0$,
Corollary \ref{cor:tap_upperbound_proj_bayes} ensures that
there exists $M:=M(\iota,\eps)>0$ where condition (b) holds with
probability at least $1-\eps/2$ for all large $n,p$.
By Lemma \ref{lm:concen_bayes} and the approximation (\ref{eq:Projapprox}),
choosing $M$ large enough ensures also that condition (c) holds with
probability approaching 1 as $n,p \to \infty$. Hence $\P[\cE]>1-\eps$ for all
large $n,p$.

Now consider the sub-level set
\[\cS=\{(\bbm,\bs) \in
\Gamma^p:p^{-1}\cF_\TAP(\bbm,\bs)<\phi(\gamma_\star)+3\iota/4\},\]
which is open
in $\Gamma^p$. We define $\cE'$ as the following event: $(\bbmB^M,\bsB^M) \in \cS$,
and the connected component $\cS_0$ of $\cS$
that contains $(\bbmB^M,\bsB^M)$ satisfies
\begin{equation}\label{eq:pf_lm:stationary_concen_claima2}
\inf_{(\bbm,\bs) \in \cS_0} \cF_\TAP(\bbm,\bs) \text{ is attained at some }
(\bbm_\star,\bs_\star) \in \cS_0.
\end{equation}
Noting that this event $\cE'$ depends only on $(\X,\y)$ and not on $\bbeta_0$,
let us define $(\bbm_\star,\bs_\star)$ as any measurable selection of a point attaining the
infimum in (\ref{eq:pf_lm:stationary_concen_claima2}) when $\cE'$ holds,
and as an arbitrary (deterministic) point of $\Gamma^p$ when $\cE'$ does not
hold. Then $(\bbm_\star,\bs_\star)$ is a measurable function of $(\X,\y)$.

We claim that $\cE \subseteq \cE'$.
Assuming momentarily this claim, on $\cE$, we then
have that $(\bbm_\star,\bs_\star)$ is a local
minimizer of $\cF_\TAP$ by construction. Since $\cS$ is open and $\cS_0$ is
connected, it is also path connected. Any continuous path connecting
$(\bbm_\star,\bs_\star)$ and $(\bbmB^M,\bsB^M)$ in $\cS_0$ must be disjoint from
$\Gamma^p[K(\rho_0) \setminus K(\rho_1)]$, by (\ref{eq:TAPlowerstatproof1})
and the definition of the sub-level set $\cS$.
Then by condition (c) of $\cE$ which ensures
$(\bbmB^M,\bsB^M) \in \Gamma^p[K(\rho_1)]$, we must have also
$(\bbm_\star,\bs_\star) \in \Gamma^p[K(\rho_1)]$. This, together with conditions
(a) and (b), imply
\begin{align*}
\Big| p^{-1} \|\bbm_\star-\bbeta_0\|_2^2 -q_\star \Big| \le \rho_1,~~~~ \Big|
p^{-1} \|\bbs_\star-\bbm_\star^{2}\|_1 -q_\star \Big|\leq \rho_1,~~~~ \Big\vert
p^{-1} \cF_{\TAP}(\bbm_\star, \bbs_\star) - \phi(\gamma_\star) \Big\vert \le
\iota
\end{align*}
on the event $\cE$, which holds with probability at least $1-\eps$. Since
$\rho_1,\iota,\eps>0$ may be chosen arbitrarily small, the lemma follows by taking $\rho_1,\iota,\eps \rightarrow 0$ sufficiently slowly as $n,p \rightarrow \infty$.

To conclude the proof, it remains to show the claim $\cE \subseteq \cE'$. Note
that $\cE$ implies $(\bbmB^M,\bsB^M) \in \cS$, so we must show
(\ref{eq:pf_lm:stationary_concen_claima2}).
By compactness of $\overline \cS_0$, there exists a sequence
$\{(\bbm^t,\bs^t)\}_{t \geq 1} \in \cS_0$ for which
\[\inf_{(\bbm,\bs) \in \cS_0} \cF_\TAP(\bbm,\bs)
=\lim_{t \to \infty} \cF_\TAP(\bbm^t,\bs^t), \qquad
\lim_{t \to \infty} (\bbm^t,\bs^t)=(\bar \bbm,\bar \bs) \in \overline \cS_0.\]
Suppose by contradiction that $\cE$ holds, but
(\ref{eq:pf_lm:stationary_concen_claima2}) does not hold. Then, since
$\cF_\TAP$ is continuous on $\cS_0$, this implies that
$(\bar \bbm,\bar \bs)$ must belong to the boundary
$\partial \cS_0=\overline \cS_0 \setminus \cS_0$. But
$\partial \cS_0 \subseteq \partial(\Gamma^p) \cup \{(\bbm,\bs) \in
\Gamma^p:p^{-1}\cF_\TAP(\bbm,\bs)=\phi(\gamma_\star)+3\iota/4\}$, and
since condition (b) of $\cE$ ensures $\inf_{(\bbm,\bs) \in \cS_0}
\cF_\TAP(\bbm,\bs)<\phi(\gamma_\star)+\iota/2$, we must then have
$(\bar \bbm,\bar \bs) \in \partial(\Gamma^p)$.

Recall from Lemma \ref{lemma:hboundary} that
${-}\bar \hs(m,s)$ extends to a lower semi-continuous function on
the closure $\overline \Gamma$ (with values in $[0,\infty]$).
All other terms of $\cF_\TAP$ also extend
continuously to $\overline\Gamma^p$. Then, defining $\cF_\TAP$ on
$\overline\Gamma^p$ by these extensions, we have by lower semi-continuity
\begin{equation}\label{eq:Fboundarycompare}
\inf_{(\bbm,\bs) \in \cS_0} \cF_\TAP(\bbm,\bs)
=\lim_{t \to \infty} \cF_\TAP(\bbm^t,\bs^t)
\geq \cF_\TAP(\bar\bbm,\bar\bs).
\end{equation}
Let $\cJ \subseteq \{1,\ldots,p\}$ be the coordinates for which
$(\bar m_j,\bar s_j) \in \partial \Gamma$, where $\cJ$ is non-empty because
$(\bar\bbm,\bar\bs) \in \partial (\Gamma^p)$. Let us write
\[\cF_\TAP(\bbm,\bs)=\sum_{j \in \cJ} {-}\hs(m_j,s_j)+R(\bbm,\bs)\]
where $R(\bbm,\bs)$ contains all other terms including
${-}\hs(m_j,s_j)$ for $j \notin \cJ$. Fix a small convex open neighborhood $\cO$ of
$(\bar\bbm,\bar \bs)$, and observe that $R(\bbm,\bs)$ is $L$-Lipschitz on $\cO$
for some $L<\infty$ (depending on $n,p$ and $\X,\y$) as long as $\cO$ is small
enough so that $(m_j,s_j) \notin \partial \Gamma$ for all $j \notin \cJ$.
For any sufficiently
small $r>0$, Lemma \ref{lemma:hboundary} ensures that for every $j \in
\cJ$, there is a point $(m_j',s_j') \in \Gamma$ for which
$\|(\bar m_j,\bar s_j)-(m_j',s_j')\|_2=r$ and
${-}\hs(\bar m,\bar s)+\hs(m_j',s_j')>(L+1)r$. Consider the point $(\bbm',\bs')
\in \Gamma^p$ with components $(m_j',s_j')$ for $j \in \cJ$ and
$(\bar m_j,\bar s_j)$ for $j \notin \cJ$, and choose $r$ small enough so that
$(\bbm',\bs')$ belongs to the above neighborhood $\cO$. Then
$\|(\bar\bbm,\bar\bs)-(\bbm',\bs')\|_2=r\sqrt{|\cJ|}$, and
\begin{align*}
\cF_\TAP(\bar\bbm,\bar\bs)&>\cF_\TAP(\bbm',\bs')+(L+1)r|\cJ|
+R(\bar\bbm,\bar\bs)-R(\bbm',\bs')\\
&\geq \cF_\TAP(\bbm',\bs')+(L+1)r|\cJ|-Lr\sqrt{|\cJ|}
>\cF_\TAP(\bbm',\bs').
\end{align*}
Together with (\ref{eq:Fboundarycompare}), this shows
\begin{equation}\label{eq:Fboundarycomparefinal}
\inf_{(\bbm,\bs) \in \cS_0} \cF_\TAP(\bbm,\bs)>\cF_\TAP(\bbm',\bs').
\end{equation}

But for any $(\bbm^t,\bs^t) \in \cO$ and each $(\bbm,\bs)$ on the linear path
between $(\bbm',\bs')$ and $(\bbm^t,\bs^t)$, we have
\begin{align}
\cF_\TAP(\bbm,\bs)&=\sum_{j \in \cJ} {-}\hs(m_j,s_j)+R(\bbm,\bs)
\leq \max\left(\sum_{j \in \cJ} {-}\hs(m_j',s_j'),
\sum_{j \in \cJ} {-}\hs(m_j^t,s_j^t)\right)+R(\bbm,\bs)\nonumber\\
&\leq \max\Big(\cF_\TAP(\bbm',\bs'),\cF_\TAP(\bbm^t,\bs^t)\Big)
+L\|(\bbm^t,\bs^t)-(\bbm',\bs')\|_2
\label{eq:middleFbound}
\end{align}
where the first inequality applies convexity of ${-}\hs(m,s)$ and the second
applies Lipschitz continuity of $R(\bbm,\bs)$ on $\cO$. By
(\ref{eq:Fboundarycomparefinal}), the condition $(\bbm^t,\bs^t) \in \cS$, and the
definition of the sub-level set $\cS$, we must have
$\max(\cF_\TAP(\bbm',\bs'),\cF_\TAP(\bbm^t,\bs^t))<\phi(\gamma_\star)+3\iota/4$
strictly. Then, choosing large enough $t \geq 1$ and small enough $r>0$,
(\ref{eq:middleFbound})
shows that also $\cF_\TAP(\bbm,\bs)<\phi(\gamma_\star)+3\iota/4$, so
$(\bbm,\bs) \in \cS$. Since this holds for all $(\bbm,\bs)$ on the line segment between
$(\bbm^t,\bs^t)$ and $(\bbm',\bs')$, we must have $(\bbm',\bs') \in \cS_0$, the
same connected component of $\cS$ as $(\bbm^t,\bs^t)$. This contradicts
(\ref{eq:Fboundarycomparefinal}), so we cannot have $(\bar\bbm,\bar\bs) \in
\partial (\Gamma^p)$, concluding the proof.
\end{proof}

\subsection{Proof of Lemma \ref{lem:marginalposterior-in-proof}}\label{app:proof_marginalposterior-in-proof}

Let $(\bbm_\star,\bbs_\star) \in \Gamma^p$ be the given local minimizer of
$\cF_\TAP$, and let $\blambda_\star=[\lambda(m_{\star,j},s_{\star,j})]_{j=1}^p$
and $\bgamma_\star=[\gamma(m_{\star,j},s_{\star,j})]_{j=1}^p$ be defined by
(\ref{eq:lambdagammastar}). Then these vectors
$(\bbm_\star,\bbs_\star,\blambda_\star,\bgamma_\star)$ satisfy the
stationary conditions
\begin{equation}\label{eq:stationary-in-lemma-proof}
\gamma_{\star,j}=\frac{n/p}{\sigma^2+S(\bbs_\star)-Q(\bbm_\star)},
\qquad m_{\star,j}=\langle \beta \rangle_{\lambda_{\star,j},\gamma_{\star,j}},
\qquad s_{\star,j}=\langle \beta^2 \rangle_{\lambda_{\star,j},\gamma_{\star,j}}.
\end{equation}
In particular, $\gamma_{\star,j}$ is constant across
coordinates $j=1,\ldots,p$. Note that (\ref{eqn:marginalposterior-in-proof-cd2}) implies that
\[S(\bbs_\star)-Q(\bbm_\star)  = p^{-1} \|\bbs_\star-\bbm_*^2\|_1 \gotop q_\star.\]  
Then, recalling $\gamma_\star=\gamma_\stat=\delta/(\sigma^2+q_\star)$
and $n/p \to \delta$,
 the first condition of (\ref{eq:stationary-in-lemma-proof}) gives
\begin{equation}\label{eq:TAPconvergence-in-lemma-proof-gamma}
\lim_{n,p \to \infty} \E\big[\big(\gamma_{\star,j}-\gamma_\star\big)^2\big]=0.
\end{equation}
Now for any bounded function $f:\supp(\mu_0) \to \R$,
note that the mapping $\gamma \mapsto \langle f(\beta)
\rangle_{\lambda_{\star,j},\gamma}$ has derivative $(-1/2)(\langle f(\beta)\beta^2
\rangle_{\lambda_{\star,j},\gamma}-\langle f(\beta) \rangle_{\lambda_{\star,j},\gamma}
\langle \beta^2 \rangle_{\lambda_{\star,j},\gamma})$, which is uniformly bounded
because $\beta$ has compact support. Then (\ref{eq:TAPconvergence-in-lemma-proof-gamma}) implies
\begin{equation}\label{eq:gammastarcomp-in-lemma-proof}
\lim_{n,p \to \infty} \E\big[\big(\langle f(\beta)
\rangle_{\lambda_{\star,j},\gamma_{\star,j}}
-\langle f(\beta) \rangle_{\lambda_{\star,j},\gamma_\star}\big)^2\big]=0.
\end{equation}

Furthermore, note that by symmetry and boundedness, (\ref{eqn:marginalposterior-in-proof-cd1}) implies that, for every coordinate $j=1,\ldots,p$, we have
\begin{equation}\label{eq:TAPconvergence-in-lemma-proof}
\E[(m_{\star,j}-\langle \beta_j \rangle_{\X,\y})^2] = p^{-1} \E[\| \bbm_\star - \bbmB \|_2^2] \to 0. 
\end{equation}
Define the first moment map
$m(\lambda):=\langle \beta \rangle_{\lambda,\gamma_\star}$.
Combining (\ref{eq:TAPconvergence-in-lemma-proof}), the second
stationary condition of (\ref{eq:stationary-in-lemma-proof}), and
Lemma \ref{lemma:cavity} and (\ref{eq:gammastarcomp-in-lemma-proof}) both applied to the
function $f(\beta)=\beta$, we have
\begin{equation}\label{eqn:m-convergence-interpolation-in-lemma-proof}
\begin{aligned}
&~\E\big[\big(m(\lambda_{\star,j})-m(\hat{\lambda}_j)\big)^2\big] \\
&\leq 8\Big(\E\big[ \big(m(\lambda_{\star,j}) - \langle \beta
\rangle_{\lambda_{\star,j},\gamma_{\star,j}} \big)^2\big] +\E\big[\big(\langle
\beta \rangle_{\lambda_{\star,j},\gamma_{\star,j}} - \langle \beta_j
\rangle_{\X,\y} \big)^2\big] +\E\big[\big(\langle \beta_j \rangle_{\X,\y} -
m(\hat{\lambda}_j) \big)^2\big]\Big) \to 0.
\end{aligned}
\end{equation}

Moreover, note that the map $\lambda \mapsto m(\lambda)$ is differentiable with derivative
$\langle \beta^2 \rangle_{\lambda,\gamma_\star}-\langle \beta
\rangle_{\lambda,\gamma_\star}^2>0$. Hence it is a strictly increasing map from $\R$
onto its image $(a(\sP_0),b(\sP_0))$ (c.f. Proposition \ref{prop:oneparamexpfam}), the interval from minimal to maximal
point of support of $\sP_0$. We denote its functional inverse by
\[m \mapsto \lambda(m), \qquad \lambda: (a(\sP_0),b(\sP_0)) \to \R.\]
Then
\[\frac{d}{dm} \langle f(\beta) \rangle_{\lambda(m),\gamma_\star}
=\Big(\langle f(\beta)\beta \rangle-\langle f(\beta)\rangle
\langle \beta\rangle\Big) \cdot \lambda'(m)
=\frac{\langle f(\beta)\beta \rangle
-\langle f(\beta)\rangle \langle \beta\rangle}
{\langle \beta^2 \rangle-\langle \beta \rangle^2},\]
where we have abbreviated $\langle \cdot \rangle=\langle \cdot
\rangle_{\lambda(m),\gamma_\star}$ on the right side.
Introducing replicas $\beta_1,\beta_2$, when $f$ is $L$-Lipshitz, this implies
\[\left|\frac{d}{dm} \langle f(\beta) \rangle_{\lambda(m),\gamma_\star}\right|
=\left|\frac{\langle (f(\beta_1)-f(\beta_2))(\beta_1-\beta_2) \rangle}
{\langle (\beta_1-\beta_2)^2 \rangle}\right|
\leq \frac{\langle L(\beta_1-\beta_2)^2 \rangle}
{\langle (\beta_1-\beta_2)^2 \rangle}=L.\]
Hence by the above inequality and (\ref{eqn:m-convergence-interpolation-in-lemma-proof}), we get
\[\lim_{n,p \to \infty} \E\Big[\big(\langle f(\beta)
\rangle_{\lambda_{\star,j},\gamma_\star} -\langle f(\beta)
\rangle_{\hat{\lambda}_j,\gamma_\star}\big)^2\Big]
\leq L \cdot \lim_{n,p \to \infty}
\E\big[\big(m(\lambda_{\star,j})-m(\hat{\lambda}_j)\big)^2\big]=0.\]
Now combining this with Lemma \ref{lemma:cavity} and (\ref{eq:gammastarcomp-in-lemma-proof})
both applied to this function $f(\beta)$, we have as desired
\[\lim_{n,p \to \infty}
\E\big[\big(\langle f(\beta_j) \rangle_{\X,\y}-\langle f(\beta)
\rangle_{\lambda_{\star,j},\gamma_{\star,j}}\big)^2\big]=0.\]
This proves Lemma \ref{lem:marginalposterior-in-proof}.

\section{Outline of remaining proofs}
\label{app:proof-outline}

The proofs of Theorems \ref{thm:local_convexity}, \ref{thm:NGD_convergence}, and \ref{thm:local_convexity_AMP_fixed_point} and Corollary \ref{cor:AMP+NGD-convergence} are intertwined and are thus carried out in parallel.
The argument is outlined as follows.
\begin{enumerate}

	\item \label{item:landscape-step}
	\textbf{Landscape analysis around AMP iterate.} We establish local convexity and approximate stationarity of the TAP free energy in a neighborhood of the AMP iterates for large $k$. This analysis applies in both the hard and easy regimes. It is carried out in Appendix \ref{app:proof_local_convexity}. This establishes Theorem \ref{thm:local_convexity_AMP_fixed_point}\ref{item:local-convexity}, as explained in Appendix \ref{sec:finishing-local-convexity}.

	\item \label{item:bayes-local-convexity-step}
	\textbf{Local convexity around the Bayes estimate in the easy regime.} In the easy regime, the AMP iterate approximates the Bayes estimate. Combining this fact with item \ref{item:landscape-step} above implies Theorem \ref{thm:local_convexity}. This is carried out in Appendix \ref{app:Bayes-local-convexity}.

	\item \label{item:NGD-step}
	\textbf{Generic analysis of natural gradient descent.} We show that NGD converges linearly whenever it is initialized in the neighborhood of a local minimizer satisfying certain local landscape properties. This establishes Theorem \ref{thm:NGD_convergence}, and is carried out in the next section, Appendix \ref{app:proof_NGD_convergence}.

	\item \label{item:AMP+NGD-step}
	\textbf{Convergence of AMP+NGD in the easy and hard regime.} In all
regimes, item \ref{item:landscape-step} above implies that AMP arrives in a
neighborhood satisfying the conditions of item \ref{item:NGD-step}. Thus, NGD
initialized at a sufficiently late iterate of AMP converges linearly. By item
\ref{item:bayes-local-convexity-step} above, in the easy regime this corresponds
to a Bayes-optimal local minimizer. Corollary \ref{cor:AMP+NGD-convergence} and Theorem \ref{thm:local_convexity_AMP_fixed_point}\ref{item:NGD-convergence} follow. This is carried out in Appendix \ref{app:AMP+NGD-convergence}.

	\item \label{item:calibration-step}
	\textbf{Calibrated inference.} AMP state evolution allows us to study the statistical properties of the local minimizer of item \ref{item:AMP+NGD-step} in both the easy and hard regimes. This gives us Theorem \ref{thm:error-and-calibration}. This is carred out in Appendix \ref{app:calibration}.

\end{enumerate}

\section{Local convexity of the TAP free energy}\label{app:proof_local_convexity}

In this section, we carry out a landscape analysis around the AMP iterates in both the hard and easy regime, corresponding to Steps \ref{item:landscape-step} and \ref{item:bayes-local-convexity-step} of the proof outline given in Appendix \ref{app:proof-outline}.
By the end of the section, we will have proved Theorems \ref{thm:local_convexity} and \ref{thm:local_convexity_AMP_fixed_point}\ref{item:local-convexity}, and established many facts that will be useful in the developments to come.

In this and later sections, we will frequently use the following notation.
For a sequence of random variables $X_n$ indexed by $n$ and a constant $c$,
we write $\plim_{n \rightarrow \infty} X_n = c$ to mean $X_n \gotop c$,
$\pliminf_{n \rightarrow \infty} X_n= c$ to mean $\sup \{t \in \reals \mid \P(X_n \leq t) \rightarrow 0 \} = c$,
and $\plimsup_{n \rightarrow \infty} X_n = c$ to mean $\pliminf_{n \rightarrow \infty} -X_n = -c$.

Recall the AMP iterates
$\bz^k,\bbm^k$ and sequence $\gamma_k$ defined by~\eqref{eq:AMPalg}.
We restate these here for the reader's convenience:
\begin{equation}\label{eq:AMP}
\begin{aligned}
\z^0&=\bbm^1=0, \qquad &\gamma_1&=\delta(\sigma^2+\E_{\beta \sim
\sP_0}[\beta^2])^{-1},\\
	\z^k&=\y-\X\bbm^k+\frac{\sb_{k-1}}{\delta}\,\z^{k-1},
	\qquad
	&\gamma_{k+1}&=\delta\left(\sigma^2+\mmse(\gamma_k)\right)^{-1},
	\\
	\bbm^{k+1}&=\denoiser\left(\bbm^k+\frac1\delta \X^\top
		\z^k,\,\gamma_k\right), 
	\qquad
	&\bs^{k+1}
		&=
		\sS\left(
			\bbm^k + \frac1\delta \X^\top \bz^k,\,\gamma_k
		\right),
\end{aligned}
\end{equation}
where $\sb_k = \gamma_k\, \mmse(\gamma_k)$, $\denoiser(\x,\gamma)=(\langle \beta
\rangle_{\gamma x_j,\gamma})_{j=1}^p$, and 
	$\sS(\x,\gamma)=(\langle \beta^2
	\rangle_{\gamma x_j,\gamma})_{j=1}^p$.
In this section, we will not use the ``$\AMP$'' subscript in $(\bz^k,\bbm^k, \bs^k)$ for cleaner notation.
The main goal of this section is to prove the following three facts:
\begin{enumerate}

	\item \textit{Approximate stationarity at $(\bbm^k,\bs^k)$.}
	For some $C,\kappa > 0$ depending only on $(\sigma^2,\delta,\sP_0)$, 
	and for all $k \ge 0$,
	\begin{equation}
	\label{eq:approx-stationarity}
		\plimsup_{n \rightarrow \infty}\;
		\frac1p \| \nabla \cF_\TAP(\bbm^k,\bs^k) \|_2^2 \leq Ce^{-\kappa k}.
	\end{equation}

	\item \textit{Local strong convexity around $(\bbm^k,\bs^k)$.} For some
$\kappa > 0$ depending only on $(\sigma^2,\delta,\sP_0)$, 
	\begin{equation}
	\label{eq:local-convexity}
		\lim_{\epsilon \rightarrow 0}\;
		\lim_{k \rightarrow \infty}\;
		\pliminf_{n \rightarrow \infty}\;
        \min_{
			\substack{
			    \|\bl\|_2/\sqrt{p} = 1 \\
			    \|\bbm - \bbm^k\|_2/\sqrt{n} \leq \epsilon \\
			    \|\bs - \bs^k \|_2/\sqrt{n} \leq \epsilon
			}
		}\;
        \frac{1}{p}\,\< \bl, \nabla^2 \cF_\TAP(\bbm,\bs) \bl \> 
        \geq \kappa.
	\end{equation}

	\item \textit{Small sub-optimality gap at $(\bbm^k,\bs^k)$.} For some $c
> 0$ depending only on $(\sigma^2,\delta,\sP_0)$,
	\begin{equation}
	\label{eq:sub-optimality-gap}
		\lim_{\epsilon \rightarrow 0}\;
		\lim_{k \rightarrow \infty}\;
		\plimsup_{n \rightarrow \infty}\;
			\left\{
        	\frac1p \cF_{\TAP}(\bbm^k,\bs^k)
        	-
        	\min_{
			\substack{
			    \|\bbm - \bbm^k\|_2/\sqrt{n} \leq \epsilon \\
			    \|\bs - \bs^k \|_2/\sqrt{n} \leq \epsilon
				}
			}\;
			\frac1p \cF_{\TAP}(\bbm,\bs)
			-
			c\epsilon^2
			\right\}
        	\leq 0.
	\end{equation}

\end{enumerate}

The remainder of this section is dedicated to establishing these three facts,
and ends with a proof of Theorems \ref{thm:local_convexity} and \ref{thm:local_convexity_AMP_fixed_point}\ref{item:local-convexity}.

\subsection{State evolution}


Recall the quantities defined in Section \ref{sec:gordon-post-AMP}
\begin{equation}\label{eq:AMPaux}
\begin{aligned}
      \bg^k
                  &:=
                  \bbm^k + \frac1\delta \X^\top \bz^k - \bbeta_0,
            \qquad
            &\bh^k
                  &:=
                  \beps - \bz^k,
      \\
      \bnu^k
                  &:=
                  \bbm^k - \bbeta_0,
            \qquad
            &\br^k
            	&:=
            	-\bz^k.
\end{aligned}
\end{equation}
Let $\bR_k \in \R^{n \times k}$, $\bG_k \in \R^{p \times k}$, $\bV_k \in \R^{p
\times k}$, and $\bH_k \in \R^{n \times k}$ be the matrices whose columns are
$\{\br^{k'}\}_{1\leq k' \leq k}$, $\{ \bg^{k'} \}_{1 \leq k' \leq k}$,
$\{ \bnu^{k'} \}_{1 \leq k' \leq k}$, and
$\{ \bh^{k'} \}_{1 \leq k' \leq k}$ respectively.

We state here an extended AMP state evolution that describes the behavior of
these quantities, via two bi-infinite matrices $\bK_g,\bK_h \in \R^{\Z_{> 0} \times \Z_{> 0}}$ whose entries are
\begin{equation}\label{eq:KgKh}
K_{g,kk'} = \gamma_{k \vee k'}^{-1}, \qquad
K_{h,kk'} = \mmse(\gamma_{(k-1) \vee (k'-1)}) = \delta \gamma_{k \vee k'}^{-1} -
\sigma^2.
\end{equation}
Denote by $\bK_{g,k}$ and $\bK_{h,k}$ their respective upper-right $k\times k$ submatrices.
Let $\beta_0 \sim \sP_0$ be independent of $(G_k)_{k \geq 1}$ a centered
Gaussian sequence with covariance $\bK_g$,
and $\eps \sim \normal(0,\sigma^2)$ be independent of $(H_k)_{k \geq 1}$ a
centered Gaussian sequence with covariance $\bK_h$.
These random variables describe the quantities (\ref{eq:AMPaux})
in the sense of the following proposition.
\begin{proposition}[State evolution]
\label{prop:amp-se}
      We have the following:

      \begin{enumerate}[label=(\alph*)]

            \item 
            $\gamma_k$ is strictly increasing and $\gamma_k \rightarrow \gamma_\alg$.

            \item 
            For any $k$, $\bK_{g,k}$ and $\bK_{h,k}$ are strictly positive definite. 

            \item 
            For any pseudo-Lipschitz test function $\psi:\R^{k+1} \to \R$,
almost surely as $n,p \to \infty$,
            \begin{equation}\label{eq:AMPSE}
            \begin{gathered}
                  \frac1p \sum_{j=1}^p
                        \psi(\beta_{0,j},g_j^1,\ldots,g_j^k)
\to \E[\psi(\beta_0,G_1,\ldots,G_k)],
                  \\
                  \frac1n \sum_{i=1}^n
                        \psi(\eps_i,h_i^1,\ldots,h_i^k)
\to \E[\psi(\eps,H_1,\ldots,H_k)].
            \end{gathered}
            \end{equation}
            \item 
            Almost surely as $n,p \to \infty$,
            \begin{equation}\label{eq:AMPVRGH}
            \begin{aligned}
                  \frac1p \bV_k^\top \bV_k
                        &\rightarrow
                        \bK_{h,k},
                  \qquad
                  &
                  \frac1n \bR_k^\top \bR_k
                        &\rightarrow
                        \delta \bK_{g,k},
                  \\
                  \frac1p \bG_k^\top \bV_k
                        &\rightarrow
                        \bK_{g,k} \sB_k^g,
                  \qquad
                  &
                  \frac1n \bH_k^\top \bR_k
                        &\rightarrow
                        \bK_{h,k},
            \end{aligned}
            \end{equation}
            where $\sB_k^g$ is the matrix with entries $\sB_{k,ij}^g =
\1_{\{i=j-1\}} \cdot \sb_{j-1}$.
      \end{enumerate}

\end{proposition}

\begin{proof}
    
    \noindent \textbf{Part (a).} By \eqref{eq:gammaalg} and continuity of
$\mmse(\gamma)$,
    $\gamma_\alg$ is the smallest solution to $\gamma = \delta \big(\sigma^2 + \mmse(\gamma)\big)^{-1}$.
    For any $\gamma \leq \gamma_1$,
    we have
    $\gamma \leq \gamma_1 = \delta\big(\sigma^2 + \E_{\beta \sim \sP_0}[\beta_0^2]\big)^{-1} < \delta\big(\sigma^2 + \mmse(\gamma)\big)^{-1}$.
	Thus, $\gamma_1 < \gamma_\alg$.
	For any $\gamma \geq \delta/\sigma^2$,
	we have $\gamma \geq \delta/\sigma^2 > \delta\big(\sigma^2 + \mmse(\gamma)\big)^{-1}$.
	Thus, $\gamma_\alg < \infty$.
	Because $\mmse(\gamma)$ is continuous and non-increasing and $\gamma_1 < \gamma_\alg < \infty$,
	we conclude $\gamma_k$ is strictly increasing and $\gamma_k \rightarrow \gamma_\alg$.
	\\

	\noindent \textbf{Part (b).} This is a consequence of $\gamma_k$ being
strictly increasing, $\gamma_k < \delta / \sigma^2$, and the definitions
(\ref{eq:KgKh}).\\

	\noindent \textbf{Part (c).} State evolution for general AMP algorithms,
of which \eqref{eq:AMPalg} is an example, has been established multiple times in
the literature. Here, we note that (\ref{eq:AMPalg}) may be written in terms of
the quantities (\ref{eq:AMPaux}) as
\[\bh^{k+1}=\X f_k({-}\bg^k,\bbeta_0)+\sb_k e_k(\bh^k,\beps),
\qquad {-}\bg^k=\X^\top e_k(\bh^k,\beps)-f_{k-1}({-}\bg^{k-1},\bbeta_0)\]
where we identify the functions
\[f_k(x,\beta_0)=\denoiser(-x+\beta_0,\gamma_k)-\beta_0,
\qquad e_k(x,\eps)=\frac{1}{\delta}(x-\eps)\]
applied elementwise. Observe that since $\sP_0$ has bounded support, both
functions $e_k,f_k$ are Lipschitz in the first argument $x$. Then, denoting by
$e_k',f_k'$ the derivatives in $x$
and by $(\eps,H_k)$ and $(\beta_0,G_k)$ the joint laws on
the right side of (\ref{eq:AMPSE}), we have
\begin{align*}
\delta\,\E[e_k'(H_k,\eps)] &=1,\\
\E[f_k'({-}G_k,\beta_0)] &= {-}\gamma_k\E\Big[\langle \beta^2
\rangle_{\gamma_k(G_k+\beta_0),\gamma_k}-\langle \beta
\rangle_{\gamma_k(G_k+\beta_0),\gamma_k}^2\Big]={-}\gamma_k\mmse(\gamma_k)={-}\sb_k.
\end{align*}
Then \cite[Theorem 1]{berthierMontanariNguyen} (with the notational
identification $\bA=\X^\top$) establishes (\ref{eq:AMPSE})
where $(G_1,\ldots,G_k)$ and $(H_1,\ldots,H_k)$ are centered Gaussian vectors
independent of $\beta_0$ and $\eps$, with covariance matrices $\bK_{g,k}$ and
$\bK_{h,k}$ having entries
\[K_{g,kk'}=\delta\,\E[e_k(H_k,\eps)e_{k'}(H_{k'},\eps)],
\qquad K_{h,kk'}=\E[f_{k-1}(G_{k-1},\beta_0)f_{k'-1}(G_{k'-1},\beta_0)].\]
Applying the martingale identity
\begin{equation}\label{eq:martingale}
\E[f_k({-}G_k,\beta_0)f_{k'}({-}G_{k'},\beta_0)]
=\E[(\E[\beta_0 \mid \beta_0+G_k]-\beta_0) \cdot
(\E[\beta_0 \mid \beta_0+G_{k'}]-\beta_0)]=\mmse(\gamma_{k \vee k'}),
\end{equation}
a straightforward induction argument shows that these covariances $\bK_{g,k}$
and $\bK_{h,k}$ coincide with the definition (\ref{eq:KgKh}),
hence establishing part (c).\\

	\noindent \textbf{Part (d).} Identifying $\br^k=\bh^k+\beps$, the
statements for $n^{-1}\bR_k^\top \bR_k$ and $n^{-1}\bH_k^\top \bR_k$ follow from
part (c) applied with $\psi(\eps,H_1,\ldots,H_k)=(H_k+\eps)(H_{k'}+\eps)$ and
with $\psi(\eps,H_1,\ldots,H_k)=H_k(H_{k'}+\eps)$.

Identifying
$\bnu^k=f_{k-1}({-}\bg^{k-1},\bbeta_0)$, the statements
for $p^{-1}\bV_k^\top \bV_k$ and $p^{-1}\bG_k^\top \bV_k$ also
follow from part (c)
applied with $\psi(\beta_0,G_1,\ldots,G_k)
=f_{k-1}({-}G_{k-1},\beta_0)f_{k'-1}({-}G_{k'-1},\beta_0)$
and with $\psi(\beta_0,G_1,\ldots,G_k)=G_kf_{k'-1}({-}G_{k'-1},\beta_0)$,
together with the martingale identity (\ref{eq:martingale}) and Stein's lemma
\[\E[G_k f_{k'-1}({-}G_{k'-1},\beta_0)]
={-}\E[G_kG_{k'-1}] \cdot \E[f_{k'-1}'({-}G_{k'-1},\beta_0)]
=\gamma_{k \vee (k'-1)} \cdot \sb_{k'-1}.\]
\end{proof}

\subsection{Approximate stationarity at the AMP iterates}
\label{app:approx-stationarity}

We let $C,\kappa > 0$ be constants depending on $(\sigma^2,\delta,\sP_0)$ whose value may change at each appearance.

The TAP gradient at the AMP iterates is given by
\begin{equation}
\label{eq:TAP-grad}
	\nabla \cF_\TAP(\bbm^k,\bs^k)
		=
		\begin{pmatrix}
			\blambda^k - \frac1{\sigma^2}\X^\top(\y - \X \bbm^k)
- \frac{n/p}{\sigma^2 + S(\bs^k) - Q(\bbm^k)}\bbm^k
			\\
			-\frac12 \bgamma^k + \frac{n/p}{2(\sigma^2 + S(\bs^k) - Q(\bbm^k))}\ones
		\end{pmatrix}.
\end{equation}
where $\blambda^k=[\lambda(m_j^k,s_j^k)]_{j=1}^p$ and
$\bgamma^k=[\gamma(m_j^k,s_j^k)]_{j=1}^p$.
Comparing the definitions \eqref{eq:AMP} with the definitions of
$(\lambda(m,s),\gamma(m,s))$ via (\ref{eq:momentspace}--\ref{eq:lambdagammastar}),
we have
\begin{equation}\label{eq:AMPlambdagamma}
\blambda^k = \gamma_{k-1}\Big(\bbm^{k-1} + \frac1\delta \X^\top \bz^{k-1}\Big),
\qquad \bgamma^k = \gamma_{k-1} \ones.
\end{equation}
Recall also $\y - \X\bbm^k = \bz^k - \frac{\sb_{k-1}}{\delta}\bz^{k-1}$.
Plugging these expressions into the gradient expression above, using that
$\|\X\|_{\op} / \sqrt{p} \leq C$ with high probability,
and applying the consequences of Proposition \ref{prop:amp-se}
\[\begin{gathered}
	\plim_{n \rightarrow \infty}\;
		\frac1p \| \bbm^k - \bbm^{k-1} \|^2
		=
		\Big(\mmse(\gamma_{k-2}) - \mmse(\gamma_{k-1})\Big) \leq Ce^{-\kappa k},
	\\
	\plim_{n \rightarrow \infty}
		\frac1n \| \bz^k - \bz^{k-1} \|^2
		=
		\plim_{n \rightarrow \infty}
		\frac1n \| \bh^k - \bh^{k-1} \|^2
		=
		\delta (\gamma_{k-1}^{-1} - \gamma_k^{-1})
		\leq
		Ce^{-\kappa k},
	\\
		\plim_{n \rightarrow \infty}
		\Big|
			\gamma_{k-1} - \frac{n/p}{\sigma^2 + S(\bs^k) - Q(\bbm^k)}
		\Big|
		=
		\Big| \gamma_{k-1} - \frac{\delta}{\sigma^2 +
\mmse(\gamma_{k-1})} \Big|
		=
		|\gamma_{k-1} - \gamma_k| \leq Ce^{-\kappa k},
	\\
		\Big|
			\frac{\gamma_{k-1}}{\delta}
			-
			\frac1{\sigma^2}\Big(1 - \frac{\sb_{k-1}}{\delta}\Big)
		\Big| = 
		\Big|
			\frac{\gamma_{k-1}}{\delta}
			-
			\frac1{\sigma^2}\Big(1 -
		\frac{\gamma_{k-1}\mmse(\gamma_{k-1})}{\delta}\Big) \Big| 
		=
		\frac{|\mmse(\gamma_{k-2})-\mmse(\gamma_{k-1})|/\sigma^2}{\sigma^2 + \mmse(\gamma_{k-2})} \le Ce^{-\kappa k},
\end{gathered}\]
we conclude \eqref{eq:approx-stationarity}.
We have used multiple times in the previous display that $\gamma_k \to
\gamma_\alg$, $|\gamma_k - \gamma_{k-1}|\leq C^{-\kappa k}$, $n/p \to \delta$, and $\gamma_{k+1} = \delta[\sigma^2 + \mmse(\gamma_k)]^{-1}$.

\subsection{Local convexity around AMP iterates}\label{app:localconvexity}

In this section we prove \eqref{eq:local-convexity}. Computing the Hessian of
$\cF_\TAP$ in (\ref{eqn:TAP}), setting
\[V(\bbm,\bs)=\sigma^2+S(\bs)-Q(\bbm), \quad
\bD_{mm}=-\diag(\partial_m^2\sh(m_j,s_j)_{j=1}^p),\]
\[\bD_{ms}=-\diag(\partial_{ms}^2\sh(m_j,s_j)_{j=1}^p), \quad
\bD_{ss} = - \diag(\partial_s^2\sh(m_j,s_j)_{j=1}^p),\]
and applying the variational representation
$\sigma^{-2}\bl_1^\top \X^\top \X\bl_1=\max_{\bu \in \R^n} 2\bu^\top
\X\bl_1-\sigma^2\|\bu\|_2^2$, we obtain that
the smallest eigenvalue of $\nabla^2 \cF_\TAP(\bbm,\bs)$ is given by
\begin{equation}
\label{eq:local-hessian-gordon-form}
      \frac{1}{p}\,\min_{\|\bl\|_2 = \sqrt{p}}
        \< \bl, \nabla^2 \cF_\TAP(\bbm,\bs) \bl \>
        =
        2
        \min_{\|\bl\|_2 = \sqrt{p}}\;
        \max_{\bu \in \reals^n}
        \Big\{
              \frac1p \bu^\top \X \bl_1
              +
              \Theta_{\TAP}(\bu,\bl_1,\bl_2;\bbm,\bs)
            \Big\},
\end{equation}
where $\bl = (\bl_1,\bl_2)$ with $\bl_1,\bl_2 \in \reals^p$ and
\[\begin{aligned}
    \Theta_{\TAP}(\bu,\bl_1,\bl_2;\bbm,\bs)
    &:=
    \Theta_{\bulk}(\bu,\bl_1,\bl_2;\bbm,\bs)
    +
    \Theta_{\spike}(\bu,\bl_1,\bl_2;\bbm,\bs),
    \\
    \Theta_{\bulk}(\bu,\bl_1,\bl_2;\bbm,\bs)
        &:=
        \frac1{2p}
        \begin{pmatrix}
            \bu^\top & \bl_1^\top & \bl_2^\top
        \end{pmatrix}
        \begin{pmatrix}
            -\sigma^2 \id_n & 0 & 0 \\[3pt]
            0 & \bD_{mm} - (n/p)V(\bbm,\bs)^{-1} \id_n & \bD_{ms} \\[3pt]
            0 & \bD_{ms} & \bD_{ss}
        \end{pmatrix}
        \begin{pmatrix}
            \bu \\[2pt] 
            \bl_1 \\[2pt]
            \bl_2
        \end{pmatrix},
    \\
    \Theta_{\spike}(\bu,\bl_1,\bl_2;\bbm,\bs)
        &:=
        \frac{(n/p) V(\bbm,\bs)^{-2}}{p}
        \begin{pmatrix}
            \bu^\top & \bl_1^\top & \bl_2^\top
        \end{pmatrix}
        \begin{pmatrix}
            0 & 0 & 0 \\[3pt]
            0 & -\frac{\bbm\bbm^\top}{p} & \frac{\bbm\ones^\top}{2p} \\[3pt]
            0 & \frac{\ones \bbm^\top}{2p} & -\frac{\ones\ones^\top}{4p}
        \end{pmatrix}
        \begin{pmatrix}
            \bu \\[2pt] 
            \bl_1 \\[2pt]
            \bl_2
        \end{pmatrix}.
\end{aligned}\]
The maximum over $\bu$ is achieved at $\bu = \X \bl_1 / \sigma^2$.
For a sufficiently large constant $C_0>0$, with probability approaching 1,
the event $\| \X \bl_1 \|_2 / \sigma^2 \leq \sqrt{p} \| \X \|_{\op} / \sigma^2
\leq C_0 \sqrt{n} $ occurs, so
\begin{equation}
      \frac{1}{p}\,\min_{\|\bl\|_2 = \sqrt{p}}
            \< \bl, \nabla^2 \cF_\TAP(\bbm,\bs) \bl \>
            =
            2
            \min_{\|\bl\|_2/\sqrt{p} = 1}\;
            \max_{\|\bu\|_2/\sqrt{n} \leq C_0}
            \Big\{
                  \frac1p \bu^\top \X \bl_1
                  +
                  \Theta_{\TAP}(\bu,\bl_1,\bl_2;\bbm,\bs)
            \Big\}.
\end{equation}
This constant $C_0>0$ will be fixed throughout the remainder of the proof.
Thus, it suffices to show that for some $\kappa > 0$,
\begin{equation}
\label{eq:local-convexity-gordon-form}
	\lim_{\epsilon \rightarrow 0}\;
	\pliminf_{n \rightarrow \infty}\;
    \min_{
          \substack{
                \|\bl\|_2/\sqrt{p} = 1 \\
                \|\bbm - \bbm^k\|_2/\sqrt{n} \leq \epsilon \\
                \|\bs - \bs^k \|_2/\sqrt{n} \leq \epsilon
          }
    }\;
    \max_{\|\bu\|_2 \leq C_0\sqrt{n}}
    \Big\{
          \frac1p \bu^\top \X \bl_1
          +
          \Theta_{\TAP}(\bu,\bl_1,\bl_2;\bbm,\bs)
    \Big\}
    \geq \kappa.
\end{equation}

\subsubsection{Gordon post-AMP: proof of Proposition \ref{prop:gordon-post-amp}}
\label{app:proof-gordon-post-amp}

As described in Section \ref{sec:gordon-post-AMP},
we lower bound the left side of the previous display using a conditional form of Gordon's comparison inequality.
Recall that  $\proj_{\bR_k} \in \R^{n \times n}$ and $\proj_{\bV_k} \in \R^{p \times p}$ 
are the projections onto the linear spans of $\br^1,\ldots,\br^k$ and
$\bnu^1,\ldots,\bnu^k$, and let $\proj_{\bR_k}^\perp,\proj_{\bV_k}^\perp$ be the
projections onto the orthogonal complements.
Let $\bg \sim \normal(\bzero,\id_p)$ and $\bh \sim \normal(\bzero,\id_n)$ be independent of each other and everything else.
Recall we define
\begin{equation}
\begin{gathered}
    \bg^*(\bu)
          :=
          \frac{1}{\sqrt{\delta}}\,\frac{1}{n}\bG_k\bK_{g,k}^{-1}\bR_k^\top \bu
          +
          \frac1{\sqrt{n}}\|\proj_{\bR_k}^\perp \bu\|_2 \bg,
    \qquad
    \bh^*(\bl_1)
          :=
          \frac1p \bH_k \bK_{h,k}^{-1} \bV_k^\top \bl_1
          +
          \frac1{\sqrt{p}} \| \proj_{\bV_k}^\perp \bl_1 \|_2 \bh.
\end{gathered}
\end{equation}
The Gordon post-AMP objective is defined as
\begin{equation}
\label{eq:AuxObj-def}
    \AuxObj_k(\bu,\bl_1,\bl_2;\bbm,\bs;\bG_k,\bH_k,\bR_k,\bV_k,\bg,\bh)
    :=
    -\sqrt{\delta}\,\frac{\<\bg^*(\bu),\bl_1\>}{p}
    +\delta\,\frac{\<\bh^*(\bl_1),\bu\>}{n}
    +\Theta_{\TAP}(\bu,\bl_1,\bl_2;\bbm,\bs).
\end{equation}
Proposition \ref{prop:gordon-post-amp} is equivalent to showing
  \begin{equation}
  \begin{aligned}
        &\pliminf_{n \rightarrow \infty}
        \min_{
              \substack{
                    \|\bl\|_2/\sqrt{p} = 1 \\
                    \|\bbm - \bbm^k\|_2/\sqrt{n} \leq \epsilon \\
                    \|\bs - \bs^k \|_2/\sqrt{n} \leq \epsilon
              }
        }\;
        \max_{\|\bu\|_2 \leq C_0\sqrt{n}}
        \Big\{
              \frac1p \bu^\top \bX \bl_1
              +
              \Theta_{\TAP}(\bu,\bl_1,\bl_2;\bbm,\bs)
        \Big\}
  \\
        &\qquad\qquad\qquad\qquad\geq
        \pliminf_{n \rightarrow \infty}
        \min_{
              \substack{
                    \|\bl\|_2/\sqrt{p} = 1 \\
                    \|\bbm - \bbm^k\|_2/\sqrt{p} \leq \epsilon \\
                    \|\bs - \bs^k \|_2/\sqrt{p} \leq \epsilon
              }
        }\;
        \max_{\|\bu\|_2 \leq C_0\sqrt{n}}
        \AuxObj_k(\bu,\bl_1,\bl_2;\bbm,\bs;\bG_k,\bH_k,\bR_k,\bV_k,\bg,\bh).
  \end{aligned}
  \end{equation}
We prove Proposition \ref{prop:gordon-post-amp} by combining the Gordon's comparison inequality
with a Gaussian conditioning argument (see e.g.\ \cite{Bolthausen2014,bayati2011dynamics}).
The next lemma handles approximations due to the Gaussian conditioning.
\begin{lemma}
\label{lem:conditional-X}
      We have
      \begin{equation}
            \Big\| 
                  \big(\X - \proj_{\bR_k}^\perp \X \proj_{\bV_k}^\perp\big)
                  -
                  \Big( {-}
                        \frac1n \bR_k \bK_{g,k}^{-1} \bG_k^\top
                        +
                        \frac1p \bH_k \bK_{h,k}^{-1}\bV_k^\top
                  \Big)
            \Big\|_{\op}
                  \gotoas
                  0.
      \end{equation}
\end{lemma}

\begin{proof}
      We use $\bA \eqnas \bB$ to denote $\| \bA - \bB \|_{\op} \gotoas 0$.
      We have
      \begin{equation}\label{eq:Xcondapprox}
      \begin{aligned}
            &\X - \proj_{\bR_k}^\perp \X \proj_{\bV_k}^\perp\\
                  &=
                  \X \proj_{\bV_k}
                  +
                  \proj_{\bR_k} \X \proj_{\bV_k}^\perp
            \\
                  &=\frac1p \X \bV_k (\bV_k^\top \bV_k/p)^{-1} \bV_k^\top
                  +
                  \frac1n \bR_k (\bR_k^\top \bR_k/n)^{-1} \bR_k^\top \X
                  \Big(
                        \id_p
                        -
                        \frac1p \bV_k(\bV_k^\top \bV_k/p)^{-1} \bV_k^\top
                  \Big)
            \\
                  &=
                  \frac1p (\bH_k - \delta^{-1}\bR_k \sB_k^g)(\bV_k^\top
\bV_k/p)^{-1}\bV_k^\top
                  +
                  \frac1{n} \bR_k (\bR_k^\top \bR_k/n)^{-1}(\bV_k - \bG_k)^\top
                  \Big(
                        \id_p
                        -
                        \frac1p \bV_k(\bV_k^\top \bV_k/p)^{-1} \bV_k^\top
                  \Big)
            \\
                  &\eqnas
                  \frac1p \bH_k \bK_{h,k}^{-1}\bV_k^\top
                  -
                  \frac1n \bR_k \sB_k^g \bK_{h,k}^{-1}\bV_k^\top
                  -
                  \frac1n \bR_k \bK_{g,k}^{-1} \bG_k^\top
                  +
                  \frac1{np}\bR_k \bK_{g,k}^{-1} \bG_k^\top \bV_k \bK_{h,k}^{-1}\bV_k^\top
            \\
                  &\eqnas
                  \frac1p \bH_k \bK_{h,k}^{-1}\bV_k^\top
                  -
                  \frac1n \bR_k \sB_k^g \bK_{h,k}^{-1}\bV_k^\top
                  -
                  \frac1n \bR_k \bK_{g,k}^{-1} \bG_k^\top
                  +
                  \frac1n\bR_k \sB_k^g \bK_{h,k}^{-1}\bV_k^\top
            \\
                  &=
                  \frac1p \bH_k \bK_{h,k}^{-1}\bV_k^\top
                  -
                  \frac1n \bR_k \bK_{g,k}^{-1} \bG_k^\top,
      \end{aligned}
      \end{equation}
      where in the fourth line we used $\X\bV_k=\bH_k - \delta^{-1}\bR_k \sB_k^g$ and $\delta^{-1}\X^\top \bR_k=\bV_k-\bG_k$ by the definitions (\ref{eq:AMP})
and (\ref{eq:AMPaux}), and in the fifth and sixth lines we used
the state evolution of Proposition \ref{prop:amp-se} and the fact that
$\bK_{g,k}$ and $\bK_{h,k}$ are positive definite.
\end{proof}

\begin{proof}[Proof of Proposition \ref{prop:gordon-post-amp}]
      For any fixed $\bR \in \reals^{n \times k}$, $\bV \in \reals^{p \times k}$,
      we have that $\X - \proj_{\bR}^\perp \X \proj_{\bV}^\perp$ is
independent of $\proj_{\bR}^\perp \X \proj_{\bV}^\perp$ since $\X$ has i.i.d.\
Gaussian entries.
Then conditional on $\bV_k,\bR_k,\bG_k,\bH_k$, observe that the fourth line of
(\ref{eq:Xcondapprox}) shows
$\X-\proj_{\bR_k}^\perp \X \proj_{\bV_k}^\perp$ is deterministic, while
$\proj_{\bR_k}^\perp \X \proj_{\bV_k}^\perp$ is equal in conditional law to
$\proj_{\bR_k}^\perp \Xtilde \proj_{\bV_k}^\perp$ for $\Xtilde$ an independent
copy of $\X$ that is also independent of $\bV_k,\bR_k,\bG_k,\bH_k$. Thus
we conclude that
      \begin{equation}
            \X
                  \eqnd
                  \X - \proj_{\bR_k}^\perp \X \proj_{\bV_k}^\perp +
\proj_{\bR_k}^\perp \Xtilde \proj_{\bV_k}^\perp.
      \end{equation}
First conditioning on all quantities but $\Xtilde$,
applying Gordon's inequality in the form of Lemma \ref{lemma:gordon}
with $\bG \equiv \sqrt{p}\Xtilde$, and then marginalizing over the randomness
of the remaining quantities, we obtain for any $t \in \R$ that
      \[\begin{aligned}
            &\P\Big(
                  \min_{
                        \substack{
                              \|\bl\|_2/\sqrt{p} = 1 \\
                              \|\bbm - \bbm^k\|_2/\sqrt{n} \leq \epsilon \\
                              \|\bs - \bs^k \|_2/\sqrt{n} \leq \epsilon
                        }
                  }\;
                  \max_{\|\bu\|_2/\sqrt{n} \leq C_0}
                  \Big\{
                        \frac1p \bu^\top \X \bl_1
                        +
                        \Theta_{\TAP}(\bu,\bl_1,\bl_2;\bbm,\bs)
                  \Big\}
                  \leq t
            \Big)
            \\
                  &\quad
                  \leq 
                  2\,\P\Big(
                        \min_{
                              \substack{
                                    \|\bl\|_2/\sqrt{n} = 1 \\
                                    \|\bbm - \bbm^k\|_2/\sqrt{n} \leq \epsilon \\
                                    \|\bs - \bs^k \|_2/\sqrt{n} \leq \epsilon
                              }
                        }\;
                        \max_{\|\bu\|_2/\sqrt{n} \leq C_0}
                        \Big\{
                              \frac1p\bu^\top\big(\X - \proj_{\bR_k}^\perp \X \proj_{\bV_k}^\perp\big)\bl_1
                              -
                              \frac1{p^{3/2}}
                              \|\proj_{\bR_k}^\perp\bu\| \< \proj_{\bV_k}^\perp \bg ,  \bl_1 \>
                                          \\
                  &\hspace{2in} +\frac1{p^{3/2}}
                              \|\proj_{\bV_k}^\perp \bl_1\| \< \proj_{\bR_k}^\perp \bh ,  \bu \>
                              +
                              \Theta_{\TAP}(\bu,\bl_1,\bl_2;\bbm,\bs)
                        \Big\}
                        \leq t
                  \Big),
      \end{aligned}\]
      where $\bg \sim \normal(\bzero,\id_p)$ and $\bh \sim \normal(\bzero,\id_n)$ are independent of each other and everything else.
      For all $\bl$ satisfying $\|\bl\|_2 = \sqrt{p}$,
      we have $|\< \proj_{\bV_k}^\perp \bg ,  \bl_1 \> - \< \bg , \bl_1\>|/p \leq \| \proj_{\bV_k} \bg \| / \sqrt{p} \gotop 0$.
      Likewise,
      for all $\|\bu \| \leq C_0\sqrt{n}$,
      we have $|\< \proj_{\bR_k}^\perp \bh ,  \bu \> - \< \bh , \bu\>|/p \leq
C_0 \sqrt{n/p}\,\| \proj_{\bR_k} \bh \| / \sqrt{p} \gotop 0$.
      Combining these observations with Lemma \ref{lem:conditional-X} and the fact that $n/p \rightarrow \delta$,
      the difference between the objective in the preceding display and $\AuxObj_k(\bu,\bl_1,\bl_2;\bbm,\bs;\bG_k,\bH_k,\bR_k,\bV_k,\bg,\bh)$ is uniformly bounded over the domain of optimization by a quantity which converges in probability to 0. The claim of the lemma follows.
\end{proof}

\subsubsection{Reduction to optimization on Wasserstein space}

Define
\begin{equation}
\begin{gathered}
      \EmpLaw(\bu,\bH_k,\bR_k,\bh)
            :=
            \frac1n
            \sum_{i=1}^n
            \delta_{u_i,h_i^1,\ldots,h_i^k,r_i^1,\ldots,r_i^k,h_i},
      \\
      \EmpLaw(\bl_1,\bl_2,\bbm,\bs,\bG_k,\bV_k,\bg)
            :=
            \frac1p
            \sum_{i=1}^p
            \delta_{l_{1i},l_{2i},m_i,s_i,g_i^1,\ldots,g_i^k,\nu_i^1,\ldots,\nu_i^k,g_i}.
\end{gathered}
\end{equation}
From the definition \eqref{eq:AuxObj-def}, it is evident that $\AuxObj_k$
is invariant under permutations of coordinates of its arguments, and hence
 is a function of its arguments only via the distributions $\EmpLaw(\bu,\bH_k,\bR_k,\bh)$ and $\EmpLaw(\bl_1,\bl_2,\bbm,\bs,\bG_k,\bV_k,\bg)$.

Let $\LL_2$ be the space of random variables on $\R$ with finite second moment,
equipped with the usual inner-product $\langle U,V \rangle_{\LL_2}=\E[UV]$, and let
$\WW_2(\R^k)$ denote the Wasserstein-2 space on $\R^k$, i.e.\ the space of
joint laws of $k$ variables belonging to $\LL_2$.
We will define a function on $\WW_2(\R^{2k+2}) \times \WW_2(\R^{2k+5})$,
which we denote by $\AuxObj_k^{(1)}$,
which is an extension of $\AuxObj_k$ in the sense that
\begin{equation}\label{eq:AuxObj_extension}
\begin{aligned}
    &\AuxObj_k^{(1)}(\EmpLaw(\bu,\bH_k,\bR_k,\bh),\EmpLaw(\bl_1,\bl_2,\bbm,\bs,\bG_k,\bV_k,\bg))
    \\
        &\qquad\qquad\qquad\qquad\qquad=
        \AuxObj_k(\bu,\bl_1,\bl_2;\bbm,\bs;\bG_k,\bH_k,\bR_k,\bV_k,\bg,\bh).
\end{aligned}
\end{equation}
To define the extension, denote the vector of random variables $(G_1,\ldots,G_k)$ by $G_1^k$, and likewise for $R_1^k$, $H_1^k$, $V_1^k$, etc.
Define functions $G^{(1)},H^{(1)}: (\LL_2)^{2k+2} \rightarrow \LL_2$ by
\begin{equation}
\begin{gathered}
    G^{(1)}(U,G_1^k,R_1^k,G)
      =
      \frac1{\sqrt{\delta}}
      \begin{pmatrix}
            G_1 
            &
            \cdots
            &
            G_k
      \end{pmatrix}
      \bK_{g,k}^{-1}
      \begin{pmatrix}
            \< R_1 , U \>_{\LL_2}
            \\
            \vdots
            \\
            \< R_k , U \>_{\LL_2}
      \end{pmatrix}
      +
      \| \proj_{R_1^k}^\perp U \|_{\LL_2} G,
      \\
      H^{(1)}(L_1,H_1^k,V_1^k,H)
            =
            \begin{pmatrix}
                  H_1
                  &
                  \cdots
                  &
                  H_k
            \end{pmatrix}
            \bK_{h,k}^{-1}
            \begin{pmatrix}
                  \< V_1 , L_1 \>_{\LL_2}
                  \\
                  \vdots
                  \\
                  \< V_k , L_1 \>_{\LL_2}
            \end{pmatrix}
            +
            \| \proj_{V_1^k}^\perp L_1 \|_{\LL_2} H,
\end{gathered}
\end{equation}
where $\proj_{R_1^k}^\perp$ denotes the projection in $\LL_2$ orthogonal
to the linear span of $(R_1,\ldots,R_k)$,
and likewise for $\proj_{V_1^k}^\perp$.
For any $(\mu_u,\mu_l) \in \WW_2(\R^{2k+2}) \times \WW_2(\R^{2k+5})$,
the extension is defined as
\[\begin{aligned}
    \AuxObj_k^{(1)}(\mu_u,\mu_l)
        &:=
        -\sqrt{\delta}\<G^{(1)}(U,G_1^k,R_1^k,G),L_1\>_{\LL_2}
        +\delta\< H^{(1)}(L_1,H_1^k,V_1^k,H) , U \>_{\LL_2}
        +\Theta_\TAP^{(1)}(\mu_u,\mu_l),
\end{aligned}\]
where
\[\begin{aligned}
      \Theta_\TAP^{(1)}(\mu_u,\mu_l)
            &:=
            \Theta_\bulk^{(1)}(\mu_u,\mu_l)
            +
            \Theta_\spike^{(1)}(\mu_l),
      \\
      \Theta_\bulk^{(1)}(\mu_u,\mu_l)
            &:=
            -\frac{\sigma^2\delta}2 \| U \|_{\LL_2}^2
            -\frac{\delta }{2V(M,S)}L_1^2
            +\frac{1}{2}(L_1,L_2)^\top \nabla^2[-\hs(M,S)] (L_1,L_2),
      \\
      \Theta_\spike^{(1)}(\mu_l)
            &:=
            -\frac{\delta
            \big(
                  \<M,L_1\>_{\LL_2}
                  -
                  \<1,L_2\>_{\LL_2}/2
            \big)^2}{V(M,S)^2},
      \\
      V(M,S) &:= \sigma^2 + \<1,S\>_{\LL_2} - \| M \|_{\LL_2}^2,
\end{aligned}\]
and $(U,H_1^k,R_1^k,H) \sim \mu_u$ and $(L_1,L_2,M,S,G_1^k,V_1^k,G) \sim \mu_l$.
It is then clear from the definition
(\ref{eq:AuxObj-def}) that (\ref{eq:AuxObj_extension}) holds.

Our goal in this section is to lower bound the right side of the final display in Proposition \ref{prop:gordon-post-amp} by the value of a saddle-point problem involving $\AuxObj_k^{(1)}$.
Recall the functions $\denoiser(\cdot)$ and $\sS(\cdot)$ from
(\ref{eq:denoiser}).
Define the distributions
\[\begin{aligned}
      \SE_{u,k}
            &:=
            \Law(H_1^k,R_1^k,H)
            \;\; 
            &\text{where}
            \;\;
            & H_1^k \sim \normal(0,\bK_{h,k})
            \indep
           	E \sim \normal(0,\sigma^2)
            \indep
            H \sim \normal(0,1),
      \\
            &&&R_i = H_i - E,
            \;\;
            i \leq k,
      \\
      \SE_{l,k}
            &:=
            \Law(G_1^k,V_1^k,G)
            \;\;
            &\text{where}
            \;\;
            &G_1^k \sim \normal(0,\bK_{g,k})
            \indep
            -V_1 \sim \sP_0
            \indep
            G \sim \normal(0,1),
      \\
            &&&V_{i+1} = \denoiser(-V_1 + G_i,\,\gamma_i) + V_1,
            \;\;
            i \leq k-1,
\end{aligned}\]
For any $\eps>0$, define the deterministic subsets of Wasserstein space
\[\begin{aligned}
      \bbSE_{u,k}^{(1)}
      :=
      \Big\{
            \Law(U,H_1^k,R_1^k,H)
            \Bigm|
            &\| U \|_{\LL_2} \leq C_0,\,
            (H_1^k,R_1^k,H) \sim \SE_{u,k}
      \Big\},
      \\
      \bbSE_{l,k}^{(1)}(\epsilon)
      :=
      \Big\{
            \Law(L_1,L_2,M,S,G_1^k,V_1^k,G)
            \Bigm|
            &\| L_1 \|_{\LL_2}^2 + \| L_2 \|_{\LL_2}^2 = 1,\,
            \| M - V_k + V_1 \|_{\LL_2} \leq \epsilon,\,
            \\
            &\| S - \sS(G_{k-1}-V_1,\,\gamma_{k-1}) \|_{\LL_2} \leq \epsilon,\,
            \\
            &(G_1^k,V_1^k,G) \sim \SE_{l,k}
            \Big\},
\end{aligned}\]
and the random subsets of Wasserstein space
\[\begin{aligned}
      \bbSEhat_{u,k}^{(1)}
      :=
      \Big\{
            \Law(U,\Hhat_1^k,\Rhat_1^k,\Hhat)
            \Bigm|
            &\| U \|_{\LL_2} \leq C_0,\,
            (\Hhat_1^k,\Rhat_1^k,\Hhat) \sim \EmpLaw(\bH_k,\bR_k,\bh),
            \\
            &\text{law consists of $n$ equal mass atoms}
      \Big\},
      \\
      \bbSEhat_{l,k}^{(1)}(\epsilon)
      :=
      \Big\{
            \Law(\Lhat_1,\Lhat_2,\Mhat,\Shat,\Ghat_1^k,\Vhat_1^k,\Ghat)
            \Bigm|
            &\| \Lhat_1 \|_{\LL_2}^2 + \| \Lhat_2 \|_{\LL_2}^2 = 1,\,
            \| \Mhat - \Vhat_k + \Vhat_1 \|_{\LL_2} \leq \epsilon,\,
            \\
            &\| \Shat - \sS(\Ghat_{k-1}-\Vhat_1,\,\gamma_{k-1}) \|_{\LL_2} \leq \epsilon,\,
            \\
            &(\Ghat_1^k,\Vhat_1^k,\Ghat) \sim \EmpLaw(\bG_k,\bV_k,\bg),
            \\
            &\text{law consists of $p$ equal mass atoms}
            \Big\}.
\end{aligned}\]
Note that since $\bbm^k=\bv^k+\bbeta_0=\bv^k-\bv^1$ and
$\bs^k=\sS(\bg^{k-1}+\bbeta_0)=\sS(\bg^{k-1}-\bv^1)$, we have that
$\bbSEhat_{l,k}^{(1)}(\epsilon)$ is the image of the set
    $\{\bl,\bbm,\bs:\|\bl\|_2/\sqrt{p} = 1,\,\|\bbm - \bbm^k\|_2/\sqrt{p} \leq
\epsilon,\,\|\bs - \bs^k \|_2/\sqrt{p} \leq \epsilon\}$ under the map
    $(\bl,\bbm,\bs) \mapsto \EmpLaw(\bl_1,\bl_2,\bbm,\bs,\bG_k,\bV_k,\bg)$.
Similarly $\bbSEhat_{u,k}^{(1)}$ is the image of the set $\{\bu:\| \bu
\|_2/\sqrt{n} \leq C_0\}$ under the map $\bu\mapsto \EmpLaw(\bu,\bH_k,\bR_k,\bh)$.
Then, by the equality (\ref{eq:AuxObj_extension}),
the optimization of $\AuxObj_k$ is equivalently expressed as
    \begin{equation}\label{eq:AuxObjkequiv}
    \begin{aligned}
        \min_{
            \substack{
                \|\bl\|_2/\sqrt{p} = 1 \\
                \|\bbm - \bbm^k\|_2/\sqrt{p} \leq \epsilon \\
                \|\bs - \bs^k \|_2/\sqrt{p} \leq \epsilon
            }
        }\;
        &\max_{\|\bu\|_2 \leq C_0\sqrt{n}}\;
        \AuxObj_k(\bu,\bl_1,\bl_2;\bbm,\bs;\bG_k,\bH_k,\bR_k,\bV_k,\bg,\bh)
        \\
        &\qquad\qquad=
        \min_{\muhat_l \in \bbSEhat_{l,k}^{(1)}(\epsilon)}\;\;
        \max_{\muhat_u \in \bbSEhat_{u,k}^{(1)}}\;
            \AuxObj_k^{(1)}(\muhat_u,\muhat_l).
    \end{aligned}
    \end{equation}

The main result of this section is the following lower bound.
\begin{lemma}[Reduction to optimization on Wasserstein space]
\label{lem:reduction-to-wasserstein-opt}
    We have
    \begin{equation}\label{eq:AuxObjsaddle}
    \begin{aligned}
        &\liminf_{\epsilon \rightarrow 0}\;\;
            \liminf_{k \rightarrow \infty}\;\;
            \pliminf_{n \rightarrow \infty}\;\;
            \min_{
                  \substack{
                        \|\bl\|_2/\sqrt{p} = 1 \\
                        \|\bbm - \bbm^k\|_2/\sqrt{p} \leq \epsilon \\
                        \|\bs - \bs^k \|_2/\sqrt{p} \leq \epsilon
                  }
            }\;
            \max_{\|\bu\|_2 \leq C_0\sqrt{n}}\;
            \AuxObj_k(\bu,\bl_1,\bl_2;\bbm,\bs;\bG_k,\bH_k,\bR_k,\bV_k,\bg,\bh)
        \\
        &\qquad\qquad\qquad\qquad\qquad\qquad
            \geq
            \liminf_{\epsilon \rightarrow 0}\;\;
            \liminf_{k \rightarrow \infty}\;\;
            \min_{\mu_l \in \bbSE_{l,k}^{(1)}(\epsilon)}\;\;
            \max_{\mu_u \in \bbSE_{u,k}^{(1)}}\;
            \AuxObj_k^{(1)}(\mu_u,\mu_l).
    \end{aligned}
    \end{equation}
\end{lemma}
\noindent The benefit of working with $\AuxObj_k^{(1)}$ in place of $\AuxObj_k$ is that $\AuxObj_k^{(1)}$ is deterministic,
and the saddle-point problem on the right side of (\ref{eq:AuxObjsaddle})
involves the exact state evolution distributions rather than empirical approximations to them.

The proof of Lemma \ref{lem:reduction-to-wasserstein-opt} relies on
(\ref{eq:AuxObjkequiv}) and the following two lemmas.
First, we bound the random maximization with a deterministic maximization.
\begin{lemma}
\label{lem:U-emp-to-pop}
Fix any $k \geq 1$ and constant $C>0$.
Then there exists a function $\delta:(0,\infty) \rightarrow (0,\infty]$ with
$\delta(x) \rightarrow 0$ as $x \rightarrow 0$, depending on $k,C$ but
not on $n,p$, such that the following holds: For any
$\mu_l=\Law(L_1,L_2,M,S,G_1^k,V_1^k,G)$ satisfying
$\|G\|_{\LL_2},\|G_i\|_{\LL_2},\|V_i\|_{\LL^2} \leq C$,
$\|L_1\|_{\LL_2} \leq 1$, and $\< G,L_1\>_{\LL_2} \geq 0$, we have
        \[\sup_{\muhat_u \in \bbSEhat_{u,k}^{(1)}}
            \AuxObj_k^{(1)}(\muhat_u,\mu_l)
            \geq
            \sup_{\mu_u \in \bbSE_{u,k}^{(1)}}
            \AuxObj_k^{(1)}(\mu_u,\mu_l)
            -
            \delta\big(W_2(\EmpLaw(\bH_k,\bR_k,\bh),\SE_{u,k})\big)\]
where $W_2(\cdot,\cdot)$ is the Wasserstein-2 distance on $\WW_2(\R^{2k+1})$.
\end{lemma}
\noindent Second, we make a replacement of the minimization variable that incurs only a small approximation error.
\begin{lemma}
\label{lem:muhat-l-to-mu-l}
Fix any $k \geq 1$ and constant $C>0$. Then
there exists a function $\delta:(0,\infty) \rightarrow (0,\infty]$ with
$\delta(x) \rightarrow 0$ as $x \rightarrow 0$, depending on $k,C$ but not
on $n,p$, such that the following holds:
For any $\muhat_l = \Law(\Lhat_1,\Lhat_2,\Mhat,\Shat,\Ghat_1^k,\Vhat_1^k,\Ghat)$
satisfying $\| \Lhat_1 \|_{\LL_2}^2 + \|\Lhat_2\|_{\LL_2}^2 = 1$,
$\| \Mhat - \Vhat_k + \Vhat_1 \|_{\LL_2} \leq \epsilon$,
and
$\| \Shat - \sS(\Ghat_{k-1}-\Vhat_1,\,\gamma_{k-1}) \|_{\LL_2} \leq \epsilon$,
    there exists $\mu_l=\Law(L_1,L_2,M,S,G_1^k,V_1^k,G)$ such that
    $\Law(\Lhat_1,\Lhat_2,\Mhat,\Shat)=\Law(L_1,L_2,M,S)$,
    $\mu_l \in \bbSE_{l,k}^{(1)}\big(\epsilon+\delta\big(
                W_2(\Law(\Ghat_1^k,\Vhat_1^k,\Ghat),\SE_{l,k})
            \big) \big)$, and
    \begin{equation}\label{eq:mulswap}
        \Big|
            \sup_{\mu_u \in \bbSE_{u,k}(C)}
            \AuxObj_k^{(1)}(\mu_u,\muhat_l)
            -
            \sup_{\mu_u \in \bbSE_{u,k}(C)}
            \AuxObj_k^{(1)}(\mu_u,\mu_l)
        \Big|
        \leq
        \delta\big(
            W_2(\Law(\Ghat_1^k,\Vhat_1^k,\Ghat),\SE_{l,k})
        \big).
    \end{equation}
\end{lemma}
\noindent The proofs of these two lemmas are given below.
First, we show they imply Lemma \ref{lem:reduction-to-wasserstein-opt}.
\begin{proof}[Proof of Lemma \ref{lem:reduction-to-wasserstein-opt}]
Define $\bbSEhat_{l,k}^{(1)}(\eps,\eps')$ by
\begin{equation}
\begin{aligned}
  \bbSEhat_{l,k}^{(1)}(\epsilon,\eps')
  :=
  \Big\{
        \Law(\Lhat_1,\Lhat_2,\Mhat,\Shat,\Ghat_1^k,\Vhat_1^k,\Ghat)
        \Bigm|
        &\| \Lhat_1 \|_{\LL_2}^2 + \| \Lhat_2 \|_{\LL_2}^2 = 1,\,
        \| \Mhat - \Vhat_k + \Vhat_1 \|_{\LL_2} \leq \epsilon,\,
        \\
        &\| \Shat - \sS(\Ghat_{k-1}-\Vhat_1,\,\gamma_{k-1}) \|_{\LL_2} \leq \epsilon,\,
        \\
        &
        W_2\big((\Ghat_1^k,\Vhat_1^k,\Ghat),\SE_{l,k}\big) 
        \leq \eps'
        \Big\}.
\end{aligned}
\end{equation}
Because
$W_2\big(\EmpLaw(\bG_k,\bV_k,\bg),\SE_{l,k}\big) \gotop 0$ by
Proposition \ref{prop:amp-se}, for fixed $\eps' > 0$,
$\bbSEhat_{l,k}^{(1)}(\epsilon) \subseteq \bbSEhat_{l,k}^{(1)}(\epsilon,\eps')$
with probability going to 1 as $n,p \rightarrow \infty$.
Thus, using \eqref{eq:AuxObjkequiv},
\begin{equation}
\begin{aligned}
    &\pliminf_{n \rightarrow \infty}\;\;
        \min_{
              \substack{
                    \|\bl\|_2/\sqrt{p} = 1 \\
                    \|\bbm - \bbm^k\|_2/\sqrt{p} \leq \epsilon \\
                    \|\bs - \bs^k \|_2/\sqrt{p} \leq \epsilon
              }
        }\;
        \max_{\|\bu\|_2 \leq C_0\sqrt{n}}\;
        \AuxObj_k(\bu,\bl_1,\bl_2;\bbm,\bs;\bG_k,\bH_k,\bR_k,\bV_k,\bg,\bh)
    \\
    &\qquad\qquad\qquad\qquad\qquad\qquad
        \geq
        \pliminf_{n \rightarrow \infty}
        \min_{\hmu_l \in \bbSEhat_{l,k}^{(1)}(\eps,\eps')}\;\;
        \max_{\hmu_u \in \bbSEhat_{u,k}^{(1)}}\;
        \AuxObj_k^{(1)}(\hmu_u,\hmu_l).
\end{aligned}
\end{equation}

We claim $W_2\big(\Law(\Ghat_1^k,\Vhat_1^k,-\Ghat),\SE_{l,k}\big) =
W_2\big(\Law(\Ghat_1^k,\Vhat_1^k,\Ghat),\SE_{l,k}\big)$.
Indeed, let $(\Ghat_1^k,\Vhat_1^k,\Ghat,G_1^k,V_1^k,G)$ be a coupling between
$\Law(\Ghat_1^k,\Vhat_1^k,\Ghat)$ and $\SE_{l,k}$.
Then $(\Ghat_1^k,\Vhat_1^k,-\Ghat,G_1^k,V_1^k,-G)$ is a coupling between
$\Law(\Ghat_1^k,\Vhat_1^k,-\Ghat)$ and $\SE_{l,k}$ with the same $\ell_2$
distance, where $(G_1^k,V_1^k,-G) \sim \SE_{l,k}$ because
$G\sim\normal(0,1)$ independently of $G_1^k,V_1^k$. So the claim follows.
Thus, if $\hmu_l := \Law(\Lhat_1,\Lhat_2,\Mhat,\Shat,\Ghat_1^k,\Vhat_1^k,\Ghat) \in \bbSEhat_{l,k}^{(1)}(\eps,\eps')$,
then so too is $\hmu_l^- := \Law(\Lhat_1,\Lhat_2,\Mhat,\Shat,\Ghat_1^k,\Vhat_1^k,-\Ghat)$.
Moreover, if $\< \Ghat , \Lhat_1 \>_{\LL_2} < 0$,
then $-\| \proj_{\Rhat_1^k}^\perp \Uhat \|_{\LL_2} \< \Ghat , \Lhat_1 \>_{\LL_2}
\geq -\| \proj_{\Rhat_1^k}^\perp \Uhat \|_{\LL_2} \< -\Ghat , \Lhat_1
\>_{\LL_2}$. This is the only term of $\AuxObj_k(\hmu_u , \hmu_l)$ that depends
on $\Ghat$, so we conclude that
$\AuxObj_k(\hmu_u , \hmu_l) \geq \AuxObj_k(\hmu_u , \hmu_l^-)$. Thus,
\begin{equation}\label{eq:Gsignflip}
    \min_{\hmu_l \in \bbSEhat_{l,k}^{(1)}(\eps,\eps')}\;\;
    \max_{\hmu_u \in \bbSEhat_{u,k}^{(1)}}\;
    \AuxObj_k^{(1)}(\hmu_u,\hmu_l)
		=
        \min_{\substack{\hmu_l \in \bbSEhat_{l,k}^{(1)}(\eps,\eps') \\ \< \Ghat , \Lhat_1 \>_{\LL_2} \geq 0}}\;\;
        \max_{\hmu_u \in \bbSEhat_{u,k}^{(1)}}\;
        \AuxObj_k^{(1)}(\hmu_u,\hmu_l).
\end{equation}
Because $\| \Ghat_i \|_{\LL_2}$, $\|\Ghat\|_{\LL_2}$, $\|\Vhat_i \|_{\LL_2} \leq
C$ on $\bbSEhat_{l,k}^{(1)}(\eps,\eps')$ for a constant $C>0$, Lemma \ref{lem:U-emp-to-pop} gives
    \[\begin{aligned}
        &\min_{\substack{\hmu_l \in \bbSEhat_{l,k}^{(1)}(\eps,\eps') \\ \< \Ghat , \Lhat_1 \>_{\LL_2} \geq 0}}\;\;
        \max_{\hmu_u \in \bbSEhat_{u,k}^{(1)}}\;
        \AuxObj_k^{(1)}(\hmu_u,\hmu_l)
    \\
        &\qquad\qquad\geq
        \min_{\substack{\hmu_l \in \bbSEhat_{l,k}^{(1)}(\eps,\eps') \\ \< \Ghat , \Lhat_1 \>_{\LL_2} \geq 0}}\;
        \max_{\mu_u \in \bbSE_{u,k}^{(1)}} \AuxObj_k^{(1)}(\mu_u,\muhat_l)
        -
        \delta\big(W_2(\EmpLaw(\bH_k,\bR_k,\bh),\SE_{u,k})\big)
    \\
    	&\qquad\qquad=
    	\min_{\hmu_l \in \bbSEhat_{l,k}^{(1)}(\eps,\eps')}\;
        \max_{\mu_u \in \bbSE_{u,k}^{(1)}} \AuxObj_k^{(1)}(\mu_u,\muhat_l)
        -
        \delta\big(W_2(\EmpLaw(\bH_k,\bR_k,\bh),\SE_{u,k})\big),
    \end{aligned}\]
    where the last equality holds by the same argument as in (\ref{eq:Gsignflip}).
    For each $\muhat_l \in \bbSEhat_{l,k}^{(1)}(\epsilon,\eps')$,
    we can select $\mu_l \in \bbSE_{l,k}^{(1)}\big(\epsilon + \delta(\eps')\big)$ as in Lemma \ref{lem:muhat-l-to-mu-l}.
    Thus, denoting the right side of the previous display by $\mathrm{RHS}$,
    \begin{equation}
    \begin{aligned}
        &\mathrm{RHS}\geq
        \min_{\mu_l \in \bbSE_{l,k}^{(1)}(\,\epsilon+ \delta(\eps')\,)}
        \;\;
        \max_{\mu_u \in \bbSE_{u,k}^{(1)}} \AuxObj_k^{(1)}(\mu_u,\mu_l)
    	-
        \delta\big(\eps'\big)
        -
        \delta\big(W_2(\EmpLaw(\bH_k,\bR_k,\bh),\SE_{u,k})\big).
    \end{aligned}
    \end{equation}
    Because $\delta(x) \rightarrow 0$ as $x \to 0$,
    and from Proposition \ref{prop:amp-se}
	$W_2(\EmpLaw(\bH_k,\bR_k,\bh),\SE_{u,k}) \gotop 0$ as $n,p \rightarrow \infty$,
    taking $n \rightarrow \infty$ followed by $\eps' \rightarrow 0$ gives
    \begin{equation}
    \begin{aligned}
        &\pliminf_{n \rightarrow \infty}\;\;
            \min_{
                  \substack{
                        \|\bl\|_2/\sqrt{p} = 1 \\
                        \|\bbm - \bbm^k\|_2/\sqrt{p} \leq \epsilon \\
                        \|\bs - \bs^k \|_2/\sqrt{p} \leq \epsilon
                  }
            }\;
            \max_{\|\bu\|_2 \leq C_0\sqrt{n}}\;
            \AuxObj_k(\bu,\bl_1,\bl_2;\bbm,\bs;\bG_k,\bH_k,\bR_k,\bV_k,\bg,\bh)
        \\
        &\qquad\qquad\qquad\qquad\qquad\qquad
            \geq
            \min_{\mu_l \in \bbSE_{l,k}^{(1)}(2\epsilon)}\;\;
            \max_{\mu_u \in \bbSE_{u,k}^{(1)}}\;
            \AuxObj_k^{(1)}(\mu_u,\mu_l).
    \end{aligned}
    \end{equation}
    The result follows by taking $k \rightarrow \infty$ followed by $\epsilon
\rightarrow 0$.
\end{proof}

\noindent We now prove Lemmas \ref{lem:U-emp-to-pop} and \ref{lem:muhat-l-to-mu-l}.


\begin{proof}[Proof of Lemma \ref{lem:U-emp-to-pop}]
    Let
        \[\AA
        =
        \Big\{
            \Law(U,\Hhat_1^k,\Rhat_1^k,\Hhat)
            \Bigm|
            \| U \|_{L_2} \leq C_0,
            \\
            (\Hhat_1^k,\Rhat_1^k,\Hhat) \sim \EmpLaw(\bH_k,\bR_k,\bh)
        \Big\},\]
which removes the constraint of $\bbSEhat_{u,k}^{(1)}$ that the law of $U$ is
atomic. First we show that 
    \begin{equation}
    \label{eq:SE-to-SEhat-intermediate}
        \sup_{\muhat_u \in \AA}
            \AuxObj_k^{(1)}(\muhat_u,\mu_l)
            \geq
            \sup_{\mu_u \in \bbSE_{u,k}^{(1)}}
            \AuxObj_k^{(1)}(\mu_u,\mu_l)
            -
            \delta\big(W_2(\EmpLaw(\bH_k,\bR_k,\bh),\SE_{u,k})\big).
    \end{equation}
    Consider $\mu_u \in \bbSE_{u,k}^{(1)}$.
    Let $\Pi_u^{(1)}$ be an optimal $\ell_2$ coupling between $\EmpLaw(\bH_k,\bR_k,\bh)$ and $\SE_{u,k}$ (its existence is guaranteed by \cite[Theorem 4.1]{villani2008optimal}).
    Then, by the gluing lemma (see, e.g., \cite[pg.~23]{villani2008optimal}),
    there exists a joint law
$\Pi_u^{(2)}=\Law(U,H_1^k,R_1^k,H,\Hhat_1^k,\Rhat_1^k,\Hhat)$ including $U$
such that $(U,H_1^k,R_1^k,H) \sim \mu_u$ and $(H_1^k,R_1^k,H,\Hhat_1^k,\Rhat_1^k,\Hhat) \sim \Pi_u^{(1)}$.
    Because $\| L_1 \|_{\LL_2} \leq 1$, $\|U\|_{\LL_2} \leq C_0$, and $\|V_i\|_{\LL_2} \leq C$ for all $i$,
    the function $(\tilde H_1^k,\tilde H) \mapsto \<H^{(1)}(L_1,\tilde
H_1^k,V_1^k,\tilde H),U\>$ (as a function on $\LL_2$) is Lipschitz
    with Lipschitz constant depending on $C_0,C > 0$.
    Likewise,
    because $\| L_1 \|_{\LL_2} \leq 1$, $\|U\|_{\LL_2} \leq C_0$, $\|G\|_{\LL_2}
\leq C$, $\|G_i\|_{\LL_2} \leq C$ for all $i$, and
$(\E[R_iR_j])_{i,j \in [k]} \succ c\,\id_k$ and $C'>\max_{i=1}^k \E[R_i^2]$ 
for some constants $C',c>0$ when $(H_1^k,R_1^k,H) \sim \SE_{u,k}$,
    the function $\tilde R_1^k \mapsto \<G^{(1)}(U,G_1^k,\tilde
R_1^k,G),L_1\>_{\LL_2}$ is Lipschitz on an $\LL_2$-neighborhood of $R_1^k$,
    with Lipschitz constant depending on $C_0,C,C',c > 0$.
(Here,
    the conditions $(\E[R_iR_j])_{i,j \in [k]} \succ c\,\id_k$ 
and $C'>\max_{i=1}^k \E[R_i^2]$ are used to check that
$\|\proj_{R_1^k}^\perp U \|_{\LL_2}$ is Lipschitz in $R_1^k$.)
Then letting $\muhat_u' = \Law(U,\Hhat_1^k,\Rhat_1^k,\Hhat) \in \AA$,
    the above Lipschitz properties and the identity
$\|(H_1^k,R_1^k,H)-(\Hhat_1^k,\Rhat_1^k,\Hhat)\|_{\LL_2}=W_2(\EmpLaw(\bH_1,\bR_k,\bh),\SE_{u,k})$ by definition of the optimal coupling $\Pi_u^{(1)}$ imply that
        \[\sup_{\muhat_u \in \AA} 
        \AuxObj_k^{(1)}(\muhat_u,\mu_l)
        \geq
        \AuxObj_k^{(1)}(\muhat_u',\mu_l)
        \geq
        \AuxObj_k^{(1)}(\mu_u,\mu_l)
        -
        \delta\big(W_2(\EmpLaw(\bH_k,\bR_k,\bh),\SE_{u,k})\big),\]
    where $\delta(\,\cdot\,)$ is as in the statement of the lemma.
Here, the error $\delta\big(W_2(\EmpLaw(\bH_k,\bR_k,\bh),\SE_{u,k})\big)$ is the
same for all $\mu_u \in \bbSE_{u,k}^{(1)}$, so taking the supremum over $\mu_u$ gives \eqref{eq:SE-to-SEhat-intermediate}.

    Next, we show
    \begin{equation}\label{eq:SEhat-intermediate-to-SEhat}
        \sup_{\muhat_u \in \bbSEhat_{u,k}^{(1)}}
            \AuxObj_k^{(1)}(\muhat_u,\mu_l)
        \geq
        \sup_{\muhat_u \in \AA}
            \AuxObj_k^{(1)}(\muhat_u,\mu_l)
    \end{equation}
    Consider $\muhat_u=\Law(U,\Hhat_1^k,\Rhat_1^k,\Hhat) \in \AA$.  
    Because $\EmpLaw(\bH_k,\bR_k,\bh)$ consists of $n$ equal sized atoms,
    there exist disjoint events $(\Omega_i)_{i \in [n]}$ with $\P(\Omega_i) = 1/n$ and $(\Hhat_1,\ldots,\Hhat_k,\Rhat_1,\ldots,\Rhat_k,\Hhat)(\omega) = \allowbreak\sum_{i=1}^n \allowbreak(h_i^1,\ldots,h_i^k,\allowbreak r_i^1,\ldots,r_i^k,h_i) \allowbreak\mathbb{I}_{\Omega_i}(\omega)$,
    where $\omega$ denotes an element of the probability space and
    $\mathbb{I}_{\Omega_i}$ denotes the indicator of the event $\Omega_i$.
Let $\cI$ be the sigma-field generated by $\Omega_1,\ldots,\Omega_n$,
and define $\Uhat=\E[U \mid \cI]$.
    By Jensen's inequality, $\|\Uhat\|_{\LL_2} \leq \| U \|_{\LL_2} \leq C_0$.
    Moreover, $\<\Rhat_i,\Uhat\>_{\LL_2}=\< \Rhat_i, U \>_{\LL_2}$,
    because $\Rhat_i$ is $\cI$-measurable, and similarly
    $\< \Hhat_i,\Uhat\>_{\LL_2} = \<\Hhat_i,U\>_{\LL_2}$ and $\<
\Hhat,\Uhat\>_{\LL_2}=\<\Hhat,U\>_{\LL_2}$.
    Combining these statements,
    we then have $\| \proj_{R_1^k}^\perp \Uhat \|_{\LL_2} \leq \|
\proj_{R_1^k}^\perp U \|_{\LL_2}$.
    Because $\<G,L_1\>_{\LL_2} \geq 0$ by assumption of the lemma,
    these equalities and inequalities imply $\AuxObj_k^{(1)}(\Law(\Uhat,\Hhat_1^k,\Rhat_1^k,\Hhat),\mu_l)\geq\AuxObj_k^{(1)}(\muhat_u,\mu_l)$.
    Finally, because $\Law(\Uhat,\Hhat_1^k,\Rhat_1^k,\Hhat) \in
\bbSEhat_{u,k}^{(1)}$, (\ref{eq:SEhat-intermediate-to-SEhat}) follows.
The lemma follows from (\ref{eq:SE-to-SEhat-intermediate}) and
(\ref{eq:SEhat-intermediate-to-SEhat}).
\end{proof}

\begin{proof}[Proof of Lemma \ref{lem:muhat-l-to-mu-l}]
Let $\muhat_l=\Law(\Lhat_1,\Lhat_2,\Mhat,\Shat,\Ghat_1^k,\Vhat_1^k,\Ghat)$
be as in the statement of the lemma, and
    let $\Pi^{(1)}$ be an optimal $\ell_2$ coupling between
$\Law(\Ghat_1^k,\Vhat_1^k,\Ghat)$ and $\SE_{l,k}$. By the gluing lemma,
    there exists a joint law
$\Pi^{(2)}=\Law(\Lhat_1,\Lhat_2,\Mhat,\Shat,\Ghat_1^k,\Vhat_1^k,\Ghat,G_1^k,V_1^k,G)$
such that
$(\Lhat_1,\Lhat_2,\Mhat,\Shat,\Ghat_1^k,\Vhat_1^k,\Ghat) \sim \muhat_l$ and $(\Ghat_1^k,\Vhat_1^k,\Ghat,G_1^k,V_1^k,G) \sim \Pi^{(1)}$.
    We then define $\mu_l = \Law(\Lhat_1,\Lhat_2,\Mhat,\Shat,G_1^k,V_1^k,G)$.
    Because $\Law(\Lhat_1,\Lhat_2,\Mhat,\Shat)$ is the same under $\muhat_l$ and $\mu_l$,
    for any fixed $\mu_u$ 
    we have $\Theta_{\bulk}^{(1)}(\mu_u,\mu_l) = \Theta_{\bulk}^{(1)}(\mu_u,\muhat_l)$
    and $\Theta_{\spike}^{(1)}(\mu_l) = \Theta_{\spike}^{(1)}(\muhat_l)$,
    whence $\Theta_{\TAP}^{(1)}(\mu_u,\mu_l) = \Theta_{\TAP}^{(1)}(\mu_u,\muhat_l)$.

    For any $\mu_u \in \bbSE_{u,k}^{(1)}$, if $(U,H_1^k,R_1^k,H) \sim \mu_u$,
    then $(\tilde G_1^k,\tilde G) \mapsto \< G^{(1)}(U,\tilde G_1^k,R_1^k,\tilde
G),\Lhat_1\>_{\LL_2}$ is $C$-Lipschitz for a constant $C>0$,
    because $R_i$, $U$, and $\Lhat_1$ all have bounded $\LL_2$-norms.
    In an $\LL_2$-neighborhood around $V_1^k$, the mapping $\tilde V_1^k \mapsto
\< H^{(1)}(\Lhat_1,H_1^k,\tilde V_1^k,H) , U \>_{\LL_2}$ is also $C$-Lipschitz
because $H_i$, $H$, and $\Lhat_1$ have bounded $\LL_2$-norms, and
$(\E[V_iV_j])_{i,j \in [k]} \succ c\,\id_k$ and $C'>\max_{i=1}^k \E[V_i^2]$ for
some constants $C',c>0$.  Thus, since
$\|(G_1^k,V_1^k,G)-(\Ghat_1^k,\Vhat_1^k,\Ghat)\|_{\LL_2}
=W_2(\Law(\Ghat_1^k,\Vhat_1^k,\Ghat),\SE_{l,k})$, we have
        \begin{align*}
        \big|
            \< G^{(1)}(U,G_1^k,R_1^k,G),\Lhat_1 \>_{\LL_2}
            -
            \< G^{(1)}(U,\Ghat_1^k,R_1^k,\Ghat),\Lhat_1 \>_{\LL_2}
        \big|
        &\leq 
        \delta\big(
            W_2(\Law(\Ghat_1^k,\Vhat_1^k,\Ghat),\SE_{l,k})
        \big),\\
        \big|
            \< H^{(1)}(\Lhat_1,H_1^k,V_1^k,H),U \>_{\LL_2}
            -
            \allowbreak
            \< H^{(1)}(\Lhat_1,\allowbreak
H_1^k,\allowbreak\Vhat_1^k,\allowbreak H),\allowbreak U \>_{\LL_2}
        \big|
        &\leq 
        \delta\big(
            W_2(\Law(\Ghat_1^k,\Vhat_1^k,\Ghat),\SE_{l,k})
        \big).
\end{align*}
    Because this holds uniformly over $\mu_u \in \bbSE_{u,k}^{(1)}$,
    we conclude (\ref{eq:mulswap}).
    Moreover, 
    $\| \Mhat - V_k + V_1 \|_{\LL_2} \leq \| \Mhat - \Vhat_k + \Vhat_1
\|_{\LL_2} + \delta\big(W_2(\Law(\Ghat_1^k,\Vhat_1^k,\Ghat),\SE_{l,k})\big)$
    and 
    $\| \Shat - \sS(G_{k-1}-V_1,\,\gamma_{k-1}) \|_{\LL_2} \leq \| \Shat -
\sS(\Ghat_{k-1}-\Vhat_1,\,\gamma_{k-1}) \|_{\LL_2} +
\delta\big(W_2(\Law(\Ghat_1^k,\Vhat_1^k,\Ghat),\SE_{l,k})\big)$, the latter holding
because $\sS(\,\cdot\,,\gamma)$ is Lipschitz. Then we obtain also
    $
        \mu_l \in \bbSE_{l,k}^{(1)}\big(
            \epsilon
            + 
            \delta\big(
                W_2(\Law(\Ghat_1^k,\Vhat_1^k,\Ghat),\SE_{l,k})
            \big)
        \big) 
    $
    as desired.
\end{proof}

\subsubsection{Dimensionality reduction in Wasserstein space}

In what remains,
we will be working with laws $\mu_u$ and $\mu_l$. We will refer to random
variables $U,H_1^k,R_1^k,H$ etc., which will be implicitly
distributed according to these laws.

Define the deterministic subsets of Wasserstein space
\begin{equation}
\begin{aligned}
      \bbSE_{u,k}^{(2)}
      :=
      \Big\{
            \Law(U,H_k,R_k,H)
            \Bigm|
            &\| U \|_{\LL_2} \leq C_0,\,
      \\
            &H_k \sim \normal(0,\delta \gamma_k^{-1}-\sigma^2)
            \indep
            H \sim \normal(0,1)
            \indep
            E \sim \normal(0,\sigma^2),
      \\
            &
            R_k = H_k - E
      \Big\},
      \\
      \bbSE_{l,k}^{(2)}(\epsilon)
      :=
      \Big\{
            \Law(L_1,L_2,M,S,G_{k-1},V_1,V_k,G)
            \Bigm|
            &\| L_1 \|_{\LL_2}^2 + \| L_2 \|_{\LL_2}^2 = 1,\,
            \| M - V_k + V_1 \|_{\LL_2} \leq \epsilon,\,
            \\
            &\| S - \sS(-V_1 + G_{k-1},\,\gamma_{k-1}) \|_{\LL_2} \leq \epsilon,\,
            \\
            &G_{k-1} \sim \normal(0,\gamma_{k-1}^{-1})
            \indep
            G \sim \normal(0,1) \indep
            -V_1 \sim \sP_0, 
            \\
            & V_k = \denoiser(-V_1 + G_{k-1},\,\gamma_{k-1}) + V_1
            \Big\}.
\end{aligned}
\end{equation}
These replace the random vectors $H_1^k,R_1^k,G_1^k,V_1^k$ defining
$\bbSE_{u,k}^{(1)}$ and $\bbSE_{l,k}^{(1)}(\epsilon)$ by scalar random
variables $H_k,R_k,G_{k-1},V_1,V_k$ having the same marginal laws, thus reducing
the dimension of the optimization over Wasserstein space.

For $\mu_u \in \bbSE_{u,k}^{(2)}$ and $\mu_l \in \bbSE_{l,k}^{(2)}(\epsilon)$,
define
\[\begin{aligned}
      G^{(2)}(U,G_{k-1},R_k,G)
            &=
            \frac{\gamma_{k-1}}{\sqrt{\delta}} G_{k-1} \< R_k , U \>_{\LL_2}
            +
            \|\proj_{R_k}^\perp U \|_{\LL_2}G,
      \\
      H^{(2)}(L_1,H_k,V_k,H)
            &=
            \frac1{\delta \gamma_k^{-1} - \sigma^2} H_k\< V_k, L_1 \>_{\LL_2} +
\| \proj_{V_k}^\perp L_1 \|_{\LL_2} H,
\end{aligned}\]
and 
\[
      \AuxObj_k^{(2)}(\mu_u,\mu_l)
            =
            -\sqrt{\delta}\<G^{(2)}(U,G_{k-1},R_k,G),L_1\>_{\LL_2}
            +\delta\< H^{(2)}(L_1,H_k,V_k,H) , U \>_{\LL_2}
            +\Theta_\TAP^{(2)}(\mu_u,\mu_l),
\]
where $\Theta_\TAP^{(2)}$ is defined identically to $\Theta_\TAP^{(1)}$ and
depends only on $\Law(U)$ and $\Law(L_1,L_2,M,S)$.
(We use the different notation $\Theta_\TAP^{(2)}$ because its arguments
$(\mu_u,\mu_l)$ are now joint laws on lower-dimensional spaces.)
\begin{lemma}[Dimensionality reduction in Wasserstein space]
\label{lem:wass-dim-reduce}
      We have
      \[
            \liminf_{\epsilon \rightarrow 0}\;\;
            \liminf_{k \rightarrow \infty}\;\;
            \min_{\mu_l \in \bbSE_{l,k}^{(1)}(\epsilon)}\;\;
            \max_{\mu_u \in \bbSE_{u,k}^{(1)}}
            \AuxObj_k^{(1)}(\mu_u,\mu_l)
            \geq
            \liminf_{\epsilon \rightarrow 0}\;\;
            \liminf_{k \rightarrow \infty}\;\;
            \min_{\mu_l \in \bbSE_{l,k}^{(2)}(\epsilon)}\;\;
            \max_{\mu_u \in \bbSE_{u,k}^{(2)}}\;
            \AuxObj_k^{(2)}(\mu_u,\mu_l).
       \]
\end{lemma}

\begin{proof}
      The key idea is to find maps
      \begin{equation}
            \Reduce_{l,k}^{(2)}: \bbSE_{l,k}^{(1)}(\epsilon) \rightarrow \bbSE_{l,k}^{(2)}(\epsilon),
            \qquad
            \Reduce_{u,k}^{(2)}: \bbSE_{u,k}^{(1)} \times
\bbSE_{l,k}^{(1)}(\epsilon) \rightarrow \bbSE_{u,k}^{(2)},
      \end{equation}
      such that
      \begin{enumerate}[label=(\alph*)]

            \item 
            For any $\mu_u \in \bbSE_{u,k}^{(1)}$ and $\mu_l \in \bbSE_{l,k}^{(1)}(\epsilon)$,
            we have
            \begin{equation}
                  \AuxObj_k^{(1)}(\mu_u,\mu_l)
                  \geq
                  \AuxObj_k^{(2)}\big(\Reduce_{u,k}^{(2)}(\mu_u),\Reduce_{l,k}^{(2)}(\mu_u,\mu_l)\big)
                  -
                  \delta(k),
            \end{equation}
            for some $\delta(k) \rightarrow 0$ as $k \rightarrow \infty$.

            \item 
            For any fixed $\mu_l \in \bbSE_{l,k}^{(1)}(\epsilon)$,
            the map $\mu_u \mapsto \Reduce_{u,k}^{(2)}(\mu_u,\mu_l)$ maps $\bbSE_{u,k}^{(1)}$ surjectively onto $\bbSE_{u,k}^{(2)}$.

      \end{enumerate}
      To see why this suffices,
      observe that we then have that for any fixed $\mu_l \in \bbSE_{l,k}^{(1)}(\epsilon)$,
      \begin{equation}
      \begin{aligned}
            \max_{\mu_u \in \bbSE_{u,k}^{(1)}}\;
                  \AuxObj_k^{(1)}(\mu_u,\mu_l)
                  &\geq
                  \max_{\mu_u \in \bbSE_{u,k}^{(1)}}
                  \AuxObj_k^{(2)}\big(\Reduce_{u,k}^{(2)}(\mu_u,\mu_l),\Reduce_{l,k}^{(2)}(\mu_l)\big)
                  -
                  \delta(k)
            \\
                  &=
                  \max_{\mu_u \in \bbSE_{u,k}^{(2)}}
                  \AuxObj_k^{(2)}\big(\mu_u,\Reduce_{l,k}^{(2)}(\mu_l)\big)
                  -
                  \delta(k),
      \end{aligned}
      \end{equation}
      where the inequality holds by item \textit{(a)},
      and the equality holds by item \textit{(b)}.
      Then, taking the minimum of both sides of the previous display over $\mu_l \in \bbSE_{l,k}^{(1)}(\epsilon)$,
      we have
      \begin{equation}
      \begin{aligned}
            \min_{\mu_l \in \bbSE_{l,k}^{(1)}}\;\;
            \max_{\mu_u \in \bbSE_{u,k}^{(1)}}\;
                  \AuxObj_k^{(1)}(\mu_u,\mu_l)
                  &\geq
                  \min_{\mu_l \in \bbSE_{l,k}^{(1)}}\;\;
                  \max_{\mu_u \in \bbSE_{u,k}^{(2)}}
                  \AuxObj_k^{(2)}\big(\mu_u,\Reduce_{l,k}^{(2)}(\mu_l)\big)
                  -
                  \delta(k)
            \\
                  &\geq
                  \min_{\mu_l \in \bbSE_{l,k}^{(2)}}\;\;
                  \max_{\mu_u \in \bbSE_{u,k}^{(2)}}
                  \AuxObj_k^{(2)}\big(\mu_u,\mu_l\big)
                  -
                  \delta(k),
      \end{aligned}
      \end{equation}
      where the second inequality holds because $\Reduce_{l,k}^{(2)}(\mu_l) \in \bbSE_{l,k}^{(2)}(\epsilon)$ for all $\mu_l \in \bbSE_{l,k}^{(1)}(\epsilon)$.

      We now define the maps $\Reduce_{l,k}^{(2)}$ and $\Reduce_{u,k}^{(2)}$,
      and show they satsify \textit{(a)} and \textit{(b)}. With respect to the
inner-product on $\LL_2$,
      let $G_{k-1}^\perp,G_1^\perp,G_2^\perp,\ldots,G_{k-2}^\perp,G_k^\perp$ be the result of Gram-Schmidt orthonormalization applied to $G_{k-1},G_1,G_2,\ldots,G_{k-2},G_k$ (in that order),
and similarly let $R_{k-1}^\perp,R_1^\perp,R_2^\perp,\ldots,R_{k-2}^\perp,R_k^\perp$
and $H_k^\perp,H_1^\perp,H_2^\perp,\ldots,H_{k-1}^\perp$
and $V_k^\perp,V_1^\perp,V_2^\perp,\ldots,V_{k-1}^\perp$ be the results of
Gram-Schmidt orthonormalization applied to $R_{k-1},R_1,R_2,\ldots,R_{k-2},R_k$
and $H_k,H_1,H_2,\ldots,H_{k-1}$ and $V_k,V_1,V_2,\ldots,V_{k-1}$ (in those orders). Then $G_1^\perp,G_2^\perp,\ldots,G_{k-2}^\perp,G_k^\perp$ are independent of $G_{k-1}$ and $H_1^\perp,H_2^\perp,\ldots,H_{k-1}^\perp$ are independent of $H_k$,
      and they are standard normal.
Observe that by definition, $\bK_{h,k}=\Cov(H_1^k)$ and
$\bK_{g,k}=\Cov(G_1^k)$, while by 
(\ref{eq:AMPVRGH}) also $\bK_{h,k}=\Cov(V_1^k)$ and $\bK_{g,k}=\delta^{-1}
\Cov(R_1^k)$. Then 
\begin{align}
\begin{pmatrix} \<H_1,U\>_{\LL_2} & \cdots & \<H_k,U\>_{\LL_2} \end{pmatrix}
\bK_{h,k}^{-1}\begin{pmatrix} \<V_1,L_1\>_{\LL_2} \\ \vdots \\
\<V_k,L_1\>_{\LL_2} \end{pmatrix}
&=\sum_{j=1}^k
\<H_j^\perp,U\>_{\LL_2}\<V_j^\perp,L_1\>_{\LL_2},\label{eq:HVGramSchmidt}\\
\begin{pmatrix} \<G_1,L_1\>_{\LL_2} & \cdots & \<G_k,L_1\>_{\LL_2} \end{pmatrix}
\bK_{g,k}^{-1}\begin{pmatrix} \<R_1,U\>_{\LL_2} \\ \vdots \\
\<R_k,U\>_{\LL_2} \end{pmatrix}
&=\sqrt{\delta}\sum_{j=1}^k
\<R_j^\perp,U\>_{\LL_2}\<G_j^\perp,L_1\>_{\LL_2}.\label{eq:GRGramSchmidt}
\end{align}
We record here that for each $k'<k$,
we can write $H_{k'} = \sum_{\ell = 1}^k c_{h,k'\ell} H_\ell^\perp$ for some
coefficients $c_{h,k'\ell}$, which will be used at the end of the proof.

      We define
      \[
      \begin{aligned}
            \Reduce_{l,k}^{(2)}(\mu_l)
                  &=
                  \Law\big(
                        L_1,L_2,M,S,G_{k-1},V_1,V_k,
                        \Gtilde
                  \big),
            \\
            \Reduce_{u,k}^{(2)}(\mu_u,\mu_l)
                  &=
                  \Law\big(
                        U,H_k,R_k,
                        \Htilde
                  \big),
      \end{aligned}
      \]
      where
\begin{equation}\label{eq:GHtilde}
            \Gtilde
                  =
                  \frac{
                        \sum_{j\neq k-1}G_j^\perp \< G_j^\perp , L_1 \>_{\LL_2}
                        + 
                        G\< G , L_1 \>_{\LL_2}
                  }{
                        \sqrt{\sum_{j\neq k-1}\< G_j^\perp , L_1 \>_{\LL_2}^2 
                        +
                        \< G , L_1 \>_{\LL_2}^2}
                  },
            \quad
            \Htilde
                  :=
                  \frac{
                        \sum_{j=1}^{k-1} H_j^\perp \< V_j^\perp , L_1 \>_{\LL_2} 
                        +
                        H \| \proj_{V_1^k}^\perp L_1 \|_{\LL_2}
                  }{
                        \| \proj_{V_k}^\perp L_1 \|_{\LL_2}
                  }.
\end{equation}
If $\mu_l$ is such that the denominator in the first and/or second expression of
(\ref{eq:GHtilde}) is
0, we may set $\Gtilde = G$ and/or $\Htilde = H$ respectively.
      By inspection, $\Reduce_{l,k}^{(2)}(\mu_l)$ depends only on $\mu_l$, and
$\Reduce_{u,k}^{(2)}(\mu_u,\mu_l)$ depends on both $\mu_u$ and $\mu_l$. 

      First,
      we establish that $\Reduce_{l,k}^{(2)}(\mu_l) \in
\bbSE_{l,k}^{(2)}(\epsilon)$ and $\Reduce_{u,k}^{(2)}(\mu_u,\mu_l) \in \bbSE_{u,k}^{(2)}$.
      Since $(G_j^\perp)_{j \neq k-1}$ and $G$ are standard normal variables
      independent of each other and of $(G_{k-1},V_1)$,
      we have that $\Gtilde \sim \normal(0,1)$ and is independent of $(G_{k-1},V_1)$. 
The remaining conditions of $\bbSE_{l,k}^{(2)}(\epsilon)$ evidently hold for
$\Reduce_{l,k}^{(2)}(\mu_l)$, by the corresponding conditions of
$\bbSE_{l,k}^{(1)}(\epsilon)$ for $\mu_l$.
      Thus, $\Reduce_{l,k}^{(2)}(\mu_l) \in \bbSE_{l,k}^{(2)}(\epsilon)$.
      Likewise,
      since $(H_j^\perp)_{j<k}$ and $H$ are standard normal variables
independent of each other and of $(H_k,E)$ (where $E=H_k - R_k$),
      and $\| \proj_{V^k}^\perp L_1 \|_{\LL_2}^2 = \sum_{j=1}^{k-1} \< V_j^\perp
, L_1 \>_{\LL_2}^2 + \| \proj_{V_1^k}^\perp L_1 \|_{\LL_2}^2$ by definition of
the above Gram-Schmidt procedure,
      we have that $\Htilde \sim \normal(0,1)$ and is independent of $(H_k,E)$.
The remaining conditions of $\bbSE_{u,k}^{(2)}$ evidently hold for
$\Reduce_{u,k}^{(2)}(\mu_u,\mu_l)$, by the corresponding conditions of
$\bbSE_{u,k}^{(1)}$ for $\mu_u$, so
      $\Reduce_{u,k}^{(2)}(\mu_u,\mu_l) \in \bbSE_{u,k}^{(2)}$.

      Now we establish item \textit{(a)}.
      The terms $\Theta_\TAP^{(1)}$ and
$\Theta_\TAP^{(2)}$ coincide by definition and depend only on
$\Law(L_1,L_2,M,S)$ and $\Law(U)$.
      Thus, 
    \begin{equation}\label{eq:Theta1Theta2}
            \Theta_\TAP^{(1)}(\mu_u,\mu_l)
                  =
                  \Theta_\TAP^{(2)}
                  \big(
                        \Reduce_{u,k}^{(2)}(\mu_u),
                        \Reduce_{l,k}^{(2)}(\mu_u,\mu_l)
                  \big).
\end{equation}
Applying (\ref{eq:HVGramSchmidt}) and the definitions of $H^{(1)}$ and
$\Htilde$, we also have under $\mu_u,\mu_l$ that
\begin{equation}\label{eq:H1H2}
      \begin{aligned}
            &\< H^{(1)}(L_1,H_1^k,V_1^k,H) , U \>_{\LL_2}
                  =
                  \sum_{j=1}^k \<H_j^\perp,U\>_{\LL_2} \< V_j^\perp ,
L_1 \>_{\LL_2}
                  +
                  \| \proj_{V_1^k}^\perp L_1 \|_{\LL_2}
                  \< H,U \>_{\LL_2}
            \\
                  &\hspace{0.2in}=
                  \frac1{\delta\gamma_k^{-1}-\sigma^2}
                  \<H_k , U \>_{\LL_2} \<V_k,L_1\>_{\LL_2}
                  +
                  \| \proj_{V_k}^\perp L_1 \|_{\LL_2} \< \Htilde, U \>_{\LL_2}
                  =
                  \< H^{(2)}(L_1,H_k,V_k,\Htilde),U\>_{\LL_2},
      \end{aligned}
\end{equation}
      where the last expression is equivalent when evaluated under the laws $\Reduce_{u,k}^{(2)}(\mu_u)$ and $\Reduce_{l,k}^{(2)}(\mu_u,\mu_l)$.
(If $\mu_l$ is such that $\|P_{V_k}^\perp L_1\|_{\LL_2}=0$ so $\Htilde=H$,
then $L_1$ is in the span of $V_k$ and this identity holds also.)
Finally, applying (\ref{eq:GRGramSchmidt}) and the definitions of $G^{(1)}$
and $\Gtilde$,
\begin{equation}\label{eq:G1G2}
      \begin{aligned}
            &\<G^{(1)}(U,G_1^k,R_1^k,G),L_1\>_{\LL_2}
                  =
                  \sum_{j=1}^k \<R_j^\perp,U\>_{\LL_2}\<G_j^\perp,L_1\>_{\LL_2}
                  +
                  \|\proj_{R_1^k}^\perp U \|_{\LL_2} \< G,L_1\>_{\LL_2}
            \\
            &\qquad\leq
                  \frac{\gamma_{k-1}}{\sqrt{\delta}}
                  \<R_{k-1},U\>_{\LL_2}\<G_{k-1},L_1\>_{\LL_2}
                  +
                  \sqrt{\sum_{j\neq k-1}\< R_j^\perp , U \>_{\LL_2}^2
                  +
                  \| \proj_{R_1^k}^\perp U \|_{\LL_2}^2}
                  \sqrt{\sum_{j\neq k-1}\< G_j^\perp , L_1 \>_{\LL_2}^2 
                  +
                  \< G , L_1 \>_{\LL_2}^2}
            \\
            &\qquad=
                  \frac{\gamma_{k-1}}{\sqrt{\delta}}
                  \<R_{k-1},U\>_{\LL_2}\<G_{k-1},L_1\>_{\LL_2}
                  +
                  \| \proj_{R_{k-1}}^\perp U \|_{\LL_2}
                  \< \Gtilde , L_1 \>_{\LL_2}
            \\
                  &\qquad=
                  \< G^{(2)}(U,G_{k-1},R_k,\Gtilde) , L_1 \>_{\LL_2}
                  +
                  \frac{\gamma_{k-1}}{\sqrt{\delta}}\<G_{k-1},L_1\>_{\LL_2}\<
R_{k-1} - R_k , U \>_{\LL_2}
+(\|\proj_{R_{k-1}}^\perp U\|_{\LL_2}-\|\proj_{R_k}^\perp U\|_{\LL_2})
\<\Gtilde,L_1\>_{\LL_2}.
      \end{aligned}
\end{equation}
(If $\mu_l$ is such that $\<G,L_1\>_{\LL_2}=0$ and $\<G_j^\perp,L_1\>_{\LL_2}=0$
for all $j\neq k-1$ so $\Gtilde=G$, then also $\<\Gtilde,L_1\>_{\LL_2}=0$ so
this statement holds with equality.)
      Note that
\[\big| \gamma_{k-1} \<G_{k-1},L_1\>_{\LL_2}\< R_k - R_{k-1} , U
\>_{\LL_2}\big| \leq C_0\gamma_{k-1} \| G_{k-1} \|_{\LL_2} \| R_k - R_{k-1}
\|_{\LL_2} \rightarrow 0\]
as $k \rightarrow \infty$, because by Proposition \ref{prop:amp-se},
      $\gamma_k \rightarrow \gamma_\alg$,
      $\| G_{k-1} \|_{\LL_2}^2 = \gamma_{k-1}^{-1} \rightarrow
\gamma_\alg^{-1}$, and $\| R_k - R_{k-1} \|_{\LL_2}^2 = \delta(\gamma_{k-1}^{-1} - \gamma_k^{-1}) \rightarrow 0 $ as $k \rightarrow \infty$.
Similarly,
\begin{align*}
\big|(\|\proj_{R_{k-1}}^\perp U\|_{\LL_2}-\|\proj_{R_k}^\perp U\|_{\LL_2})
\<\Gtilde,L_1\>_{\LL_2}\big|
&\leq \|\Gtilde\|_{\LL_2}\|\proj_{R_{k-1}} U-\proj_{R_k}U\|_{\LL_2}\\
&\leq C_0\|\Gtilde\|_{\LL_2}\|R_{k-1}-R_k\|_{\LL_2}(\|R_{k-1}\|_{\LL_2}^{-1}
+\|R_k\|_{\LL_2}^{-1}) \to 0
\end{align*}
as $k \to \infty$. 
      Combining these bounds with (\ref{eq:Theta1Theta2}), (\ref{eq:H1H2}),
and (\ref{eq:G1G2}) gives item \textit{(a)}.

      Lastly, we establish item \textit{(b)}.
      Fix $\mu_l \in \bbSE_{l,k}^{(1)}(\epsilon)$,
      and consider any
$\mu_u^{(2)}=\Law(U,H_k,R_k,\Htilde) \in \bbSE_{u,k}^{(2)}$.
Let $E=H_k-R_k$, so that by definition of $\bbSE_{u,k}^{(2)}$ we have
$E \indep \Htilde \indep H_k$. Now consider a multivariate normal distribution
$\nu=\Law(H_k,E,\Htilde,H,H_1^\perp,\ldots,H_{k-1}^\perp)$ where
$(H_k,E)$ have the same law as in $\mu_u^{(2)}$,
$H,H_1^\perp,\ldots,H_{k-1}^\perp \sim \normal(0,1)$
are independent of each other and of $(H_k,E)$, and $\Htilde$ is defined by
the second expression of (\ref{eq:GHtilde}), or by $\Htilde=H$ if the
denominator of this expression is 0. Then, by the gluing lemma, there exists a
joint law of $(U,H_k,E,R_k,\Htilde,H,H_1^\perp,\ldots,H_{k-1}^\perp)$
where $E=H_k-R_k$, $(U,H_k,R_k,\Htilde) \sim \mu_u^{(2)}$,
and $(H_k,E,\Htilde,H,H_1^\perp,\ldots,H_{k-1}^\perp) \sim \nu$.
From this joint law, define
      $H_k^\perp = H_k/\sqrt{\delta \gamma_k^{-1} - \sigma^2}$, and for
      $k' < k$
      define $H_{k'} = \sum_{\ell = 1}^k c_{h,k'\ell} H_\ell^\perp$ and
$R_{k'}=E-H_{k'}$,
      where $c_{h,k'\ell}$ are the coefficients in the preceding Gram-Schmidt
procedure. Finally, let $\mu_u=\Law(U,H_1^k,R_1^k,H)$. By this construction,
$H_1^k$ is a function of $H_1^\perp,\ldots,H_{k-1}^\perp,H_k$,
hence $H_1^k \indep E \indep H$. Furthermore, the joint law of $H_1^k$ is
precisely $\normal(0,\bK_{h,k})$ by reversal of the Gram-Schmidt procedure.
Thus $\mu_u \in \bbSE_{u,k}^{(1)}$, and by construction
$\Reduce_{u,k}^{(2)}(\mu_u,\mu_l)=\mu_u^{(2)}$.
      Thus, $\mu_u \mapsto \Reduce_{u,k}^{(2)}(\mu_u,\mu_l)$ is surjective.
\end{proof}

\subsubsection{Reduction to analysis at fixed point}

Define the distributions
\begin{equation}
\begin{aligned}
      \SE_{u,\infty}
            &:=
            \Law(H_\infty,R_\infty,H)
            \;\; 
            &\text{where}
            \;\;
            & H_\infty \sim \normal(0,\delta\gamma_\alg^{-1}-\sigma^2)
            \indep
            E \sim \normal(0,\sigma^2)
            \indep
            H \sim \normal(0,1),
      \\
            &&&R_\infty = H_\infty - E,
      \\
      \SE_{l,\infty}
            &:=
            \Law(G_\infty,V_1,V_\infty,G)
            \;\;
            &\text{where}
            \;\;
            &G_\infty \sim \normal(0,\gamma_\alg^{-1})
            \indep
            -V_1 \sim \sP_0
            \indep
            G \sim \normal(0,1),
      \\
            &&&V_\infty = \denoiser(-V_1 + G_\infty,\,\gamma_\alg) + V_1.
\end{aligned}
\end{equation}
Define the deterministic subsets of Wasserstein space
\[
\begin{aligned}
      \bbSE_u^{(3)}
      :=
      \Big\{
            \Law(U,H_\infty,R_\infty,H)
            \Bigm|
            &\| U \|_{\LL_2} \leq C_0,\,
            (H_\infty,R_\infty,H) \sim \SE_{u,\infty}
      \Big\},
      \\
      \bbSE_l^{(3)}(\epsilon)
      :=
      \Big\{
            \Law(L_1,L_2,M,S,G_\infty,V_1,V_\infty,G)
            \Bigm|
            &\| L_1 \|_{\LL_2}^2 + \| L_2 \|_{\LL_2}^2 = 1,\,
            \| M - V_\infty + V_1 \|_{\LL_2} \leq \epsilon,\,
            \\
            &\| S - \sS(-V_1 + G_\infty,\,\gamma_\infty) \|_{\LL_2} \leq \epsilon,\,
            \\
            &(G_\infty,V_1,V_\infty,G) \sim \SE_{l,\infty}
            \Big\}.
\end{aligned}
\]
For $\mu_u \in \bbSE_u^{(3)}$ and $\mu_l \in \bbSE_l^{(3)}(\epsilon)$,
define
\begin{align*}
      G^{(3)}(U,G_\infty,R_\infty,G)
            &=
            \frac{\gamma_\alg}{\sqrt{\delta}} G_\infty \< R_\infty , U \>_{\LL_2}
            +
            \|\proj_{R_\infty}^\perp U \|_{\LL_2}G,
      \\
      H^{(3)}(L_1,H_\infty,V_\infty,H)
            &=
            \frac1{\delta \gamma_\alg^{-1} - \sigma^2} H_\infty\< V_\infty, L_1 \>_{L_2} + \| \proj_{V_\infty}^\perp L_1 \|_{L_2} H
\end{align*}
and 
\[
      \AuxObj_\infty^{(3)}(\mu_u,\mu_l)
            =
            -\sqrt{\delta}\<G^{(3)}(U,G_\infty,R_\infty,G),L_1\>_{\LL_2}
            +\delta\< H^{(3)}(L_1,H_\infty,V_\infty,H) , U \>_{\LL_2}
            +\Theta_\TAP^{(3)}(\mu_u,\mu_l),
\]
where
\[
\begin{aligned}
      \Theta_\TAP^{(3)}(\mu_u,\mu_l)
            &:=
            \Theta_\bulk^{(3)}(\mu_u,\mu_l)
            +
            \Theta_\spike^{(3)}(\mu_u,\mu_l),
      \\
      \Theta_\bulk^{(3)}(\mu_u,\mu_l)
            &:=
            -\frac{\sigma^2\delta}2 \| U \|_{\LL_2}^2
            -\frac{\delta}{2V_\star}L_1^2
            +\frac{1}{2}(L_1,L_2)^\top \nabla^2 [{-}\hs(M,S)] (L_1,L_2),
      \\
      \Theta_\spike^{(3)}(\mu_u,\mu_l)
            &:=
            -\frac{\delta
            \big(
                  \<M_\infty,L_1\>_{\LL_2}
                  -
                  \<1,L_2\>_{\LL_2}/2
            \big)^2}{V_\star^2},\\
M_\infty &:= V_\infty - V_1, \qquad S_\infty := \sS(-V_1 +
G_\infty,\,\gamma_\infty), \qquad V_\star=V(M_\infty,S_\infty).
\end{aligned}
\]
Here, we have substituted
$(M_\infty,S_\infty)$ for $(M,S)$ in all terms of $\Theta_\TAP^{(2)}$ except
for $\nabla^2[-\hs(M,S)]$. We will make this final substitution
at a later step of the proof.

We have the following reduction.
\begin{lemma}[Reduction to analysis at fixed point]
\label{lem:wass-reduce-to-fixed-pt}
      We have
      \begin{equation}
      \label{eq:Aux3-lower-bound}
            \liminf_{\epsilon \rightarrow 0}\;\;
            \liminf_{k \rightarrow \infty}\;\;
            \min_{\mu_l \in \bbSE_{l,k}^{(2)}(\epsilon)}\;\;
            \max_{\mu_u \in \bbSE_{u,k}^{(2)}}\;
            \AuxObj_k^{(2)}(\mu_u,\mu_l)
            =
            \liminf_{\epsilon \rightarrow 0}\;\;
            \min_{\mu_l \in \bbSE_l^{(3)}(\epsilon)}\;\;
            \max_{\mu_u \in \bbSE_u^{(3)}}\;
            \AuxObj_\infty^{(3)}(\mu_u,\mu_l).
      \end{equation}
\end{lemma}

\begin{proof}
      Consider $\mu_u \in \bbSE_{u,k}^{(2)}$ and $\mu_l \in \bbSE_{l,k}^{(2)}(\epsilon)$.
      For random variables distributed according to these laws,
      define
\[
      \begin{gathered}
            H_\infty
                  :=
                  \sqrt{\frac{\delta \gamma_\alg^{-1}-\sigma^2}{\delta \gamma_k^{-1} - \sigma^2}}\, H_k,
            \qquad
            R_\infty
                  :=
                  H_\infty - E,
            \\
            G_\infty
                  :=
                  \sqrt{\frac{\gamma_{k-1}}{\gamma_\alg}}\, G_{k-1},
            \qquad
            V_\infty
                  :=
                  \denoiser(-V_1 + G_\infty,\,\gamma_\alg) + V_1.
      \end{gathered}
\]
      Define
\[
            \Reduce_{u,k}^{(3)}(\mu_u)
                  :=
                  \Law(U,H_\infty,R_\infty,H),
            \qquad
            \Reduce_{l,k}^{(3)}(\mu_l)
                  :=
                  \Law(L_1,L_2,M,S,G_\infty,V_1,V_\infty,G).
\]
      By construction, $\Reduce_{u,k}^{(3)}(\mu_u) \in \bbSE_u^{(3)}$
      and $\Reduce_{u,k}^{(3)}$ is surjective.
      Because $\gamma_k \rightarrow \gamma_\alg$ and
$\denoiser(\,\cdot\,,\,\gamma)$ and $\sS(\,\cdot\,,\,\gamma)$ are uniformly
Lipschitz continuous for $\gamma$ in compact sets (by boundedness of the support
of $\sP_0$),
\[
      \begin{gathered}
            \| H_k - H_\infty \|_{\LL_2} \leq \delta(k),
                  \qquad
                  \| R_k - R_\infty \|_{\LL_2} \leq \delta(k),
            \\
            \| G_k - G_\infty \|_{\LL_2} \leq \delta(k),
                  \qquad
                  \| V_k - V_\infty \|_{\LL_2} \leq \delta(k),
            \\
            \| M - M_\infty \|_{L_2} \leq \delta(k,\epsilon),
                  \qquad 
                  \| S - S_\infty \|_{L_2} \leq \delta(k,\epsilon),
      \end{gathered}
\]
      where $\delta(k)$ represents a quantity which satisfies $\delta(k) \rightarrow 0$ as $k \rightarrow \infty$,
      and $\delta(k,\epsilon)$ represents a quantity which satisfies
$\delta(k,\epsilon) \rightarrow 0$ as $k \rightarrow \infty$ and $\epsilon \rightarrow 0$ (in any order).
      Thus, $\Reduce_{l,k}^{(3)}(\mu_l) \in \bbSE_l^{(3)}(\delta(k,\epsilon))$.
      Note that $G^{(2)},H^{(2)}$ and $G^{(3)},H^{(3)}$ are continuous in their
arguments at $(U,G_\infty,R_\infty,G)$ and $(L_1,H_\infty,V_\infty,H)$, because
$\| R_\infty \|_{L_2} > 0$ and $\| V_\infty \|_{L_2} > 0$ (where these
conditions ensure continuity of the terms $\|\proj_{R_\infty}^\perp U \|_{L_2}$ and $\| \proj_{V_\infty}^\perp L_1 \|_{L_2}$).
      Also, $E(M,S)$ and $\< M , L_1 \>_{\LL_2}$ are continuous in $M,S$.
      Thus,
\[
      \begin{aligned}
            \min_{\mu_l \in \bbSE_{l,k}^{(2)}(\epsilon)}
            \max_{\mu_u \in \bbSE_{u,k}^{(2)}}
            \AuxObj_k^{(2)}(\mu_u,\mu_l)
                  &\geq
                  \min_{\mu_l \in \bbSE_{l,k}^{(2)}(\epsilon)}
                  \max_{\mu_u \in \bbSE_{u,k}^{(2)}}
                  \AuxObj_\infty^{(3)}(
                        \Reduce_{u,k}^{(3)}(\mu_u),
                        \Reduce_{l,k}^{(3)}(\mu_l)
                  ) - \delta(k,\epsilon)
            \\
                  &\geq
                  \min_{\mu_l \in \bbSE_l^{(3)}(\delta(k,\epsilon))}
                  \max_{\mu_u \in \bbSE_u^{(3)}}
                  \AuxObj_\infty^{(3)}(
                        \mu_u,
                        \mu_l
                  ) - \delta(k,\epsilon),
      \end{aligned}
\]
      where the second inequality holds by surjectivity of $\Reduce_{u,k}^{(3)}$ and because for $\mu_l \in \bbSE_{l,k}^{(2)}(\epsilon)$ we have $\Reduce_{l,k}^{(3)}(\mu_l) \in \bbSE_l^{(3)}(\delta(k,\epsilon))$.
      The result follows by taking $k \rightarrow \infty$ followed by $\epsilon \rightarrow 0$.
\end{proof}

\subsubsection{Analysis of conditional lower bound}
\label{sec:cond-lower-bound}

Note that the min-max problem $\min_{\mu_l \in \bbSE_l^{(3)}(\epsilon)} \max_{\mu_u \in \bbSE_u^{(3)}} \AuxObj_\infty^{(3)}(\mu_u,\mu_l)$ can be equivalently written as the min-max problem on $\LL_2$
\begin{equation}\label{eq:L2minmax}
\begin{aligned}
      &
      \min_{
      	\substack{
      		\| M - V_\infty + V_1 \|_{\LL_2} \leq \epsilon
      		\\
      		\| S - \sS(-V_1 + G_\infty,\,\gamma_\infty) \|_{\LL_2} \leq \epsilon
      	}
      }
      \min_{
  		\| L_1 \|_{\LL_2}^2 + \| L_2 \|_{\LL_2}^2 = 1
      }\;
      \max_{\| U \|_{\LL_2} \leq C_0}
            -\gamma_\alg \< R_\infty , U \>_{\LL_2} \< G_\infty , L_1 \>_{\LL_2}
            -
            \sqrt{\delta} \| \proj_{R_\infty}^\perp U \|_{\LL_2} \< G , L_1
\>_{\LL_2}
      \\
            &\qquad\qquad\qquad
            +
            \frac{\delta}{\delta \gamma_\alg^{-1}-\sigma^2}
            \< H_\infty , U \>_{\LL_2} \< V_\infty , L_1 \>_{\LL_2}
            +
            \delta \| \proj_{V_\infty}^\perp L_1 \|_{\LL_2} \< H , U \>_{\LL_2}
            -\frac{\sigma^2\delta}2 \| U \|_{\LL_2}^2
            +
            \sI,
\end{aligned}
\end{equation}
where
\[
    \sI =
        {-}\frac{\delta}{2V_\star}L_1^2
        +\frac{1}{2}(L_1,L_2)\nabla^2 [{-}\hs(M,S)](L_1,L_2)
            -\frac{\delta
            \big(
                  \<M_\infty,L_1\>_{\LL_2}
                  -
                  \<1,L_2\>_{\LL_2}/2
            \big)^2}{V_\star^2}.
\]
Here, the max is over couplings of $U$ to
$(H_\infty,R_\infty,E,H) \sim \SE_{u,\infty}$, and the
min is over couplings of $(L_1,L_2,M,S)$ to
$(G_\infty,V_1,V_\infty,G) \sim \SE_{l,\infty}$. We now
proceed to reduce this to a min-max optimization over scalars in $\R$.

Introduce the change of variables $\sqrt{\delta} U = \alpha_\| R_\infty / \|
R_\infty \|_{\LL_2} + \alpha_\perp Z_\perp$,
where $\< Z_\perp , R_\infty \>_{\LL_2} = 0$, $\| Z_\perp \|_{\LL_2} = 1$,
and $\alpha_\perp=\|\proj_{R_\infty}^\perp \sqrt{\delta} U\|_{\LL_2} \geq 0$.
Then, letting $\nu_0^2 = \delta \gamma_\alg^{-1} - \sigma^2$,
applying $\|R_\infty\|_{\LL_2}^2=\delta\gamma_{\alg}^{-1}$,
and (for now) fixing $M,S$,
the optimization in (\ref{eq:L2minmax}) may be rewritten as
\begin{align*}   
      &\min_{
      		\| L_1 \|_{\LL_2}^2 + \| L_2 \|_{\LL_2}^2 = 1
      }\;
      \max_{\substack{\alpha_\perp^2+\alpha_\|^2 \leq C_0^2/\delta \\ \alpha_\perp \geq 0}}
      \max_{\substack{\| Z_\perp \|_{\LL_2} =1\\\<Z_\perp,R_\infty\>_{\LL_2}=0}}
            -\sqrt{\gamma_\alg} \,\alpha_\| \< G_\infty , L_1 \>_{\LL_2}
            -
            \alpha_\perp \< G , L_1 \>_{\LL_2}
      \\
            &\qquad\qquad
            +
            \sqrt{\delta}
            \Big\<
                  \frac{\<V_\infty,L_1\>_{\LL_2}}{\nu_0^2} H_\infty
                  +
                  \|\proj_{V_\infty}^\perp L_1 \|_{\LL_2}H,\,
                  \alpha_\| \frac{R_\infty}{\|R_\infty\|_{\LL_2}}
                  +
                  \alpha_\perp Z_\perp
            \Big\>_{\LL_2}
            -\frac{\sigma^2}2 (\alpha_\perp^2 + \alpha_\|^2)
            + \sI.
\end{align*}
Applying independence of $(R_\infty,H_\infty)$ with $H$,
and also evaluating explicitly the maximum over $Z_\perp$, this is equal to
\begin{align*}
      \\
      &\min_{
      		\| L_1 \|_{\LL_2}^2 + \| L_2 \|_{\LL_2}^2 = 1
      }\;
      \max_{\substack{\alpha_\perp^2+\alpha_\|^2 \leq C_0^2/\delta \\ \alpha_\perp \geq 0}}
            -\sqrt{\gamma_\alg} \,\alpha_\| \< G_\infty , L_1 \>_{\LL_2}
            -
            \alpha_\perp \< G , L_1 \>_{\LL_2}
            +
            \sqrt{\delta}
            \Big\<
                  \frac{\<V_\infty,L_1\>_{\LL_2}}{\nu_0^2} H_\infty,\,
                  \alpha_\| \frac{R_\infty}{\|R_\infty\|_{\LL_2}}
            \Big\>_{\LL_2}
      \\
            &\qquad\qquad
            +
            \sqrt{\delta}\alpha_\perp
            \left\|
                  \frac{\<V_\infty,L_1\>_{\LL_2}}{\nu_0^2} \proj_{R_\infty}^\perp H_\infty
                  +
                  \|\proj_{V_\infty}^\perp L_1 \|_{\LL_2}\proj_{R_\infty}^\perp H
            \right\|_{\LL_2}
            -
            \frac{\sigma^2}2 (\alpha_\perp^2 + \alpha_\|^2)
            +
            \sI
    	\\
        &=
      \min_{
      		\| L_1 \|_{\LL_2}^2 + \| L_2 \|_{\LL_2}^2 = 1
      }\;
      \max_{\substack{\alpha_\perp^2+\alpha_\|^2 \leq C_0^2/\delta \\ \alpha_\perp \geq 0}}
            -\sqrt{\gamma_\alg} \,\alpha_\| \< G_\infty , L_1 \>_{\LL_2}
            -
            \alpha_\perp \< G , L_1 \>_{\LL_2}
            +
            \sqrt{\delta}
            \Big\<
                  \frac{\<V_\infty,L_1\>_{\LL_2}}{\nu_0^2} H_\infty,\,
                  \alpha_\| \frac{R_\infty}{\|R_\infty\|_{\LL_2}}
            \Big\>_{\LL_2}
        \\
            &\qquad\qquad
            +
            \sqrt{\delta}\alpha_\perp
            \sqrt{
                  \frac{\<V_\infty,L_1\>_{\LL_2}^2}{\nu_0^4}
\|\proj_{R_\infty}^\perp H_\infty\|_{\LL_2}^2
                  +
                  \|\proj_{V_\infty}^\perp L_1 \|_{\LL_2}^2
            }
            -\frac{\sigma^2}2 (\alpha_\perp^2 + \alpha_\|^2)
            +\sI.
\end{align*}
Thus we have replaced the maximization over $U$ by a maximization over scalar
variables $\alpha_\parallel,\alpha_\perp$.

Next, we introduce the variables and Lagrange multipliers
\[
\begin{aligned}
\< V_\infty, L_1  \>_{\LL_2} &= l_v  &&(\text{with multiplier } \lambda_{l_v}/2), 
&\<M_\infty,L_1\>_{\LL_2}&-\<1,L_2\>_{\LL_2}/2 = l_m && (\text{with multiplier }
\lambda_{l_m}/2), \\
\| L_1 \|_{\LL_2}^2 &= \barl_1^2 && (\text{with multiplier } \alpha_1/2), 
&\| L_2 \|_{\LL_2}^2 &= (1 -\barl_1^2) && (\text{with multiplier } \alpha_2/2).
\end{aligned}
\]
Fix a constant $C(\sP_0)>\max(\|V_\infty\|_{\LL_2},\|M_\infty\|_{\LL_2}+1/2)$
so that $|l_v|,|l_m|<C(\sP_0)$. Then,
recalling $\|V_\infty\|_{\LL_2}^2=\nu_0^2$ and the definition of $\sI$,
for any compact subset $K \subset \reals^4$,
the previous display is lower bounded by
\begin{align*}
      &\min_{\substack{|l_v|,|l_m| \leq C(\sP_0) \\ \barl_1\in[0,1]}}
\min_{L_1,L_2 \in \LL_2}
\max_{(\lambda_{l_v},\lambda_{l_m},\alpha_1,\alpha_2) \in K}
      \max_{\substack{\alpha_\perp^2+\alpha_\|^2 \leq C_0^2/\delta \\ \alpha_\perp \geq 0}}
      -\sqrt{\gamma_\alg} \,\alpha_\| \< G_\infty , L_1 \>_{\LL_2}
            -
            \alpha_\perp \< G , L_1 \>_{\LL_2}
            +
            \frac{\sqrt{\delta} \alpha_\| l_v}{\nu_0^2 \| R_\infty \|_{\LL_2}}
            \< H_\infty , R_\infty \>_{\LL_2}
      \\
            &\;\;           
            +
            \sqrt{\delta}\alpha_\perp
            \sqrt{
                  \frac{l_v^2}{\nu_0^4} \|\proj_{R_\infty}^\perp
H_\infty\|_{\LL_2}^2
                  +
                  \Big(
                        \barl_1^2 - \frac{l_v^2}{\nu_0^2}
                  \Big)
            }
            -\frac{\sigma^2}2 (\alpha_\perp^2 + \alpha_\|^2)
            -\frac{\delta}{V_\star^2} l_m^2-\frac{\delta}{2V_\star}L_1^2
+\frac{1}{2}(L_1,L_2)[{-}\nabla^2 \hs(M,S)](L_1,L_2)
      \\
            &\;\;
            +
            \frac{\lambda_{l_v}}{2}(\< V_\infty,L_1\>_{\LL_2} - l_v)
            +
            \frac{\lambda_{l_m}}{2}(\<M,L_1\>_{\LL_2} - \< 1 , L_2 \>_{\LL_2}/2 - l_m)
            +
            \frac{\alpha_1}{2}(\|L_1\|_{\LL_2}^2 - \barl_1^2)
      		+
            \frac{\alpha_2}{2}(\|L_2\|_{\LL_2}^2 - (1 - \barl_1^2)).
\end{align*}
As a further lower bound, we may exchange $\min_{L_1,L_2\in\LL_2}$ with the
following two maxes, and then explicitly minimize the quadratic function of
$(L_1,L_2)$, to obtain
\begin{align*}
      &\min_{\substack{|l_v|,|l_m| \leq C(\sP_0) \\ \barl_1\in[0,1]}}\;
      \max_{(\lambda_{l_v},\lambda_{l_m},\alpha_1,\alpha_2)\in K}\;
      \max_{\substack{\alpha_\perp^2+\alpha_\|^2 \leq C_0^2/\delta \\ \alpha_\perp \geq 0}}
            \frac{\sqrt{\delta} \alpha_\| l_v}{\nu_0^2 \| R_\infty \|_{\LL_2}}
            \< H_\infty , R_\infty \>_{\LL_2}
            +
            \sqrt{\delta}\alpha_\perp
            \sqrt{
                  \frac{l_v^2}{\nu_0^4} \|\proj_{R_\infty}^\perp
H_\infty\|_{\LL_2}^2
                  +
                  \Big(
                        \barl_1^2 - \frac{l_v^2}{\nu_0^2}
                  \Big)
            }
      \\
            &\;\;
            -\frac{\sigma^2}2 (\alpha_\perp^2 + \alpha_\|^2)
            -\frac{\delta}{V_\star^2} l_m^2
            -
            \frac{\lambda_{l_v}l_v}{2}
            -
            \frac{\lambda_{l_m}l_m}{2}
            -
            \frac{\alpha_1 \barl_1^2}{2}
            -
            \frac{\alpha_2(1 - \barl_1^2)}{2}
            -
            \frac\delta{2 V_\star} \barl_1^2
-\frac{1}{2}\E\left[\ba(M)^\top\bD(M,S)^{-1}\ba(M)\right],
\end{align*}
where
\begin{align*}
      \ba(m) = \begin{pmatrix}
			    \frac{\lambda_{l_v}V_\infty  + \lambda_{l_m} m}{2} - (\sqrt{ \gamma_{\alg}} \alpha_\| G_\infty + \alpha_\perp G)
			    \\
			    -\frac14 \lambda_{l_m}
			\end{pmatrix},
\quad \bD(m,s)=\begin{pmatrix} \alpha_1 & 0 \\ 0 & \alpha_2
\end{pmatrix}+\nabla^2[{-}\hs(m,s)].
\end{align*}
By (\ref{eq:hesshlowerbound}),
$\nabla^2[-\sh(m,s)] \succeq 2c(\sP_0) \id_2$ for an appropriately chosen
constant $c(\sP_0) > 0$.
Restricting to the range $\alpha_1,\alpha_2 \geq {-}c(\sP_0)$,
the operator norm of $\bD(m,s)^{-1}$ is then
bounded above by a constant uniformly over $m,s$.
Thus, $\bD(m,s)^{-1}$ is also uniformly Lipschitz in $m,s$.
Since $\ba(M)$ is also uniformly bounded in $\LL_2$,
we may then take the minimum of the previous display over $(M,S)$ satisfying
$\| M - V_\infty + V_1 \|_{\LL_2} \leq \epsilon$, $\| S - \sS(-V_1 +
G_\infty,\,\gamma_\infty)\|_{\LL_2} \leq \epsilon$ and combine the above bounds and Lipschitz properties to get
\[
\begin{aligned}
	&\liminf_{\epsilon \rightarrow 0}\;\;
        \min_{\mu_l \in \bbSE_l^{(3)}(\epsilon)}\;\;
        \max_{\mu_u \in \bbSE_u^{(3)}}\;
        \AuxObj_\infty^{(3)}(\mu_u,\mu_l)
    \\
      &\geq \min_{\substack{|l_v|,|l_m| \leq C(\sP_0) \\ \barl_1\in[0,1]}}\;
      \sup_{\substack{(\lambda_{l_v},\lambda_{l_m},\alpha_1,\alpha_2)\in K \\
\alpha_1,\alpha_2 \geq {-}c(\sP_0)}}\;
      \max_{\substack{\alpha_\perp^2+\alpha_\|^2 \leq C_0^2/\delta \\ \alpha_\perp \geq 0}}
            \frac{\sqrt{\delta} \alpha_\| l_v}{\nu_0^2 \| R_\infty \|_{\LL_2}}
            \< H_\infty , R_\infty \>_{\LL_2}
            +
            \sqrt{\delta}\alpha_\perp
            \sqrt{
                  \frac{l_v^2}{\nu_0^4} \|\proj_{R_\infty}^\perp
H_\infty\|_{\LL_2}^2
                  +
                  \Big(
                        \barl_1^2 - \frac{l_v^2}{\nu_0^2}
                  \Big)
            }
      \\
            &\;\;
            -\frac{\sigma^2}2 (\alpha_\perp^2 + \alpha_\|^2)
            -\frac{\delta}{V_\star^2} l_m^2
            -
            \frac{\lambda_{l_v}l_v}{2}
            -
            \frac{\lambda_{l_m}l_m}{2}
            -
            \frac{\alpha_1 \barl_1^2}{2}
            -
            \frac{\alpha_2(1 - \barl_1^2)}{2}
            -
            \frac{\delta}{2 V_\star} \barl_1^2
            -
      		\frac12 \E[\ba(M_\infty)^\top \bD(M_\infty,S_\infty)^{-1} \ba(M_\infty)].
\end{aligned}
\]
We have now replaced all optimizations over $\LL_2$ by optimizations of scalar
quantities.

Finally, let us further simplify this scalar optimization.
Note that $\langle H_\infty,R_\infty
\rangle_{\LL_2}=\|H_\infty\|_{\LL_2}^2=\nu_0^2$,
$\|R_\infty\|_{\LL_2}^2=\delta\gamma_\alg^{-1}=\sigma^2+\nu_0^2$,
and $\|\proj_{R_\infty}^\perp
H_\infty\|_{\LL_2}^2=\sigma^2\nu_0^2/(\sigma^2+\nu_0^2)$.
Make a change of variables $\bar l_1 \to l_p$ given by
\[
      \sqrt{
            \frac{l_v^2}{\nu_0^4} \|\proj_{R_\infty}^\perp H_\infty\|_{\LL_2}^2
            +
            \Big(
                  \barl_1^2 - \frac{l_v^2}{\nu_0^2}
            \Big)
      }
      =
      \sqrt{
            \barl_1^2 - \frac{\gamma_\alg}{\delta} l_v^2
      }
      =:
      l_p,
      \quad
      \text{where $\barl_1^2 = l_p^2 + \frac{\gamma_\alg}{\delta}l_v^2$,}
\]
Then the previous optimization objective is
\begin{align*}
&\sqrt{\gamma_\alg} \alpha_\| l_v
            +
            \sqrt{\delta}\alpha_\perp l_p
            -\frac{\delta}{V_\star^2} l_m^2
            -
            \frac{\lambda_{l_v}l_v}{2}
            -
            \frac{\lambda_{l_m}l_m}{2}
    		-
            \frac{\alpha_1}{2}\Big(l_p^2 + \frac{\gamma_\alg}{\delta} l_v^2\Big)
        \\
            &\;\;    
            -
            \frac{\sigma^2}2 (\alpha_\perp^2 + \alpha_\|^2)
            -
            \frac{\alpha_2}{2}\Big(1 - l_p^2 - \frac{\gamma_\alg}{\delta} l_v^2\Big)
            -
            \frac{\delta}{2V_\star}\Big(l_p^2 + \frac{\gamma_\alg}{\delta} l_v^2\Big)
      		-
            \frac12
            \E\left[
                  \ba(M_\infty)\bD(M_\infty,S_\infty)^{-1}\ba(M_\infty)
            \right].
\end{align*}
Apart from the linear term ${-}\alpha_2/2$, this is quadratic in all
optimization variables. For the last term $\E[\ba(M_\infty)^\top
\bD(M_\infty,S_\infty)^{-1}\ba(M_\infty)]$,
we may first evaluate the expectation over $G$, noting that it is
independent of $(G_\infty,V_\infty,M_\infty,S_\infty)$. Then, introducing the
shorthand
\begin{equation}
\begin{gathered}
\bz = (l_v,l_m,\lambda_{l_v},\lambda_{l_m},\alpha_\|)^\top, \quad
      \bA
            =
            \begin{pmatrix}
                  -\frac{\gamma_\alg}{\delta}(\alpha_1 - \alpha_2 +
\frac{\delta}{V_\star})
                        &
                  0
                        \\
                  0
                        &
                  -\frac{2\delta}{V_\star^2}
            \end{pmatrix}, \quad
      \bB
            =
            \begin{pmatrix}
                  -\frac12 & 0 & \sqrt{\gamma_{\alg}}
                  \\
                  0 & -\frac12 & 0
            \end{pmatrix},
      \\
D_\star=\bD(M_\infty,S_\infty), \quad
      \bC
            =
            -
            \E
            \left[
                  \begin{pmatrix}
                        \frac{V_\infty}{2} & 0
                        \\
                        \frac{M_\infty}{2} & - \frac14
                        \\
                        -\sqrt{\gamma_{\alg}}G_\infty & 0
                  \end{pmatrix}
                  D_\star^{-1}
                  \begin{pmatrix}
                        \frac{V_\infty}{2} & \frac{M_\infty}{2} & -\sqrt{\gamma_{\alg}}G_\infty
                        \\
                        0 & - \frac14 & 0
                  \end{pmatrix}
            \right]-\sigma^2\begin{bmatrix}0&0&0\\0&0&0\\0&0&1
            \end{bmatrix},
\end{gathered}\label{eq:spike_quadratic_matrices}
\end{equation}
this optimization objective may be re-written as
\[-\frac{\alpha_2}{2}+
		\frac12
		\begin{pmatrix}
			l_p & \alpha_\perp
		\end{pmatrix}
		\begin{pmatrix}
			-\big(\alpha_1 - \alpha_2 + \frac{\delta}{V_\star}\big) 
				&
				\sqrt{\delta}
				\\
				\sqrt{\delta}
				&
				-\sigma^2 - \E[D_\star^{-1}]_{11}
		\end{pmatrix}
		\begin{pmatrix}
			l_p \\ \alpha_\perp
		\end{pmatrix}
                +\frac12
	\bz^\top
	\begin{pmatrix}
	    \bA & \bB
	    \\
	    \bB^\top & \bC
	\end{pmatrix}
	\bz.
\]

%
Finally, take the compact domain $K \subset \R^4$ so that
\[K=K_\lambda \times K_\alpha \subset \R^2 \times \R^2,
\qquad \alpha_1,\alpha_2 \geq {-}c(\sP_0) \text{ for all }
(\alpha_1,\alpha_2) \in K_\alpha.\]
We may exchange the minimization over $(l_v,l_m,l_p)$ with the maximization
over $(\alpha_1,\alpha_2)$, expand the domain of $l_p$ to $[0,1]$ (noting that
$l_p \in [0,1]$ when $\bar l_1 \in [0,1]$), and
restrict the domain $\alpha_\perp^2 + \alpha_\|^2 \leq C_0^2/\delta$ further to
$\alpha_\perp^2 \leq C'$ and $\alpha_\|^2 \leq C'$ where $C'=C_0^2/(2\delta)$.
All three operations yield a further lower bound of the optimization problem,
so we arrive at
\begin{align}
	&\liminf_{\epsilon \rightarrow 0}\;\;
        \min_{\mu_l \in \bbSE_l^{(3)}(\epsilon)}\;\;
        \max_{\mu_u \in \bbSE_u^{(3)}}\;
        \AuxObj_\infty^{(3)}(\mu_u,\mu_l)\label{eq:finalAuxObj3lower}\\
	&\geq \sup_{(\alpha_1,\alpha_2) \in K_\alpha}
\min_{\substack{|l_v|,|l_m| \leq C(\sP_0) \\ l_p\in[0,1]}}\;
	\sup_{(\lambda_{l_v},\lambda_{l_m}) \in K_\lambda}\;
	\max_{\substack{\alpha_\|^2 \leq C'\\ \alpha_\perp^2 \leq C',\,\alpha_\perp \geq 0}}\;
	\bulk(\alpha_1,\alpha_2;l_p,\alpha_\perp) 
	+ 
	\spiked(\alpha_1,\alpha_2;l_v,l_m,\lambda_{l_v},\lambda_{l_m},\alpha_\|)
\nonumber
\end{align}
where we set
\begin{align*}
\bulk(\alpha_1,\alpha_2;l_p,\alpha_\perp)
		&:={-}\frac{\alpha_2}{2}
                +\frac12
		\begin{pmatrix}
			l_p & \alpha_\perp
		\end{pmatrix}
		\begin{pmatrix}
			-\big(\alpha_1 - \alpha_2 + \frac{\delta}{V_\star}\big) 
				&
				\sqrt{\delta}
				\\
				\sqrt{\delta}
				&
				-\sigma^2 - \E[D_\star^{-1}]_{11}
		\end{pmatrix}
		\begin{pmatrix}
			l_p \\ \alpha_\perp
		\end{pmatrix},\\
	\spiked(\alpha_1,\alpha_2;l_v,l_m,\lambda_{l_v},\lambda_{l_m},\alpha_\|)
	&:=
	\frac12
	\bz^\top
	\begin{pmatrix}
	    \bA & \bB
	    \\
	    \bB^\top & \bC
	\end{pmatrix}
	\bz.
\end{align*}
Observe that for fixed $(\alpha_1,\alpha_2) \in K_\alpha$, the resulting
min-max problem becomes separable in the remaining optimization variables,
with $l_p,\alpha_\perp$ participating in one part of the objective (``bulk''), and $l_v,l_m,\lambda_{l_v},\lambda_{l_m},\alpha_\|$ participating in another part of the objective (``spike''). 
Thus, we can analyze these min-max problems separately.

\subsubsection{The bulk min-max problem is positive}
\label{sec:bulk}

Here we show that, choosing $C'=C_0^2/2\delta$ sufficiently large,
for any $\eps>0$ there exist $\alpha_1,\alpha_2$ with
$|\alpha_1|,|\alpha_2|<\eps$ such that
\begin{equation}\label{eq:bulkminmax}
	\bulk(\alpha_1,\alpha_2):=\min_{l_p \in [0,1]}\;
	\max_{\substack{\alpha_\perp^2 \leq C',\;\alpha_\perp \geq 0} } 
	\bulk(\alpha_1,\alpha_2;l_p,\alpha_\perp)
	> 0
\end{equation}
strictly.

For fixed $l_p$,
the unconstrained supremum over $\alpha_\perp \geq 0$ is achieved at 
$\sqrt{\delta} l_p/(\sigma^2+\E[D_\star^{-1}]_{11})$,
which is smaller than ${C'}^2$ uniformly over $l_p \in [0,1]$ for sufficiently large $C'$. Thus
\begin{align*}
\bulk(\alpha_1,\alpha_2)=\inf_{l_p \in [0,1]} L(l_p;\alpha_1,\alpha_2):=
\inf_{l_p\in[0,1]}\frac{l_p^2\delta}{2}\Big[\frac{1}
{\sigma^2+\E[D_\star^{-1}]_{11}}-\frac{1}{V_\star}\Big]
-
\frac12[\alpha_2(1-l_p^2)+\alpha_1l_p^2].
\end{align*}
Let us introduce the simplified notation $(\beta_0,z) \sim \sP_0 \times \normal(0,1)$
and $\lambda=\gamma_\alg \beta_0+\sqrt{\gamma_\alg}z$, so that
\begin{equation}\label{eq:GMSVinftylaw}
(\sqrt{\gamma_\alg}G_\infty,M_\infty,S_\infty,V_\infty) \overset{d}{=}
\Big(z,\langle \beta \rangle_{\lambda,\gamma_\alg},
\langle \beta^2 \rangle_{\lambda,\gamma_\alg},
\langle \beta \rangle_{\lambda,\gamma_\alg}-\beta_0\Big).
\end{equation}
Then recall from (\ref{eq:hesshinverse}) that at $\alpha_1=\alpha_2=0$, we have
\begin{equation}\label{eq:Dinvalpha0}
D_\star^{-1}|_{\alpha_1=\alpha_2=0}=\nabla^2[{-}\hs(M_\infty,S_\infty)]^{-1}
\overset{d}{=}\begin{pmatrix} \Var[\beta \mid \lambda] &
\Cov[\beta,\beta^2 \mid \lambda] \\ \Cov[\beta,\beta^2 \mid \lambda]
& \Var[\beta^2 \mid \lambda] \end{pmatrix}.
\end{equation}
In particular, $\sigma^2+\E[D_\star^{-1}]_{11}
=\sigma^2+\E[\Var[\beta \mid \lambda]]=\sigma^2+\langle 1,S_\infty
\rangle_{\LL_2}-\|M_\infty\|_{\LL_2}^2=V_\star$, so $L(l_p;0,0)=0$
for any $l_p \in [0,1]$.
Calculating the derivative using $\partial_\alpha
D(\alpha)^{-1}={-}D(\alpha)^{-1}[\partial_\alpha D(\alpha)]D(\alpha)^{-1}$, we have 
\begin{align}
\partial_{\alpha_1} L(l_p; 0, 0) =&~ \frac{l_p^2}{2}\Big( \delta \frac{
\E[\Var[\beta \mid \lambda]^2]}{(\sigma^2 + \E[ \Var[\beta \mid \lambda]])^2} -
1\Big)=: \frac{l_p^2}{2} \chi_1, \\
\partial_{\alpha_2} L(l_p; 0, 0) =&~ \frac{l_p^2}{2}\Big( \delta \frac{
\E[\Cov[\beta, \beta^2 \mid \lambda]^2]}{(\sigma^2 + \E[ \Var[\beta \mid
\lambda]])^2} + 1\Big) - \frac12 =: \frac{l_p^2}{2}\chi_2 - \frac12.
\end{align}
Observe that $\sigma^2+\E[\Var[\beta \mid
\lambda]]=\sigma^2+\mmse(\gamma_\alg)=\delta/\gamma_\alg$,
and the condition $\phi''(\gamma_\alg)>0$ implies, by Lemma
\ref{lem:second_derivative_phi}, that $\delta/\gamma_\alg^2>\E[\Var[\beta \mid
\lambda]^2]$ strictly. Then $\chi_1<0$.
%

Now define $\bar L(l_p; \zeta) = L(l_p; -\chi_2 \zeta,  \chi_1 \zeta)$. Then we
have $\bar L(l_p; 0) = 0$ and 
\[
\partial_\zeta \bar L(l_p; 0) =  -\chi_1/2 >0
\]
where this bound holds uniformly over $l_p \in [0,1]$.
Moreover, by the smoothness of $\bar L$ and the boundedness of the support of $\sP_0$ (hence the boundedness of the posterior variance and covariance),  we have
\[
\sup_{l_p \in [0, 1]} \sup_{0 < \zeta < \eps} | \partial_\zeta^2 \bar L(l_p; \zeta) | < \infty 
\]
for some $\eps>0$.
This immediately implies that, for any $\eps>0$, there exists
some $\zeta \in (0,\eps)$ such that $\inf_{l_p \in [0, 1]} \bar L(l_p; \zeta) > 0$
strictly. Therefore, for any $\eps>0$, there exist $(\alpha_1,\alpha_2)$
with $|\alpha_1|,|\alpha_2|<\eps$ such that (\ref{eq:bulkminmax}) holds, as claimed.

\subsubsection{The spike min-max problem is non-negative}
\label{sec:spike}

\newcommand{\posta}{{b_\star}}
\newcommand{\postb}{{r_\star}}
\newcommand{\postc}{{s_\star}}
\newcommand{\postd}{{w_\star}}

In this section,
we show that for sufficiently large choices of $C'=C_0^2/2\delta$ and the
compact domain $K_\lambda$, and for any $\eps>0$ sufficiently
small and all $\alpha_1,\alpha_2$ with $|\alpha_1|,|\alpha_2|<\eps$,
\begin{equation}\label{eq:spikeminmax}
	\spiked(\alpha_1,\alpha_2):=
	\min_{|l_v|,|l_m| \leq C(\sP_0)}\;
	\sup_{(\lambda_{l_v},\lambda_{l_m}) \in K}\;
	\max_{\alpha_\parallel^2 \leq C'}\;
	\spiked(\alpha_1,\alpha_2;l_v,l_m,\lambda_{l_v},\lambda_{l_m},\alpha_\|) \geq 0.
\end{equation}
Define $\bz=(\bz_1^\top,\bz_2^\top)^\top$ where $\bz_1=(l_v,l_m)^\top,\bz_2=(\lambda_{l_v},\lambda_{l_m},\alpha_{\|})$.
The unconstrained supremum over $\bz_2$ is achieved at $(-\bC)^{-1}\bB\bz_1$,
and it is direct to check that this is uniformly bounded for
$|l_v|,|l_m| \leq C(\sP_0)$. Then
for sufficiently large $K_\lambda$ and $C'>0$,
\begin{align*}
\spike(\alpha_1,\alpha_2)=
  \frac12
  \inf_{\bz_1\in[-C(\sP_0),C(\sP_0)]^2}\sup_{\bz_2\in K \times [-C',C']}
      \bz^\top
      \begin{pmatrix}
            \bA & \bB
            \\
            \bB^\top & \bC
      \end{pmatrix}
      \bz=
\inf_{\bz_1\in[-C(\sP_0),C(\sP_0)]^2}
\frac12\bz_1^\top (\bA-\bB\bC^{-1}\bB^\top)\bz_1
\end{align*}
where $\bA,\bB,\bC$ are defined in \eqref{eq:spike_quadratic_matrices}.
Let us introduce again $(\beta_0,z) \sim \sP_0 \times \normal(0,1)$ and
$\lambda=\gamma_\alg\beta_0+\sqrt{\gamma_\alg} z$, and recall
(\ref{eq:GMSVinftylaw}) and (\ref{eq:Dinvalpha0}). Let us further define
\[k_\star=\E[\kappa_4[\beta \mid \lambda]],
\qquad b_\star=\E[\Var[\beta \mid \lambda]^2]\]
where $\kappa_4[\beta \mid \lambda]$ is the fourth cumulant of the posterior law
of $\beta$. Applying $V_\star=\sigma^2+\E[\Var[\beta \mid
\lambda]]=\sigma^2+\mmse(\gamma_\alg)=\delta/\gamma_\alg$,
the cumulant relations (\ref{eq:cumulants}), and calculations similar to
(\ref{eq:hessiancalculations}) whose details we omit here for brevity,
we obtain at $\alpha_1=\alpha_2=0$ that
\[
\begin{aligned}
\bA|_{\alpha_1=\alpha_2=0} =&~ \frac{\gamma_\alg^2}{\delta}\begin{bmatrix}
-1 & 0   \\
 0& -2
\end{bmatrix},~~~~ \bB = \begin{bmatrix} -1/2  & 0 & \sqrt{\gamma_{\alg}} \\
0 & -1/2 & 0
\end{bmatrix},\\
\bC|_{\alpha_1=\alpha_2=0} =&~ -\frac{1}{4}
\begin{pmatrix}
\posta & 0  &-2\sqrt{\gamma_{\alg}}\posta\\
 0& \frac{1}{2}\posta+\frac14k_\star &\sqrt{\gamma_{\alg}}k_\star
\\
-2\sqrt{\gamma_{\alg}}\posta&\sqrt{\gamma_{\alg}}k_\star &
4(\gamma_\alg k_\star+\delta/\gamma_\alg)
\end{pmatrix}.
\end{aligned}
\]
Therefore
\begin{align*}
&\qquad
\bA - \bB \bC^{-1} \bB^\top|_{\alpha_1=\alpha_2=0}\\
&=
\frac{\gamma_\alg^2}{\delta}\begin{bmatrix}
-1 & 0   \\
 0& -2
\end{bmatrix}\\
&\quad
+4
\begin{bmatrix} 
-1/2  & 0 & \sqrt{\gamma_\alg} \\
0 & -1/2 & 0
\end{bmatrix}
\begin{pmatrix}
\posta & 0  &-2\sqrt{\gamma_{\alg}}\posta\\
 0& \frac{1}{2}\posta+\frac14k_\star &\sqrt{\gamma_{\alg}}k_\star
\\
-2\sqrt{\gamma_{\alg}}\posta&\sqrt{\gamma_{\alg}}k_\star &
4(\gamma_\alg k_\star+\delta/\gamma_\alg)
\end{pmatrix}^{-1}
\begin{bmatrix}
 -1/2  & 0 \\
 0 & -1/2 \\
 \sqrt{\gamma_\alg}  & 0
\end{bmatrix}\\
&=
\frac{\gamma_\alg^2}{\delta}\begin{bmatrix}
-1 & 0   \\
 0& -2
\end{bmatrix}+\begin{bmatrix} 
1 & 0 & 1 \\
0 & 1 & 0   
\end{bmatrix}
\begin{pmatrix}
\posta & 0  &\posta\\
 0& \frac{1}{2}\posta+\frac14k_\star &-\frac12k_\star
\\
\posta&-\frac12k_\star &
 k_\star+\delta/\gamma^2_\alg
\end{pmatrix}^{-1}
\begin{bmatrix}
 1 & 0\\
 0 & 1 \\
 1 & 0  
\end{bmatrix}\\
&=
\frac{\gamma_\alg^2}{\delta}\begin{bmatrix}
-1 & 0   \\
 0& -2
\end{bmatrix}+\begin{bmatrix} 
1&0&0  \\
0& 1   & 0  
\end{bmatrix}
\begin{pmatrix}
\posta & 0  &0\\
 0& \frac{1}{2}\posta+\frac14k_\star &-\frac12k_\star
\\
0&-\frac12k_\star &
 k_\star-\posta+\delta/\gamma^2_\alg
\end{pmatrix}^{-1}
\begin{bmatrix}
 1  &0 \\
0& 1 \\
 0  &0   
\end{bmatrix}\\
&=
\frac{\gamma_\alg^2}{\delta}\begin{bmatrix}
-1 & 0   \\
 0& -2
\end{bmatrix}+
\begin{bmatrix}
1/\posta&0\\
0&\frac{2k_\star+2(\delta/\gamma^2_\alg-\posta)}{\posta k_\star+(\posta+k_\star/2)(\delta/\gamma_{\alg}^2-\posta)}
\end{bmatrix}.
\end{align*}
By Lemma \ref{lem:second_derivative_phi}, the condition
$\phi''(\gamma_\alg)>0$ implies $\delta/\gamma^2_\alg>\posta$ strictly.
It follows from some basic algebra that 
\begin{align*}
\bA-\bB\bC^{-1}\bB^\top\succ 0
\end{align*} when $\alpha_1=\alpha_2=0$. Moreover,  it can be verified that
$\bA-\bB\bC^{-1}\bB^\top$ is continuous in $(\alpha_1,\alpha_2)$ in a
neighborhood around the origin. Therefore, there exists some $\eps>0$ such that
$\bA-\bB\bC^{-1}\bB^\top\succ 0$ when $|\alpha_1|,|\alpha_2|<\eps$,
which implies $\spiked(\alpha_1,\alpha_2) \geq 0$.

Applying (\ref{eq:bulkminmax}) with (\ref{eq:spikeminmax}) back to
(\ref{eq:finalAuxObj3lower}), we have shown that for any $K_\alpha$ containing
an open neighborhood of 0,
\[\liminf_{\epsilon \rightarrow 0}\;\;
        \min_{\mu_l \in \bbSE_l^{(3)}(\epsilon)}\;\;
        \max_{\mu_u \in \bbSE_u^{(3)}}\;
        \AuxObj_\infty^{(3)}(\mu_u,\mu_l)
\geq \sup_{(\alpha_1,\alpha_2) \in K_\alpha}
\bulk(\alpha_1,\alpha_2)+\spiked(\alpha_1,\alpha_2)>0\]
strictly.

Combining Proposition \ref{prop:gordon-post-amp} with Lemmas \ref{lem:reduction-to-wasserstein-opt} and \ref{lem:wass-reduce-to-fixed-pt}, and the lower-bound on 
\eqref{eq:Aux3-lower-bound} established in Sections \ref{sec:cond-lower-bound},
\ref{sec:bulk}, and \ref{sec:spike}, we obtain
\eqref{eq:local-convexity-gordon-form}. This then concludes the proof of the
local strong convexity around the AMP iterates as stated in \eqref{eq:local-convexity}.

\subsubsection{Local convexity around the AMP iterate: proof of Theorem \ref{thm:local_convexity_AMP_fixed_point}\ref{item:local-convexity}}
\label{sec:finishing-local-convexity}

The claim \eqref{eq:approx-stationarity} was established in Appendix
\ref{app:approx-stationarity}, and \eqref{eq:local-convexity}
in Appendix \ref{app:localconvexity}. The last claim \eqref{eq:sub-optimality-gap} holds because,
if $\tfrac1p\cF_{\TAP}$ is $\kappa$-strongly convex on $\{\|\bbm - \bbm^k\|_2/\sqrt{n} \leq \epsilon,\,\|\bs - \bs^k \|_2/\sqrt{n} \leq \epsilon\}$, 
then
\begin{equation}
	\min_{
	\substack{
	    \|\bbm - \bbm^k\|_2/\sqrt{n} \leq \epsilon \\
	    \|\bs - \bs^k \|_2/\sqrt{n} \leq \epsilon
		}
	}\;
	\frac1p \cF_{\TAP}(\bbm,\bs)
	\geq 
	\frac1p \cF_{\TAP}(\bbm,\bs) - \frac1{2p\kappa} \| \nabla \cF_{\TAP}(\bbm^k,\bs^k) \|_2^2.
\end{equation}
Moreover, in this case, there exists a local minimizer $(\bbm_\star,\bs_\star)$ satisfying 
\begin{equation}
	\frac1p
	\Big(
		\| \bbm^k - \bbm_\star \|_2^2
		+
		\| \bs^k - \bs_\star \|_2^2
	\Big)
	\leq
	\frac{\| \nabla \cF_{\TAP}(\bbm^k,\bs^k) \|_2^2}{p\kappa^2},
\end{equation}
provided the right-hand side is smaller than $\epsilon^2$.
By \eqref{eq:approx-stationarity},
for any fixed $\kappa,\epsilon > 0$,
there exists sufficiently large $k$ that that the right-hand side is smaller than $\epsilon^2$ with high probability.
Combining these facts implies Theorem \ref{thm:local_convexity_AMP_fixed_point}\ref{item:local-convexity}.

\subsubsection{Local convexity in the easy regime: proof of Theorem \ref{thm:local_convexity}}
\label{app:Bayes-local-convexity}

Theorem \ref{thm:local_convexity} now follows from combining Theorem \ref{thm:local_convexity_AMP_fixed_point}\ref{item:local-convexity}
with the following claim: Under Assumptions \ref{ass:Bayesian_linear_model},
\ref{ass:uniquemin}, and \ref{ass:bayesoptimal},
\[\lim_{k\rightarrow \infty}
\plim_{n \rightarrow \infty} \| \bbm^k - \bbmB \|_2^2/p = 0, \qquad \lim_{k
\rightarrow \infty} \plim_{n \rightarrow \infty} \| \bs^k - \bsB \|_2^2/p=0\]
where $\bbm^k,\bs^k$ denote the AMP iterates. I.e.,
in the easy regime, the AMP iterates approximate the Bayes
estimates well for large $k$.

To show this claim, observe that in
this setting $\gamma_\alg = \gamma_\stat$, whence $\lim_{k \rightarrow \infty} \plim_{n \rightarrow \infty} \| \bbm^k - \bbeta_0 \|^2 / p = \mmse(\gamma_\stat)$.
Because the entries of $\bbm^k$ are bounded,
we in fact have $\lim_{k \rightarrow \infty} \lim_{n \rightarrow \infty} \E[\| \bbm^k - \bbeta_0 \|^2] / p = \mmse(\gamma_\stat)$.
Since $\bbm^k$ is a $(\bX,\by)$-measurable estimator of $\bbeta_0$ and $\bbmB=\E[\bbeta_0
\mid \bX,\by]$, we have $\E[\| \bbm^k - \bbeta_0 \|^2 ] = \E[\| \bbmB - \bbeta_0 \|^2] + \E[\| \bbm^k - \bbmB \|^2]$.
Because $\E[ \| \bbmB - \bbeta_0 \|^2 ]/p \rightarrow \mmse(\gamma_\stat)$ by Theorem \ref{thm:freeenergy},
we conclude that $\lim_{k \rightarrow \infty} \plim_{n \rightarrow \infty} \| \bbm^k - \bbmB \|^2/p = 0$.
Now taking $(\bbm_\star,\bs_\star)$ to be any local minimizer as in Theorem \ref{thm:marginalposterior},
we have also that $\lim_{k \rightarrow \infty} \plim_{n \rightarrow \infty} \| \bbm_\star - \bbmB \|^2 / p = 0$.

For the second statement, for $m \in (a(\sP_0),b(\sP_0))$,
define the map $s(m,\gamma) = \< \beta^2 \>_{\lambda_\gamma(m),\gamma}$, where
$\lambda_\gamma(m)$ is the unique value of $\lambda$ such that $\< \beta
\>_{\lambda_\gamma(m),\gamma} = m$ as guaranteed by Proposition
\ref{prop:oneparamexpfam}.
The function $s(\,\cdot\,,\,\cdot\,)$ is continuous and bounded because
$\lambda_\gamma$ is continuously differentiable and $\sP_0$ has bounded support.
Note that $s_j^k = s(m_j^k,\gamma_{k-1})$ and $s_{\star,j} = s(m_{\star,j},\gamma_{\star,j})$.
Because $\gamma_k \rightarrow \gamma_\stat$ by Proposition \ref{prop:amp-se},
$\gamma_{\star,j} \gotop \gamma_\stat$ 
by (\ref{eq:TAPconvergence-in-lemma-proof-gamma}),
and $\lim_{k \rightarrow \infty} \plim_{n \rightarrow \infty} \|\bbm^k -
\bbm_{\star}\|_2^2/p = 0 $,
we conclude $\|\bs_\star - s(\bbm_{\star},\gamma_\stat)\|_2^2/p \gotop 0$ and $\lim_{k
\rightarrow \infty} \plim_{n \rightarrow \infty} \|\bs^k -
s(\bbm_{\star},\gamma_\stat)\|_2^2/p = 0$.
By Theorem \ref{thm:marginalposterior},
we also have $\|\bs_{\star} - \bs_{\sB}\|_2^2/p \gotop 0$.
Thus
$\lim_{k \rightarrow \infty} \plim_{n \rightarrow \infty} \| \bs^k - \bsB
\|^2/p=0$, as claimed.

\section{Local convergence of NGD}\label{app:proof_NGD_convergence}

In this section, we carry out a generic analysis of natural gradient descent,
corresponding to Step \ref{item:NGD-step} of the proof outline given in Appendix
\ref{app:proof-outline}, and we prove Theorem \ref{thm:NGD_convergence}.

In this section, we will not use ``$\NGD$'' subscripts to denote the iterates of $\NGD$ for notational compactness.
To avoid notational confusion, let us denote the relative
entropy function and its Bregman divergence here by
\begin{align}
H(\bbm, \bbs) =&~ \sum_{j = 1}^p - \sh(m_j, s_j), \label{eqn:entropy-NGD-proof}\\
D_H((\bbm, \bbs), (\bbm', \bbs')) =&~ H(\bbm, \bbs) - H(\bbm', \bbs') - (\bbm -
\bbm', \bbs - \bbs')^\top \nabla H(\bbm', \bbs').
\end{align}
Then the NGD algorithm (\ref{eqn:NGD}) is equivalent to the Bregman gradient algorithm
\begin{equation}\label{eq:NGDbregmanform}
(\bbm^{k+1}, \bbs^{k+1}) = \argmin_{\bbm, \bbs}\;(\bbm, \bbs)^\top \nabla \cF_{\TAP}(\bbm^k, \bbs^k) + \frac{1}{\eta} D_H((\bbm, \bbs), (\bbm^k, \bbs^k)). 
\end{equation}
Indeed, this minimization over $(\bbm,\bs)$ is convex, and the first-order
condition for its minimizer is precisely the update for $(\bbm^{k+1},\bs^{k+1})$
in (\ref{eqn:NGD}) because $\nabla H(\bbm,\bs)=(\blambda,{-}\tfrac{1}{2}\bgamma)$
where $(\blambda,\bgamma)$ solves $m_j=\langle \beta \rangle_{\lambda_j,\gamma_j}$
and $s_j=\langle \beta^2 \rangle_{\lambda_j,\gamma_j}$.

We first state a lemma that ensures local convergence of an abstract form of
this Bregman gradient algorithm under certain assumptions.
\begin{lemma}\label{lem:NGD-local-analysis}
Let $\cX \subseteq \R^p$ be a convex open domain,
let $F,H,E:\cX \to \R$ be twice-continuously differentiable
with $F(\bx) = H(\bx) + E(\bx)$, and
let $\bx_\star \in \cX$ be such that $\nabla F(\bx_\star) = \bzero$. Suppose
that for some constants $\mu,\nu,L,\eps>0$,
\begin{align}
\nabla^2 H(\bx) \succeq&~ \nu \cdot \id_p \text{ for all } \bx \in \cX, \label{eqn:strong_convex_H_lem}\\
\nabla^2 F(\bx) \succeq &~ \mu \cdot \nabla^2 H(\bx) \text{ for all }
\bx \in \cX \text{ with } \|\bx-\bx_\star\|_2 \leq \eps \sqrt{p},\label{eqn:relative_convexity_lem}\\
\nabla^2 F(\bx) \preceq&~ L \cdot \nabla^2 H(\bx) \text{ for all } \bx \in \cX.\label{eqn:relative_smoothness_lem}
\end{align}
Suppose further that for all $\bx \in \cX$ and some constant $C_0>0$,
\begin{align}\label{eqn:gradient_energy_bound_lem}
\| \nabla E(\bx) \|_2 &\le C_0\sqrt{p}, \quad
E(\bx) \ge {-}C_0p, \quad H(\bx) \ge {-}C_0p, \quad F(\bx_\star) \le C_0p.
\end{align}

Define the Bregman distance $D_H(\y,\bx)=H(\y)-H(\bx)-(\y-\bx)^\top \nabla
H(\bx)$, consider the Bregman gradient algorithm 
\begin{equation}\label{eqn:NGD_x_in_lemma}
\bx^{k+1}=\argmin_{\bx \in \cX}\,\bx^\top \nabla F(\bx^k)+\frac{1}{\eta}
D_H(\bx,\bx^k),
\end{equation}
and suppose that this minimizer exists and is unique for each $k \geq 0$.
Then there exist $\eta_0,C>0$ depending only on $(\nu, \mu, L, \eps, C_0)$ such that, for any
stepsize $\eta \in (0, \eta_0]$ and any
initialization $\bx^0$ satisfying $\| \bx^0 - \bx_\star \|_2 \le \eps\sqrt{ p}$
and $F(\bx^0) - F(\bx_\star) \le p \cdot \mu\nu \eps^2 /8$, we have for all $k \geq
0$ that
\begin{align}
F(\bx^k) - F(\bx_\star) \le&~ C (1 - \mu \eta)^k \cdot D_H(\bx_\star, \bx^0), \label{eqn:function_value_convergence_NGD_in_lem} \\
\| \bx^k - \bx_\star \|_2 \le&~ \sqrt{C (1 - \mu \eta)^k \cdot D_H(\bx_\star, \bx^0)}. \label{eqn:distance_convergence_NGD_in_lem} 
\end{align}
\end{lemma}

We prove Lemma \ref{lem:NGD-local-analysis} in Appendix \ref{sec:proof-NGD-local-analysis} below.
First, we use it to prove Theorem \ref{thm:NGD_convergence}.

\begin{proof}[Proof of Theorem \ref{thm:NGD_convergence}]
	We will prove that \eqref{eqn:strong_convex_H_lem}, \eqref{eqn:relative_convexity_lem}, \eqref{eqn:relative_smoothness_lem}, and \eqref{eqn:gradient_energy_bound_lem} hold for $H$ as in \eqref{eqn:entropy-NGD-proof},
	\[
		E(\bbm,\bbs)=\frac{n}{2} \log 2\pi\sigma^2 +
		\frac{1}{2\sigma^2} \| \by - \bX \bbm \|_2^2
		+\frac{n}{2}\log\left(1+\frac{S(\bbs)-Q(\bbm)}{\sigma^2}\right),
	\]	
	and $\bx_\star = (\bbm_\star,\bs_\star)$.
	Note that $\cF_\TAP(\bbm,\bs)=H(\bbm,\bs)+E(\bbm,\bs)$.

The statement \eqref{eqn:strong_convex_H_lem} was shown in
(\ref{eq:hesshlowerbound}).
%
	For the remaining statements, we compute
	\[\begin{aligned}
	\nabla E(\bbm, \bs)&=\Big[ \frac{1}{\sigma^2} (\X^\top \X\bbm - \X^\top \y) -
	\frac{n}{p} \frac{\bbm}{\sigma^2 + S(\bs)-Q(\bbm)};\; \frac{n}{2 p}
	\frac{\ones}{\sigma^2 + S(\bs)-Q(\bbm)}\Big],\\
	\nabla_{\bbm}^2 E(\bbm, \bs)&= \frac{1}{\sigma^2} \X^\top \X - \frac{n}{p} \frac{1}{\sigma^2 + S(\bs)-Q(\bbm)} \id_p - \frac{2 n}{p^2} \frac{\bbm \bbm^\top}{(\sigma^2 + S(\bs)-Q(\bbm))^2},\\
	\nabla_{\bbm}\nabla_{\bs} E(\bbm, \bs)&=\frac{n}{p^2} \frac{1}{(\sigma^2 + S(\bs)-Q(\bbm))^2} \bbm \ones^\top, \\
	\nabla_{\bs}^2 E(\bbm, \bs) &={-}\frac{n}{2 p^2} \frac{1}{(\sigma^2 +
	S(\bs)-Q(\bbm))^2} \ones \ones^\top.
	\end{aligned}\]
	On the event $\|\bX\|_{\op} \leq C_0$ which holds with probability approaching 1
	for a sufficiently large constant $C_0>0$, it is then direct to check that
	\begin{equation}\label{eq:Ebounds}
	|E(\bbm,\bs)| \leq C, \qquad
	\|\nabla E(\bbm,\bs)\|_2 \leq C\sqrt{p}, \qquad
	\|\nabla^2 E(\bbm,\bs)\|_{\op} \leq \overline E
	\end{equation}
	for constants $C,\overline E>0$.

	To prove (\ref{eqn:relative_convexity_lem}), on the event (\ref{eq:Ebounds})
	note that 
	\[
	\nabla^2 \cF_{\TAP}(\bbm, \bs)  \succeq \nabla^2 H(\bbm, \bs) - \| \nabla^2
	E(\bbm, \bs) \|_{\op} \id_{2p} \succeq \nabla^2 H(\bbm, \bs) - \overline E \cdot \id_{2p}. 
	\]
	Furthermore, by Theorem \ref{thm:local_convexity}, for all
	$(\bbm,\bs) \in \ball((\bbm_\star,\bs_\star),\eps\sqrt{p})$ we have 
	$\nabla^2 \cF_{\TAP}(\bbm, \bs)  \succeq \kappa \id_{2p}$.
	Combining these gives
	$(\overline E + \kappa) \nabla^2 \cF_{\TAP}(\bbm, \bs) \succeq \kappa \nabla^2
	H(\bbm, \bs)$ which proves (\ref{eqn:relative_convexity_lem}). 
	Also, on the event (\ref{eq:Ebounds}) we have
	\[
	\nabla^2 \cF_{\TAP}(\bbm, \bs) =  \nabla^2 H(\bbm, \bs) + \nabla^2 E(\bbm, \bs) \preceq (1 + \overline E/\nu) \nabla^2 H(\bbm, \bs).
	\]
	which shows (\ref{eqn:relative_smoothness_lem}).
	The first condition of (\ref{eqn:gradient_energy_bound_lem}) is shown in
	(\ref{eq:Ebounds}). The second and third conditions are trivial by the
	bounds $E(\bbm,\bs) \geq \frac{n}{2}\log 2\pi\sigma^2$ and $H(\bbm,\bs) \geq 0$.
	The last condition is implied by \eqref{eqn:TAP_energy_convergence}.
	Thus, \eqref{eqn:gradient_energy_bound_lem} holds.

	The assumptions of the theorem imply  $\| \bx^0 - \bx_\star \|_2 \le \eps\sqrt{ p}$ and $F(\bx^0) - F(\bx_\star) \le p \cdot \mu\nu \eps^2 /8$.
	Thus, the statement of the theorem follows from Lemma \ref{lem:NGD-local-analysis}.
\end{proof}

\subsection{Proof of Lemma~\ref{lem:NGD-local-analysis}}\label{sec:proof-NGD-local-analysis}

The proof adapts the analysis of \cite{lu2018relatively}.\\

\noindent
{\bf Step 1. Show that $F(\bx^{k+1}) \le F(\bx^k)$.} First, by \cite[Lemma
3.2]{chen1993convergence}, the iterates $\bx^k$ satisfy the three-point
inequality, for any $\bx \in \cX$,
\begin{equation}\label{eqn:three_point_in_lemma_proof}
\nabla F(\bx^{k})^\top (\bx^{k+1}-\bx) + \eta^{-1} \cdot D_H( \bx^{k+1}, \bx^k) +
\eta^{-1} \cdot D_H(\bx, \bx^{k+1}) \le \eta^{-1} \cdot D_H(\bx, \bx^{k}).
\end{equation}
(The proof in \cite{chen1993convergence} extends directly to our setting of an
open domain $\cX \subseteq \R^p$, under the given assumption that the infimum in
(\ref{eqn:NGD_x_in_lemma}) is always attained at some $\bx \in \cX$.)
Furthermore, the relative smoothness condition (\ref{eqn:relative_smoothness_lem}) implies that
$\eta^{-1}H-F$ is convex when $\eta \le 1/L$, and hence for all $\bx \in \cX$,
\begin{equation}\label{eqn:relative_smoothness_in_lemma_proof}
F(\bx) \le F(\bx^k) + \nabla F(\bx^k)^\top(\bx - \bx^k) + \eta^{-1} \cdot
D_H(\bx, \bx^k).
\end{equation}
Combining the inequalities above and noting that $D_H(\y,\bx) \geq 0$ by
convexity of $H$, we obtain
\[
\begin{aligned}
F(\bx^{k+1}) \le F(\bx^k) + \nabla F(\bx^k)^\top(\bx^{k+1} - \bx^k) + \eta^{-1} \cdot D_H(\bx^{k+1}, \bx^k)
\le F(\bx^k)  - \eta^{-1} \cdot D_H(\bx^k, \bx^{k+1}) \le F(\bx^k). 
\end{aligned}
\]

\noindent
{\bf Step 2. Show that $\| \bx^k - \bx_\star \|_2 \le \eps \sqrt{p}$ and
$F(\bx^{k}) - F(\bx_\star) \le p \cdot \mu\nu \eps^2 / 8$.} We prove this by
induction. Note that for $k = 0$, this is satisfied by the assumption of the
lemma. Assume this holds for $k$. Then the result of Step 1 above
implies that $F(\bx^{k+1}) \le F(\bx^k) \le F(\bx_\star) + p \cdot \nu \eps^2 /
8$. Observe that for any $\rho < \eps$, by the strong convexity $\nabla^2 F(\bx)
\succeq \mu \cdot \nabla^2 H(\bx) \succeq \mu\nu \cdot \id_p$ when
$\|\bx-\bx_\star\|_2 \leq \eps \sqrt{p}$, we have the implication 
\begin{equation}\label{eqn:growth_condition_in_proof}
\| \bx - \bx_\star \|_2/\sqrt{p} \le \eps~ \text{and}~ F(\bx) - F(\bx_\star)
\leq p \cdot \mu\nu \rho^2 / 2~~~ \Rightarrow~~~ \| \bx - \bx_\star \|_2
/\sqrt{p} \le \rho.
\end{equation}
As a consequence, $F(\bx^k) - F(\bx_\star) < p \cdot \mu\nu \eps^2 / 8$ implies that $\| \bx^k - \bx_\star \|_2 /\sqrt{p} \le \eps/2$. 

By the optimality condition of (\ref{eqn:NGD_x_in_lemma}) for $\bx^{k+1}$
compared to the value at $\bx^k$, we have 
\begin{equation}\label{eqn:three_point_special_in_proof_lem}
\nabla F(\bx^{k})^\top(\bx^{k+1} - \bx^k) + \eta^{-1} \cdot D_H( \bx^{k+1}, \bx^k)  \le 0. 
\end{equation}
From the definitions $D_H(\y, \bx)=H(\y)-H(\bx)-(\y-\bx)^\top \nabla H(\bx)$
and $F(\bx)=H(\bx)+E(\bx)$, we have 
\[
\begin{aligned}
\Big\vert \nabla F(\bx^k)^\top(\bx^{k+1} - \bx^k) + D_H(\bx^{k+1}, \bx_k) \Big\vert 
&= \Big\vert \nabla E(\bx^k)^\top(\bx^{k+1} - \bx^k) + H(\bx^{k+1}) - H(\bx^k) \Big\vert \\
&\le C_0\sqrt{p} \cdot \| \bx^{k+1} - \bx^k \|_2 + \vert H(\bx^{k+1}) \vert + \vert H(\bx^k) \vert
\end{aligned}
\]
where the last inequality uses assumption (\ref{eqn:gradient_energy_bound_lem}). 
Also $H(\bx^k) \geq -C_0p$ and $H(\bx^k)=F(\bx^k)-E(\bx^k)
\leq F(\bx_\star)+p \cdot \mu\nu\eps^2/8+C_0p$, and similarly for $H(\bx^{k+1})$,
so
\begin{equation}\label{eqn:nabla_F_inner_product_x_estimate_in_proof_lem}
\begin{aligned}
\Big\vert \nabla F(\bx^k)^\top(\bx^{k+1} - \bx^k) + D_H(\bx^{k+1}, \bx_k) \Big\vert 
\le&~ C\sqrt{p} \cdot \| \bx^{k+1} - \bx^k \|_2 + Cp
\end{aligned}
\end{equation}
for some constant $C>0$ that depends on $(\mu,\nu,\eps,C_0)$. Furthermore, by
the definition of $D_H$ and by (\ref{eqn:strong_convex_H_lem}), we have 
\begin{equation}\label{eqn:D_H_lower_strong_convex_in_proof_lem}
D_H( \bx^{k+1}, \bx^k) \ge \nu\cdot \| \bx^{k+1} - \bx^k \|_2^2. 
\end{equation}
Combining (\ref{eqn:three_point_special_in_proof_lem}), (\ref{eqn:nabla_F_inner_product_x_estimate_in_proof_lem}), and (\ref{eqn:D_H_lower_strong_convex_in_proof_lem}), we get 
\[
- C\sqrt{p} \cdot \| \bx^{k+1} - \bx^k \|_2 - Cp +  (\eta^{-1} - 1) \nu\cdot \| \bx^{k+1} - \bx^k \|_2^2 \le 0, 
\]
which implies when $\eta<1$ that
\[
\| \bx^{k+1} - \bx^k \|_2 \le \frac{C'\sqrt{p}}{\sqrt{\eta^{-1} -1}}
\]
for some $(\mu,\nu,\eps,C_0)$ dependent constant $C'>0$. Then, for $\eta$ small enough so that $C'/\sqrt{\eta^{-1} -1 }\le \eps/2$, we obtain 
\[
\| \bx^{k+1} - \bx_\star \|_2 \le \| \bx^k - \bx_\star \|_2 + \| \bx^{k+1} - \bx^k \|_2 \le \eps \sqrt{p}, 
\]
completing the induction.\\

\noindent
{\bf Step 3. Finish the proof.} Since $\|\bx^k-\bx_\star\|_2 \leq \eps\sqrt{
p}$, the relative convexity condition (\ref{eqn:relative_convexity_lem}) implies
that $F-\mu H$ is convex on the line segment between $\bx^k$ and $\bx_\star$, so
\[
F(\bx^k) + \nabla F(\bx^k)^\top(\bx_\star - \bx^k) \le F(\bx_\star) - \mu \cdot D_H(\bx_\star, \bx^k). 
\]
By the three point inequality (\ref{eqn:three_point_in_lemma_proof}) applied
with $\bx=\bx_\star$ and relative smoothness (\ref{eqn:relative_smoothness_in_lemma_proof}), we obtain 
\[
F(\bx^{k+1}) \le F(\bx^k) + \nabla F(\bx^k)^\top(\bx_\star - \bx^k) + \eta^{-1} \cdot D_H(\bx_\star, \bx^k) -\eta^{-1} \cdot D_H(\bx_\star, \bx^{k+1}). 
\]
Combining the inequalities above, we obtain 
\[
F(\bx^{k+1}) \le F(\bx_\star) + (\eta^{-1} - \mu) \cdot D_H(\bx_\star, \bx^k) - \eta^{-1} \cdot D_H(\bx_\star, \bx^{k+1}). 
\]
Multiplying by $[1 / (1 - \mu \eta)]^{k+1}$ and summing over $k$ to telescope the sums of the last two terms, 
\[
\sum_{j=0}^{k-1} \Big( \frac{1}{1 - \mu \eta} \Big)^{j+1} F(\bx^{j+1}) \le
\sum_{j=0}^{k-1} \Big( \frac{1}{1 - \mu \eta} \Big)^{j+1}  F(\bx_\star) +\eta^{-1} D_H(\bx_\star, \bx^0). 
\]
Now applying $F(\bx^{j+1}) \ge F(\bx^k)$ for all $j \le k - 1$ to the left side, we obtain
\[
F(\bx^k) \le F(\bx_\star) + \eta^{-1} ( 1 - \mu \eta)^k D_H(\bx_\star, \bx^0). 
\]
This proves (\ref{eqn:function_value_convergence_NGD_in_lem}). Then
(\ref{eqn:distance_convergence_NGD_in_lem}) follows from the conclusion of Step
2 and (\ref{eqn:growth_condition_in_proof}). This completes the proof of Lemma \ref{lem:NGD-local-analysis}.

\section{AMP+NGD in the easy and hard regimes}

In this section, we establish the convergence of AMP+NGD in the easy and hard regime,
and demonstrate that it achieves calibrated inference in both regimes. This corresponds to Steps \ref{item:AMP+NGD-step} and \ref{item:calibration-step} of the proof outline given in Appendix \ref{app:proof-outline}.
We will carry out Step \ref{item:AMP+NGD-step} in Appendx \ref{app:AMP+NGD-convergence} and Step \ref{item:calibration-step} in Appendix \ref{app:calibration}.

\subsection{Convergence of AMP+NGD: proofs of Corollary \ref{cor:AMP+NGD-convergence} and Theorem \ref{thm:local_convexity_AMP_fixed_point}\ref{item:NGD-convergence}}
\label{app:AMP+NGD-convergence}

\begin{proof}[Proof of Corollary \ref{cor:AMP+NGD-convergence}]
	Under Assumptions \ref{ass:Bayesian_linear_model}, \ref{ass:uniquemin}, and \ref{ass:bayesoptimal} (i.e., in the easy regime),
	with high probability the local minimizers $(\bbm_\star,\bs_\star)$ of
Theorems \ref{thm:local_convexity} and \ref{thm:local_convexity_AMP_fixed_point}\ref{item:local-convexity} coincide.
	Thus, combining Theorem
\ref{thm:local_convexity_AMP_fixed_point}\ref{item:local-convexity} with
\eqref{eq:approx-stationarity} and \eqref{eq:sub-optimality-gap} gives that for
large enough $T_0$, initializing NGD at $(\bbm_{\NGD}^0,\bs_{\NGD}^0) =
(\bbm_{\AMP}^{T_0},\bs_{\AMP}^{T_0})$ satisfies the conditions of Theorem \ref{thm:NGD_convergence}.
	Corollary \ref{cor:AMP+NGD-convergence} follows.
\end{proof}

\begin{proof}[Proof of Theorem \ref{thm:local_convexity_AMP_fixed_point}\ref{item:NGD-convergence}]
	Although Theorem \ref{thm:NGD_convergence} is stated for the Bayes-optimal local minimizer,
	by inspecting its proof, we see that it applies to any local minimizer around which \eqref{eqn:strong_convex_H_lem}, \eqref{eqn:relative_convexity_lem}, \eqref{eqn:relative_smoothness_lem}, and \eqref{eqn:gradient_energy_bound_lem} can be established.
All of these conditions hold
for the local minimizer $(\bbm_\star,\bs_\star)$ of Theorem
\ref{thm:local_convexity_AMP_fixed_point}\ref{item:local-convexity} by the exact
same argument as in the proof of Theorem \ref{thm:NGD_convergence}, except for
the condition $\cF_{\TAP}(\bbm_\star,\bs_\star) \leq C_0 p$ in \eqref{eqn:gradient_energy_bound_lem}.
	For this, we use that, for sufficiently large but fixed $k$, with high probability
	\begin{equation}
	\begin{aligned}
		&\frac1p \cF_\TAP(\bbm_\star,\bs_\star)
			\leq
			\frac1p \cF_\TAP(\bbm_{\AMP}^k,\bs_{\AMP}^k)
		\\
			&=
			\frac{n/p}{2}\log(2\pi\sigma^2)
			+
			\frac1pD_0(\bbm_{\AMP}^k,\bs_{\AMP}^k)
			+
			\frac1{2p\sigma^2}\|\by-\bX\bbm_{\AMP}^k\|_2^2
			+
			\frac{n/p}{2}\log\Big( 1 + \frac{S(\bs_{\AMP}^k) - Q(\bbm_{\AMP}^k)}{\sigma^2}\Big).
	\end{aligned}
	\end{equation}
	Letting $C > 0$ denote a constant depending only on $(\sigma^2,\delta,\sP_0)$ whose value can change at each appearance,
	with high probability $\| \by - \bX \bm_{\AMP}^k\|_2^2/ \leq C p$
because $\| \bX \|_{\op} \leq C$, $\| \beps \|_2^2/p \leq C$, and $\| \bbeta_0 -
\bbm_{\AMP}^k \|_2^2 /p \leq C$ as its entries are bounded.
	By state evolution (Proposition \ref{prop:amp-se}), $\log\Big( 1 + \frac{S(\bs_{\AMP}^k) - Q(\bbm_{\AMP}^k)}{\sigma^2}\Big) \leq C$ with high probability.
	Finally, using the definition of $\bbm_{\AMP}^k$ and $\bs_{\AMP}^k$ in \eqref{eq:AMP} and the definition of $D_0$ (see \eqref{eq:hdef}),
	we have
	\begin{equation}
		D_0(\bbm_{\AMP}^k,\bs_{\AMP}^k)
		=
		\sum_{j=1}^p 
		\Big(
			-\frac12\gamma_{k-1} s_{\AMP,j}^k
			+
			\lambda_{\AMP,j}^{k-1} m_{\AMP,j}^k
			-
			\log \E_{\beta \sim \sP_0}
			\Big[
				e^{-(\gamma_{k-1}/2)\beta^2 + \lambda_{\AMP,j}^{k-1}\beta}
			\Big]
		\Big),
	\end{equation}
	where $\blambda_{\AMP}^{k-1}$ is as in \eqref{eq:AMPlambdagamma}.
	Because $(\gamma,\lambda) \mapsto \log \E_{\beta \sim \sP_0}\big[e^{-(\gamma/2)\beta^2 + \lambda\beta}\big]$ is pseudo-Lipschitz (a consequence of the boundedness of the support of $\sP_0$),
	we conclude that $D_0(\bbm_{\AMP}^k,\bs_{\AMP}^k) \leq C$ with high probability by state evolution (Proposition \ref{prop:amp-se}).
	Having established \eqref{eqn:strong_convex_H_lem}, \eqref{eqn:relative_convexity_lem}, \eqref{eqn:relative_smoothness_lem}, and \eqref{eqn:gradient_energy_bound_lem},
	Theorem \ref{thm:local_convexity_AMP_fixed_point}\ref{item:NGD-convergence} follows.
\end{proof}

\subsection{Risk and calibration of local minimizer: proof of Theorem \ref{thm:error-and-calibration}}
\label{app:calibration}

We prove Theorem \ref{thm:error-and-calibration}.
The statement (\ref{eq:localminerror}) follows by combining
$\lim_{n,p \to \infty} p^{-1}\|\bbm^k-\bbeta_0\|_2^2=\mmse(\gamma_{k-1})$ and
$\lim_{k \to \infty} \gamma_{k-1}=\gamma_\alg$ from AMP
state evolution (Proposition \ref{prop:amp-se}(a,d)) with \ref{thm:local_convexity_AMP_fixed_point}\ref{item:NGD-convergence}.

For the remaining statements of Theorem \ref{thm:error-and-calibration},
denote $\blambda_\star=[\lambda(m_{j,\star},s_{j,\star})]_{j=1}^p$
and $\bgamma_\star=[\gamma(m_{j,\star},s_{j,\star})]_{j=1}^p$.
By the stationary condition
$\nabla \cF_\TAP(\bbm_\star,\bs_\star) = \bzero$,
\begin{equation}\label{eq:stationaryrestatement}
\begin{gathered}
	\blambda_\star=
		\frac1{\sigma^2}\bX^\top(\by - \bX \bbm_\star)
		+\frac{n/p}{\sigma^2 + S(\bs_\star) - Q(\bbm_\star)}
\cdot\bbm_\star, 
	\qquad
	\bgamma_\star=\frac{\delta}{\sigma^2 + S(\bs_\star) -
Q(\bbm_\star)} \cdot \ones.
\end{gathered}
\end{equation}
In particular, $\gamma_{j,\star}$ is constant across coordinates $j$.
Define $\blambda^k=[\lambda(m_j^k,s_j^k)]_{j=1}^p$ and
$\bgamma^k=[\gamma(m_j^k,s_j^k)]_{j=1}^p$ where $(\bbm^k,\bs^k)$ are
the AMP iterates of (\ref{eq:AMP}). Recall
that $\lim_{k \to \infty} \plimsup_{n,p \to \infty} p^{-1}\|\nabla
\cF_\TAP(\bbm^k,\bs^k)\|_2^2=0$ by (\ref{eq:approx-stationarity}),
where $\nabla \cF_\TAP(\bbm^k,\bs^k)$ is given
explicitly by (\ref{eq:TAP-grad}). Comparing this with
(\ref{eq:stationaryrestatement}) and using Theorem \ref{thm:local_convexity_AMP_fixed_point}\ref{item:NGD-convergence},
this implies
\begin{equation}\label{eq:lambdagammaconvergence}
	\lim_{k \rightarrow \infty} \plimsup_{n,p \rightarrow \infty} \frac1p \|
\blambda^k - \blambda_\star \|_2^2 = 0
	\qquad
	\text{and}
	\qquad
	\lim_{k \rightarrow \infty} \plimsup_{n,p \rightarrow \infty} \frac1p
\|\bgamma^k - \bgamma_\star\|_2^2 = 0.
\end{equation}
Recalling from (\ref{eq:AMPlambdagamma}) that $\bgamma^k$ is deterministic with
all coordinates equal to $\gamma_{k-1}$, and that
$\lim_{k \to \infty} \gamma_{k-1}=\gamma_\alg$, the second limit of
(\ref{eq:lambdagammaconvergence}) proves the statement
$\gamma_{j,\star}=\gamma(m_{j,\star},s_{j,\star}) \gotop \gamma_\alg$.

Finally, recall from (\ref{eq:AMPlambdagamma}) that
$\blambda^k=\gamma_{k-1}(\bbm^{k-1}+\frac{1}{\delta}\bX^\top \bz^{k-1})$,
whose state evolution is implied by that for
$(\bbeta^0,\bg^k)$ in Proposition \ref{prop:amp-se}(c).
Applying this state evolution, the first limit of (\ref{eq:lambdagammaconvergence}),
and the convergence $\lim_{k \to \infty} \gamma_{k-1}=\gamma_\alg$, we
obtain for any bounded and Lipschitz function $\varphi: \reals^2 \rightarrow [0,1]$
that
	\[\frac1p \sum_{j=1}^p \varphi(\beta_{0,j},\lambda_{j,\star})
		\gotop
		\E_{(\beta_0,z) \sim \sP_0 \times
\normal(0,1)}\big[\varphi(\beta_0,\gamma_\alg \beta_0 +
\sqrt{\gamma_\alg}\,z)\big].\]
By symmetry,
$\E[\varphi(\beta_{0,j},\lambda_{j,\star})]$ does not depend on $j$.
Because $\varphi$ is bounded, we conclude 
	\[\E[\varphi(\beta_{0,j},\lambda_{j,\star})] \rightarrow \E_{\beta_0,z
\sim \sP_0 \times \normal(0,1)}[\varphi(\beta_0,\gamma_\alg \beta_0 +
\sqrt{\gamma_\alg}\,z)].\]
Because the distribution of $\gamma_\alg \beta_0 + \sqrt{\gamma_\alg}\,z$ is absolutely continuous with respect to Lebesgue measure,
for any $f$ that is Lipschitz and bounded and for any non-empty open set $A$, we
may use the above convergence together with
standard lower and upper Lipschitz approximations of $\lambda \mapsto \mathbf{1}\{\lambda \in A\}$ to conclude that
\[\begin{gathered}
	\E[f(\beta_j)\mathbf{1}\{\lambda_{j,\star} \in A\}] \rightarrow \E_{(\beta_0,z) \sim \sP_0 \times \normal(0,1)}\big[f(\beta)\mathbf{1}\{\gamma_\alg \beta_0 + \sqrt{\gamma_\alg}\,z \in A\}\big],
	\\
	\E[\mathbf{1}\{\lambda_{j,\star} \in A\}] \rightarrow \E_{(\beta_0,z) \sim \sP_0 \times \normal(0,1)}\big[\mathbf{1}\{\gamma_\alg \beta_0 + \sqrt{\gamma_\alg}\,z \in A\}\big] > 0,
\end{gathered}\]
whence
$\E[f(\beta_j) \mid \lambda_{j,\star} \in A] \rightarrow \E_{(\beta_0,z) \sim \sP_0 \times \normal(0,1)}\big[f(\beta) \mid \gamma_\alg \beta_0 + \sqrt{\gamma_\alg}\,z \in A \big]$.

\bibliographystyle{myalpha}
\bibliography{main.bbl}

\newcommand{\etalchar}[1]{$^{#1}$}
\begin{thebibliography}{LTBS{\etalchar{+}}15}
\expandafter\ifx\csname url\endcsname\relax
  \def\url#1{\texttt{#1}}\fi
\expandafter\ifx\csname doi\endcsname\relax
  \def\doi#1{\burlalt{doi:#1}{http://dx.doi.org/#1}}\fi
\expandafter\ifx\csname urlprefix\endcsname\relax\def\urlprefix{URL }\fi
\expandafter\ifx\csname href\endcsname\relax
  \def\href#1#2{#2}\fi
\expandafter\ifx\csname burlalt\endcsname\relax
  \def\burlalt#1#2{\href{#2}{#1}}\fi

\bibitem[Ama98]{Amari1998}
Shun-ichi Amari.
\newblock Natural gradient works efficiently in learning.
\newblock {\em Neural Computation}, 10(2):251--276, 1998.

\bibitem[BDMK16]{barbier2016mutual}
Jean Barbier, Mohamad Dia, Nicolas Macris, and Florent Krzakala.
\newblock The mutual information in random linear estimation.
\newblock In {\em 2016 54th Annual Allerton Conference on Communication,
  Control, and Computing (Allerton)}, pages 625--632. IEEE, 2016.

\bibitem[BKM17]{blei2017variational}
David~M Blei, Alp Kucukelbir, and Jon~D McAuliffe.
\newblock Variational inference: {A} review for statisticians.
\newblock {\em Journal of the American Statistical Association},
  112(518):859--877, 2017.

\bibitem[BKM{\etalchar{+}}19]{barbier2019optimal}
Jean Barbier, Florent Krzakala, Nicolas Macris, L{\'e}o Miolane, and Lenka
  Zdeborov{\'a}.
\newblock Optimal errors and phase transitions in high-dimensional generalized
  linear models.
\newblock {\em Proceedings of the National Academy of Sciences},
  116(12):5451--5460, 2019.

\bibitem[Bla85]{Blair1985}
Charles Blair.
\newblock Problem complexity and method efficiency in optimization.
\newblock {\em SIAM Review}, 27(2):264--265, 1985.

\bibitem[BLM15]{bayati2015universality}
Mohsen Bayati, Marc Lelarge, and Andrea Montanari.
\newblock Universality in polytope phase transitions and message passing
  algorithms.
\newblock {\em Annals of Applied Probability}, 25(2):753--822, 2015.

\bibitem[BM11]{bayati2011dynamics}
Mohsen Bayati and Andrea Montanari.
\newblock The dynamics of message passing on dense graphs, with applications to
  compressed sensing.
\newblock {\em IEEE Transactions on Information Theory}, 57(2):764--785, 2011.

\bibitem[BM19]{barbier2019adaptive}
Jean Barbier and Nicolas Macris.
\newblock The adaptive interpolation method: a simple scheme to prove replica
  formulas in {B}ayesian inference.
\newblock {\em Probability Theory and Related Fields}, 174:1133--1185, 2019.

\bibitem[BMDK20]{barbier2020mutual}
Jean Barbier, Nicolas Macris, Mohamad Dia, and Florent Krzakala.
\newblock Mutual information and optimality of approximate message-passing in
  random linear estimation.
\newblock {\em IEEE Transactions on Information Theory}, 66(7):4270--4303,
  2020.

\bibitem[BMN19]{berthierMontanariNguyen}
Rapha{\"e}l Berthier, Andrea Montanari, and Phan-Minh Nguyen.
\newblock State evolution for approximate message passing with non-separable
  functions.
\newblock {\em Information and Inference: A Journal of the IMA}, 9(1):33--79,
  01 2019.

\bibitem[Bol14]{Bolthausen2014}
Erwin Bolthausen.
\newblock An iterative construction of solutions of the {TAP} equations for the
  {S}herrington--{K}irkpatrick model.
\newblock {\em Communications in Mathematical Physics}, 325(1):333--366, Jan
  2014.

\bibitem[BT03]{BECK2003167}
Amir Beck and Marc Teboulle.
\newblock Mirror descent and nonlinear projected subgradient methods for convex
  optimization.
\newblock {\em Operations Research Letters}, 31(3):167--175, 2003.

\bibitem[Cel22]{celentano2022sudakov}
Michael Celentano.
\newblock Sudakov-{F}ernique post-{AMP}, and a new proof of the local convexity
  of the {TAP} free energy, 2022, arXiv:2208.09550 [math.PR].

\bibitem[CFM23]{celentano2023local}
Michael Celentano, Zhou Fan, and Song Mei.
\newblock Local convexity of the {TAP} free energy and {AMP} convergence for
  $\mathbb{Z}_2$-synchronization.
\newblock {\em The Annals of Statistics}, 51(2):519--546, 2023.

\bibitem[CL21]{chen2021universality}
Wei~Kuo Chen and Wai-Kit Lam.
\newblock Universality of approximate message passing algorithms.
\newblock {\em Electronic Journal of Probability}, 26:1 -- 44, 2021.

\bibitem[CM22]{celentano2022fundamental}
Michael Celentano and Andrea Montanari.
\newblock Fundamental barriers to high-dimensional regression with convex
  penalties.
\newblock {\em The Annals of Statistics}, 50(1):170--196, 2022.

\bibitem[CS12]{carbonetto2012scalable}
Peter Carbonetto and Matthew Stephens.
\newblock Scalable variational inference for {B}ayesian variable selection in
  regression, and its accuracy in genetic association studies.
\newblock {\em Bayesian Analysis}, 7(1):73--108, 2012.

\bibitem[CT93]{chen1993convergence}
Gong Chen and Marc Teboulle.
\newblock Convergence analysis of a proximal-like minimization algorithm using
  {B}regman functions.
\newblock {\em SIAM Journal on Optimization}, 3(3):538--543, 1993.

\bibitem[DMM09]{donoho2009message}
David~L Donoho, Arian Maleki, and Andrea Montanari.
\newblock Message-passing algorithms for compressed sensing.
\newblock {\em Proceedings of the National Academy of Sciences},
  106(45):18914--18919, 2009.

\bibitem[DMM10]{donoho2010message}
David~L Donoho, Arian Maleki, and Andrea Montanari.
\newblock Message passing algorithms for compressed sensing: {I}. motivation
  and construction.
\newblock In {\em 2010 IEEE information theory workshop on information theory
  (ITW)}, pages 1--5. IEEE, 2010.

\bibitem[EAMS22]{el2022sampling}
Ahmed El~Alaoui, Andrea Montanari, and Mark Sellke.
\newblock Sampling from the {S}herrington-{K}irkpatrick {G}ibbs measure via
  algorithmic stochastic localization.
\newblock In {\em 2022 IEEE 63rd Annual Symposium on Foundations of Computer
  Science (FOCS)}, pages 323--334. IEEE, 2022.

\bibitem[EAMS23]{el2023sampling}
Ahmed El~Alaoui, Andrea Montanari, and Mark Sellke.
\newblock Sampling from mean-field {G}ibbs measures via diffusion processes,
  2023, arXiv:2310.08912 [math.PR].

\bibitem[FMM21]{fan2021tap}
Zhou Fan, Song Mei, and Andrea Montanari.
\newblock {TAP} free energy, spin glasses and variational inference.
\newblock {\em The Annals of Probability}, 49(1):1 -- 45, 2021.

\bibitem[GBJ18]{GiordanoBroderickJordan2018}
Ryan Giordano, Tamara Broderick, and Michael~I. Jordan.
\newblock Covariances, robustness, and variational {B}ayes.
\newblock {\em Journal of Machine Learning Research}, 19(51):1--49, 2018.

\bibitem[GDKZ23]{ghio2023sampling}
Davide Ghio, Yatin Dandi, Florent Krzakala, and Lenka Zdeborov{\'a}.
\newblock Sampling with flows, diffusion and autoregressive neural networks: A
  spin-glass perspective, 2023, arXiv:2308.14085 [cond-mat.dis-nn].

\bibitem[Gor85]{Gordon1985}
Yehoram Gordon.
\newblock Some inequalities for {G}aussian processes and applications.
\newblock {\em Israel Journal of Mathematics}, 50(4):265--289, Dec 1985.

\bibitem[Gor88]{gordon1988}
Yerhoram Gordon.
\newblock On {M}ilman's inequality and random subspaces which escape through a
  mesh in $\mathbb{R}^n$.
\newblock In {\em Geometric Aspects of Functional Analysis}, pages 84--106,
  Berlin, Heidelberg, 1988. Springer Berlin Heidelberg.

\bibitem[GP10]{guillemin2010differential}
Victor Guillemin and Alan Pollack.
\newblock {\em Differential topology}, volume 370.
\newblock American Mathematical Society, 2010.

\bibitem[GWSV11]{guo2011estimation}
Dongning Guo, Yihong Wu, Shlomo~S Shitz, and Sergio Verd{\'u}.
\newblock Estimation in {G}aussian noise: Properties of the minimum mean-square
  error.
\newblock {\em IEEE Transactions on Information Theory}, 57(4):2371--2385,
  2011.

\bibitem[HS23]{han2022universality}
Qiyang Han and Yandi Shen.
\newblock {Universality of regularized regression estimators in high
  dimensions}.
\newblock {\em The Annals of Statistics}, 51(4):1799 -- 1823, 2023.

\bibitem[Kab03]{kabashima2003cdma}
Yoshiyuki Kabashima.
\newblock A {CDMA} multiuser detection algorithm on the basis of belief
  propagation.
\newblock {\em Journal of Physics A: Mathematical and General}, 36(43):11111,
  2003.

\bibitem[KMS{\etalchar{+}}12]{krzakala2012probabilistic}
Florent Krzakala, Marc M{\'e}zard, Francois Sausset, Yifan Sun, and Lenka
  Zdeborov{\'a}.
\newblock Probabilistic reconstruction in compressed sensing: algorithms, phase
  diagrams, and threshold achieving matrices.
\newblock {\em Journal of Statistical Mechanics: Theory and Experiment},
  2012(08):P08009, 2012.

\bibitem[KMTZ14]{krzakala2014variational}
Florent Krzakala, Andre Manoel, Eric~W Tramel, and Lenka Zdeborov{\'a}.
\newblock Variational free energies for compressed sensing.
\newblock In {\em 2014 IEEE International Symposium on Information Theory},
  pages 1499--1503. IEEE, 2014.

\bibitem[LFN18]{lu2018relatively}
Haihao Lu, Robert~M Freund, and Yurii Nesterov.
\newblock Relatively smooth convex optimization by first-order methods, and
  applications.
\newblock {\em SIAM Journal on Optimization}, 28(1):333--354, 2018.

\bibitem[LFW23]{li2023approximate}
Gen Li, Wei Fan, and Yuting Wei.
\newblock Approximate message passing from random initialization with
  applications to $\mathbb{Z}_2$ synchronization.
\newblock {\em Proceedings of the National Academy of Sciences},
  120(31):e2302930120, 2023.

\bibitem[LHM10]{logsdon2010variational}
Benjamin~A Logsdon, Gabriel~E Hoffman, and Jason~G Mezey.
\newblock A variational {B}ayes algorithm for fast and accurate multiple locus
  genome-wide association analysis.
\newblock {\em BMC Bioinformatics}, 11:1--13, 2010.

\bibitem[LM19]{lelarge2019fundamental}
Marc Lelarge and L{\'e}o Miolane.
\newblock Fundamental limits of symmetric low-rank matrix estimation.
\newblock {\em Probability Theory and Related Fields}, 173:859--929, 2019.

\bibitem[LTBS{\etalchar{+}}15]{loh2015efficient}
Po-Ru Loh, George Tucker, Brendan~K Bulik-Sullivan, Bjarni~J Vilhj{\'a}lmsson,
  Hilary~K Finucane, Rany~M Salem, Daniel~I Chasman, Paul~M Ridker, Benjamin~M
  Neale, Bonnie Berger, et~al.
\newblock Efficient {B}ayesian mixed-model analysis increases association power
  in large cohorts.
\newblock {\em Nature Genetics}, 47(3):284--290, 2015.

\bibitem[LW22]{li2022non}
Gen Li and Yuting Wei.
\newblock A non-asymptotic framework for approximate message passing in spiked
  models, 2022, arXiv:2208.03313 [math.ST].

\bibitem[MFC{\etalchar{+}}19]{maillard2019high}
Antoine Maillard, Laura Foini, Alejandro~Lage Castellanos, Florent Krzakala,
  Marc M{\'e}zard, and Lenka Zdeborov{\'a}.
\newblock High-temperature expansions and message passing algorithms.
\newblock {\em Journal of Statistical Mechanics: Theory and Experiment},
  2019(11):113301, 2019.

\bibitem[Mon23]{montanari2023sampling}
Andrea Montanari.
\newblock Sampling, diffusions, and stochastic localization, 2023,
  arXiv:2305.10690 [cs.LG].

\bibitem[MS22]{mukherjee2022variational}
Sumit Mukherjee and Subhabrata Sen.
\newblock Variational inference in high-dimensional linear regression.
\newblock {\em The Journal of Machine Learning Research}, 23(1):13703--13758,
  2022.

\bibitem[MT06]{montanari2006analysis}
Andrea Montanari and David Tse.
\newblock Analysis of belief propagation for non-linear problems: The example
  of {CDMA} (or: {H}ow to prove {T}anaka's formula).
\newblock In {\em 2006 IEEE Information Theory Workshop (ITW)}, pages 160--164.
  IEEE, 2006.

\bibitem[MW23a]{mei2023deep}
Song Mei and Yuchen Wu.
\newblock Deep networks as denoising algorithms: Sample-efficient learning of
  diffusion models in high-dimensional graphical models, 2023, arXiv:2309.11420
  [cs.LG].

\bibitem[MW23b]{montanari2023posterior}
Andrea Montanari and Yuchen Wu.
\newblock Posterior sampling from the spiked models via diffusion processes,
  2023, arXiv:2304.11449 [math.ST].

\bibitem[OYM17]{ormerod2017variational}
John~T Ormerod, Chong You, and Samuel M{\"u}ller.
\newblock A variational {B}ayes approach to variable selection.
\newblock {\em Electronic Journal of Statistics}, 11:3549--3594, 2017.

\bibitem[Pan13]{panchenko2013sherrington}
Dmitry Panchenko.
\newblock {\em The {S}herrington-{K}irkpatrick model}.
\newblock Springer Science \& Business Media, 2013.

\bibitem[PP09]{payaro2009hessian}
Miquel Payar{\'o} and Daniel~P Palomar.
\newblock Hessian and concavity of mutual information, differential entropy,
  and entropy power in linear vector gaussian channels.
\newblock {\em IEEE Transactions on Information Theory}, 55(8):3613--3628,
  2009.

\bibitem[QS23]{qiu2022tap}
Jiaze Qiu and Subhabrata Sen.
\newblock The tap free energy for high-dimensional linear regression.
\newblock {\em The Annals of Applied Probability}, 33(4):2643 -- 2680, 2023.

\bibitem[Roc97]{rockafellar1997convex}
R~Tyrrell Rockafellar.
\newblock {\em Convex analysis}, volume~11.
\newblock Princeton University Press, 1997.

\bibitem[RP16]{reeves2016replica}
Galen Reeves and Henry~D Pfister.
\newblock The replica-symmetric prediction for compressed sensing with
  {G}aussian matrices is exact.
\newblock In {\em 2016 IEEE International Symposium on Information Theory
  (ISIT)}, pages 665--669. IEEE, 2016.

\bibitem[RS22]{ray2022variational}
Kolyan Ray and Botond Szab{\'o}.
\newblock Variational {B}ayes for high-dimensional linear regression with
  sparse priors.
\newblock {\em Journal of the American Statistical Association},
  117(539):1270--1281, 2022.

\bibitem[Sto13]{stojnic2013framework}
Mihailo Stojnic.
\newblock A framework to characterize performance of {LASSO} algorithms, 2013,
  ].

\bibitem[TAH18]{thrampoulidis2018precise}
Christos Thrampoulidis, Ehsan Abbasi, and Babak Hassibi.
\newblock Precise error analysis of regularized {M}-estimators in high
  dimensions.
\newblock {\em IEEE Transactions on Information Theory}, 64(8):5592--5628,
  2018.

\bibitem[Tal10]{talagrand2010mean}
Michel Talagrand.
\newblock {\em Mean field models for spin glasses: Volume I: Basic examples},
  volume~54.
\newblock Springer Science \& Business Media, 2010.

\bibitem[Tan02]{tanaka2002statistical}
Toshiyuki Tanaka.
\newblock A statistical-mechanics approach to large-system analysis of {CDMA}
  multiuser detectors.
\newblock {\em IEEE Transactions on Information theory}, 48(11):2888--2910,
  2002.

\bibitem[TAP77]{thouless1977solution}
David~J Thouless, Philip~W Anderson, and Robert~G Palmer.
\newblock Solution of `{S}olvable model of a spin glass'.
\newblock {\em Philosophical Magazine}, 35(3):593--601, 1977.

\bibitem[TOH15]{thrampoulidis2015regularized}
Christos Thrampoulidis, Samet Oymak, and Babak Hassibi.
\newblock Regularized linear regression: A precise analysis of the estimation
  error.
\newblock In {\em Proceedings of The 28th Conference on Learning Theory},
  volume~40 of {\em Proceedings of Machine Learning Research}, pages
  1683--1709. PMLR, Jul 2015.

\bibitem[TS11]{turner_sahani_2011}
Richard~Eric Turner and Maneesh Sahani.
\newblock {\em Two problems with variational expectation maximisation for time
  series models}, pages 104--124.
\newblock Cambridge University Press, 2011.

\bibitem[TV04]{tulino2004random}
Antonia~M Tulino and Sergio Verd{\'u}.
\newblock Random matrix theory and wireless communications.
\newblock {\em Foundations and Trends in Communications and Information
  Theory}, 1(1):1--182, 2004.

\bibitem[Ver18]{vershynin2018high}
Roman Vershynin.
\newblock {\em High-dimensional probability: An introduction with applications
  in data science}, volume~47.
\newblock Cambridge university press, 2018.

\bibitem[Vil08]{villani2008optimal}
C.~Villani.
\newblock {\em Optimal Transport: Old and New}.
\newblock Springer Berlin Heidelberg, 2008.

\bibitem[WJ08]{wainwright2008graphical}
Martin~J Wainwright and Michael~I Jordan.
\newblock Graphical models, exponential families, and variational inference.
\newblock {\em Foundations and Trends in Machine Learning}, 1(1--2):1--305,
  2008.

\bibitem[WT05]{WangTitterington2005}
Bo~Wang and D.~M. Titterington.
\newblock Inadequacy of interval estimates corresponding to variational
  {B}ayesian approximations.
\newblock In {\em Proceedings of the Tenth International Workshop on Artificial
  Intelligence and Statistics}, volume~R5 of {\em Proceedings of Machine
  Learning Research}, pages 373--380. PMLR, Jan 2005.

\bibitem[YPB20]{yang2020alpha}
Yun Yang, Debdeep Pati, and Anirban Bhattacharya.
\newblock $\alpha$-variational inference with statistical guarantees.
\newblock {\em Annals of Statistics}, 48(2):886--905, 2020.

\end{thebibliography}

\end{document}